\UseRawInputEncoding
\documentclass[12equipped]{article}
\usepackage{amsfonts,amssymb,amsmath}
\usepackage{setspace}

\usepackage{tikz-cd} 
\input epsf.tex

\usepackage{color}
\usepackage{amssymb}
\usepackage{hyperref}

\def\r#1#2{{\color{red}{#1\over #2}}} 
\def\b#1#2{{\color{blue}{#1\over #2}}}
\def\hh{\tilde{h}}
\def\gg#1#2{{\color{green}{#1\over #2}}}

\title{Holomorphic Floer theory I: exponential integrals in finite and infinite dimensions}

\author {Maxim Kontsevich, Yan Soibelman}

\begin{document}

\maketitle

\newtheorem{thm}{Theorem}[subsection]
\newtheorem{defn}[thm]{Definition}
\newtheorem{lmm}[thm]{Lemma}
\newtheorem{rmk}[thm]{Remark}
\newtheorem{prp}[thm]{Proposition}
\newtheorem{conj}[thm]{Conjecture}
\newtheorem{exa}[thm]{Example}
\newtheorem{cor}[thm]{Corollary}
\newtheorem{que}[thm]{Question}
\newtheorem{ack}{Acknowledgements}
\newtheorem{exe}{Exercise}
\newtheorem{thm-constr}[thm]{Theorem-construction}
\newtheorem{prp-constr}[thm]{Proposition-construction}

\newcommand{\C}{{\bf C}}
\newcommand{\K}{{\bf k}}
\newcommand{\R}{{\bf R}}
\newcommand{\N}{{\bf N}}
\newcommand{\Z}{{\bf Z}}
\newcommand{\Q}{{\bf Q}}
\newcommand{\GG}{\Gamma}
\newcommand{\A}{A_{\infty}}
\newcommand{\g}{{\bf g}}
\newcommand{\p}{{\bf p}}
\newcommand{\CC}{{\mathcal C}}
\newcommand{\MM}{{\mathcal M}}
\newcommand{\LL}{{\mathcal L}}
\newcommand{\FF}{{\mathcal F}}
\renewcommand{\H}{{\mathcal{H}}}
\newcommand{\EE}{{\mathcal E}}
\newcommand{\Aa}{{\mathcal A}}
\newcommand{\DD}{{\mathcal D}}
\newcommand{\UU}{{\mathcal U}}
\newcommand{\ZZ}{{\mathcal Z}}
\newcommand{\epi}{\twoheadrightarrow}
\newcommand{\mono}{\hookrightarrow}
\newcommand\ra{\rightarrow}
\newcommand\uhom{{\underline{Hom}}}
\newcommand\OO{{\mathcal O}}
\newcommand{\WP}{\widehat{P}}
\newcommand{\x}{\widehat{x}}
\newcommand{\y}{\widehat{y}}
\newcommand{\HH}{{\bf H}}
\newcommand{\epp}{\varepsilon}

\newcommand{\XX}{\mathcal X}

\newcommand{\M}{{\mathsf{M}}}

\newcommand{\G}{{\mathsf{G}}}

\newcommand{\Gr}{{\mathsf{Gr}}}
\newcommand{\T}{{\mathsf{T}}}
\newcommand{\h}{{\hbar}}


\begin{abstract}

In the first of the series of papers devoted to our project ``Holomorphic Floer Theory" we discuss exponential integrals and related wall-crossing structures. We emphasize two points of view on the subject: the one based on the ideas of deformation quantization and the one based on the ideas of Floer theory. Their equivalence is a corollary of our generalized Riemann-Hilbert correspondence. In the case of exponential integrals this amounts to several comparison isomorphisms between local and global versions of de Rham and Betti cohomology. We develop the corresponding theories in particular generalizing Morse-Novikov theory to the holomorphic case. We prove that arising wall-crossing structures are analytic. As a corollary, perturbative expansions of exponential integrals are resurgent. Based on a careful study of finite-dimensional exponential integrals we propose a conjectural approach to infinite-dimensional exponential integrals. We illustrate this approach in the case of Feynman path integral with holomorphic Lagrangian boundary conditions as well as in the case of the complexified Chern-Simons theory. We discuss the arising perverse sheaf of infinite rank as well as analyticity of the corresponding ``Chern-Simons wall-crossing structure". We develop a general theory of quantum wave functions and show that in the case of Chern-Simons theory it gives an alternative description of the Chern-Simons wall-crossing structure based on the notion of generalized Nahm sum. We propose several conjectures about analyticity and resurgence of the corresponding perturbative series.

\end{abstract}

\setcounter{tocdepth}{4}

\tableofcontents

\section{Introduction}

This is the first paper in the series devoted  to our project ``Holomorphic Floer Theory".\footnote{Technically, the paper [KoSo12] is also part of this series, but it is devoted to the topic which can be treated independently.}
We started the project in 2014 on our trip to the workshop on  wall-crossing formulas at the Simons Center for Geometry and Physics. Currently various  files on the project  have in total several hundred pages. We have taught  many lecture courses and have given many talks on the subject all over the world, some of which can be found online. Finally it has become clear for us that the project is too big for a single paper or even a book. It is rather a program with many ramifications and  subprojects, each  being interesting for its own sake. For that reason we have decided to publish  the material as a series of papers.  

\subsection{Warming up: Morse theory in holomorphic setting}\label{motivation}

Let $X$ be a complex manifold, $f:X\to \C$ a  holomorphic function with finitely many isolated Morse critical points $\{x_1,...,x_k\}$ (in complex sense) and different critical values $\{z_1,...,z_k\}$. For what follows it will be convenient to assume that $f$ is proper, although this condition can be weakened. Then the points $x_i$  are Morse critical points of the $C^\infty$ function $Re(f)$ having the same Morse index $n=dim_\C X$. Therefore for a general Riemannian metric  the Morse complex of $Re(f)$ has trivial differential. Indeed generically there are no gradient lines between critical points of $Re(f)$.

Equivalently, the differential is trivial on the Floer complex of the pair of complex  Lagrangian submanifolds $L_0=X$ and $L_1=graph(df)$ considered as {\it real} Lagrangian submanifolds of the {\it real} symplectic manifold $T^\ast X$ (the $C^\infty$ symplectic form $\omega_{T^\ast X}$ is the real part of the holomorphic symplectic form $\omega^{2,0}_{T^\ast X}$). Triviality of the Floer differential follows from the fact that for a generic compatible with $\omega_{T^\ast X}$ almost complex structure $J$ on $T^\ast X$ there are no $J$-holomorphic discs\footnote{Such discs are also called pseudo-holomorphic.} with the boundary on $L_0\cup L_1$. 

On the other hand the gradient lines and  the pseudo-holomorphic discs can  appear for  special values of $t\in \C^\ast$ if we replace  $f$ by $f/t$.  The virtual number of gradient lines or pseudo-holomorphic discs can change as $t$ varies. This change is controlled by the wall-crossing formulas (see [KoSo1] for a general discussion). {\it This observation indicates that wall-crossing structures introduced in [KoSo7] as the structures underlying the wall-crossing formulas can be thought of as a replacement of the Floer and Morse complexes in the holomorphic setting}. 

Categorically,  Floer complexes are $Hom's$ in an appropriate Fukaya category $\FF(T^\ast X,  Re(\omega^{2,0}_{T^\ast X})+i\, Im(\omega^{2,0}_{T^\ast X}))$. The above discussion implies that from the categorical perspective one has to consider ``the family of Fukaya categories" $\FF_t(T^\ast X)=\FF(T^\ast X, {1\over{t}}Re(\omega^{2,0}_{T^\ast X})+{i\over{t}} Im(\omega^{2,0}_{T^\ast X}))$ in order to make the wall-crossing phenomenon visible.

 Let us now consider the irregular $D$-module on the punctured complex line  $\C^\ast_t$ \footnote{We will use this kind of notation to stress that $t$ is the standard coordinate on $\C^\ast$.} which is the direct image of the $D$-module $e^{f/t}\OO_{X\times \C^\ast}$ under the natural projection to the second factor. Then the above-mentioned pseudo-holomorphic discs appear iff  $t$ belongs to a Stokes ray of the $D$-module, where a Stokes  ray is the ray $\R_{>0}\cdot(z_i-z_j)$ where $z_i\ne z_j$ are two critical values of $f$.
 Stokes isomorphisms of the spaces of solutions of the corresponding differential equation can be written in terms of the virtual numbers of the discs. The relation between pseudo-holomorphic discs and Stokes isomorphisms is a simple incarnation of our generalized Riemann-Hilbert correspondence, which relates deformation quantization and the Floer theory in the framework of complex symplectic manifolds.  Wall-crossing structures appear in this framework quite  naturally. 
 
 In this paper we will discuss (among other things) wall-crossing  structures in the  case of exponential integrals. 
 By the exponential integral we mean the function in a variable $t\in \C^\ast$ given by the formula $I(t):=I_C(t):=I_C(f,t)=\int_C e^{f/t}vol$ (sometimes we use this term for $I_C(1)$). Here $vol$ is a fixed holomorphic volume form on  $X$ and $C$ is an appropriate integration cycle. One can express $I(t)$  as a linear combination of exponential integrals over specific integration cycles known as {\it thimbles}. They are real non-compact integration cycles formed by the gradient lines of $Re(f/t)$ with respect to an auxiliary Hermitian metric (say K\"ahler metric in case if $X$ is a K\"ahler manifold). As one changes $Arg(t)$ and crosses a Stokes ray the linear combination of the integrals over thimbles gets changed by a Stokes automorphism. Collection of all Stokes automorphisms encodes the corresponding wall-crossing structure which coincides with the one coming from the above-mentioned holomorphic version of Morse theory. 
 
 Roughly, this paper is devoted to a discussion of various aspects and generalizations of this phenomenon. Besides  of the  connection with holomorphic Morse theory, exponential integrals can be put into the framework of exponential Hodge theory as exponential periods. In particular there will be de Rham and Betti aspects of the story. Furthermore, asymptotic expansion of $I_C(t)$ as $t\to 0$ is typically a divergent series whose Borel transform admits endless analytic continuation in the sense of \'Ecalle. In other words, it is {\it resurgent}. Based on the results obtained in the finite-dimensional case we will make several proposals for infinite-dimensional exponential integrals. But before discussing the contents of the paper in detail we would like to say few general words about the Holomorphic Floer Theory and our subsequent papers on the project.

\subsection{What is Holomorphic Floer Theory?}



By {\it Holomorphic Floer Theory}  (HFT for short) we understand an unspecified collection of results and conjectures of Floer-theoretical nature in the framework of {\it complex} symplectic manifolds.\footnote{There is a $2$-categorical upgrade of the HFT in which people discuss a quaternionic version of the Floer complex where instead of pseudo-holomorphic discs one counts certain maps of $3$-dimensional balls satisfying the so-called Fueter equations (see e.g. [DoRez]). We will not discuss this theory in current paper.} We remark that  Lagrangian submanifolds which appears in HFT do not have to be complex (a.k.a. holomorphic) although the ambient symplectic manifold always does. This framework appears in many questions including the one about similarity of  Fukaya categories of symplectic manifolds and  categories of $C^\infty$ holonomic $DQ$-modules (see e.g. [BreSo], [Kap] for the early attempts to solve this puzzle). 

 In the case of complex symplectic manifolds we propose a  {\it generalized Riemann-Hilbert correspondence} (RH-correspondence for short):

{\bf The category of holonomic $DQ$-modules on a complex symplectic manifold $M$ is derived equivalent to the Fukaya category of $M$.}

We will call the $DQ$-module side of the RH-correspondence the ``de Rham side", while the Fukaya side of the RH-correspondence will be called the
``Betti side". This terminology corresponds to the one of Simpson in the non-abelian Hodge theory, where one has algebraic bundles with connections on the de Rham side and locally constant sheaves of their flat sections on the Betti side. We will explain in subsequent papers on the project how ``Dolbeault side" corresponding to Higgs bundles as well as the analog of harmonic bundles appear in HFT. Then we will have the full package of the {\it generalized non-abelian Hodge theory}. The role of the twistor parameter is played by the coordinate $t$. More precisely the geometry which appears  in this paper is related to the $1$-parameter family of {\it holomorphic} rescaled symplectic forms: $\omega^{2,0}\mapsto \omega^{2,0}/t$.  For the generalized non-abelian Hodge theory one considers the family of Fukaya categories corresponding to the family of {\it non-holomorphic} symplectic forms ${\omega^{2,0}\over t}+t\overline{\omega^{2,0}}, t\in \C^\ast$. We will use this family  in one of our forthcoming papers. We will also explain that in the case of complex symplectic surfaces the generalized non-abelian Hodge theory is related to the theory of periodic monopoles.

The above-mentioned Riemann-Hilbert correspondence should be called {\it global} since it deals with the categories associated with the whole symplectic manifold $M$.
 As we will see there are also {\it local} versions of the categories associated with a  neighborhood of a fixed Lagrangian subvariety of $M$. Furthermore there is a local version of the RH-correspondence. In  the simplest case of exponential integrals  discussed in this paper we will be interested only in the corollaries of the  categorical equivalences, namely in the comparison isomorphisms between de Rham and Betti cohomologies, both local and global.

A precise statement of the generalized RH-correspondence requires  more structures, and will be discussed in the subsequent papers of our project. The above vague statement is sufficient for the purposes of the current paper. We make only one remark here. Since supports of holonomic $DQ$-modules can be singular, the RH-correspondence requires a more general notion of the Fukaya category in which objects can have supports which are {\it singular} Lagrangian subvarieties. We will  propose such a theory in one of the subsequent papers on the project.

One can also say that the part of the project which is devoted to the interplay between Fukaya categories and categories of $DQ$-modules  can be called {\it Holomorphic Floer Quantization} (HFQ for short), as it ``quantizes'' complex symplectic manifolds in two different ways: via deformation quantization and via Fukaya categories. Then the RH-correspondence says that these two ways of quantization are equivalent.
As a side remark we mention that HFQ  also gives a natural point of view on the notion of {\it coisotropic brane}. Indeed supports of general $DQ$-modules are coisotropic, and the RH-correspondence claims that there is a subcategory of the category of all $DQ$-modules which is equivalent to the category of branes with Lagrangian support.

The  simplest version of wall-crossing formulas which appear in  relation to  exponential integrals are known in physics as Cecotti-Vafa wall-crossing formulas 
(or $2d$ wall-crossing formulas).  The comparison isomorphisms between de Rham and Betti cohomology associated with the pair $(X,f)$ 
explain these wall-crossing formulas. Exponential integrals can be understood as  {\it exponential periods}, i.e. as  pairings of de Rham cocycles with Betti cycles.
 Such  pairings are incarnations of the {\it global Betti-to-global de Rham isomorphism} in the exponential Hodge theory. On the other hand  this isomorphism 
is a corollary  of the global Riemann-Hilbert correspondence. This observation supports  our terminology of Betti and de Rham sides of the RH-correspondence.

Let us list without comments some  topics which will be discussed in the subsequent papers in the series:

\begin{itemize}

\item{Exponential integrals} (this paper).

\item{Universal Fukaya category as a sheaf on an analytic $2$-stack}.

\item{Fukaya categories which include objects with singular Lagrangian supports}.

\item{Fukaya categories  with parabolic structure}.

\item{Coisotropic branes and non-commutative Hodge structures}.

\item{Deformation quantization as a sheaf on an analytic $2$-stack and the corresponding category of holonomic $DQ$-modules}.

\item{Universal generalized Riemann-Hilbert correspondence for sheaves of categories on analytic $2$-stacks}.

In relation to our paper [KoSo12] devoted to analytic wall-crossing structures we mention the (essentially finished)  paper on

\item{Algebraic wall-crossing structures}

in which algebraicity of the generating series arising in HFT is discussed.

\item{Compactifications of $M$  and Calabi-Yau categories with corners}.

\item{ Holonomic $DQ$-modules with parabolic structure}.

\item{Generalized Riemann-Hilbert correspondence in dimension one: rational, trigonometric
and elliptic cases}.

\item{Generalized Riemann-Hilbert correspondence for higher-dimensional quantum tori}.

\item{Twistor families of $\A$-categories, periodic monopoles and generalized non-abelian Hodge theory in dimension one}.

\end{itemize}

Each of the above topics has many ramifications which we  plan to discuss in the corresponding papers as well.

Rest of the Introduction is devoted to the brief description of the contents of the paper. We remark that in most cases each section of the paper starts with its own short description of the  contents.  Also, the paper contains many remarks which put  the discussion in a bigger context and  explain ``what is going on" from the perspective of future sections or future  papers.

\subsection{De Rham and Betti cohomology in the case of functions}

Sections \ref{Betti and de Rham cohomology} and \ref{case of 1-forms} are devoted to exponential integrals where the exponent is a regular function or a closed $1$-form (i.e. multivalued function).
The approaches in both cases are similar: we define local and global Betti and de Rham cohomology and interpret the integrals as periods, i.e. an isomorphism between global Betti and global de Rham cohomology.  Exponential integrals in the case of functions are defined for non-compact complex manifolds (e.g. complex smooth affine algebraic varieties). An important role in the story is played by the choice of a ``good" compactification $\overline{X}$. 

In Sections \ref{de Rham cohomology},  \ref{Betti cohomology} we discuss definitions of the de Rham and Betti cohomology in the case of functions. Assume for simplicity that $X$ is a smooth complex affine algebraic variety.
Then we can compactify $X$ to a smooth projective variety $\overline{X}$ such that $f$ extends to a regular map $\overline{f}: \overline{X}\to {\bf P}^1$. Let us denote $\overline{X}-X=D_h\sqcup D_v$, where $D_v=\overline{f}^{-1}(\infty)$. In the theory of exponential integrals $I_C(f,t)$  we need to introduce three types of divisors: the divisor $D_0\subset X$ which contains the boundary of the integration cycle $C$, the ``vertical" divisor $D_v$ where $f$ has infinite values, and the ``horizontal" divisor $D_h$ where $f$ has finite limit.


Global de Rham cohomology are basically the hypercohomology of the complex of differential forms on $X$ with respect to the twisted de Rham differential $d+df\wedge (\bullet)$. Then the role of the compatification is to control ``critical points at infinity". Global Betti cohomology are basically the relative integer cohomology of ${X}$  with respect to the union of $D_0$ and the fiber $f^{-1}(z)$ where $Re(z)\ll 0$.

From the perspective of the RH-correspondence  compactifications are needed in order to impose restrictions on the supports of objects of the categories on the de Rham and Betti sides. In  subsequent papers on the project we will explain the role of {\it log extensions} which are partial compactifications of complex symplectic manifolds such that on the compactifying divisors the symplectic form has poles of order one. That will give a more conceptual perspective on a rather ad hoc compactifications which appear in this paper.

In addition to the global de Rham and Betti cohomology one can define their local versions.
This is done by replacing the compactification $\overline{X}$ by the $f$-preimage of the union of small discs about critical values of $f$. In particular, if $f$ is Morse, one can choose a local basis in the relevant relative homology groups consisting of local thimbles and integrate over them  closed (in the twisted sense) de Rham forms. Since in the case of functions local thimbles can be extended to the global ones we obtain the relationship between local and global Betti (co)homology. Furthermore, rescaling $f\mapsto f/t$ we obtain a vector bundle over $\C^\ast_t$ or a vector bundle over the formal punctured or non-punctured disc depending on the choice of cohomology theory.

Altogether one has four comparison isomorphisms between the bundles arising from local and global de Rham and Betti cohomology. They are discussed in Sections \ref{De Rham local to global and Hodge}, \ref{Betti local to global}, \ref{local de Rham to Betti}. In those sections we also make comments about the relation of the comparison isomorphisms with local and global versions of the Riemann-Hilbert correspondence. In particular Stokes isomorphisms which appear in the Deligne-Malgrange irregular Riemann-Hilbert correspondence can be described in terms of the comparison isomorphisms of local and global Betti cohomology.

In Section \ref{perverse sheaf} we also discuss an alternative description of the Betti data in terms of certain perverse sheaves on the line. This description will be useful later when we start to discuss infinite-dimensional exponential integrals and  resurgence of their perturbative expansions. In that case the language of perverse sheaves of infinite rank will allow us to speak about Betti data without working out the foundations of semi-infinite integration cycles.

Wall-crossing formulas and wall-crossing structures in the case of functions are discussed in Section \ref{WCF}. For simplicity we concentrate mostly on the case of Morse functions. Then a chosen homology class is an integer combination of the standard ones   given by thimbles, and the exponential integral is an integer linear combination of the integrals over thimbles. These data give rise to a section of  a finite rank holomorphic vector bundle over the one-dimensional complex manifold, which is $\C_t^\ast$ with finitely many Stokes rays removed. One can easily write down the wall-crossing formula for the jump of the section as $Arg(t)$ crosses the direction of a Stokes ray. 
At the level of individual integrals over thimbles the wall-crossing formulas look  as follows:
$$I_i(t)\mapsto I_i(t)+n_{ij}I_j(t),$$
where $n_{ij}\in \Z$ is the intersection index of thimbles emanating from the critical points $x_i$ and $x_j$ (see Section \ref{WCF} for a more precise statement).

Collection of such wall-crossing formulas is equivalent to a certain Riemann-Hilbert problem. The latter is in turn can be interpreted as the gluing conditions for a new  holomorphic vector bundle over $\C^\ast_t$ which can be extended to the whole line $\C_t$.  

Section \ref{dependence on parameters} is devoted to a generalization of Section \ref{WCF} to the case of families of functions. In this framework one cannot work just with wall-crossing formulas as in Section \ref{WCF}.  Instead in Section \ref{dependence on parameters} one has to use a more  general notion of wall-crossing structure, which is roughly speaking a local system of stability data on a graded Lie algebra (see [KoSo1], [KoSo 7], [KoSo12] about the terminology).

\subsection{The case of closed $1$-forms}

Suppose that instead of a pair $(X,f)$ as before we have a pair $(X,\alpha)$, where $\alpha$ is a general closed holomorphic $1$-form on $X$. The case of functions corresponds to  $\alpha=df$. The generalization to $1$-forms is discussed in mainly Section\ref{case of 1-forms}, while some related questions appear in Sections \ref{applications}, \ref{generalizations}. The case of closed $1$-forms has many subtle points which did not play a role in the case of functions.

For example now $X$ does not have to be non-compact. A non-trivial example is when $X$ is a compact Riemann surface endowed with a holomorphic $1$-form.
Furthermore, definitions of the twisted de Rham and Betti cohomology involve more complicated compactifications. They are discussed in Sections \ref{de Rham for 1-forms} and Section \ref{Betti for 1-forms} respectively. Roughly, the  reason for this more complicated picture is the possibility for a closed $1$-form to have a first order pole at a divisor (i.e. to be a logarithmic form near the divisor).

One can try to define de Rham  and Betti cohomology in a way similar to the one in the case of functions, but replacing the trivial local system by the one associated with the form $\alpha$. Then one faces several problems, since in the global case one should take into account the above-mentioned compactifications. This complicates enormously the comparison of local and global cohomology.
E.g. differently from the case of functions the behavior of global thimbles can be chaotic. As a result the Betti global-to-local isomorphism becomes a non-trivial statement. In order to prove it we revisit some parts of the
Morse-Novikov theory  in the real case. This is done in Section \ref{Betti cohomology for real one-forms}. The comparison theorem of global and local Betti cohomology is proved in Section \ref{non-archimedean Betti local to global} for the cohomology with values in the non-archimedean local system. The comparison of real and holomorphic stories is done in Sections \ref{real and complex cases}, \ref{Stokes rays for 1-forms}. In Section \ref{Stokes rays for 1-forms} we also generalize the Betti global-to-local isomorphism for the cohomology with coefficients in the ``universal local system". Roughly this means that we consider not just a single non-archimedean rank one local system, but the whole neighborhood of it. This neighborhood is given by a non-archimedean tube domain. Stokes isomorphisms also have meaning in this universal setting. Having all these constructions we define the wall-crossing structure corresponding to holomorphic closed $1$-forms in Section \ref{WCS for 1-forms}. We claim its analyticity in Section \ref{analyticity of WCS} in Theorem \ref{analyticity of WCS for 1-forms}. 

The proof of Theorem \ref{analyticity of WCS for 1-forms} is based on several intermediate constructions and results, some of them are postponed until Section \ref{applications}. We also use [KoSo12] where it was proven that analytic wall-crossing structures form a connected component in the space of all wall-crossing structures. We deform the central charge to a rational one (here one has to go beyond the framework of holomorphic $1$-forms to the one of pairs of $C^\infty$ $1$-forms which is discussed in more detail in Section  \ref{Morse-Novikov for pair of 1-forms}). Finally we use a result of pure topological nature proved  in Section \ref{rationality of WCS}. In fact in Section \ref{rationality of WCS} we prove    {\it rationality}  of the corresponding WCS. It implies analyticity. The proof is non-trivial and uses constructible sheaves on the torus which are local systems outside of a point.

Rest of Section \ref{applications} as well as Section \ref{generalizations} is devoted to examples and generalizations of some definitions, constructions  and results from Section \ref{case of 1-forms}.

One of the examples (see  Section \ref{gamma-function}) is the wall-crossing structure associated with a holomorphic $1$-form on a compact complex curve. In particular we construct a WCS associated with the $\Gamma$-function, where all general definitions are illustrated in a one-dimensional example.

An example  of square-tiled Riemann surfaces discussed in Section \ref{square-tiled surfaces} illustrates difficulties of the notion of thimble in the case of $1$-forms. Even the convergence of integrals is not obvious (they can possibly diverge for the higher-dimensional analogs of the square-tiled surfaces). Furthermore, although analyticity of the arising WCS is not difficult to show in this case, we notice a stronger result of rationality of the matrix elements of Stokes isomorphisms.

In Sections \ref{de Rham for pair of 1-forms}, \ref{Morse-Novikov for pair of 1-forms} we discuss a generalization of the story with holomorphic closed one-forms to the case of a pair of compatible closed real-valued $C^\infty$ one-forms. The pair $(Re(\alpha), Im(\alpha))$ associated with a holomorphic closed $1$-form $\alpha$ is an example of such a pair. On the Betti side we obtain another  generalization of the Morse-Novikov theory. Combining  with the de Rham side we obtain the full package of comparison isomorphisms between local and global Betti and de Rham cohomology.

Another generalization is discussed in Section \ref{allowable symplectic forms}. There we consider a pair of real symplectic forms satisfying certain compatibility conditions (the prototype is of course a pair consisting of the real and imaginary parts of a holomorphic symplectic form). In order to understand  properties of the corresponding WCS one should use Fukaya categories. We only sketch this approach.

Comparison isomorphisms between local and global de Rham and Betti cohomology are discussed in Section \ref{comparison isomorphisms for 1-forms}. Similarly to the case of functions the diagram of comparison isomorphisms is a special case of our generalized Riemann-Hilbert correspondence. The generalized RH-correspondence is not discussed in this paper in maximal generality, although some related material is presented in Section \ref{Fukaya categories}. Nevertheless some of its features can be observed in this simple example. For example, all relevant cohomology theories are defined over different ground rings. In order to compare them one should take an appropriate formal expansions sometimes corrected by the Stokes isomorphisms.
The comparison diagram is still conjectural in the case of closed $1$-forms (see Conjecture \ref{comparison of isomorphisms diagram}).

Section \ref{comparison isomorphisms and wheels of lines} is devoted to an interpretation of the Betti global-to-local isomorphisms in terms of the wheels of projective lines. This language was developed in [KoSo12] in order to formulate the property of a wall-crossing structure to be analytic. Instead of coherent sheaves on $\C^\ast$ we now work with their ``universal" versions, which are coherent sheaves over an analytic neighborhood of a wheel of projective lines in a toric variety. 

In Section  \ref{categorical analogs} we discuss the categorification of  the comparison isomorphisms. This is also done  from the perspective of  generalized Riemann-Hilbert correspondence, which will be discussed elsewhere. Nevertheless the reader which has some basic knowledge of Fukaya categories and deformation quantization will be able to understand how to categorify comparison isomorphisms from previous sections. In particular we explain that the relation between local and global Betti cohomology is encoded in the difference between a small Liouville neighborhood of a Lagrangian submanifold (no pseudo-holomorphic discs with the boundary on the Lagrangian submanifold) and the whole symplectic manifold where such discs can appear. The details about local and global versions of the generalized Riemann-Hilbert correspondence will be discussed in subsequent papers of our HFT program.

\subsection{Infinite-dimensional case and quantum wave functions}

Starting with Section \ref{infinite-dimensional case} we discuss infinite-dimensional  exponential integrals including their generalization to the case of closed $1$-forms. Main examples are: Feynman path integral with boundary conditions on two complex Lagrangian submanifolds and complexified Chern-Simons functional integral over the space of connections on a real oriented $3$-manifold. 

Let us consider for example the space $P(L_0, L_1)=\{\phi| \phi:[0,1]\to M\}$ of smooth paths in the complex symplectic manifold $(M, \omega^{2,0}) $ with the endpoints  $\phi(0)\in L_0, \phi(1)\in L_1$, where $L_i, i=0, 1$ are complex Lagrangian submanifolds of $M$. Space $P(L_0, L_1)$ carries a natural holomorphic closed $1$-form given by the integration of $\omega^{2,0}$ over paths in $P(L_0, L_1)$. Going to the universal abelian cover of $P(L_0, L_1)$ we obtain the infinite-dimensional complex manifold endowed with a holomorphic function $S$.

The corresponding exponential integral is not well-defined, but its formal expansions at critical points of $S$ (which are in bijection with $L_0\cap L_1$) are well-defined series. Then we can mimic finite-dimensional considerations and obtain the corresponding wall-crossing structure. Informally local perturbative expansions at the smooth transversal intersection points should be thought of as integrals over local thimbles.
But in the infinite-dimensional case we do not have a rigorously defined Betti side of the story as well as the top degree volume form (``Feynman measure"). One can speculate that a complexified action functional defines a locally-trivial bundle outside of the set of its critical values, but it is not easy to make mathematical sense of this type of statements. In some cases the set of critical values can be dense, thus even the conjecture about such a bundle can be problematic. 

Instead of developing in detail the Betti approach we develop another one based  on the de Rham considerations. The latter is expected to be equivalent to the ill-defined infinite-dimensional exponential integrals, similarly to the relation between Feynman path integrals and the Hamiltonian formalism in quantum mechanics. In our case we develop a theory of quantum wave functions. Roughly, such a theory should produce a holonomic cyclic $DQ$-module associated with a complex Lagrangian submanifold or (conjecturally) with a possibly singular Lagrangian subvariety. Furthermore, adding a non-trivial  Hamiltonian to the story one should be able not only  transport the Lagrangian submanifolds but also the corresponding quantum wave functions.

Development of the theory of quantum wave functions requires more data and more careful considerations. It is based on the notion of Harish-Chandra pair, which is one of the standard techniques in deformation quantization of symplectic manifolds. Harish-Chandra pairs appear naturally in the framework  of formal differential geometry of Gelfand and Kazhdan (see [GeKazh]). The theory takes care about the change of local coordinates (symplectic or quantum in our case) as long as there is a standard local model. We review all that in Section \ref{Harish-Chandra pairs}. Next we add Lagrangian submanifolds to the story and introduce the above-mentioned additional data which we call {\it quantum wave function structure} in Section \ref{QWFS}.\footnote{In the present paper we discuss only {\it smooth} Lagrangian subvarieties for which the language of Harish-Chandra pairs is adequate. The generalization to the singular case will be discussed elsewhere.}

For a pair of quantum wave functions one can define their pairing which belongs to $\C[[t]]$ in the simplest case of transversal intersection. More generally it contains fractional powers of $t$ and polynomials in $\log t$.
The arising series are typically divergent.

We expect that for a large class of quantum wave functions discussed in Section \ref{resurgence} under the name of {\it normalizing objects}, the Borel transforms of the pairings admit analytic continuation along any path avoiding a countable set of singularities. In other words, they are resurgent. 
The normalization problem of quantum wave functions is addressed in Section \ref{normalization}.

Another question discussed in Section \ref{sums and polytopes} is the resurgence properties of sums over integer points of rational polytopes generalizing Euler-Maclaurin formula. These results are used later in Section \ref{Nahm sums} in the discussion of Nahm sums. The unifying concept is the one of {\it resurgence package} which we discuss in the next subsection.

\subsection{Resurgence packages}\label{resurgence package}

By a {\it simple} resurgence package we understand the following data 1)-3) satisfying the following condition 4):

1) At most countable discrete subset $\{z_i\}_{i\in I}\subset \C$.

2) For each $i\in I$ a power series $\phi_i=c_{i,0}+tc_{i,1}+t^2c_{i,2}+...\in \C[[t]]$.

3) A collection of integers $(n_{ij})_{i,j\in I, i\ne j}, n_{ij}\in \Z$.

4) Consider the Borel transforms $B(\phi_i):=B(\phi_i)(s)=\sum_{k\ge 0}{c_{i,k}\over{k!}}(s-z_i)^k$. Then each $B(\phi_i)$ is required to be a  germ of analytic function at $z_i$, which admits the endless analytic continuation in the sense of \'Ecalle. Keeping the same notation for the analytic continuation, we require that near the point $z_j, j\ne i$ we have $B(\phi_i)=n_{ij}B(\phi_j)log(s-z_j)+f_{ij}(s)$, where $f_{ij}$ is a holomorphic function in the punctured neighborhood of $z_j$.

Simple resurgence packages with finite sets $I$ appear naturally in the study of exponential integrals for holomorphic  functions  with  holomorphic Morse critical points and pairwise different critical values and a given volume element. 
In this case $\phi_i$ will be the formal expansions of the  modified exponential integrals $e^{-z_i/t}(2\pi t)^{-dim_\C X/2}I_i(t)$.

For non-Morse critical points one has a more general notion of {\it resurgence package} whose definition is similar, but  we do not discuss it here.

Resurgence packages appear quite naturally in HFT.  One can hope that any analytic wall-crossing structure gives rise to a resurgence package (cf. [KoSo12]). In this paper our main source of analytic wall-crossing structures and resurgence packages will be a pair of complex Lagrangian subvarieties of a complex symplectic manifold.

Furthermore, we will also  discuss various generalizations of 1)-4) to the infinite-dimensional case. This will include the heuristic picture from physics and precise mathematical prediction for the resurgence package. We will also illustrate this discussion in  examples.

\subsection{Chern-Simons theory and generalized Nahm sums}
Section \ref{CS theory} is devoted to  several conjectures in Chern-Simons theory for compact closed oriented $3$-dimensional manifolds. The Chern-Simons theory depends on a compact simple  group $G_c$ (gauge group of the theory) and on the positive integer $k$ (level of the theory). One can replace $G_c$ by its complexification $G$ and introduce a small parameter $\hbar={2\pi i\over k}$. Then the {\it complexified} Chern-Simons theory can be approached perturbatively by looking at the formal expansions of the partition function as $\hbar\to 0$.  In this way one obtains a collection of formal series parametrized by connected components of the critical locus of the Chern-Simons functional
$$CS(A)=\int_{M^3}Tr\left({1\over{2}}A\wedge dA+{1\over{3}}A\wedge A\wedge A\right).$$
These series together with the Stokes isomorphisms describing the ``interaction" of the series associated with different connected components should give full information about the partition function of the Chern-Simons theory. By analogy with the finite-dimensional exponential integrals, the same data should describe the ``Chern-Simons wall-crossing structure". We propose a different description of the latter. The Betti version of it is based on the hypothetical perverse sheaf of infinite rank, while the de Rham version is based on the  theory of quantum wave functions. Conjecturally these two versions should agree. Furthermore, we propose a couple of conjectures which relate this wall-crossing structure (or rather the Borel resummation of the local perturbative expansions) with the conventional Chern-Simons theory associated with the group $G_c$ and integer values of $k$.

In relation to the complexified Chern-Simons theory we will also discuss in Section \ref{Nahm sums}  the  {\it generalized Nahm sums}. The key role in the construction of the related resurgence package is played by the dilogarithm function. Considerations of Section \ref{Nahm sums} make sense beyond their applications to Chern-Simons theory.

Generalized Nahm sum is an expression
$$Z_{N}=\sum_{\substack{0\le j_1,...,j_{d}\le N-1\\
\overline{j}/N:=(j_1/N,...,j_{d}/N)\in P}}\chi(\overline{j})\prod_{1\le i\le d}((j_i)_q!)^{a_{i}}q^{{1\over{2}}\sum_{1\le i_1,i_2\le d}b_{i_1i_2}j_{i_1}j_{i_2}+\sum_i c_ij_i}$$
which depends on a positive integer $N$, a rational polyhedron $P\subset [0,1)^d$, a rational symmetric matrix $(b_{ij})$, collection of integers $(a_i)$, and a finite order character $\chi$ of the group $\Z^d$. Here $q=e^{2\pi i\over N}$.
As $N\to \infty$ the properly defined critical points of this expression approach to the critical points of the multivalued function
$$f_{a,b,\chi}({x})=-\sum_{\beta}{a_i}Li_2(x_i)-\sum_{1\le i\le d} log(\chi_i)log(x_i)+{1\over {2}}\sum_{i_1,i_2}b_{i_1i_2}log(x_{i_1})log(x_{i_2}), $$
where $x=(x_1,...,x_d)\in (\C^\ast)^d,$
and the classical dilogarithm function  is defined by the analytic continuation of the series $Li_2(z)=\sum_{k\ge 1}{z^k\over{k^2}}$.

We will also discuss the associated resurgence package. For that we interpret the critical points of $f_{a,b,\chi}$ as intersection points of two complex Lagrangian submanifolds in $(\C^\ast)^{2d}$. 
 The first one is essentially $graph(df_{a,b,\chi})$ and the second one is a union of Lagrangian subtori corresponding to the faces of $P$. We endow each Lagrangian subvariety with quantum wave functions. Then the series $\phi_a$ in the resurgence package are given by the pairings of these quantum wave functions. 
 


The resurgence package gives via Borel resummation a collection of analytic functions in sectors in $\C_\hbar$. We conjecture that an integer linear combination of these functions evaluated at $\hbar={2\pi i\over N}$ coincides with $Z_N$.
Taking $N=k+c$ where $c$ is the dual Coxeter number of $G$ one should recover the
Witten-Turaev-Reshetikhin invariants of $3$-dimensional manifolds as well as their resurgence properties. 

{\it Acknowledgements.}
The work of M.K. and Y.S.  was partially supported by ERC--SyG project, Recursive and Exact New Quantum Theory (ReNewQuantum), 
which received funding from the European Research Council (ERC) under the European Union's Horizon 2020 research and innovation program
under grant agreement No 810573. The work of Y.S.  was in addition partially  supported  by Simons Foundation  grant MP-TSM-00001658.  
The authors  thank to IHES for
excellent research conditions.

\part{Finite-dimensional exponential integrals}

\section{Betti and de Rham cohomology of a variety with function}\label{Betti and de Rham cohomology}

Some results of this section are either known to experts or belong to the folklore. We apologize to the reader for those places where the references are inadequate. 

Although exponential integrals  can be considered in the category of complex manifolds with holomorphic functions on them, we will restrict ourselves to the algebra-geometric framework. This makes a clear parallel with the de Rham theory of a single variety (no function), which makes sense over an arbitrary field of characteristic zero. Furthermore  for a general holomorphic function $f$ there is a subtle issue of imposing the growth restrictions of $f$ at infinity.

\subsection{Data for exponential integrals}\label{data for exp int}

We are going to consider exponential integrals over  chains with boundary allowing the boundary to belong to a given closed algebraic subset of an  algebraic variety. More precisely the exponential integral will depend on the following data:

1) Smooth  $n$-dimensional algebraic variety $X$ over the field $\C$ of complex numbers.

2) Closed algebraic subset $D_0\subset X$, $dim_\C D_0<n$.

3) Regular function $f\in \mathcal{O}(X):=\Gamma(X,\mathcal{O}_X)$.

4) Algebraic volume form $vol:=vol_X\in \Gamma(X,\Omega^n_X)$, which we assume for simplicity to be nowhere vanishing.

5) A singular  Borel-Moore $n$-chain $C$ (``chain of integration") given by an integer linear combination of locally closed real semialgebraic oriented $n$-dimensional submanifolds of $X(\C)$ (see e.g. [KoSo13], Section 8 for the notion of real semialgebraic Borel-Moore chain) and satisfying the following properties:

5a) $Supp \,\partial C\subset D_0(\C)$;

5b) $Re(f)_{|Supp(C)}: Supp(C)\to \R$ is a proper map, bounded from above.

\begin{defn} Having the data 1)-5) we  define the exponential integral as $I:=I_C(f)=\int_C e^f vol$ (the integral is absolutely convergent by 5b)).

\end{defn}

\begin{defn} Given the data 1)-4) we define the bifurcation set $Bif:=Bif((X,D_0),f)\subset \C$ as the minimal finite set of points such that for any $z\in \C-Bif((X,D_0),f)$ there exists an open (in analytic topology) neighborhood $U\subset \C$ and homeomorphism $f^{-1}(U)\simeq U\times f^{-1}(z)$ which is compatible with natural projections of both spaces to $U$, and such that it induces a homeomorphism $f^{-1}(U)\cap D_0\simeq U\times (f^{-1}(z)\cap D_0)$.\footnote{The definition makes sense. Indeed, the condition is satisfied on non-empty Zariski open subset of $\C$. Its  complement is a finite set.}

\end{defn}

Suppose that we can compactify $X$ to a smooth projective variety $\overline{X}$ such that $f$ extends to a regular map $\overline{f}: \overline{X}\to {\bf P}^1$. Let us denote $\overline{X}-X=D_h\sqcup D_v$, where $D_v=\overline{f}^{-1}(\infty)$.

{\it In what follows unless we say otherwise, we will make the following simplifying assumptions:}

i) $D_0\cup D_h\cup D_v$ is a divisor with normal crossings.

ii) The restriction of $\overline{f}$ to any smooth open stratum of $D_h-D_v$ or  $(D_h-D_v)\cap \overline{D}_0$ does not have critical points.

Notice that the assumptions i) and ii) are different from the assumptions of non-degeneracy in [Moch5], Def. 2.6.
The assumption i) can be always achieved by choosing an appropriate compactification of $X$, whereas the assumption ii) is imposed in order to simplify the exposition.

The role of the divisors in the intuitive picture of exponential integrals is the following: $D_0$ contains ``endpoints" of integration cycles, the ``vertical divisor" $D_v$ is the ``divisor at infinity" where the function $f$ grows indefinitely, and the ``horizontal divisor"  $D_h$ is the ``divisor at infinity" where the function $f$ has finite limit. The assumption ii) implies that the set $Bif$ consists of points $f(x), x\in X$ such that $x$ is either a critical point of $f$ or a critical point of the restriction $f_{|D_0}$. Thus ii) can be called the  {\bf absence of critical points at infinity assumption}.

Generalizing the notion of exponential integral let us consider $I(t):=I_C(t)=\int_C e^{f/t} vol$ as a function of $t\in \C^\ast$, and which we also call the exponential integral. In order to make sense of $I_C(t)$ for a given $t$ we assume that $Re(f/t)_{|Supp(C)}$ satisfies 5b). It follows that the set of allowed $t$ belongs to an open sector in $\C^\ast$. If $Supp(C)$ is not compact then the sector is strict. For a compact $Supp(C)$ the set of allowed $t$ is the whole $\C^\ast$.

\begin{rmk}\label{exp int and comparison isomorphism}
As we will see later the exponential integral $I_C(t)$  is an example of a matrix element of the global  de Rham-to-Betti isomorphism. In what follows we will discuss different types of cohomology theories: global and local de Rham and Betti cohomology, as well as four comparison isomorphisms between them.
\end{rmk}

\subsection{Twisted global de  Rham cohomology}\label{de Rham cohomology}

In this subsection we can replace the field $\C$ by an arbitrary field of characteristic zero.

We assume the set up of the Section \ref{data for exp int}. Then the holomorphic (in fact algebraic) volume form $vol_X$ is closed with respect to the differential $d_f=d+df\wedge (\bullet)$\footnote{We can equally well to use the differential  $d_f=d-df\wedge (\bullet)$. In fact as we will see, it makes sense to consider the whole family of differentials $td+df\wedge (\bullet), t\in \C$, and the above choices correspond to $t=1$ and $t=-1$ respectively.}Notice that the differential $d_f$ gives rise to a complex of sheaves (in Zariski topology)
$$\Omega^{\bullet}_X:=\Omega^0_X\to \Omega^1_X\to...\to \Omega^n_X.$$
Recall that $D_0$ is a divisor with normal crossings,
we denote by $\Omega^{\bullet}_{X,D_0}$ a subcomplex consisting of differential forms whose restriction to $D_0$ is equal to $0$. 

\begin{defn} We define the twisted global de Rham cohomology as the graded abelian group
$$H_{DR,glob}^{\bullet}((X,D_0),f)={\mathbb H}^{\bullet}(X_{Zar},(\Omega^{\bullet}_{X,D_0},d_f)).$$

\end{defn}

Sometimes we will skip either of the words ``twisted", ``global" or both, if it will not lead to a confusion.

\begin{rmk} \label{analytic exponent}
Notice that in the complex analytic  case the expression $e^f$ makes sense. In particular, if $\alpha$ is a holomorphic de Rham differential form which is closed with respect to the de Rham differential then the analytic differential form $e^{-f}\alpha$ is $d_f$-closed. The converse is also true.

\end{rmk}

\begin{defn}\label{rotated de Rham cohomology}
Let us fix $t\in \C^\ast$. We define the graded $\C$-vector space $H^\bullet_{DR,glob, t}(X, D_0,f)$ as $H_{DR,glob}^{\bullet}((X,D_0),f/t)$.
\end{defn}

{\it In order to alleviate the notation we will often assume that $D_0=\emptyset$ and  drop  it  from the notation}. 

The following result  was conjectured by the first author and recently proved by Mochizuki (see [Moch6]). We will need it later when considering a more general case of closed $1$-forms.
 
 \begin{thm}\label{degeneration in Kahler case}
 Let $X$ be a non-compact complex manifold which admits a K\"ahler metric, and $f: X\to \C$ a holomorphic function such that the set of critical points $Crit(f)$ is compact. Then $R\Gamma(X, (\Omega_X^{\bullet}[[t]], td+ df\wedge(\bullet)))$ as well as its hypercohomology $\mathbb{H}^{\bullet}(X,(\Omega^{\bullet}_X[[t]],td+df\wedge(\bullet)))$ is a free $\Z$-graded $\C[[t]]$-module.
 
 \end{thm}

\subsection{Local de Rham  cohomology}\label{De Rham local to global and Hodge}

In the case $f=0$ the triples $(X,D_0,f)$ reduce to  pairs $(X,D_0)$. The latter are studied in  Hodge theory. In particular, the cohomology groups $H^{\bullet}(X,D_0)$ carry a mixed Hodge structure. In the case when $X$ is projective and $D_0=\emptyset$ the Hodge structure is pure. The analog of the latter story in the case of exponential Hodge theory (see e.g. [KoSo5], [FreJos]) corresponds to the case when $D_0=\emptyset$ and $f: X\to \C$ is a proper map. 

\begin{rmk} \label{non-proper function} In case when $f$ is not a proper map, we  can sometimes extend it to a  proper map $\overline{f}: \overline{X}\to \C$, such that $D_h:=\overline{X}-X$ is a normal crossing divisor, $D_h=\cup_i D_i$, such that  any intersection $D_{i_1}\cap...\cap D_{i_k}$ of its irreducible components does not contain critical points of $\overline{f}$, and the restriction of $\overline{f}$ to such an intersection is a locally trivial analytic fibration over its image in $\C$. In that case $f$ is topologically isotrivial at infinity in  the sense of [KoSo5].\footnote{A function $f$ on {smooth}  algebraic variety $X/\C$ is called   topologically isotrivial at infinity if there exists
a  $C^\infty$-manifold with boundary $U\subset X(\C),\,\,\dim_\R U=2\dim X$ such that $f_{|U}:U\to \C$ is proper, and a homeomorphism
 $h:\partial U \times [0,+\infty)\simeq X(\C)- {int}\, U\, ,\,\,\,h(x,0)=x\mbox{ for }x\in \partial U\,,$
such that $f(h(x,t))=f(x)$ for any $x\in \partial U$.} 
 Then the de Rham  cohomology for the pairs $(X,f)$ and $(\overline{X}, \overline{f})$ coincide.

\end{rmk}

\begin{rmk}\label{spectral sequence}

In this remark we allow  $D_0\ne \emptyset$. We will not use this remark anywhere in the paper, so we omit here the details.

Similarly to the usual Hodge theory  there  is a  spectral sequence converging to  $H^\bullet_{DR}(X, D_0,f)$. Assume that the pair $(X, f)$ as well as all pairs $(D_{0, j},f)$ for each smooth irreducible component $D_{0,j}\subset D_0$ satisfy the condition from Remark \ref{non-proper function}. Then the spectral sequence degenerates at the term $E_2$ (this happens e.g. when $f$ is proper). The degeneration follows from the weight theory for exponential mixed Hodge structures developed in [KoSo5].
\end{rmk}

Let $Z_i:=f^{-1}(z_i)\cap Crit(f)$ be the critical fiber of $f$ over the point $z_i$. 

\begin{defn}\label{local de Rham}
We define the local de Rham cohomology $H^\bullet_{DR,loc}(X,f)$ as the  $\C[[t]]$-module $\oplus_{i\in S}\mathbb{H}^{\bullet}(U(Z_i),(\Omega^{\bullet}_X[[t]],td+df\wedge(\bullet))),$
where $U(Z_i)$ is any open neighborhood of the critical locus $Z_i$ in $X$ in analytic topology. Each summand in the RHS is called the local de Rham cohomology associated with $Z_i$ (or $z_i$) and is denoted by $H^\bullet_{DR,loc,z_i}(X,f)$.

\end{defn}

\begin{rmk}\label{formal neighborhood}
In the Definition \ref{local de Rham} the hypercohomology does not depend on the choice of the open neighborhood, since the complex of sheaves is acyclic outside of $Crit(f)$. One can show that $U(Z_i)$ in the definition can be replaced by the completion $U_{form}(Z_i)$ of $X$ along $Z_i$ (i.e. by the ``formal neighborhood\rq\rq{} of $Z_i$ in $X$).

\end{rmk}
It follows from the Proposition below  that the local de Rham cohomology is a free $\C[[t]]$-module. Moreover, it carries an algebraic regular singular connection $\nabla_{loc}$  with pole of  order $2$ at $t=0$.

Global de Rham cohomology gives rise to a $D$-module on $\C$. The Proposition below summarize some  of its properties (see e.g. [KaKoPa1], [Sab2], [Sab5] for more details and proofs). They form a part of the  non-commutative Hodge theory package (see [KaKoPa1]).

\begin{prp}\label{bundle with meromorphic connection}
Assume that $f$ is proper. Then:

a) The coherent sheaf $\mathcal{H}_{DR,glob}^\bullet(X,f)$ on $\C$ defined by a finite type $\C[t]$-module  
$${H}^\bullet_{DR,glob}(X,f)[t]:=\mathbb{H}_{Zar}^\bullet(X\times {\bf  C}^1_t, (pr_X^\ast(\Omega^\bullet_X, td_X+df\wedge(\bullet))))$$
is in fact a graded vector bundle. Specialization ${H}^\bullet_{DR,glob}(X,f)[t]$ to any $t=t_0\in \C^\ast$ is naturally isomorphic  to $H_{DR, glob, t_0}^{\bullet}(X,D_0, f)$ from Definition \ref{rotated de Rham cohomology}.  The  restriction of the graded vector bundle $\mathcal{H}_{DR,glob}^\bullet(X,f)$ to $\C^\ast$  carries a flat connection $\nabla$.

b) The connection $\nabla$ has regular singularity at $t=\infty$ and second order pole at $t=0$. Over the field Laurent series $\C((t))$  we have an isomorphism $$(\mathcal{H}_{DR,glob}^\bullet(X,f),\nabla)\simeq \oplus_{i\in S}e^{z_i/t}\otimes (E_i,\nabla_i),$$
where $\{z_i\}_{i\in S}$ is the set of critical values of the function $f$, 
and $e^{z_i/t}$ denotes the irregular $D$-module on $\C$ corresponding to this exponential function, and $(E_i,\nabla_i)$ is a regular singular connection (a.k.a. Hukuhara-Levelt-Turrittin, or HLT decomposition).

c) The fiber of $\mathcal{H}_{DR,glob}^\bullet(X,f)$ at $t=0$ is isomorphic to the hypercohomology of the complex of sheaves $(\Omega_X^{\bullet}, df\wedge (\bullet))$ in either Zariski or analytic topology. It can be also computed as
$$\oplus_{i\in S}\mathbb{H}^\bullet(U(Z_i), (\Omega^{\bullet}_X,  df\wedge (\bullet))),$$
where $U(Z_i)$ is a sufficiently small neighborhood of the critical locus $Z_i$ in analytic topology.

d)
The following {\bf de Rham  global-to-local  isomorphism} of $\C[[t]]$-modules endowed with irregular  connections holds:
$$iso_{DR,loc}: {H}^\bullet_{DR,glob}(X,f)[[t]]:={H}^\bullet_{DR,glob}(X,f)[t]\otimes_{\C[t]}\C[[t]]\simeq$$
$$ \oplus_{i\in S}\mathbb{H}^{\bullet}(U(Z_i),(\Omega^{\bullet}_X[[t]],td+df\wedge(\bullet)))=H^\bullet_{DR,loc}(X,f).$$
After tensoring each summand with $\C((t))$ and the rank one $D$-module $\mathcal{M}=e^{z_i/t}\cdot \C((t))$ we obtain the HLT decomposition b).

e) For any $t\in \C^\ast$ there is a non-degenerate pairing
$$H_{DR,glob,-t}^{\bullet}(X,f)\otimes H_{DR,glob,t}^{\bullet}(X,f)\to \C[-2\,dim_\C X],$$
which extends to a non-degenerate pairing at $t=0$.

\end{prp}

More generally, let us assume that we can compactify $X$ by a normal crossing divisor $D_v$, i.e. $D_v=\overline{f}^{-1}(\infty)$ for the corresponding extended map $\overline{f}:\overline{X}\to {\bf P}^1$.  Let us consider the sheaf $\Omega^{\bullet}_{\overline{X},f}$ on $\overline{X}$ such that for any Zariski open set $U$ we have: $\Omega^{\bullet}_{\overline{X},f}(U)$ consists of log-differential forms on $U-D_v$ which are  log-forms along $D_v\cup D_h$ and remain such after wedge product with $df$. Then the rank of the hypercohomology
$${\mathbb H}^{\bullet}(\overline{X},(\Omega^{\bullet}_{\overline{X},f},t_1d+t_2 df\wedge (\bullet)))$$
 does not depend on $(t_1,t_2)\in \C^2$ (see e.g. [KaKoPa2], Theorem 2.18, where the above sheaf of forms is denoted by $\Omega^{\bullet}_{\overline{X}}(log(D_v\cup D_h),f)$).


I

\subsection{Global Betti cohomology}\label{Betti cohomology}

Let us assume until the end of this section the framework of  Section \ref{data for exp int} unless we say otherwise.
Let us fix a  real number $c>0$ and consider the following singular homology of a pair:
$$H_{\bullet}(X,D_0\cup f^{-1}(Re(z)\le -c),\Z)\simeq H_{\bullet}(X,D_0\cup f^{-1}(-c),\Z).$$
Clearly we obtain an inductive system with respect to the natural order on the numbers $c$.
If $-c<Re(Bif((X,D_0),f))$ (i.e. $-c$ is smaller than real parts of finitely many points from the bifurcation set) than the relative homology stabilizes. We remark that simplifying assumptions i) and ii) from Section \ref{data for exp int} are not necessary for that.

\begin{defn} We define the global Betti homology of $(X, D_0, f)$ as
$$H_{\bullet}^{Betti, glob}((X,D_0),f,\Z):=H_{\bullet}((X,D_0),f^{-1}(-\infty),\Z),$$
where the notation in the RHS means the stabilized relative homology groups for sufficiently large $c>0$.

We define global Betti cohomology in a similar way:
$$H^{\bullet}_{Betti,glob}((X,D_0),f,\Z):=H^{\bullet}((X,D_0),f^{-1}(-\infty),\Z).$$

\end{defn}
Sometimes we will omit the word ``global" in the text. Hopefully it will not lead to a confusion.

In what follows it will be useful to have a description of the global (and later the local) Betti cohomology groups in terms of the following data.

Let $B\subset \R^2$ denote a topological submanifold with the boundary, such that $B$ is homeomorphic to a closed unit disc. We assume that $\partial B\cap Bif=\emptyset$, and we fix a  point $b\in \partial B$.

With a pair $(B,b)$ as above and $k\in \Z_{\ge 0}$ we associate the abelian group
$$V(B,b):=H^k(f^{-1}(B), (D_0\cap f^{-1}(B))\cup f^{-1}(b),\Z).$$
If $(B^{\prime},b)\subset (B,b)$ in the sense that $B^{\prime}$ is contained in $B$ and the marked point $b$ is the only common point of the boundaries,
there is a natural restriction homomorphism $\rho_{B,B^{\prime}}:V(B,b)\to V(B^{\prime},b)$. Moreover if $(B,b)$ can be retracted to a ``bouquet" of $(B^{\prime},b)$ and $(B^{\prime\prime},b)$ in such a way that there are no critical values of $f$ on $B-(B^\prime\cup B^{\prime\prime})$ (see the figure)

\vspace{2mm}

FIGURE OF THE BOUQUET

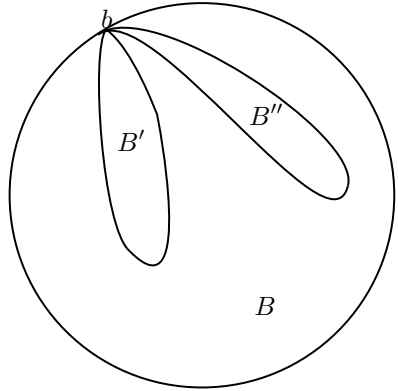
\tikzset{every picture/.style={line width=0.75pt}} 

\begin{tikzpicture}[x=0.5pt,y=0.5pt,yscale=-1,xscale=1]

\draw   (95,151.5) .. controls (95,71.14) and (160.14,6) .. (240.5,6) .. controls (320.86,6) and (386,71.14) .. (386,151.5) .. controls (386,231.86) and (320.86,297) .. (240.5,297) .. controls (160.14,297) and (95,231.86) .. (95,151.5) -- cycle ;
\draw    (161.96,29.89) .. controls (201.96,-0.11) and (361,106) .. (351,144) ;
\draw    (161.96,29.89) .. controls (201.96,-0.11) and (338,202) .. (351,144) ;
\draw    (167.96,25.89) .. controls (188,40) and (205.92,89) .. (206.46,90) .. controls (207,91) and (237,246) .. (185.96,193.89) ;
\draw    (167.96,25.89) .. controls (156,42) and (165,177) .. (185.96,193.89) ;

\draw (174,100) node [anchor=north west][inner sep=0.75pt]   [align=left] {$\displaystyle B'$};
\draw (274,80) node [anchor=north west][inner sep=0.75pt]   [align=left] {$\displaystyle B''$};
\draw (278,226) node [anchor=north west][inner sep=0.75pt]   [align=left] {$\displaystyle B$};
\draw (162,9) node [anchor=north west][inner sep=0.75pt]   [align=left] {$\displaystyle b$};

\end{tikzpicture}

\vspace{2mm}

then the restriction maps give rise to an isomorphism (additivity property)
$$V(B,b)\simeq V(B^{\prime},b)\oplus V(B^{\prime\prime},b).$$
Notice that in this case the bifurcation points inside of $B$ are split between $B^{\prime}$ and $B^{\prime\prime}$.
The above isomorphism allows us to reduce the description of the data $V(B,b)$ to the case when $B$ contains a single bifurcation point inside.

Now we see that  $H^{\bullet}_{Betti,glob}((X,D_0),f,\Z)\simeq V(B,b)$ where $B=D(r)$ is the standard disc $|z|\le r$ of a sufficiently large radius $r$, so that $B\cap Bif=\emptyset$ and $b=\partial B\cap \R_{<0}$.

\begin{defn}\label{rotated Betti cohomology}
Let $t\in \C^\ast$. We define the graded abelian group
$$H^{\bullet}_{Betti,glob,t}((X,D_0),f,\Z):=H^{\bullet}((X,D_0),(f/t)^{-1}(-\infty),\Z).$$

\end{defn}

Clearly Betti homology and Betti cohomology are dual to each other after extension of scalars to $\Q$. Moreover, one has the following result.

\begin{prp}(Poincar\'e duality).
In the above notation let $X^{\prime}=\overline{X}-D_v-\overline{D}_0$ and $D^{\prime}_0=D_h-(D_h\cap D_v)$.Then  we have the following isomorphism\footnote{Using the formalism of six functors one can show that the result holds without the assumptions i) and ii from Section \ref{data for exp int}).}

$$H_\bullet^{Betti,glob}((X,D_0),f)\simeq H_{Betti,glob}^\bullet((X^{\prime},D^{\prime}_0),-f)[2\,dim_\C X].$$
\end{prp}

{\it Proof.}
Let us define $X^0=X-D_0$. We can add a boundary (with corners) to $X^0$ by making real blow-ups and thus embedding $X^0$ into the
real compact manifold with corners $\overline{X^0}$. The non-intersecting components $\partial_\pm X^0$ of the boundary $\overline{X^0}-X^0$ can be chosen in such a way that $f^{-1}(-\infty)\in \partial_+X^0$ and $(-f)^{-1}(-\infty)\in \partial_-X^0$. We add the boundary $S^1$ at infinity to $\C=\R^2$ in such a way that left and right half-circles $S^1_\pm$ are images of $\partial_\pm X^0$. Then the result follows from the usual Poincar\'e duality in the compact case.
$\blacksquare$

As immediately follows from the definition, the family of abelian groups  $(H^{\bullet}_{Betti,glob, t}((X,D_0),f,\Z))_{t\in \C^\ast}$ gives rise to a {\it local system} (i.e. locally-constant sheaf) on $\C^\ast$ which we denote by   $\mathcal{H}^{\bullet}_{Betti,glob}((X,D_0),f)$.

\subsection{Local Betti cohomology and Betti global-to-local isomorphism}\label{Betti local to global}

Let us assume for simplicity that $D_0=\emptyset$ and drop it from the notation. Recall that  $S=\{z_1,...,z_k\}$ denote the finite  set of critical values of $f$.

\begin{defn} \label{Betti local cohomology}
For each critical value $z_i$, a sufficiently small positive $\varepsilon$ and $t\in \C^\ast$ we define the local Betti cohomology  associated with the point $z_i$ and $t$ as a graded abelian group 
$$H^\bullet_{Betti,loc,z_i,t}(X,f)=V(D(z_i/t,\varepsilon),b_t)=H^{\bullet}((f/t)^{-1}(D(z_i/t,\varepsilon)),f^{-1}(b_t),\Z),$$ 
where $D(z_i/t,\varepsilon)$ is a closed disc of radius $\varepsilon$ with the center at $z_i/t$, and $b_t:=b_\theta=z_i+\varepsilon\cdot e^{i\theta}, \theta=Arg(t)\in \R/2\pi i\Z$ is a point on the boundary circle. The direct sum of these groups over all $z_i$ is called the local Betti cohomology $H^\bullet_{Betti, loc,t}(X,f)$.
\end{defn}
Similarly to the global case, the  family of the local  Betti cohomology forms a local system $\mathcal{H}^\bullet_{Betti, loc}(X,f)$
on $\C^\ast_t$. 

Since both global and local Betti cohomology groups are described topologically in terms of the vector spaces $V(B,b)$, one can use this description in order to relate global and local Betti cohomology.
Indeed,  choose a sufficiently large disc, containing the set $S$. From each point $b_t^{(i)}=b_{\theta_t+\pi}^{(i)}$ let us shoot a ray $l_{\theta_t+\pi}^{(i)}$ in the direction $\theta_t+\pi$. Since $\theta_t$ is the same for all discs, we obtain a finite collection of parallel rays.

We can construct a homotopy of the big disc with portions of rays $l_{\theta_t+\pi}^{(i)}$ inside to the same disc but with rays deformed near the boundary circle in such a way that they all intersect at the same point denoted by $b_t$, which belongs to  the boundary of the big disc. We will denote deformed rays (they are curved paths now) by $p_{\theta_t+\pi}^{(i)}$.

\vspace{2mm}

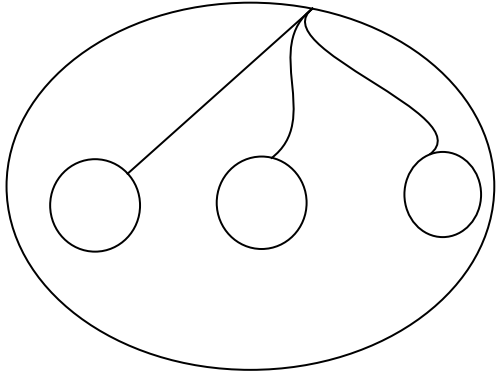
\tikzset{every picture/.style={line width=0.75pt}} 

\begin{tikzpicture}[x=0.5pt,y=0.5pt,yscale=-1,xscale=1]

\draw   (87,143.5) .. controls (87,66.73) and (169.6,4.5) .. (271.5,4.5) .. controls (373.4,4.5) and (456,66.73) .. (456,143.5) .. controls (456,220.27) and (373.4,282.5) .. (271.5,282.5) .. controls (169.6,282.5) and (87,220.27) .. (87,143.5) -- cycle ;
\draw   (120,158) .. controls (120,138.67) and (135.22,123) .. (154,123) .. controls (172.78,123) and (188,138.67) .. (188,158) .. controls (188,177.33) and (172.78,193) .. (154,193) .. controls (135.22,193) and (120,177.33) .. (120,158) -- cycle ;
\draw   (246,156) .. controls (246,136.67) and (261.22,121) .. (280,121) .. controls (298.78,121) and (314,136.67) .. (314,156) .. controls (314,175.33) and (298.78,191) .. (280,191) .. controls (261.22,191) and (246,175.33) .. (246,156) -- cycle ;
\draw   (388,149.75) .. controls (388,131.94) and (400.98,117.5) .. (417,117.5) .. controls (433.02,117.5) and (446,131.94) .. (446,149.75) .. controls (446,167.56) and (433.02,182) .. (417,182) .. controls (400.98,182) and (388,167.56) .. (388,149.75) -- cycle ;
\draw    (319,8.5) -- (178,134.5) ;
\draw    (287,122.5) .. controls (327,92.5) and (279,38.5) .. (319,8.5) ;
\draw    (407,119.5) .. controls (447,89.5) and (279,38.5) .. (319,8.5) ;

\end{tikzpicture}

\vspace{2mm}

The complement of the big disc to the union of all  paths $p_{\theta_t+\pi}^{(i)}$ and all discs $D(z_i/t,\varepsilon)$ can be homotopically retracted to $\cup_i(D(z_i/t,\varepsilon)\cup p_{\theta_t}^{(i)})$. This gives rise to an isomorphism which is called 
{\it Betti global-to-local isomorphism}. We summarize this discussion in the form of the following statement.

\begin{prp}\label{Betti local to global isomorphism}  Consider the open subset of $\C^\ast_t$ consisting of such $t\in \C^\ast$ which do not belong to the union of rays  $Arg(t)=Arg(z_i-z_j)\in \R/2\pi \Z$ for all  $i\ne j$. Then the above construction gives rise to a  canonical isomorphism of local systems restricted to this open subset:
$$iso_{Betti}: \mathcal{H}^\bullet_{Betti,glob}(X,f)\simeq \mathcal{H}^\bullet_{Betti, loc}(X,f).$$
\end{prp}
The restriction of $iso_{Betti}$ to the ray  $Arg(t)=\theta$ will be denoted by $iso_\theta$.

Notice that local Betti cohomology can be also computed via the sheaf of vanishing cycles:
$$\oplus_{z\in S}H^{\bullet}((f/t)^{-1}(z/t),\phi_{{f-z}\over{t}}(\underline{\Z}_X)$$
$$\simeq\oplus_{z\in S}H^{\bullet}(U_{\varepsilon}((f/t)^{-1}(z/t)),(f/t)^{-1}(z/t+\varepsilon\cdot e^{i\theta_t}, \Z)),$$
where $S$ is the set of critical values of $f$, $\theta_t=Arg(t)$, and $U_\epsilon(\bullet)$ denote sufficiently small open neighborhood of the set $\bullet$, and $\phi_g(F)$ denote the functor of vanishing cycles of the function $g$ applied to the local system $F$. The RHS can be identified with  $H_{Betti, glob, t}^\bullet(X,f)$ by the Proposition \ref{Betti local to global isomorphism}.

Let $\theta_{ij}:=Arg(z_i-z_j)\in \R/2\pi \Z$, where $ i\ne j$. 
\begin{defn}\label{Stokes rays}
We call the ray $\{t|Arg(t)=\theta_{ij}\}=\R_{>0}\cdot (z_i-z_j)\subset \C_t$ a {\it Stokes ray}.  

\end{defn}

 We will denote the Stokes ray with the slope $\theta$ by $s_\theta$.
Notice that points of $S$ which belong to a Stokes ray have natural order (by the distance to the vertex).
All rays with the vertex at $0\in \C$ which are not Stokes will be called {\it generic}.

For every Stokes ray $s_{\theta}$ we have a {\it Stokes  isomorphism} $T_{\theta}$ of the graded abelian group $H^\bullet_{Betti,loc, t}(X,f)$ for $Arg(t)$ sufficiently close to $\theta$. 
More precisely we can 
choose two generic rays $l_{\theta^{\pm}}$
with slightly bigger (resp. smaller) argument than 
$\theta$.
Then the desired isomorphism is $iso_{\theta^-}^{-1}\circ iso_{\theta^+}$. This isomorphism has the form $T_\theta=id+\sum_{i\ne j, Arg(z_i-z_j)=\theta}T_{ij}$, where 
$$T_{ij}: H^{\bullet}(D(z_i,\varepsilon),b_{\theta_{ij}+\pi},\Z)\to  H^{\bullet}(D(z_j,\varepsilon),b_{\theta_{ij}+\pi},\Z)$$
are homomorphisms of graded abelian groups for any $i\ne j$ and sufficiently small $\varepsilon$.

Restricting to the standard circle of directions $S^1_\theta$ we obtain a local system  of Betti cohomology over disjoint union of 
Stokes sectors (i.e. sectors bounded by consecutive Stokes rays).
A  fiber of this local system, admits a direct sum decomposition outside of the intersections of $S^1_\theta$ with the union of Stokes rays. Furthermore, it
is endowed with the lower-triangular integer isomorphism $T_\theta$  for 
each Stokes ray $s_\theta$.
 
We will call the above local system on $S^1_\theta$ endowed  with Stokes isomorphisms the {\it Betti data} associated with the pair $(X,f)$.
 This Betti data is a simple example of an  {\it analytic wall-crossing  structure} from [KoSo12].



\subsection{Equivalent descriptions of the Betti data}\label{perverse sheaf}

For each $k\in \Z_{\ge 0}$ and sufficiently small $\varepsilon >0$ we define a family of abelian groups $(\FF_z)_{z\in \C}$ by the formula:
$$\FF_z=H^k(X, D_0\cup f^{-1}(|z^{\prime}-z|\le \varepsilon),\Z).$$
It is easy to see that $\FF_z$ are stalks of a constructible sheaf $\FF$ of abelian groups on $\C$ which is a local system  on $\C-Bif$.

\begin{prp} The sheaf $\FF$ is perverse and $\R\Gamma(\C,\FF)=0$ i.e.
$$H^0(\C,\FF)\simeq H^1(\C,\FF)=0.$$

\end{prp}

{\it Proof.} Proof is completely analogous to the one in  [KaKoPa1] (see e.g. Remark 3.9 (b) and the arguments preceding it). $\blacksquare$

Let us summarize several equivalent descriptions of the Betti data, similarly to those we gave in the de Rham case. Taking into account the above discussion most of the proofs either immediate or can be derived from  [FreJos] (see Chapter 2) or  [KaKoPa1] (see Section 2).

\begin{prp}\label{descriptions of Betti data}

The following six descriptions of the Betti data of the pair $(X,f)$ are equivalent:

a) For each critical value $z_i$ a local system $E_i$ of abelian groups on $S^1_\theta$ and for any Stokes ray $\theta=\theta_{ij}=Arg(z_i-z_j), i\ne j$ a morphism of fibers $T_{\theta}: E_{i,\theta}\to E_{j,\theta}$ such that $T_{ii}$ are invertible.

b) Constructible sheaf $\FF$ on $\C$ with the set of singularities $\{z_i\}$ satisfying the property $\R\Gamma(\C,\FF)=0$.

c) Perverse sheaf $\mathcal{G}^\bullet$ on $\C$ with the set of singularities $\{z_i\}$ satisfying the property $\R\Gamma(\C,\mathcal{G}^\bullet)=0$

d) Stokes sheaf (see e.g.  [Sab4]) associated with a meromorphic connection on $\C^\ast_t$ with singular terms $\{e^{(z-z_i)^{-1}/t}\}$ (see [KoSo7] about the notion ``singular term"). The corresponding filtration is  called Deligne-Malgrange-Stokes filtration in [KaKoPa1], Section 2.3.2.

e) Constructible sheaf on $\C_t$ with the singular support which belongs to the union of $\C_t$ without small disc about $t=0$  and the positive conormal bundles to the small circles about points $z_i$ associated with singular terms $e^{z_i/t}$ (see [KaKoPa1], [KoSo7], and more details and generalizations will be given in [KoSo9]).

f) Collection of vector spaces $V(B,b)$ (see Section \ref{Betti cohomology}) such that $z_i\notin \partial B$ subject to the retraction-direct sum property from Section \ref{Betti cohomology}.

\end{prp}

In particular  Betti data form a symmetric monoidal category for which the tensor product amounts to the addition on the set of critical values and to the convolution at the level of constructible sheaves.

 \subsection{Comparison isomorphisms}\label{local de Rham to Betti}
 
  Let us assume for simplicity that $D_0=\emptyset$ and drop it from the notation. The results hold in more general case.
 
 The integration over cycles defines a non-degenerate pairing
$$H_{\bullet}^{Betti,glob}(X,f)\otimes H^{\bullet}_{DR,glob}(X,f)\to \C.$$
This pairing should be thought of as {\it exponential period map}. Consequently exponential integrals can be interpreted as exponential periods
of the cohomology class in $H^\bullet_{DR,glob}(X,f)$ of the volume form $vol_X$. Adding the parameter $t$ to the story we obtain a holomorphic section $t\mapsto [vol_X(t)]$ of the bundle $\mathcal{H}_{DR,glob}(X,f)$.

If the function $f$ is equal to zero, then the exponential period becomes the usual period of the volume form.
It follows from definitions that the following {\it comparison isomorphism between global Betti and de Rham cohomology} holds.

\begin{prp}\label{global comparison iso}
If the pair $(X,f)$ is defined over the field of complex numbers and $t\in \C^\ast$ then we have a comparison isomorphism $iso_t: H_{DR,glob,t}^{\bullet}(X,f)\simeq H_{Betti,glob,t}^{\bullet}(X,f)\otimes \C$.
\end{prp}
 
 Assuming that $f$ is a proper map, recall a vector bundle $\mathcal{H}^{\bullet}_{DR,glob}(X,f)$ endowed with a meromorphic connection $\nabla$ (see Proposition \ref{bundle with meromorphic connection}).
The bundle carries a covariantly constant lattice $\Gamma:=(\Gamma_t)_{t\in \C^{\ast}}:=iso^{-1}_t(H^\bullet_{Betti,glob,t}(X, f)$. 
 One can prove (see [KaKoPa1]) that the meromorphic connection described in Proposition \ref{bundle with meromorphic connection} satisfies the following property:
 {\it Stokes filtration at $t=0$ is compatible with the lattice $\Gamma$}.





Finally one has the following  local version of the Proposition \ref{global comparison iso}. First, let us recall that there exists a functor ${RH}_{loc}^{-1}$ from the category of local systems of $\C$-vector spaces to the category of regular singular connections of $\C((t))$-vector spaces. It is obtained as the composition of the inverse $RH^{-1}$ to the Riemann-Hilbert functor $RH$ and the functor of taking the formal completion at $t=0$ of the corresponding $D$-module on $\C^\ast$. Then the local version can be formulated such as follows.

\begin{prp} 
We have the following isomorphism of vector bundles over the punctured formal disc, which are endowed with regular singular connections:

$${RH}^{-1}_{loc}\!(\mathcal{H}^\bullet_{Betti,loc}(X,f)\otimes \C)\simeq \mathcal{H}^\bullet_{DR,loc}(X,f).$$

Same is true for a summand of the LHS and RHS corresponding to each critical value $z_i\in S$.

\end{prp}
{\it Proof.} Follows from [KaKoPa1], Lemma 3.11. $\blacksquare$

\subsection{Wall-crossing structure for exponential integrals}\label{WCF}

Given a pair $(X,f)$ consisting of a smooth complex affine algebraic variety $X, dim_\C X=n<\infty$ and a regular function $f$ on it let us choose an algebraic volume form $vol$ on $X$. Recall the exponential integral $I_C(t)=\int_C e^{f/t}vol$ which is an analytic function in $t$, where $t$ belongs to a sector depending on $C$.  


More precisely, if we do not fix the integration cycle $C$ but keep the volume form fixed we can interpret the exponential integral $I(t)$ as a morphism of sheaves of abelian groups 
$$\H^{Betti,glob}_\bullet(X,f)\to {\OO}_{\C^\ast_t}^{an}$$
given by $\gamma\mapsto \int_\gamma e^{f/t}vol$. 
Let us assume that $t_0\in \C^\ast$ does not belong to a Stokes ray. Then there is a canonical isomorphism of Betti global and local homology $H^{Betti,loc,t}_\bullet(X,f)\simeq H^{Betti,glob,t}_\bullet(X,f)$ for any $t$ which belongs to a sufficiently small sector $V$ containing the ray  $\R_{>0}\cdot t_0$. Hence we have the corresponding morphism of sheaves of abelian groups $\oplus_i \H^{Betti,loc,z_i}_\bullet(f^{-1}(V),f)\to {\OO}_{\C^\ast_t}^{an}(V)$ in the obvious notation. 

Let us now consider a Stokes ray with the slope $\theta$ and choose a basis of local sections $(\gamma^{(i),\pm}_a)$ of the direct sum of local Betti homology for $Arg(t)=\theta\pm \epsilon$ for all sufficiently small $\epsilon>0$. In other words we take the restriction of the local system of Betti local homology to the small sectors on the left and on the right of the Stokes ray. Then we have two different vectors  of analytic functions $(\int_{\gamma^{(i),-}_a}e^{f/t}vol)_{i,a}$ and 
$(\int_{\gamma^{(i),+}_a}e^{f/t}vol)_{i,a}$. Formulas which relate them are called {\it wall-crossing formulas}.  
More precisely, let us fix a critical value $z_i$ as well as a section $\gamma\in H^{Betti,loc,z_i,t}_\bullet(X,F)$ of the local system of Betti local homology trivialized in a small neighborhood of the Stokes ray $s_\theta$. Then the function $t\mapsto \int_\gamma e^{f/t}vol$ considered in a small sector $0<\theta-Arg(t)<\epsilon$ admits analytic continuation to the small sector $0<Arg(t)-\theta<\epsilon$ (i.e. it is the analytic continuation from the sector on the right of $s_\theta$ to the sector on the left of $s_\theta$, and this analytic continuation is equal at the ray with the slope $\theta_+>\theta$ to 
$$\int_{iso_{\theta_+}^\ast(\gamma)}e^{f/t}vol+\sum_{j\ne i, Arg(z_i-z_j)=\theta}\int_{(iso_{\theta_+}^\ast\circ T_{ji})\gamma}e^{f/t}vol,$$
where $iso_{\theta}^\ast$ denote the dual to the Betti global-to-local isomorphism.

Let us now illustrate the above general considerations in a special case  in which
the bases in local and global Betti homology have explicit geometric description. This is the case when $f$ is a Morse function with different critical values.
Let us choose a Hermitian metric on the compactification $\overline{X}$ and $\theta\in \R/2\pi \Z$.
Then there is a nice basis in the global Betti homology associated with the direction $\theta+\pi$  consisting of {\it Lefschetz thimbles} $th_{i,\theta+\pi}$. Let us recall this notion.

Fix a critical point $x_i\in X$. Then the thimble $th_{i,\theta+\pi}$ is by definition the union of gradient lines  of the function $Re(e^{-i\theta}f)$ emerging from the critical point $x_i$. \footnote{If our $\overline{X}$ was K\"ahler, we could choose the metric to be the K\"ahler metric and the same gradient line is also an integral curve for the Hamiltonian function $Im(e^{-i\theta}f)$  with respect to the {\it symplectic} structure.} Hence $f(th_{i,\theta+\pi})$ is a ray $Arg(t)=\theta+\pi$ emerging from the critical value $z_i=f(x_i)$. 
Choice of the Hermitian metric on the compactification of $X$ ensures that the gradient lines cannot go to infinity in finite time, hence each local thimble gives rise to a global one. 

We remark that  the orientation of the thimble $th_{i,\theta+\pi}$  is not defined canonically, hence its homology class is defined up to a sign only. Moreover for odd $n$ it is not possible to choose the orientation in a covariantly constant way in $\theta$.

Consider  the following collection of  exponential integrals, where $t$ does not belong to the union of Stokes rays and $\theta=Arg(t)$:
$$I_i(t)=\int_{th_{i,\theta+\pi}}e^{f/t}vol.$$

Assume that the set $S=\{z_1,...,z_k\}$ of critical values of $f$ is in generic position in the sense that no straight line contains at least three points from $S$.  Then a Stokes ray  with the vertex at $0\in \C$ contains exactly two different critical values  which can be  ordered by their proximity to the vertex.

It  is easy to see that if we cross a Stokes ray
$s_{ij}:=s_{\theta_{ij}}$ containing critical values $z_i,z_j, i<j$, then the exponential integral $I_i(t)$ changes such as follows:
$$I_i(t)\mapsto I_i(t)+n_{ij}I_j(t),$$
where $n_{ij}\in \Z$ is the number of gradient trajectories of the function 
$$Re(e^{i(Arg(z_i-z_j))}f)$$
 joining critical values $x_i$ and $x_j$ (equivalently the intersection index of the opposite thimbles outcoming from the critical points $x_i, x_j$). These are the wall-crossing formulas for the pair $(X,f)$ satisfying the above-mentioned assumptions. \footnote{Physicists call the above wall-crossing formulas ``Cecotti-Vafa wall-crossing formulas" or ``$2d$ wall-crossing formulas". Mathematicians  simply call them ``Picard-Lefschetz formulas".}
 
They can be encoded in a different way, which is useful for studying the resurgence of perturbative expansions of exponential integrals.
In order to explain that let us modify our exponential integrals such as follows:
$$I_i^{mod}(t):=\left({1\over{2\pi t}}\right)^{n/2}e^{-z_i/t}I_i(t).$$

As $t\to 0$ the stationary phase expansion ensures that the modified integral has a formal power series expansion for which we keep the same notation:
$$I_i^{mod}(t)=c_{i,0}+c_{i,1}t+....\in \C[[t]],$$
where $c_{i,0}\ne 0$.
The jump of the modified exponential integral across the Stokes ray $s_{ij}$ is given by $\Delta(I_i^{mod}(t))=n_{ij}I_j^{mod}(t)e^{-(z_i-z_j)/t}$. Therefore the vector $\overline{I}^{mod}(t)=(I_1^{mod}(t),...,I_k^{mod}(t)), k=|S|$ satisfies the Riemann-Hilbert problem on $\C$ with known jumps across the Stokes rays and known asymptotic expansion as $t\to 0$ (notice that because of our ordering of the points in $S$, the function $e^{-{(z_i-z_j)/t}}$ has trivial Taylor expansion as $t\to 0$ along the Stokes ray $s_{ij}$).

In abstract terms, we consider a Riemann-Hilbert problem for a sequence of $\C^k$-valued functions (here $k$ is the rank of the Betti cohomology, which is under our assumptions is equal to the cardinality $|S|=k$) $\Psi_1(t),...,\Psi_k(t)$ on $\C^{\ast}-\cup (Stokes\,rays)$ each of which has a formal power asymptotic expansion in $\C[[t]]$ as $t\to 0$, and which satisfy the following jumping conditions along the Stokes rays $s_{ij}$:

$$ \Psi_j\mapsto \Psi_j,$$
$$\Psi_i\mapsto \Psi_i+n_{ij}e^{-{{{z_i-z_j}\over{t}}}}\Psi_j.$$

This collection $(\Psi_i)_{1\le i\le k}$ gives rise to a holomorphic vector bundle. Each $\Psi_i(t)$ is a vector $(\Psi_{ij}(t))_{1\le j\le k}$.

The above  Riemann-Hilbert problem and  the resurgence properties of Taylor expansions of its solutions   are discussed in [KoSo12]  from the point of view of introduced in the loc.cit. notion of analytic wall-crossing structure. 

\begin{rmk}\label{resurgence and volume of fiber}
One can rewrite the exponential integral over the thimble such as follows:
$$\int_{th_{i,\theta+\pi}}e^{f/t}vol=\int_{l_{\theta+\pi}}e^{s/t}vol_f(s)ds,$$
where $vol_f(s)$ is the volume of the nearby cycle in the fiber $f^{-1}(s)$ defined by the thimble with respect to the Gelfand-Leray form ${vol\over{df}}$,  and $l_{\theta+\pi}$ is the admissible ray of the slope $\theta+\pi$. Changing variables $t\mapsto 1/\lambda$ we can think of the exponential integral over a thimble as of the  Laplace transform  of the function $vol_f(s)$. Hence resurgence properties of the exponential integral can be deduced from the analytic behavior of this function.

More precisely, the shifted Borel transform $\mathcal{B}(I_i^{mod}(t))(s)=\sum_{k\ge 0}c_{i,k}{(s-z_i)^k\over {k!}}$ is equal up to a universal constant to $({d\over {ds}})^{n/2}(\int_{\gamma(s)}vol)$ provided $n$ is even. Here $\gamma(s)\in H_n(X, f^{-1}(s), \Z)$ for $s$ sufficiently close to $z_i$ is the relative homology class defined the thimble $\theta+\pi, \theta=Arg(s-z_i)$. For $n$ odd we understand the fractional derivative using the standard integral representation. 

\end{rmk}

\begin{rmk}\label{relation to Hodge structure}
This remark is in a sense a continuation of the previous one.
Alternatively the same Borel transform is equal to $({d\over {ds}})^{n/2-1}(\int_{\partial \gamma(s)}{vol\over{df}})$, 
where $\partial \gamma(s)\in H_{n-1}(f^{-1}(s),\Z)$ is the boundary of $\gamma(s)$.
In case if the fibers $f^{-1}(s)$ endowed with the volume forms $vol/df$ are Zariski open in compact Calabi-Yau varieties $\overline{f^{-1}(s)}$ one can consider the pairing of $\partial \gamma(s)$ with the lowest term 
$F_{n-1}H^{n-1}_{DR}(\overline{f^{-1}(s)})$, where $F_\bullet$ denote the Hodge filtration. Then by the Griffiths transversality the above $n/2-1$ derivative gives us the pairing of the same cycle with the class in the middle term 
$F_{n/2}H^{n-1}_{DR}(\overline{f^{-1}(s)}$ of the Hodge filtration. 
\end{rmk}

\begin{rmk}\label{volume form and global Betti-to-de Rham}
  
 
 Assume that $f$ is Morse with different critical values. Let us fix a top degree holomorphic form $vol_X\in \Omega^{n,0}(X)$.  Since $(d+(df/t)\wedge\bullet)(vol_X)=0$ the top degree form defines the class $[vol_X]\in H_{DR}^n(X, f/t)$.
We denote by $J_l(t), 1\le l\le k$ the number $e^{-z_l/t}\langle[vol_X], th_{l, t}\rangle$, where $k$ is the cardinality of the set of critical points and $z_l$ is the critical value of $f$ corresponding to the critical point $x_l$. Then under the global isomorphism of de Rham and Betti cohomology the cohomology class $[vol_X]$ is mapped to the Betti cohomology class $iso_t([vol_X])$ such that its pairing with the relative homology class of $th_{l,t}$ is equal to $J_l(t)$.

 \end{rmk}

Wall-crossing formulas for the exponential integral are encoded in the corresponding {\it wall-crossing structure} (WCS for short) on $\C^\ast_t$ or on the circle of directions $S^1_\theta$, as described in [KoSo12], Section 6.2. The reader can find in the loc.cit. or in the [KoSo7]  the precise definition and some results concerning WCS. In a few words the WCS  is a local system on a topological space (which is $\C^\ast$ or $S^1$ in our case) of more basic structures called {\it stability data on a graded Lie algebra}. In the case of exponential integral of a Morse function one considers stability data on the Lie algebra $\mathfrak{gl}(k,\C)$, where $k$ is the number of critical values of $f$. This is a graded Lie algebra with respect to its root system.

In order to define the WCS  in the Morse case the only remaining piece of information is the collection of integers $n_{ij}$ which appeared in the wall-crossing formulas. In  the non-Morse case one should use their generalizations, which are the morphisms $T_{ij}$. The morphisms $T_{ij}$ are derived from  Betti global-to-local isomorphisms. Hence the WCS for exponential integrals  is also determined by those isomorphisms.

It is also mentioned in [KoSo12] that similar results should hold in Quantum Field Theory, with exponential integrals being replaced by (renormalized) Feynman integrals.





\subsection{Wall-crossing structure in the case of dependence on parameters}\label{dependence on parameters}

In this subsection we discuss a generalization in the case when the function $f$ in the exponential integral depends on parameters.

Let $\pi: {\mathcal X}\to U$ be a smooth fiber bundle, such that each fiber and the total space
are complex algebraic varieties. Suppose that $f: {\mathcal X}\to \C$ is a regular
function such that for any compact subset $K\subset U$ the restriction $f_{|\pi^{-1}(K)}$
is a proper map. In other words, we consider a family of proper regular functions
$f_u: X_u:=\pi^{-1}(u)\to \C$. We denote by $S_u\subset \C$ the set of critical values of $f_u$ and make the following

{\bf Continuity Assumption:}

{\it The map $u\mapsto max\{|z|\in \R_{\ge 0}|z\in S_u\}$ is locally bounded.}

Since each $S_u$ is a finite set the assumption means that the above map is continuous
(although the number of critical values of $f_u$ can change when we vary $u\in U$).

Let $S_u:=(z_i(u))_{1\le i\le m}$. We denote by $U^0\subset U$ an open domain where
the number of critical values is constant (and hence automatically maximal), so  $u\mapsto z_i(u)$ is a continuous
function for each $i\in I$. Let ${\mathcal S}\to U^0$ be
the corresponding family of sets with the fiber $S_u$ over $u\in U^0$.

The Betti cohomology $H^{\bullet}_{Betti,glob,t}:=H^{\bullet}_{Betti,glob,t}(X_u,f_u,\Z)$ gives rise to a local system of abelian groups on
$\C^{\ast}_t\times {\mathcal S}$.

The morphisms $T_{ij}$  can jump across the {\it  walls of first kind} (see [KoSo1] ). The latter are union of real
subvarieties in $U^0$ for which at least three points $z_i=z_i(u)$ align. Let us denote by $T_{ab}^{\pm}$ the corresponding morphisms
on two sides of the wall and assume that only three points $z_i,z_j,z_k$ are aligned on the wall in such a way that $z_j$ lies between $z_i$ and $z_k$ we have the following wall-crossing formulas (see [KoSo1], Section 2.9 for the details):
$$T_{ij}^+=T_{ij}^-, T_{jk}^+=T_{jk}^-, T_{ik}^+=T_{ik}^-+T_{ij}^-T_{jk}^-.$$

It follows that given a continuous path $\phi:[0,1]\to U^0$ and using the above wall-crossing formula one can recalculate the Betti data at $\phi(1)$ from the Betti data at $\phi(0)$.
An appropriate framework for the exponential integrals depending on parameters is the one of the {\it wall-crossing structures} and {\it analytic wall-crossing structures} (see [KoSo7], [KoSo12]).

Consider for any
$1\le i\ne j\le n$, the following variety
$${\mathcal W}_{ij}=\{(t,u)\in \C^{\ast}\times U^0| Im(t^{-1}(z_i(u)-z_j(u)))=0\}.$$

The varieties ${\mathcal W}_{ij}$ are  {\it walls of second kind} in the sense of [KoSo1].

In order to define the desired WCS we  specify first a local system of lattices ${\Gamma}\to \C^{\ast}\times U^0$, a local system of central charges, and a local system
${\g}\to \C^{\ast}\times U^0$ of $\underline{\Gamma}$-graded Lie algebras (say, over ${\bf Q}$).

In our case the fiber ${\Gamma}_{t,u}$ does not depend on $t\in \C^{\ast}$ and coincide with the lattice
$\Gamma_u=Ker (\Z^{S_u}\to \Z)$, where the map is given by $(z_1(u),...,z_m(u))\mapsto \sum_{1\le i\le m}z_i(u)$. Thus $\Gamma_u$ is
isomorphic to the root lattice $A_{m-1}$.

The fiber of the local system of graded Lie algebras is
$$\g_{t,u}=End(\oplus_{1\le i\le m}H^{\bullet}_{Betti,glob,z_i(u),t})\otimes {\bf Q}$$
endowed with the natural structure of Lie algebra.

The local system of central charges is a local system of homomorphisms from the local system of lattices to $\C$.
In our case it is given by
$$Z_u:\Gamma_{t,u}\to \C, e_{ij}\mapsto t^{-1}(z_j(u)-z_i(u)),$$
where $e_{ij}, i<j$ is the standard basis of the Lie algebra $\mathfrak{sl}(m)$.

\begin{rmk} \label{grading and walls}

a)The grading combined with the central charge gives rise to a semisimple derivation of the Lie algebra $\g_{t,u}$.
It would be interesting to generalize this observation to the non-semisimple case.

b) The above-mentioned walls of second kind arise in this framework as the set $\{(t,u)\in \C^{\ast}\times U^0| Z_{t,u}(\gamma)\in \R_{>0}\}$ for some $\gamma\ne 0\in \Gamma_u$ such that $\g_{t,u,\gamma}\ne 0$ for the corresponding graded component.

\end{rmk}

Finally, in order to define WCS we need to specify a collection of elements $a_{t,u,\gamma}\in \g_{t,u,\gamma}$ attached to walls
of second kind, such that
for any small loop in $\C^{\ast}\times U^0$ the following ``triviality of monodromy'' condition (wall-crossing formula) is satisfied:
$$\prod^{\to}exp(a_{t,u,\gamma})=id,$$
where the product is taken in the clockwise order over all (possibly countably many) intersections of the loop with the walls (see [KoSo7]).
The elements $a_{t,u,\gamma}$ are  derived from  the morphisms $T_{ij}$. Summarizing, we obtain a WCS on $\C^{\ast}\times U^0$.
We expect that after a modification of the local system of Betti cohomology along
the walls of second kind we can  extend the WCS to $\C^{\ast}\times U$. 


\begin{exa}\label{focus-focus singularity}

Let
$U=\C, {\mathcal X}=\C^{n+1}$, and $f(x_1,...,x_n,u)={x_1^3\over{3}}-ux_1+\sum_{i\ge 2}x_i^2$.
Then for each $u$ the set $S_u$ of critical points of $f_u$ is $(\pm\sqrt{u},0,0,...,0)$. They are all Morse
critical points. There are two critical values
equal to $\pm 2/3 u^{3/2}$.

Let us fix for simplicity the parameter $t=1$. Then the walls of second kind are derived from the equation $Im(u^{3/2})=0$ or, equivalently,
$u^3\in \R_{\ge 0}$. This condition gives rise to three rays emerging from each critical value and having consecutive angles $2\pi/3$.
One can check that the intersection indices of the corresponding thimbles are given by $n_{ij}=1$. The triviality of monodromy
condition amounts to the identity for $2\times 2$ matrices $(XY)^3=id$, where $X(e_1)=-e_2, X(e_2)=-e_1$ corresponds to exchange of two branches
and take care about orientation and $Y(e_1)=e_1, Y(e_2)=e_1+e_2$ corresponds to the monodromy of the focus-focus singularity
of the $\Z$-affine structure on $\R^2-\{(0,0)\}$. Having this ``initial data'' (see [KoSo7]) we can construct the WCS on
$\C^{\ast}\times \C$. 
\end{exa}

\begin{rmk} \label{root A local system}
In fact the above Proposition (as well as the expectation about the extension of the local system) are not specific for
the WCS arising from the family of holomorphic functions with Morse critical points. It is true for any WCS associated with
the local system of $A_{m-1}$ root lattices on $U^0$ and the graded Lie algebra $\mathfrak{gl}(m)$ as long as the wall-crossing formulas
are the Cecotti-Vafa wall-crossing formulas. This remark can be used for a construction of holonomic $D$-module on
a curve which corresponds to a given spectral curve.

\end{rmk}

\section{Betti and de Rham cohomology of a variety with  closed $1$-form}\label{case of 1-forms}

\subsection{Compactification of a variety with closed $1$-form}\label{compactification in case of 1-forms}

In the previous section to a triple $(X,D_0,f)$ we associated families over $\C_t^\ast$ of the (twisted) de Rham cohomology $H_{DR,t}^\bullet(X,D_0,f)$ and the Betti cohomology $H_{Betti,t}^\bullet(X, D_0,f)$.  Exponential integrals were formulated in terms of the pairing between the de Rham and Betti cocycles.
In this section we are going to generalize these results from the case of exact algebraic (or holomorphic) $1$-forms $\alpha=df$ to the case of arbitrary algebraic (or holomorphic) closed $1$-forms.
 
Thus, let $X$ be a complex smooth algebraic variety, $\alpha$ a closed algebraic $1$-form on $X$, and $D_0\subset X$ is a normal crossing divisor. The following result holds.

\begin{prp}\label{compactification for 1-forms}There is a smooth complex projective algebraic variety $\overline{X}\supset X$ containing  a normal crossing divisor $D:=\overline{D}_0\cup D_h\cup D_v\cup D_{log}$ such that:

a) among all the divisors $\overline{D}_0, D_h, D_v, D_{log}$ only $D_v$ and $D_{log}$ can have common irreducible components;

b) $X=\overline{X}-( D_h\cup D_v\cup D_{log})$,

c) divisor $D_0$ is the intersection $\overline{D}_0\cap X$,

d) for any point  $x\in \overline{X}$ there exists an analytic neighborhood $U$ and meromorphic in $U$ closed $1$-forms $\alpha_{reg},\alpha_{log}, \alpha_\infty$ such that the $1$-form $\alpha$ is represented in $U$ as a sum:

$$\alpha=\alpha_{reg}+\alpha_{log}+\alpha_{\infty},$$
and  furthermore the summands satisfy the following properties:

1) the form $\alpha_{reg}$ is regular on $U$,

2) the form $\alpha_{log}$ can be expressed  in the local coordinates  such as follows $\alpha_{log}=\sum_i c_id\,log\,z_i$, where $c_i\in \C-\{0\}$,
and $\prod_i z_i=0$ is a local equation for $D_{log}$,

3) the form $\alpha_{\infty}$ can be locally near $D_v\cap U$ written as $\alpha_\infty=df$, where
$f$ is an analytic function in a neighborhood of $D_v$, which has the form
$$f={c\over{\prod_j z_j^{k_j}}}(1+o(1)),k_j\ge 1$$
where $c\in \C-\{0\}$ and $z_j$ are local coordinates in $U$ in which the divisor $D_v$ is written as
$\prod_jz_j=0$. This condition is equivalent to the condition that $f={c(z)\over{\prod_j z_j^{k_j}}}$, where $c(z)$ is invertible holomorphic function on $U$.

\end{prp}

{\it Proof.} First let us choose an arbitrary simple normal crossing compactification $\overline{X}_{norm}$ of $X$ and denote by $D_{norm}$ the corresponding simple normal crossing divisor. Thus $D_{norm}=\overline{X}_{norm}-X$. Let $\overline{D}_0^0$ denote the closure of $D_0$ in $\overline{X}_{norm}$ and $D_{h}^0$ denote the union of those components of $D_{norm}$ at which $\alpha$ has finite limit. For each other irreducible component $D_{norm,j}$ of $D_{norm}$ we define the period of $\alpha$ as $c_j={1\over{2\pi i}}\int_{S^1_j}\alpha$, where $S^1_j$ is a small circle centered at a smooth point of $D_{j}$. We denote by $D_{log}^0$ the union of those irreducible components for which $c_j\ne 0$. Then in a small analytic neighborhood of a point of $D_{norm}$ we can write $\alpha-\sum_jc_j dlog(z_j)=dF$, where $F$ is a meromorphic function in the neighborhood having a pole of order at least $1$ on $D_{norm}-\overline{D}_0^0-D_h^0$, and $z_j$ are coordinates near $D_{log}^0$ in which $D_{log}^0$ can be written as  $\prod_j z_j=0$. We remark that the function $F$ is not defined canonically since it depends on the choice of coordinates $\{z_j\}$ as well as on the overall constant.
Nevertheless $F$ defines  a section of the sheaf $\OO(\ast(D_{norm}))/\OO_{\overline{X}_{norm}}$, where $\ast(divisor)$ means poles of an arbitrary  order at the divisor.

Thus $F$ admits locally a  non-canonical presentation $F=F_{sing}+F_{reg}$, where $F_{reg}$ is holomorphic. Analogously we  define the divisor $D_{v}^0$ by the condition that $\alpha$ has poles of order bigger than $1$ at $D_{v}^0$. Then we see that $D_{v}^0$ and $D_{log}^0$ can intersect (this is true even in the case when $\overline{X}$ is a curve).

If $F$ is well-defined globally on $X$ one can consider a possibly singular algebraic subset $\overline{graph(F)}\subset \overline{X}_{norm}\times {\bf P}^1$, the closure of the graph of $F$
 considered as an algebraic subset of $X\times \C$. 
 
 Otherwise we proceed such as follows.
 From the short exact sequence of sheaves
$$0\to \OO_{\overline{X}_{norm}}\to \OO(\ast(D_{norm}))\to  \OO(\ast(D_{norm}))/\OO_{\overline{X}_{norm}}\to 0$$
we derive the morphism of cohomology groups $H^0(\overline{X}_{norm},  \OO(\ast(D_{norm}))\to H^1(\overline{X}_{norm}, \OO_{\overline{X}_{norm}})$.
Since $F$ always defines a class in $H^1(\overline{X}_{norm}, \OO_{\overline{X}_{norm}})$ we can use the latter in order to ``twist" the trivial ${\bf P}^1$-bundle over $\overline{X}_{norm}$ and obtain a possibly non-trivial ${\bf P}^1$-bundle $Y$ over the same base, which in addition is trivialized over $X$. For that we use the isomorphism $H^1(\overline{X}_{norm},\OO_{\overline{X}_{norm}})\simeq H^1(\overline{X}_{norm}, {\bf G}_a)$ and the inclusion ${\bf G}_a\to PGL(2)=Aut({\bf P}^1)$ as the group of shifts $\{x\mapsto x+const\}$. The trivialization $Y_{|X}\simeq X\times {\bf P}^1$ gives rise to a divisor $X\times\{0\}\subset Y$. We denote its closure by $\overline{graph(F)}$. This divisor coincides with the above-discussed divisor $\overline{graph(F)}$ in the case when $F$ is defined globally. The divisor $\overline{graph(F)}$ is irreducible, possibly singular, and it
projects one-to-one onto an open part of $\overline{X}_{norm}$ outside of the intersection of $\overline{graph}$ with the divisor $D_\infty=\{\infty_{{\bf P}^1}\}\times \overline{X}_{norm}$.

Let us introduce the set $Z:=\partial_\infty(\overline{graph}):=\overline{graph}-graph$. This is a closed subset of $Y$ of $codim\ge 2$.

Notice that  $(\overline{graph}\, \cup D_\infty)-\partial_\infty(\overline{graph})$ is a closed smooth  divisor in  $Y-Z$ containing a closed normal crossing divisor isomorphic to $D_0\subset X\simeq graph$. By the embedded resolution of singularities there exists a smooth variety $Y_1$ together with a morphism $\pi_1: Y_1\to Y$ which induces an isomorphism $Y_1-\pi^{-1}_1(Z)\to Y-Z$, and such that
$\pi^{-1}_1(Z)$ as well as $\pi^{-1}_1(Z)\cup graph \cup \pi^{-1}_1(D_\infty-Z)$ are simple normal crossing divisors. 

Notice that $D_{graph}=\overline{\pi_1^{-1}(graph)}$ is a smooth divisor in $Y_1$ containing a s.n.c. divisor
$A= \overline{\pi_1^{-1}(graph)}-\pi_1^{-1}(D_\infty)$ as well as the divisor $B=\overline{\pi_1^{-1}(D_0)}$.
Let $Z_1=A\cap B$. Then $A\cup B$ is a n.c.d  outside $Z_1\subset D_{graph}$. By the embedded resolution of singularities there is a sequence of blow-ups of $D_{graph}$ with smooth centers contained in $Z_1$ or in its proper transforms such that in the resulting smooth variety the proper transforms of $A,B$  and exceptional divisors form a s.n.c. divisor. We can repeat the sequence of blow-ups with the same centers in the ambient variety $Y_1$. The resulting resolution of singularities $\pi_2:Y_2\to Y$ enjoy the following properties:

a) $\pi^{-1}_2(Z)$ as well as $\pi^{-1}_2(Z)\cup graph \cup \pi^{-1}_2(D_\infty-Z)$ are simple normal crossing divisors. 

b) the union of $\overline{\pi_2^{-1}(D_0\subset graph)}$ with the other divisors in $\overline{\pi_2^{-1}(graph)}$ is a s.n.c. divisor.

Finally we take  $\overline{X}:=\overline{\pi_2^{-1}(graph)}$. Then $\overline{X}$ is a smooth subvariety of $Y_2$ which contains an open subset isomorphic to $\overline{graph}-Z\supset X$.

We can consider the closure  $\overline{D}_0$  as a subset of $\overline{X}$. We define $D_v=\pi_2^{-1}(D_\infty)\cap \overline{X}$.
The remaining divisors are of two types. We denote by $D_{log}$ the union of those for which the residue of $\alpha$ is not equal to zero. The union of others are denoted by $D_h$.
Finally we denote by $f$ a holomorphic map $f:\overline{X}\to {\bf P}^1$ corresponding to $F$. The Proposition is proved.
$\blacksquare$

\begin{defn}\label{log-extension} 
Let us fix a compactification $\overline{X}$ as in the Proposition \ref{compactification for 1-forms}.

Let $X_{log}=X\cup D_{log}$. We call $X_{log}$ the logarithmic extension (log-extension for short) of a variety $X$ associated with the compactification $\overline{X}$. 

\end{defn}

 \subsection{Twisted de Rham cohomology in the case of $1$-forms}\label{de Rham for 1-forms}

We will use the conventions of the  Section \ref{compactification in case of 1-forms}.
Using local  computations in coordinates one shows that the sheaf  of de Rham differential forms on $X$, vanishing on each smooth component of the divisor $D_0$  is isomorphic  to $\Omega^{\bullet}_X(log\,D_0)(-D_0)$.  
Using e.g.  local computations in coordinates one can show that the sheaf  of de Rham differential forms on $X$, vanishing on each smooth component of the divisor $D_0$  is isomorphic  to $\Omega^{\bullet}_X(log\,D_0)(-D_0)$.  

If a smooth component of $D_{log}$ is locally in analytic coordinates given by the equation $x_i=0$ and does not belong to $D_v$ then near this component we can  express the form $\alpha$ as $\alpha={c_i dx_i\over{x_i}}+regular\, terms$, where $c_i\ne 0$ is a residue of $\alpha$ at this component. Notice that this is true in Zariski topology as well, since the notion of log-form is algebro-geometric. Clearly the change of $\alpha$ by $\alpha/t$ does not affect $D_{log}$.

We denote by $\Omega^{\bullet}_{\overline{X},D}$ the sheaf of de Rham forms on $\overline{X}$ having logarithmic poles on $(D_{log}-D_v)\cup D_h$ 
and whose restriction to $D_0$ vanishes.
 
\begin{defn}\label{De Rham cohomology for 1-forms}Let  $t\in \C^{\ast}$. We define  the (twisted) global de Rham cohomology  by the formula 
$$H^{\bullet}_{DR, glob, t}(X,D_0,\alpha)={\mathbb H}^{\bullet}(\overline{X}, (\Omega^\bullet_{\overline{X}, D}, td+{\alpha}\wedge (\bullet))).$$
The hypercohomology can be taken either in Zariski or in analytic topology.

\end{defn}

Varying the parameter $t$ we can think of  de Rham cohomology for $1$-forms as a $\C[t]$-module, similarly to the case of functions.

The global de Rham cohomology defined in this way is different from the naive one given by $\R\Gamma(X, (\Omega_X^\bullet, td+\alpha\wedge(\bullet)))$. The relationship between them will be explained in
Proposition \ref{naive de Rham}. In the case when $\alpha=df$ both definitions are equivalent. 

\begin{rmk}\label{independence on compactification}
a) We are not going to discuss here the question whether the de Rham cohomology for $1$-forms  depends on the choice of  compactification. 
One can hope that the answer is negative.

b) In what follows in order to save the notation we will often (but not always) assume that $D_0=\emptyset$ and omit $D_0$ from the notation. In this case we will denote the de Rham cohomology by $H^{\bullet}_{DR,t}(X,\alpha)$. All the results hold in the case of non-trivial $D_0$ as well. 
\end{rmk}

Notice that we have 
$$(X\to \overline{X})_\ast(\Omega^{\bullet}_X(log\,D_0)(-D_0))=\Omega^{\bullet}_{\overline{X}}(log(\overline{D}_0))(-\overline{D}_0)(\ast(\overline{D}_0+D_h+D_v+D_{log})).$$
As before, the notation  $\ast(divisor)$ means that we  consider  differential forms with poles of arbitrary order on the divisor.

\begin{defn}\label{logarithmic subsheaf}
Define a subsheaf $\Omega^{\bullet}_{\overline{X},\alpha}$ of the sheaf
$\Omega^{\bullet}_{\overline{X}}(log(\overline{D}_0+D_h+D_v+D_{log}))(-\overline{D}_0)$ as the one consisting of   forms $\eta\in \Omega^{\bullet}_{\overline{X}}(log(\overline{D}_0+D_h+D_v+D_{log}))(-\overline{D})$  such that $\alpha\wedge \eta$ still belongs to $\Omega^{\bullet}_{\overline{X}}(log(\overline{D}_0+D_h+D_v+D_{log}))(-\overline{D}_0)$ .
\end{defn}
This subsheaf is closed with respect to the de Rham differential $d$ as well as with respect to the multiplication operator $\alpha\wedge (\bullet)$.

\begin{prp}\label{almost vector bundle for 1-forms} 

a) $\Omega^{\bullet}_{\overline{X},\alpha}$ is a vector bundle on $\overline{X}$.

b) The hypercohomology ${\mathbb H}^{\bullet}(\overline{X},(\Omega_{\overline{X},\alpha}^{\bullet},td+{\alpha}\wedge (\bullet)))$ considered as an analytic family of vector spaces on $\C^\ast_t$ is in fact a vector bundle outside of finitely many points $t\in \C^\ast$. Furthermore it naturally extends  to $t=0$ giving a vector bundle on a sufficiently small disc $|t|<\epsilon$.

\end{prp}

{\it Proof.} In order to prove a) notice that the only non-locally free part of this coherent sheaf can be supported on $D_v$. Then part a) can be checked locally, and it follows from [KaKoPa2], Section 2.4 where it was proved for $\alpha=df$.  Part b) follows from the observation that the family in $t$ of the  hypercohomology groups is a coherent sheaf on $\C^\ast_t$ as well as from the degeneration of the Hodge-to-de Rham spectral sequence.\footnote{To establish the degeneration result can considers the bundle with meromorphic connection corresponding to the meromorphic extension of $\alpha$ to $\overline{X}$. Then the desired result comes from the degeneration at $E_1$ page of the spectral sequence of the pushforward to the point of the corresponding mixed twistor $D$-module (see [Moch4], [Sab6]).}
$\blacksquare$



\begin{prp} \label{qis with log forms} For any $t\in \C^\ast$ the natural embedding induces a quasi-isomorphism of complexes of sheaves
$$\left(\Omega^{\bullet}_{\overline{X},\alpha}, td+{\alpha}\wedge (\bullet)) \right)\to \left(\Omega^\bullet_{\overline{X}, D}, td+{\alpha}\wedge (\bullet)\right).$$

\end{prp}

{\it Proof.} Same as the proof of a similar statement in [KaKoPa2], Section 2.4.$\blacksquare$

It follows that in the definition of the twisted de Rham cohomology we can use the complex $\left(\Omega^{\bullet}_{\overline{X},\alpha}, d+{\alpha\over{t}}\wedge (\bullet)) \right)$. In either case the corresponding coherent sheaf gives rise to a vector bundle over a small disc in $\C_t$.
The following proposition shows that the naive approach via the twisted de Rham complex on $X$ gives a different answer.

\begin{prp}\label{naive de Rham} The natural embedding of complexes 
$$\left(\Omega^\bullet_{\overline{X},\alpha},d+{\alpha\over{t}}\wedge (\bullet) \right)\to (X\to \overline{X})_\ast\left(\Omega^\bullet_X,  d+{\alpha\over{t}}\wedge (\bullet)\right)$$
 is a quasi-isomorphism as long as $c_l/t+n$ is not equal to zero for all $l$ and all $n\in \Z_{\ge 1}$.

\end{prp}
{\it Sketch of the proof.}
It suffices to make local computations in analytic coordinates near smooth points of the divisors.  Using  appropriate filtrations on the complexes in question, we reduce the Proposition to a $1$-dimensional problem.  Namely, in the local coordinate $x$ on the line we should compare the cohomology of the $2$-terms complex $\C[x]\to \C[x]{dx\over{x}}, x^n\mapsto (n+{c\over{t}})x^n {dx\over{x}}$ with the similar cohomology of the $2$-terms complex $\C[x,x^{-1}]\to \C[x, x^{-1}]{dx\over{x}}$. The differential here is just the differential 
$d+{c\over{t}}{dx\over{x}}$ written in the basis $\{x^n\}$. We see that the cohomology groups coincide as long as $n+{c\over{t}}\ne 0$ for $n\ge 1$. This concludes the proof. $\blacksquare$

Hence outside of the values of $t\in \C^\ast$ specified in the Proposition \ref{naive de Rham} we have $H^{\bullet}_{DR,glob, t}(X,D_0,\alpha)\simeq \mathbb{H}^\bullet(\overline{X},(X\to \overline{X})_\ast\left(\Omega^\bullet_X,  d+{\alpha\over{t}}\wedge (\bullet)\right) )$.

\begin{cor}\label{de Rham cohomology as almost vector bundle} The family of vector spaces 
$$\mathbb{H}^\bullet(\overline{X},(X\to \overline{X})_\ast\left(\Omega^\bullet_X,  d+{\alpha\over{t}}\wedge (\bullet)\right) )$$
gives rise to a holomorphic vector bundle outside of the union of finitely many arithmetic progressions and finitely many points .
\end{cor}
{\it Proof.} This follows immediately from the Propositions \ref{almost vector bundle for 1-forms},   \ref{qis with log forms}, \ref{naive de Rham}. $\blacksquare$



From now on we will make the following two assumptions:
\[
\boxed{
\begin{array}{l}
\text{A1) $D_0=\emptyset$.} \\
\text{A2) The restriction of $\alpha$ to smooth strata of $D_h-(D_v\cup D_{log})$ has no zeros.}
\end{array}
}
\]

\begin{rmk}\label{comment on the assumption}
 The  assumption A2) implies that the locus of zeros of $\alpha$ is compact. This fact will be used in both de Rham and Betti frameworks. We expect that all results formulated below can be generalized to the case when  both assumptions A1) and A2) are omitted.

\end{rmk}

\subsection{De Rham global-to-local isomorphism in the case of $1$-forms}\label{de Rham local to global for 1-forms}

Let   $X$ be a smooth complex quasiprojective algebraic variety and $\alpha$  a closed regular $1$-form on $X$. Let us fix a projective compactification $\overline{X}$ as in the Proposition \ref{compactification for 1-forms}. 

We denote by $\ZZ(\alpha)$ the closed set of zeros of $\alpha$. By  assumption   A2) this set is compact. We have a decomposition $\ZZ(\alpha)$ into the disjoint union of finitely many non-empty connected components $\ZZ_i:=\ZZ_i(\alpha)$, $\ZZ(\alpha)=\sqcup_i \ZZ_i$.

\begin{defn}\label{local de Rham for 1-forms} We define local de Rham cohomology as a $\C[[t]]$-module 
$$H^\bullet_{DR,loc}(X,\alpha)={\mathbb H}^{\bullet}(U(\ZZ(\alpha)), (\Omega_{U(\ZZ(\alpha))}^{\bullet}[[t]], td+{\alpha}\wedge(\bullet))),$$
where $U(\ZZ(\alpha))\subset X$ is an open neighborhood   of $\ZZ(\alpha)$ in analytic topology.
\end{defn}
Local de Rham cohomology does not depend on the neighborhood of $\ZZ(\alpha)$ because the complex of sheaves $(\Omega_{U(\ZZ(\alpha))}^{\bullet}[[t]], td+{\alpha}\wedge(\bullet)))$ is acyclic outside of $\ZZ(\alpha)$.
Clearly 
$$H^\bullet_{DR,loc}(X,\alpha)\simeq
\oplus_{i\in I}\mathbb{H}^\bullet(U(\ZZ_{i}),(\Omega^\bullet_{U(\ZZ_{i}),\alpha}[[t]], td+\alpha\wedge (\bullet))), $$
where $U(\ZZ_{i})\subset X$ is an open neighborhood of the component $\ZZ_{i}$.

\begin{rmk}\label{whole space}
Similarly to the case of functions in Definition \ref{local de Rham} we can replace $U(\ZZ(\alpha))$ in Definition \ref {local de Rham for 1-forms} by $X$ and use Zariski topology. Also we obtain the same answer if we replace $U(\ZZ(\alpha))$ by $\overline{X}$ and the sheaf $\Omega_{U(\ZZ(\alpha))}^{\bullet}$ by $\Omega_{\overline{X},\alpha}^{\bullet}$.
\end{rmk}



Similarly to the case of functions, the global de Rham cohomology 
$$H_{DR,glob, t}^\bullet(X,\alpha)$$
 gives rise to a coherent sheaf $\mathcal{H}^\bullet_{DR,glob}(X,\alpha)$ on $\C$, i.e. a $\Z$-graded finitely generated $\C[t]$-module $H^\bullet_{DR,glob}(X,\alpha)[t]$. 
We define 
$${H}_{DR,glob}^\bullet(X,\alpha)[[t]]={H}_{DR,glob}^\bullet(X,\alpha)[t]\otimes_{\C[t]} \C[[t]].$$

  Then we have the following {\bf de Rham global-to-local isomorphism for $1$-forms}.

\begin{prp}\label{local to global de Rham for 1-forms}

1) There is a natural isomorphism of $\C[[t]]$-modules:
$${H}_{DR,glob}^\bullet(X,\alpha)[[t]]\simeq H^\bullet_{DR,loc}(X,\alpha).$$

2) Furthermore, $H^\bullet_{DR,loc}(X,\alpha)$ is a topologically free $\C[[t]]$-module.

\end{prp}
{\it Proof.}
Isomorphism 1) follows from Remark \ref{whole space} and the acyclicity of the complex of sheaves  $(\Omega_{\overline{X},\alpha}^{\bullet}[[t]], td+{\alpha}\wedge(\bullet)))$ on the complement to $\ZZ(\alpha)$.

In order to prove 2)  we notice that  there exists an open neighborhood $U\subset X$ of $\ZZ_{red}(\alpha)$ where $\alpha=df$ for some analytic function $f$ such that $f_{|\ZZ_{red}(\alpha)}=0$. Then we use Theorem \ref{degeneration in Kahler case}.$\blacksquare$

It follows from the Proposition \ref {local to global de Rham for 1-forms} that the coherent sheaf $\mathcal{H}_{DR,glob}(X,\alpha)$ gives rise to a 
free $\C[[t]]$-module (equivalently, a vector bundle over a formal disc),
$$H^\bullet_{DR,glob}(X,\alpha)[[t]]=\Gamma(\C, \mathcal{H}_{DR,glob}(X,\alpha))\otimes_{\C[t]}\C[[t]].$$ 
 Proposition \ref{local to global de Rham for 1-forms} ensures that there is an isomorphism of  vector bundles over $Spf(\C[[t]])$ defined by local and global de Rham cohomology.

 The vector bundle coming from the global de Rham cohomology does not have a priori a natural connection because the corresponding family is not isomonodromic. \footnote{This is different from the case of function when we do have such a connection.} Nevertheless one can define a natural meromorphic connection  on $H^\bullet_{DR,glob}(X,\alpha)[[t]]$.
 For that we use the function $f$  from the proof of Proposition \ref{local to global de Rham for 1-forms} specified by the condition $f_{|\ZZ_{red}(\alpha)}=0$.
 The corresponding connection is regular singular at $t=0$, with the actual order of  pole   at most $2$. Notice that in the case of functions it does not coincide with the restriction of the global connection to the formal neighborhood of $t=0$.

\subsection{Betti global and local cohomology in the case of $1$-forms}\label{Betti for 1-forms}

Let  $X$ be as before a  complex manifold of finite dimension $n$, $\alpha$ a holomorphic
closed $1$-form on $X$. We fix a compactification as in the Proposition \ref{compactification for 1-forms} and use the notation $\ZZ(\alpha)$ for
the set of zeros of $\alpha$ as before. Recall that  $\ZZ(\alpha)=\sqcup_{i\in I}\ZZ_i(\alpha)$ is a union  of finitely many compact connected components.  


Let $\overline{X}_{cor}$ denote the oriented real  blow-up of $\overline{X}$ at $D:=\overline{X}-X$. This is a real compact manifold with corners, which is homotopy equivalent to $X$. Its boundary $\partial\overline{X}_{cor}$ is the oriented real blow-up of $\overline{X}-X$. There is a natural projection map $p:\overline{X}_{cor}\to \overline{X}$ which maps  $\partial\overline{X}_{cor}$ onto the union of divisors $D_v\cup D_{log}\cup D_h$.

The differential form $\alpha/t$ gives  rise to a flat connection on the trivial rank $1$ bundle on $X$. {\it We denote by $E:=E_{\alpha,t}$ the corresponding locally-constant sheaf of flat sections}. By homotopy reasons $E$ extends to the local system on $\overline{X}_{cor}$. We keep the same notation for this extension. Let ${D}_v^\R$ denotes the proper transform of $D_v$ under the above-mentioned oriented real blow-up.
In other words it is a set of points $y\in \overline{X}_{cor}$ such that $p(y)\in D_v$. Recall the meromorphic function $f$ from the Proposition \ref{compactification for 1-forms}, part 3). The germ of $f$ at each point $x\in D_v$ is well-defined modulo adding a germ of a holomorphic function. Let $f_x$ denote a representative of the germ of $f$ at $x$.

\begin{defn}\label{irregular real blow-up}
For any $t\in \C^\ast$ we denote by $D_{v,t}^{\R,\pm}$  closed subset of $\partial\overline{X}_{cor}$ defined by its germs such as follows:

a) at a point $y$ such that $p(y)\notin D_v$ the germ is empty;

b) for point $y\in \overline{X}_{cor}$ such that $p(y)\in D_v$, the germ is the intersection $\{\overline{Re(\pm f_{p(y)}/t)\ge 0}\}\cap \partial\overline{X}_{cor}$.
Here we use the above notation for a representative of the germ of $f$, and take the intersection with the boundary of the closure of the set of points where the real value of the corresponding function is positive or negative.
\end{defn}

We handle the log divisor similarly. We denote by $D_{log}^{\R}$ the set $\{y\in  \overline{X}_{cor}|p(y)\in D_{log}\}$

\begin{defn}
 We denote by $D_{log,t}^{\R,\pm}$  the closure of the set of such points $y\in D_{log}^{\R}$  that for any irreducible component $D_{i,log}\subset D_{log}$ containing $p(y)$ the following holds: $\mp Re(c_i/t)>0$. Here $c_i$ is the residue of $\alpha$ at $D_{i,log}$.
 \end{defn}

 Let $\overline{X}_{cor,t}^\pm$ denote the open subset of $\overline{X}_{cor}$ obtained by throwing away $(D_{log,t}^{\R}-D_{log,t}^{\R,\mp})\cup (D_v^{\R}-D_v^{\R,\mp})\cup D_h$. Recall that when working with sheaves we  systematically use the notation like $(X\to Y)_\ast(F)$ for the direct, etc. image of the sheaf $F$ under the morphism $f:X\to Y$ instead of more traditional $f_\ast(F)$, etc.

\begin{defn}\label{global Betti for 1-forms}
We define the  global Betti cohomology at $t\in \C^{\ast}$ as 
$$H_{Betti,glob, t}^\bullet(X,\alpha)=H^\bullet(\overline{X}_{cor},(X\to \overline{X}_{cor})_\ast( E_{\alpha,t})\otimes(\overline{X}_{cor,t}^+\to \overline{X}_{cor})_{!}(\underline{\Z}_{\overline{X}_{cor}^+}))$$
$$\simeq H^{\bullet}(\overline{X}_{cor}, D_{v,t}^{\R,-}\cup D_{log,t}^{\R,+}, E_{\alpha,t}).$$
\end{defn}

If $\alpha=df$ the  Betti cohomology is isomorphic to  $H^k(X, (f/t)^{-1}(-\infty), \Z)\otimes \C$, hence the new definition agrees with the  one for functions. On the other hand some  periods of  $\alpha$ can be non-trivial in the non-exact case, hence in general there is no $\Z$-lattice  in Betti cohomology. {\it Therefore the Betti cohomology for $1$-forms is defined over $\C$ only}. Notice also that we use only the compactifying divisor $D_v\cup D_{log}$ in the definition of Betti cohomology. The rest of $D$ does not play any role.
 
 {\it The family over $\C^\ast_t$ of Betti cohomology groups $H^\bullet_{Betti,glob, t}(X,\alpha)$ forms a coherent sheaf over $\C^\ast$ which we will denote by $\mathcal{H}_{Betti,glob}(X,\alpha)$.} Its restriction to a small punctured disc about $t=0$ is a holomorphic graded vector bundle.

 In what follows we will use the version of Betti cohomology with coefficients in local systems (i.e. locally constant sheaves) over   non-archimedean rings. For that we will need an appropriate notion of the Novikov ring which we discuss below.
  
 \begin{defn}\label{Novikov rings}
 Let $A$  be a commutative ring with the unit which has no zero divisors. We define the {\it Novikov ring} $Nov_A$ as the ring of series  $\sum_{i\ge i_0}a_iT^{\lambda_i}$ such that $a_i\in A, \lambda_i\in \R$ and $ \lambda_i\to +\infty$ as $i\to +\infty$. There are several important special cases  which will appear in the paper: $A=\Z, \R, \C$ and $A=P[T^{i\R}]=\sum_{j_1\le j\le j_2}a_jT^{\sqrt{-1}\mu_j}, a_j\in P$, where $P$ is either $\Z,\R$ or $\C$. If the ring $P$ is clear from the context we will sometimes denote the Novikov ring simply by $Nov$.
\end{defn}  
 
 One immediate application of the Novikov rings to the global Betti cohomology is explained in the next remark.
 
 \begin{rmk}\label{covering by two intersecting arcs}
 A slightly different but equivalent definition of the Betti cohomology can be obtained such as follows. Instead of the condition $Re(c_i/t)>0$ we will use the condition $-\pi/2-\epsilon<Arg(c_i/t)<\pi/2+\epsilon$ for sufficiently small $\epsilon>0$. Similarly, instead of the condition $Re(c_i/t)<0$ we use the condition $\pi/2-\epsilon<Arg(c_i/t)<3\pi/2+\epsilon$. This gives the covering of $S^1$ by two open overlapping arcs.  We have defined above the ring $Nov$ consisting of series $\sum_{i\ge i_0}a_iT^{b_i}$ such that $a_i\in \Z, b_i\in \C, Re(b_i)\to  +\infty$ as $i\to +\infty$. Then the map $\int\alpha:H_1(X,\Z)\to  \C\subset Nov$ induces the homomorphism $\pi_1(X)\to GL(1,Nov)$. Hence for each $t\in \C^\ast$  the rescaled $1$-form $\alpha/t$ gives rise to  a non-archimedean local system which we (abusing the notation) will  denote by $E_{\alpha,t}$. The  holonomy of this local system along each of  the two connected components  of the above intersecting arcs has the non-archimedean norm strictly less than $1$. Hence the holonomy is non-trivial. It follows that $H^1(S^1, E_{\alpha,t})=0$ (to show that one can consider the exact sequence of the corresponding pair).   To each ``sufficiently narrow'' convex cone $Cone_{\alpha/t}\subset H^1(X,\R)$ containing the ray $\R_{>0}\cdot [Re(\alpha/t)]$ one can assign a tube domain $U(Cone_{\alpha/t})$ in the space of non-archimedean local systems $Hom(\pi_1(X)\to Nov)$ such that the Betti cohomology $H_{Betti,glob, t}(X,\rho)$ does not change as long as
 $\rho\in  U(Cone_{\alpha/t})$. Here $t\in \C^\ast$ is fixed, but we can allow it to vary in a sector in such a way that the above cone does not change. 
 All the above can be said about local systems with values in any ring $Nov_A$.
 \end{rmk}

Let us now define the {\it local Betti cohomology}. Differently from the global Betti cohomology this one will live over the ring $\Z$ making the situation similar to the previously discussed case of functions. Here are more details.

For an element $i\in I$ we fix a sufficiently small positive $\varepsilon$ and consider an $\varepsilon$-neighborhood of the connected component $\ZZ_i(\alpha)$ which we denote by $U_{\varepsilon,i}:=U_{\varepsilon}(\ZZ_i(\alpha))\subset X$ (in order to define the neighborhood we choose a Riemannian metric on $X$, but the cohomology below does not depend on the choice). For each $i\in I$ let us fix a holomorphic function $W_i$ on $U_{\varepsilon,i}$ such that $dW_i=\alpha$ and $W_i=0$ on $\ZZ_i(\alpha)$. Such a function $W_i$ does exist, because the restriction of $\alpha$ to $U_{\varepsilon,i}$ is exact as long as $\varepsilon$ is sufficiently small.

For $\theta\in \R/2\pi \Z, \theta=Arg(t)$ and $j\in I$ we define a graded $\Z$-module  $H_{Betti, j, t}^\bullet(X,\alpha)$ as the relative cohomology
$H^{\bullet}(U_{\varepsilon,j},U_{\varepsilon,j}\cap W_j^{-1}(\varepsilon\cdot e^{i\theta}),\Z)$.

\begin{defn}\label{local Betti for 1-forms}
Let us fix $t\in \C^{\ast}$. The direct sum $\oplus_j H_{Betti, j, t}^\bullet(X,\alpha):=H^\bullet_{Betti,loc,t}(X,\alpha)$  is called the local Betti cohomology .
\end{defn}

Sometimes instead  we will say   ``local Betti cohomology in the direction $\theta=Arg(t)$".

Notice that the local Betti cohomology is a graded $\Z$-module. Then varying $t\in \C^\ast$ we obtain a locally constant sheaf  $\mathcal{H}_{Betti,loc,\Z}^\bullet(X,\alpha)$ on $\C^\ast$ of finite rank graded $\Z$-modules. Tensoring with $\OO_{\C^\ast}$ we get {\it a holomorphic vector bundle} $\mathcal{H}_{Betti,loc}^\bullet(X,\alpha)$. It carries a meromorphic structure at $t=0$ coming from the regular singular flat connection. {\it The space of germs of its meromorphic sections at  $t=0$ denoted by $H_{Betti,loc}^\bullet(X,\alpha)$ is a free graded $\C\{t\}[t^{-1}]$-module.} 

\begin{rmk}\label{non-zero critical values}
Our definition of local Betti cohomology uses a \\
non-canonical choice of the function $W$ normalized  in such a way that $W=0$ on 
$\ZZ(\alpha)$. More generally we may assume that the critical values
$w_i:=W_{|\ZZ_i(\alpha)}$ are different. Then we can modify the above formulas by taking $W_i=W-w_i$ near each component $\ZZ_i(\alpha)$.

\end{rmk}

\subsection{Global de Rham-to-Betti isomorphism for $1$-forms}\label{global de Rham to Betti for 1-forms}

Let $D_{log}^\pm$ denote the union of those smooth components of $D_{log}$ for which $\mp Re(c_i)>0$.
We denote by $\Omega^{\bullet}_{\overline{X}, D^\pm}$ the subsheaves of  $\Omega^{\bullet}_{\overline{X},D}$  defined in the same way as the latter, but such that the condition on logarithmic poles  is imposed at $D_{log}^\pm$ only instead of the whole $D_{log}$. We use the notation      $\Omega^{\bullet}_{\overline{X}, D^{\pm},t}$  for  similarly defined sheaves but with $\alpha$ being replaced by  $\alpha/t$.

\begin{prp}\label{global de Rham to global Betti}
1) For any $t\in \C^\ast$  such that $Re(c_i/t)>0$ for all $i$ then there is a natural isomorphism 
$$\Phi_t^+: H^\bullet_{DR,glob,t}(X,\alpha)\simeq H^\bullet_{Betti,glob,t}(X,\alpha).$$ 

2) If $Re(c_i/t)\le 0$ then there is a natural isomorphism 
$$\Phi_t^-:  H^\bullet_{DR,glob,t}(X,\alpha)\simeq H^\bullet(\overline{X}_{cor}, E_{\alpha, t}).$$

\end{prp}

{\it Sketch of the proof without technical details.} It suffices to establish a derived quasi-isomorphism of complexes of sheaves on $\overline{X}$ on the Betti and de Rham sides, namely,  $(\overline{X}_{cor}\to \overline{X})_\ast( E_{\alpha,t})\otimes(\overline{X}_{cor,t}^+\to \overline{X}_{cor})_{!}(\underline{\Z}_{\overline{X}_{cor}^+}))\simeq
(\Omega^\bullet_{\overline{X}, D}, d+{\alpha\over{t}}\wedge (\bullet))$. For the restrictions of the sheaves to $X$ this is obvious. It is left to consider their restrictions to analytic neighborhoods of points of $D_{v,t}^{\R}$ and $D_{log,t}^\R$. Then the result should follow from considerations in coordinates similar to those in the proof of  Proposition \ref{naive de Rham}. $\blacksquare$

Let us illustrate the Proposition \ref{global de Rham to global Betti} in $1$-dimensional case. Since we can work in analytic topology we take $X=\C^\ast$, so  $\overline{X}=\C{\bf P}^1$, and take $ \alpha=c{dx\over{x}}-dx$.
Then $D_{log}=\{0\}$, $D_v=\{\infty\}$ and $D=\{0\}\cup \{\infty\}$.
We may assume that $t=1$ and then consider two cases: $Re(c)>0$ and $Re(c)\le 0$. Then we have two complexes:

i) $\C[x]\to \C[x]{dx\over{x}}$;

ii) $\C[x,x^{-1}]\to \C[x,x^{-1}]{dx\over{x}}$.

In both cases the differential is $d+\alpha\wedge (\bullet)$.

There is a natural embedding $i)\to ii)$, which is a quasi-isomorphism as long as $c\notin \Z_{\ge 1}$. Notice that ${D}_{log}^\R$ is the oriented real blow-up of $\C{\bf P}^1$ at the origin.  We denote by $S^1_0$ the exceptional fiber (circle). Complex ii) is the complex of differential forms on $X$ with poles of arbitrary order on $D_{log}$. Its  cohomology is isomorphic to $H^\bullet(X, E_\alpha)$ which is in turn isomorphic to $H^\bullet(\overline{X}_{cor}, E_\alpha)$, where we set $E_\alpha:=E_{\alpha, t=1}$. The isomorphism follows from  homotopy equivalence of $X$ and $\overline{X}_{cor}$. If $Re(c)\le 0$ then i) is quasi-isomorphic to ii) and hence i) computes $H^\bullet(\overline{X}_{cor}, E_\alpha)$. This is part 2) of Proposition \ref{global de Rham to global Betti}. 

In order to demonstrate part 1) in this example  we use the Poincar\'e duality between $H^\bullet(\overline{X}_{cor}, S^1_0, E_\alpha)$ and $H^\bullet(\overline{X}_{cor}, E_{-\alpha})$. Let us now consider dual complexes to i) and ii). They are quasi-isomorphic if $Re(c)\notin \Z_{\le -1}$. Dual to i) is the complex
$x\C[x]\to \C[x]dx$ endowed with the differential  $d-\alpha\wedge (\bullet)$. We see that the cohomology of the dual to i) complex is isomorphic if $Re(c)<0$ to $H^\bullet(\overline{X}_{cor}, E_{-\alpha})$. Thus we have reduced the case 1) to case 2).


\begin{rmk} \label{non-archimedean local system} Notice that in the RHS of the Proposition \ref{global de Rham to global Betti} we can replace $E_{\alpha,t}$ by an arbitrary local system $\rho$. Then it becomes a finite type $\Z[\pi_1(X)]$-module. In this way we  obtain a coherent sheaf on the torus 
$$Hom(H_1(X,\Z)\to {\bf G}_m).$$
Let us now fix $\theta=Arg(t)$. Then the Betti cohomology $H^\bullet_{Betti,glob,t}(X,\rho):=H^\bullet_{Betti,glob,\theta}(X,\rho)$ makes sense if $\rho$ is a local system over the Novikov ring $Nov_\theta=\{\sum_{l\in \Z}n_le^{-{\lambda_l\over{t}}}\}$,  where $n_l\in \Z$, and coefficients $\lambda_l\in\C$ satisfy the condition $Re(\lambda_le^{-i\theta})\to +\infty$ as long as $l\to +\infty$. Same is true for other versions of the Novikov ring discussed previously, or for a local system over a complete valuation field $K$. We will use this  generalization later.

\end{rmk}

\subsection{Betti local cohomology in the case of $1$-forms on real manifolds}\label{Betti cohomology for real one-forms}

In this subsection we work in a bigger generality of real-valued $1$-forms on real manifolds. Although the results of this subsection are interesting by themselves, we present them here for another purpose: they will be used in the next subsection in the construction of the Betti global-to-local isomorphism for non-archimedean local systems.

Let $Y$ be a compact manifold with corners which are of codimension   less or equal than $2$.  Let $\xi$ be a vector field on $Y$  such that the set of zeros $\ZZ(\xi)$ is the disjoint union of finitely many  connected  components $\ZZ_i(\xi), i\in I$. We assume 
that the pair $(Y, \xi)$ satisfies the following properties:

1) There  exist non-intersecting open subsets  $U_i^\prime\subset Y, i\in I$ such that $\ZZ_i(\xi)\subset U_i^\prime$ and smooth real functions $f_i: U_i^\prime\to \R$ such that $Crit(f_i)=\ZZ_i(\xi)$, the restriction of $f_i$ to $\ZZ_i(\xi)$ is equal to $0$.

2) We have $\xi(f_i)>0$ in $U_i^\prime-\ZZ_i(\xi)$.

We call such $\xi$ a  {\it gradient-like} vector field.

Assume that the boundary of $Y$ (i.e. the union of non-trivial corners) has a decomposition $\partial Y=\partial_-Y\cup \partial_{h}Y\cup \partial_+Y$ (here the subscript $h$ stands for ``horizontal'') such that $\partial_-Y\cap \partial_+Y=\emptyset$ and:


3) For each point $x_+\in \partial_+ Y$ the vector $\xi(x_+)$ is non-zero and pointed outside of $Y$, and 
for each point $x_-\in \partial_- Y$ the vector $\xi(x_-)$ is non-zero and pointed inside of $Y$.

4) The vector field $\xi$ is tangent to $\partial_{h}Y$ and $\xi_{|\partial_{h}Y}\ne 0$ everywhere on $\partial_{h}Y$.\footnote{This assumption is a real analog of the holomorphic assumption that $\alpha_{|D_h}$ does not have zeros. We expect that our results hold without these assumptions.}

\vspace{2mm}


\vspace{3mm}

\tikzset{every picture/.style={line width=0.75pt}} 

\begin{tikzpicture}[x=0.5pt,y=0.5pt,yscale=-1,xscale=1]

\draw  [line width=0.75] [line join = round][line cap = round] (88.5,53) .. controls (88.5,77.12) and (81.78,101.75) .. (80.5,126) .. controls (79.85,138.35) and (80.57,152.28) .. (83.5,164) .. controls (84.31,167.22) and (87.23,171.47) .. (87.5,175) .. controls (88.72,190.9) and (88.33,207.8) .. (90.5,223) .. controls (90.9,225.83) and (90.39,232.89) .. (93.5,233) .. controls (117.19,233.82) and (130.02,230.38) .. (152.5,229) .. controls (177.78,227.45) and (206.11,226.09) .. (231.5,227) .. controls (242.15,227.38) and (253.79,231.93) .. (264.5,233) .. controls (281.73,234.72) and (298.92,235.16) .. (316.5,236) .. controls (329.8,236.63) and (343.15,237.27) .. (356.5,236) .. controls (359.75,235.69) and (361.76,232.1) .. (364.5,231) .. controls (366.39,230.25) and (375.48,231.17) .. (375.5,230) .. controls (376.22,196.17) and (377.14,185.23) .. (371.5,155) .. controls (368.32,137.96) and (362.69,121.02) .. (359.5,104) .. controls (356.53,88.15) and (355.3,74.18) .. (356.5,58) .. controls (356.68,55.62) and (360.64,49.77) .. (359.5,50) .. controls (334.13,55.07) and (304.8,49.3) .. (279.5,47) .. controls (256.1,44.87) and (229.56,46.63) .. (207.5,46) .. controls (192.49,45.57) and (177.5,43.25) .. (162.5,44) .. controls (159.75,44.14) and (157.25,45.87) .. (154.5,46) .. controls (153.04,46.07) and (91.25,48.12) .. (86.5,49) .. controls (84.86,49.3) and (86.5,52.33) .. (86.5,54) ;
\draw    (82.5,123) -- (143.5,123) ;
\draw [shift={(145.5,123)}, rotate = 180] [color={rgb, 255:red, 0; green, 0; blue, 0 }  ][line width=0.75]    (10.93,-3.29) .. controls (6.95,-1.4) and (3.31,-0.3) .. (0,0) .. controls (3.31,0.3) and (6.95,1.4) .. (10.93,3.29)   ;
\draw    (86,89) -- (143.5,90.93) ;
\draw [shift={(145.5,91)}, rotate = 181.93] [color={rgb, 255:red, 0; green, 0; blue, 0 }  ][line width=0.75]    (10.93,-3.29) .. controls (6.95,-1.4) and (3.31,-0.3) .. (0,0) .. controls (3.31,0.3) and (6.95,1.4) .. (10.93,3.29)   ;
\draw    (85,165) -- (146.5,164.03) ;
\draw [shift={(148.5,164)}, rotate = 179.1] [color={rgb, 255:red, 0; green, 0; blue, 0 }  ][line width=0.75]    (10.93,-3.29) .. controls (6.95,-1.4) and (3.31,-0.3) .. (0,0) .. controls (3.31,0.3) and (6.95,1.4) .. (10.93,3.29)   ;
\draw    (91,201) -- (148.5,200.03) ;
\draw [shift={(150.5,200)}, rotate = 179.04] [color={rgb, 255:red, 0; green, 0; blue, 0 }  ][line width=0.75]    (10.93,-3.29) .. controls (6.95,-1.4) and (3.31,-0.3) .. (0,0) .. controls (3.31,0.3) and (6.95,1.4) .. (10.93,3.29)   ;
\draw    (359,95) -- (424.5,95.97) ;
\draw [shift={(426.5,96)}, rotate = 180.85] [color={rgb, 255:red, 0; green, 0; blue, 0 }  ][line width=0.75]    (10.93,-3.29) .. controls (6.95,-1.4) and (3.31,-0.3) .. (0,0) .. controls (3.31,0.3) and (6.95,1.4) .. (10.93,3.29)   ;
\draw    (365,129) -- (426.5,129) ;
\draw [shift={(428.5,129)}, rotate = 180] [color={rgb, 255:red, 0; green, 0; blue, 0 }  ][line width=0.75]    (10.93,-3.29) .. controls (6.95,-1.4) and (3.31,-0.3) .. (0,0) .. controls (3.31,0.3) and (6.95,1.4) .. (10.93,3.29)   ;
\draw    (374,168) -- (428.5,168) ;
\draw [shift={(430.5,168)}, rotate = 180] [color={rgb, 255:red, 0; green, 0; blue, 0 }  ][line width=0.75]    (10.93,-3.29) .. controls (6.95,-1.4) and (3.31,-0.3) .. (0,0) .. controls (3.31,0.3) and (6.95,1.4) .. (10.93,3.29)   ;
\draw    (120,48) -- (163.5,45.13) ;
\draw [shift={(165.5,45)}, rotate = 176.23] [color={rgb, 255:red, 0; green, 0; blue, 0 }  ][line width=0.75]    (10.93,-3.29) .. controls (6.95,-1.4) and (3.31,-0.3) .. (0,0) .. controls (3.31,0.3) and (6.95,1.4) .. (10.93,3.29)   ;
\draw    (206,47) -- (264.5,46.03) ;
\draw [shift={(266.5,46)}, rotate = 179.05] [color={rgb, 255:red, 0; green, 0; blue, 0 }  ][line width=0.75]    (10.93,-3.29) .. controls (6.95,-1.4) and (3.31,-0.3) .. (0,0) .. controls (3.31,0.3) and (6.95,1.4) .. (10.93,3.29)   ;
\draw    (154.5,230) -- (198,229.04) ;
\draw [shift={(200,229)}, rotate = 178.74] [color={rgb, 255:red, 0; green, 0; blue, 0 }  ][line width=0.75]    (10.93,-3.29) .. controls (6.95,-1.4) and (3.31,-0.3) .. (0,0) .. controls (3.31,0.3) and (6.95,1.4) .. (10.93,3.29)   ;
\draw    (227.5,228) -- (280.53,236.68) ;
\draw [shift={(282.5,237)}, rotate = 189.29] [color={rgb, 255:red, 0; green, 0; blue, 0 }  ][line width=0.75]    (10.93,-3.29) .. controls (6.95,-1.4) and (3.31,-0.3) .. (0,0) .. controls (3.31,0.3) and (6.95,1.4) .. (10.93,3.29)   ;

\end{tikzpicture}

\vspace{3mm}

\begin{rmk}\label{compatible 1-forms}
The term ``gradient-like''  comes from the following way of construction of such vector fields.
Let $\beta$ be a closed $1$-form on $Y$ such that for the set $\ZZ(\beta)$ of zeros of $\beta$ which is a union of finitely many connected components, we have $\ZZ(\beta)=\ZZ(\xi)$, and furthermore, outside of $\ZZ(\beta)$ we have $\beta(\xi)>0$. Having a Riemannian metric on $X$ we can canonically assign to $\beta$  a gradient-like  vector field $\xi=v_\beta$ via
the isomorphism $T^\ast Y\to TY$. In general, if $\beta(\xi)>0$ outside of the set $\ZZ(\beta)$ we will say that {\it 
$\beta$ is compatible with $\xi$}. 

Here is an example of this construction which is important for this paper: if $X$ is a complex compact manifold, $\alpha$ is a holomorphic $1$-form with compact set of zeros. Let $\beta=Re(\alpha)$. Then any hermitian metric on $X$ gives rise to a gradient-like vector field.

\end{rmk}

\begin{prp}\label{neighborhood of zeros}
Let us keep the  above notation and assume the properties 1)-4).

Then for each connected component $\ZZ_i(\xi)$ and any sufficiently small $\epsilon>0$ there exists a compact smooth manifold with corners $U_i\subset U_i^\prime$ which satisfies the following properties:

i) $int(U_i)\supset \ZZ_i(\xi)$;

ii) the boundary $\partial U_i$  can be decomposed
into a union $\partial U_i=\partial_-U_i\cup \partial_{h}U_i\cup\partial_+U_i$ such that the assumptions 1)-4) are satisfied with $X$ being replaced by $U_i$;

iii) $f_i$ restricted to $\partial_\pm U_i$ is equal to $\pm \epsilon$.


\end{prp}

{\it Proof.}  
By assumption 1) for each connected component $\ZZ_i(\xi)$ we can find an open neighborhood $U_i^\prime$ and function $f_i$ on $U_i^\prime$, which is equal to zero on $\ZZ_i(\xi)$. Hence $f_i^{-1}(0)-\ZZ_i(\xi)$ is a smooth manifold. Then we can find an open relatively compact subset $V_i\subset f_i^{-1}(0)$ which contains $\ZZ_i(\xi)$ and such that the boundary of the closure $\partial\overline{V}_i$ is a smooth hypersurface in $f_i^{-1}(0)$.

Next we construct a germ $g_i$ at $\overline{V}_i$ of a smooth function  in a neighborhood of this compact set in $X$ which satisfies the following properties:

a) $g_i$ is equal to  zero at all points $x\in V_i$;

b) on the boundary $\partial\overline{V}_i$ the germ $g_i$ is determined by the following properties:

b1) $\xi(g_i)=0$;

b2) in a neighborhood of $x\in f_i^{-1}(0)$ the germ $g_i$ is equal to zero at all points of $\overline{V}_i$, and
for each point $y\notin \overline{V}_i$ the germ is equal to $exp(-{1\over dist(y,\partial\overline{V}_i)})$. Here $dist$ can be taken with respect to any Riemannian metric.

Clearly the set of such germs is a section of the \'etale space of a sheaf on $\overline{V}_i$. Since $\overline{V}_i$  is compact we obtain in this way a section of this sheaf over some small neighborhood of $\overline{V}_i$ in $X$. Then we define $U_i$ as
$U_i=f_i^{-1}([-\epsilon,\epsilon])\cap g_i^{-1}((0,\delta))$, where $\epsilon$ and $\delta$ are sufficiently small positive numbers.
By definition $U_i$ is a relatively compact manifold with corners. Finally we set $\partial_\pm U_i=f_i^{-1}(\pm\epsilon)\cap \overline{U}_i$, and $\partial_h U_i=g_i^{-1}(\delta)\cap \overline{U}_i$.
$\blacksquare$

\vspace{2mm}


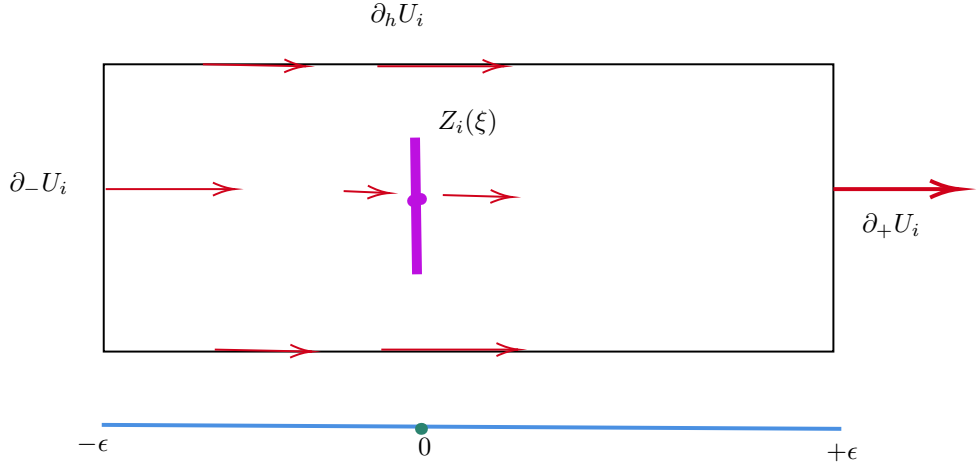
\tikzset{every picture/.style={line width=0.6pt}} 

\begin{tikzpicture}[x=0.6pt,y=0.6pt,yscale=-1,xscale=1]

\draw   (112,81) -- (480,81) -- (480,226) -- (112,226) -- cycle ;
\draw [color={rgb, 255:red, 208; green, 2; blue, 27 }  ,draw opacity=1 ]   (162,81) -- (212,81.96) ;
\draw [shift={(214,82)}, rotate = 181.1] [color={rgb, 255:red, 208; green, 2; blue, 27 }  ,draw opacity=1 ][line width=0.75]    (10.93,-3.29) .. controls (6.95,-1.4) and (3.31,-0.3) .. (0,0) .. controls (3.31,0.3) and (6.95,1.4) .. (10.93,3.29)   ;
\draw [color={rgb, 255:red, 208; green, 2; blue, 27 }  ,draw opacity=1 ]   (250,82) -- (312,82) ;
\draw [shift={(314,82)}, rotate = 180] [color={rgb, 255:red, 208; green, 2; blue, 27 }  ,draw opacity=1 ][line width=0.75]    (10.93,-3.29) .. controls (6.95,-1.4) and (3.31,-0.3) .. (0,0) .. controls (3.31,0.3) and (6.95,1.4) .. (10.93,3.29)   ;
\draw [color={rgb, 255:red, 208; green, 2; blue, 27 }  ,draw opacity=1 ]   (168,225) -- (215,225.96) ;
\draw [shift={(217,226)}, rotate = 181.17] [color={rgb, 255:red, 208; green, 2; blue, 27 }  ,draw opacity=1 ][line width=0.75]    (10.93,-3.29) .. controls (6.95,-1.4) and (3.31,-0.3) .. (0,0) .. controls (3.31,0.3) and (6.95,1.4) .. (10.93,3.29)   ;
\draw [color={rgb, 255:red, 208; green, 2; blue, 27 }  ,draw opacity=1 ]   (252,225) -- (319,225) ;
\draw [shift={(321,225)}, rotate = 180] [color={rgb, 255:red, 208; green, 2; blue, 27 }  ,draw opacity=1 ][line width=0.75]    (10.93,-3.29) .. controls (6.95,-1.4) and (3.31,-0.3) .. (0,0) .. controls (3.31,0.3) and (6.95,1.4) .. (10.93,3.29)   ;
\draw [color={rgb, 255:red, 208; green, 2; blue, 27 }  ,draw opacity=1 ]   (113,144) -- (176,144) ;
\draw [shift={(178,144)}, rotate = 180] [color={rgb, 255:red, 208; green, 2; blue, 27 }  ,draw opacity=1 ][line width=0.75]    (10.93,-3.29) .. controls (6.95,-1.4) and (3.31,-0.3) .. (0,0) .. controls (3.31,0.3) and (6.95,1.4) .. (10.93,3.29)   ;
\draw [color={rgb, 255:red, 208; green, 2; blue, 27 }  ,draw opacity=1 ][line width=1.5]    (480,144) -- (540,144) ;
\draw [shift={(543,144)}, rotate = 180] [color={rgb, 255:red, 208; green, 2; blue, 27 }  ,draw opacity=1 ][line width=1.5]    (14.21,-4.28) .. controls (9.04,-1.82) and (4.3,-0.39) .. (0,0) .. controls (4.3,0.39) and (9.04,1.82) .. (14.21,4.28)   ;
\draw [color={rgb, 255:red, 74; green, 144; blue, 226 }  ,draw opacity=1 ][line width=1.5]    (111,263) -- (484,265) ;
\draw [color={rgb, 255:red, 189; green, 16; blue, 224 }  ,draw opacity=1 ][line width=3.75]    (269,118) -- (270,187) ;
\draw [color={rgb, 255:red, 208; green, 2; blue, 27 }  ,draw opacity=1 ]   (233,145) -- (254,145.91) ;
\draw [shift={(256,146)}, rotate = 182.49] [color={rgb, 255:red, 208; green, 2; blue, 27 }  ,draw opacity=1 ][line width=0.75]    (10.93,-3.29) .. controls (6.95,-1.4) and (3.31,-0.3) .. (0,0) .. controls (3.31,0.3) and (6.95,1.4) .. (10.93,3.29)   ;
\draw [color={rgb, 255:red, 208; green, 2; blue, 27 }  ,draw opacity=1 ]   (283,147) -- (316,147.94) ;
\draw [shift={(318,148)}, rotate = 181.64] [color={rgb, 255:red, 208; green, 2; blue, 27 }  ,draw opacity=1 ][line width=0.75]    (10.93,-3.29) .. controls (6.95,-1.4) and (3.31,-0.3) .. (0,0) .. controls (3.31,0.3) and (6.95,1.4) .. (10.93,3.29)   ;
\draw  [color={rgb, 255:red, 42; green, 131; blue, 111 }  ,draw opacity=1 ][line width=4.5] [line join = round][line cap = round] (272,265) .. controls (272,265) and (272,265) .. (272,265) ;
\draw  [color={rgb, 255:red, 189; green, 16; blue, 224 }  ,draw opacity=1 ][line width=4.5] [line join = round][line cap = round] (268,150) .. controls (268.97,149.03) and (270.63,149) .. (272,149) ;

\draw (279,102) node [anchor=north west][inner sep=0.6pt]   [align=left] {$\displaystyle Z_{i}( \xi )$};
\draw (97,267) node [anchor=north west][inner sep=0.6pt]   [align=left] {$\displaystyle -\epsilon $};
\draw (475,271) node [anchor=north west][inner sep=0.6pt]   [align=left] {$\displaystyle +\epsilon $};
\draw (269,268) node [anchor=north west][inner sep=0.6pt]   [align=left] {0};
\draw (63,134) node [anchor=north west][inner sep=0.6pt]   [align=left] {$\displaystyle \partial _{-} U_{i}$};
\draw (245,49) node [anchor=north west][inner sep=0.6pt]   [align=left] {$\displaystyle \partial _{h} U_{i}$};
\draw (493,155) node [anchor=north west][inner sep=0.6pt]   [align=left] {$\displaystyle \partial _{+} U_{i}$};

\end{tikzpicture}

\vspace{3mm}

Now we take $U=\cup_{i\in I}U_i$ and define $f$ being a smooth function which coincides with $f_i$ on each $U_i$. 
Then $int(U)\supset \ZZ(\xi)$, the boundary of $U$ has a decomposition $\partial U=\partial_-U\cup \partial_hU\cup \partial_+U$ such that the properties 1)-4) are satisfied with $U$ replacing $X$ and moreover $f_{|\partial_\pm U}=\pm \epsilon$. Replacing $f$ by $f/\epsilon$ we may assume that $\epsilon=1$.

\begin{defn}\label{good neighborhood of zeros}
Let $Y$ be a manifold with corners without boundary and $\xi$ be a vector field on $Y$ such that the set of zeros $\ZZ(\xi)$ is compact.

We call a compact submanifold with corners $U\subset Y$ a good neighborhood of $\ZZ(\xi)$ if $U$ has full dimension, $\ZZ(\xi)\subset int(U)$ and there exists a smooth function
$f:U\to [-1,1]$ such that $\xi(f)>0$ on $U-\ZZ(\xi)$ and $f$ satisfies the conditions of Proposition \ref{neighborhood of zeros} with $\epsilon=1$.

\end{defn}

The Proposition  \ref{neighborhood of zeros} shows that the set of good neighborhoods of $\ZZ(\xi)$ is non-empty. Unless we say otherwise we will assume below that all appearing neighborhoods $U$ of $\ZZ(\xi)$ are good.

\begin{rmk}\label{product structure}
For a good neighborhood $U$ we see that there is a canonical diffeomorphism  $\psi: f^{-1}([-1,0))\simeq \partial_-U\times [-1,0)\subset f^{-1}(-1)\times [-1,0)$ and the function $f$ is identified under this diffeomorphism with the projection $pr_2$ to the second factor. Moreover this diffeomorphism identifies trajectories of $\xi$ with $\{x\}\times [-1,0), x\in \partial_-U$. Similar result holds for $f^{-1}((0,1])$ and $\partial_+U$.

\end{rmk}




Let $U$ be a good neighborhood of $Z(\xi)$ and $f$ be the corresponding function. We may assume that $f_{|\ZZ(\xi)}=0$.

\begin{defn}\label{stable and unstable sets}

The set of stable points $U^{st}\subset U$ is defined as $U^{st}=\{x\in U|e^{t\xi}(x)\in U, \forall t\ge 0\}=\{x\in U|lim_{t\to +\infty}f(e^{t\xi}(x))=0\}$.

Similarly, the set of unstable points $U^{unst}\subset U$ is defined as  
$U^{unst}=\{x\in U|e^{t\xi}(x)\in U, \forall t\le 0\}=\{x\in U|lim_{t\to -\infty}f(e^{t\xi}(x))=0\}$.
\end{defn}

Then $U^{st}$ and $U^{unst}$ are closed subsets of $U$ and $U^{st}\cap U^{unst}=Z(\xi)$. Moreover under the diffeomorphism $\psi$ from Remark \ref{product structure} the stable subset $U^{st}$ is identified with $\ZZ(\xi)\sqcup ((U^{st}\cap\partial_-U)\times [-1,0))$, and $U^{st}\cap \partial_-U$ is a closed subset of $int(\partial_-U)$.
Similar facts hold for $U^{unst}$ if we replace $\partial_-U$ by $\partial_+U$ and  the interval  $[-1,0)$ by $(0,1]$.

\begin{rmk}\label{thimbles in real case}
If all zeros are simple (i.e. locally are critical points of a Morse function) then both stable and unstable sets are topological cells.

\end{rmk}

Under some additional assumption on the vector field $\xi$  the embeddings $\ZZ(\xi)\subset U^{st}$ and $U^{unst}\supset \ZZ(\xi)$ become  homotopy equivalences.  Without any additional assumptions one can prove the following weaker statement, which is sufficient for our purposes.

\begin{prp}\label{isomorphism of cohomology pairs}
Let $\rho$ be a locally constant sheaf on $U$.

Then we have a natural isomorphisms of cohomology  
$$H^\bullet(U,\partial_-U,\rho)\to H^\bullet(U^{st}, U^{st}\cap \partial_-U,\rho).$$ 

Similarly we have a natural isomorphism 
$$H^\bullet(U,\partial_+U,\rho)\to H^\bullet(U^{unst}, U^{unst}\cap \partial_+U,\rho).$$

\end{prp}

{\it Proof.} We will prove the result for $U^{st}$. The case of $U^{unst}$ is similar.
The excision isomorphism ensures that $H^\bullet(U^{st}, U^{st}\cap \partial_-U,\rho)\simeq H^\bullet(U^{st}\cup \partial_-U, \partial_-U,\rho)$. Therefore we would like to prove that  $H^\bullet(U^{st}\cup \partial_-U, \partial_-U,\rho)\simeq H^\bullet(U,\partial_-U,\rho)$.

Consider the exact triangle 
$$H^\bullet(U,U^{st}\cup \partial_-U,\rho)\to  H^\bullet(U, \partial_-U,\rho)\to H^\bullet(U^{st}\cup \partial_-U,\partial_-U, \rho).$$
In order to establish the Proposition it suffices to prove the following Lemma.

\begin{lmm} \label{triviality of relative cohomology} $H^\bullet(U,U^{st}\cup \partial_-U, \rho)=0$.

\end{lmm}
{\it Proof.} Notice that $H^\bullet(U,U^{st}\cup \partial_-U, \rho)\simeq H_\bullet^{BM}(U-(U^{st}\cup \partial_-U),\rho\otimes or_Y)$, where $H_\bullet^{BM}$ is the notation for the Borel-Moore homology, and $or_Y$ denote the orientation bundle for $Y$. If we  can show  that $U-(U^{st}\cup \partial_-U)$ is homeomorphic to $\partial_+U\times [0,1)$, then using the K\"unneth formula we conclude that $H_\bullet^{BM}(U-(U^{st}\cup \partial_-U),\rho\otimes or_Y)=0$ since the Borel-Moore homology of the interval $[0,1)$ is trivial.
 
In order to prove the homeomorphism of the desired topological spaces let us observe that there exists
 a smooth function $\phi$ on $U$ such that $\phi\ge 0, \phi_{|U^{st}\cup \partial_-U}=0$, and $\phi>0$ on $U-(U^{st}\cup \partial_-U)$. Then $\partial_+U\times \R_{\le 0}$ is homeomorphic to $U-(U^{st}\cup \partial_-U)$ via the map $(x,t)\mapsto (e^{t\phi\xi}(x))$. This completes the proof.
$\blacksquare$



\begin{prp}\label{independence of good neighborhood}
The cohomology groups $H^\bullet(U,\partial_-U,\rho)$ do not depend on the good neighborhood $U$ of $\ZZ(\xi)$.

\end{prp}

{\it Proof.} By Proposition \ref{isomorphism of cohomology  pairs} it suffices to prove the result for the cohomology groups $H^\bullet(U^{st}, U^{st}\cap \partial_-U,\rho)$.  Notice that $H^\bullet(U^{st}, U^{st}\cap \partial_-U,\rho)\simeq H^\bullet(U^{st}, U^{st}-\ZZ(\xi),\rho)$. Indeed, the diffeomorphism $\psi$ from Remark \ref{product structure} identifies $U^{st}-\ZZ(\xi)$ with $(U^{st}\cap \partial_-U)\times [0,1)$. Finally we observe that the germ of $U^{st}$ near $\ZZ(\xi)$ does not depend on the choice of good neighborhood. $\blacksquare$


\begin{defn}\label{local Betti in real case}
The  local Betti cohomology of the pair $(Y,\xi)$ with coefficients in the locally constant sheaf $\rho$ is defined as 
$H^\bullet_{loc}(Y, \xi, \rho)=H^\bullet(U,\partial_-U,\rho)$ for any good neighborhood $U$ of $\ZZ(\xi)$.

The local Betti homology groups $H_{loc,\bullet}(Y, \xi, \rho)$ are defined similarly.

\end{defn}

\begin{rmk}\label{microlocal approach}


The Proposition \ref{isomorphism of cohomology pairs} can be approached via the microlocal theory of sheaves (see [KasSch4]). Indeed, it is a statement of the type of ``propagation of singularities", which stays in the origin of  the  microlocal  approach to the  Morse theory.

\end{rmk}

Main question of the next subsection is the relation between local and global Betti cohomology. As we will see, the expected global-to-local isomorphism holds in the non-archimedean framework (where it is highly non-trivial). Differently from the case of holomorphic functions, we could not prove the isomorphism over $\C$.

\subsection{Betti global-to-local isomorphism in non-archimedean framework}\label{non-archimedean Betti local to global}

In this subsection we will formulate a version of the global-to-local isomorphism of Betti cohomology with coefficients in a local system over a non-archimedean field. Since we will work with homology groups, it will be rather local-to-global isomorphism, dual to the one we are interested in.

Assume that $Y$ is a smooth compact manifold with corners and  $\xi$  a  vector field on $Y$ such that 
$\partial Y=\partial_-Y \cup \partial_h Y\cup \partial_+Y$ as in Section \ref{Betti cohomology for real one-forms}. We also assume that $\xi$ is a gradient-like near $\ZZ(\xi)$. 

Let $K$ be a  field endowed with a non-trivial non-archimedean valuation $val: K^\times\to  \R$. The  valuation map induces a norm $|\bullet|$ on $K$, and we assume that $K$ is complete with respect to the norm (i.e. $K$ is a  non-archimedean field). Let $\rho$ be a rank one $K$-local system on $Y$. Such local systems bijectively correspond to elements of $H^1(Y,K^\times)$. Via the  map $Log|\bullet|: K^\times\to \R$ we obtain a class in $H^1(Y,\R)$ which we denote by $Log(\rho)$.  Assume that $Log(\rho)=[\beta]$ where $\beta$ is a closed $1$-form such that outside of the set $\ZZ(\xi)$ we have the pointwise inequality $\beta(\xi)>0$.

 Assume that $\rho$ is metrized, i.e. its fibers carry a non-archimedean norm which continuously varies with respect to the point of $Y$.
A metrization is the same as a choice of non-archimedean norm on the fibers of $\rho$ such that the norm $|1_x|$ of a locally constant section $x\mapsto 1_x\in \rho_x\simeq K$ is a $C^0$ function on $Y$. 
In case if this function is smooth or, even weaker, of class $C^1$ we denote  by  $\beta=d(Log(|1_x|))$   the corresponding closed $1$-form which represents the  class $Log(\rho)\in H^1(Y,\R)$.

\begin{defn}\label{absence of saddle connections}
We say that our gradient-like vector field $\xi$ satisfies the absence of saddle connections assumption if there is no non-constant smooth map $\phi:(\R_\tau,\partial/\partial \tau)\to (Y-\ZZ(\xi),\xi)$ such that for a  Riemannian metric on $Y$ we have $dist(\phi(\tau), \ZZ(\xi))\to 0$ as $\tau\to +\infty$ and $\tau\to -\infty$. (This property does not depend on a choice of metric).

\end{defn}


The term is motivated by the  theory of dynamical systems where it is used in the framework of holomorphic $1$-forms on complex curves. 

\begin{rmk}
Absence of saddle connections does not hold generically in $C^\infty$ setting. E.g. consider the standard $d$-dimension sphere $S^d, d\ge 1$ and the gradient field $grad\, h$ of the height function $h$. Then $h$ is Morse and it is easy to see that a small perturbation of $grad\, h$ has saddle connections.

The situation is different in the holomorphic setting. Let $X$ be  a complex manifold  endowed with a holomorphic closed $1$-form $\alpha$ such that the set of zeros $\ZZ(\alpha)$ is compact (and hence is a union of finitely many connected components). Chose a Hermitian metric on $X$. Let $\xi$ be the gradient-like vector field corresponding to the real closed $1$-form $Re(\alpha/t)$,where $t\in \C^\ast$ is a parameter.
Then the absence of saddle connections comes from cohomological reasons, since for generic $t$ we can ensure that the integrals of $\alpha/t$ over the homology classes $\gamma\in H_1(X,\ZZ(\alpha),\Z)$ are not real positive numbers.

We remark that we do not discuss here the real manifold with corners $Y$ which appeared in considerations above. In general it is a compact submanifold with corners in $X$ (considered as a real manifold) such that $\xi$ behaves near the boundary $\partial Y$ as described in Section \ref{Betti cohomology for real one-forms}.

\end{rmk}

We start with some preliminary results which will be used in the proof of the main theorem below.

Returning to the Section \ref{Betti cohomology for real one-forms} notice that under the assumptions of the loc.cit. one can ``glue'' to $Y$ the half-stripes $\partial_+Y\times [0,+\infty)$ to $\partial_+Y$ and  $\partial_-Y\times [-\infty,0]$ to $\partial_-Y$, so that the flow of $\xi$ extends indefinitely with respect to $t\in \R$ without reaching the boundary of $Y$. We denote the obtained manifold $Y\cup \partial_+Y\times [0,+\infty)\cup \partial_-Y\times [-\infty,0]$ with corners by $Y^{ext}$ (extended $Y$). Same can be done with a good neighborhood $U$. We may assume that the vector field $\xi$ is extended to the half-stripes.


Let now $(Y,\xi)$ be as at the beginning of this subsection, and $\beta$ a smooth closed $1$-form compatible with $\xi$ in the sense that $\beta(\xi)>0$ outside of $\ZZ(\xi)$.
Before proceeding further we remark that we can replace $Y$ by $Y^{ext}$  and extend $\xi$ and $\beta$ there. For $\xi$ it is explained above. For $\beta$ we notice that replacing $\beta$ by $\beta+dg$ for an appropriate smooth function $g$ near $\partial_-Y\sqcup \partial_+Y$ we can assume that $Lie_\xi(\beta)=0$ near $\partial_-Y\cup \partial_+Y$. This guarantees an extension of $\beta$ to $Y^{ext}$ in such a way that $\beta(\xi)>0$ and $Lie_\xi(\beta)=0$ on $Y^{ext}-Y$. In the discussion below we will use the notation like $e^{(a,b]\xi}(\ZZ)$ for the result of the application of the flow $e^{t\xi}, t\in (a,b]$ to the set $\ZZ$.

\begin{prp}\label{integral estimate}
Assume that $Y,\xi,\beta$ as above and $\xi$ satisfies the absence of saddle connections  assumption.  Let $U$ be a good neighborhood of $\ZZ(\xi)$.

Then for any $C>0$ there exists $T>0$ and good neighborhood $U_C\supset \ZZ(\xi)$  such that $U_C\subset U$, $\partial_\pm U_C=\partial_\pm U\cap U_C$ and
the flow on $Y^{ext}$ generated by $\xi$ satisfies the following properties:

1) The sets $e^{[-T,0)\xi}(\partial_-U_C), U_C, e^{(0,T]\xi}(\partial_+U_C)$ are disjoint, and the natural map of their disjoint union given by $(x,t)\mapsto e^{t\xi}(x)$ on the first and the last set, and the natural embedding on $U_C$ is an embedding to $Y^{ext}$ (this implies that the image is a manifold with corners).

2)  For any $x\in \partial_-U_C$ we have $\int_{e^{[-T,0]\xi}(x)}\beta\ge C$ and similarly for any $x\in \partial_+U_C$ we have $\int_{e^{[0,T]\xi}(x)}\beta\ge C$. 

\end{prp}

{\it Proof.}  We will explain the construction of $U_C$ and the property 2) in the case of $\partial_+U_C$. The rest of the proof is similar. Main ingredient of the proof consists of the following result.
\begin{lmm}\label{infinite integral} For any $x\in \partial_+U\cap U^{unst}$ one has $\int_{e^{[0,+\infty)\xi}(x)}\beta=+\infty$.

\end{lmm}

In order to prove Lemma we start with the observation that the ray $[0,+\infty)$ can be divided into the union of intervals $[t_0=0,t_1], [t_1,t_2],...$ such that the trajectory of $\xi$ belongs to $U$ for $t$ belonging to $[t_1,t_2], [t_3,t_4],...$ and to $Y^{ext}-int(U)$ 
for $t\in [t_0,t_{1}], [t_2,t_3],...$. Let us consider first the case when there are infinitely many  intervals. 
Since there are infinitely many corresponding pieces of trajectories of $\xi$ with endpoints on $\partial_-U$ and $\partial_+U$ which stay inside $U$ and since the integral of $\beta$ over each piece is bounded below by the same constant we conclude that the total integral is infinite.

If there are finitely many such intervals then there exists time $a$ such that for all $t>a$ the trajectory of $\xi$ stays either entirely in $U$ or entirely outside of $U$.
In the latter case using the fact that $\beta(\xi)>0$ outside of $\ZZ(\xi)$ and compactness of $Y-int(U)$ we conclude that there exists a positive constant $C_2$ such that $\beta(\xi)(x)\ge C_2$ for $x\in Y-int(U)$ and hence on $Y^{ext}-U$. It follows that $\int_{e^{(a,+\infty)\xi}(x)}\beta=+\infty$ for $x\in \partial_+U\cap U^{unst}$.

Assume the former, i.e. that for $t>a$ the trajectory stays in $U$. Then it must belong to $U^{st}$, since to $t\to +\infty$ the point of trajectory has to  approach $\ZZ(\xi)$. On the other hand, since the starting point $x$ belongs to $U^{unst}$ we can take its limit along the trajectory as $t\to -\infty$. In this way we obtain a saddle connection, i.e. a trajectory with endpoints at $\ZZ(\xi)$. This contradicts to the absence of saddle connections assumption. The Lemma is proved. $\blacksquare$

We proceed with the proof of the Proposition by defining a continuous (in fact smooth) strictly increasing in second variable function $F:\partial_+U\times [0,+\infty)\to [0,+\infty)$ such that $F(x,t)=\int_{e^{[0,t]\xi(x)}}\beta$. Then Lemma \ref{infinite integral} is equivalent to the claim that $lim_{t\to +\infty}F(x,t)=+\infty$ for any $x\in \partial_+U\cap U^{unst}$.

\begin{lmm}\label{local bound}
For any $C>0$ there exists an open subset $V_C\subset \partial_+U$ which contains $\partial_+U\cap U^{unst}$ and $T>0$ such that $F(x,T)\ge C$ for any $x\in V_C$.

\end{lmm}
Lemma follows from compactness of $\partial_+U\cap U^{unst}$.

Let us now choose a subset $V^\prime\subset V_C$, where $V^\prime$ 
is a compact domain with smooth boundary and $int(V^\prime)\supset \partial_+U\cap U^{unst}$. Such $V^\prime$ always exists since $\partial_+U\cap U^{unst}$ is a compact subset of a smooth manifold. Consider now the union of intervals in $\R_{\ge 0}$ for which the trajectory $e^{-t\xi}(x), x\in V^\prime-\partial_+U\cap U^{unst}$ belongs to $U$ 
as long as $t$ belongs to one of these intervals. Let $W$ denote the union of pieces of these trajectories. We set $U_C=U^{st}\cup U^{unst}\cup W$. Shrinking $U_C$ further we obtain  a good neighborhood of $\ZZ(\xi)$ satisfying both conditions 1) and 2) of the Proposition \ref{integral estimate}. This finishes the proof. $\blacksquare$

Here is the main result of this subsection.

\begin{thm-constr}\label{isomorphism for non-archimedean local systems}
Let $Y,\xi, \rho, \beta$ be as above, and $\beta$ is compatible with $\xi$. 
We also assume that $\xi$ satisfies the absence of saddle connections assumption

Then  we will define  a non-degenerate pairing
$$H_\bullet(Y,\partial_+Y,D(\rho))\otimes H_{loc,\bullet}(Y, \xi, \rho)\to K,$$
where $D(\rho)=\rho^\ast\otimes or_Y$ denote the dual local system, and $or_Y$ is the orientation local system of $Y$. It gives rise to an isomorphism $iso_{Y,\rho}: H_{loc,\bullet}(Y, \xi, \rho)\simeq H_\bullet(Y,\partial_-Y, \rho)$.

Furthermore, if the vector field $\xi$ and the norm $|\bullet|$ vary continuously with respect to a parameter belonging to a connected topological space, but
the  set $\ZZ(\xi)$  
 does not vary,  then the pairing does not change. 
\end{thm-constr}
{\it Sketch of the proof.} Full proof with all the details will take too much space, so we will indicate main steps only.

First we need to define a pairing between a chain representing  a class
in  $H_\bullet(Y,\partial_+Y, D(\rho))$ and a representative of another class in $H_{loc,\bullet}(Y, \xi, \rho)$.
This can be thought of as a version of the dual to {\it Betti global-to-local isomorphism} in the non-archimedean setting. We can replace $Y$ by $Y^{ext}$ since they are homotopy equivalent. The advantage of that is that for $Y^{ext}$ all  trajectories of $\xi$ do exist for all $t\in \R$, while for $Y$ they can reach the boundary $\partial Y$ in finite time.

Recall that local Betti homology can be defined as $H_\bullet(U,\partial_-U,\rho)$ where $U$ is any good neighborhood of $\ZZ(\xi)$.  
Let us fix such $U$ and construct a sequence of good neighborhoods $U=U_0\supset U_1\supset U_2\supset...$ as well as the sequence $T_1<T_2<T_3<..., T_i\to +\infty$ by induction, taking $C=1,2,3...$ in Proposition \ref{integral estimate}.  Without loss of generality we may assume  that $\partial_hU_n, n\ge 0$ are disjoint. 

Starting with a class $\alpha_0\in H_\bullet(U,\partial_-U,\rho)$ we can construct a sequence of classes $\alpha_n\in H_\bullet(U_n, \partial_-U_n,\rho),n\ge 1$ using the natural isomorphisms 

$$H_\bullet(U_n,\partial_-U_n,\rho)\simeq H_\bullet(U_{n+1}\cup \partial_-U_n,\partial_-U_n,\rho)\simeq H_\bullet(U_{n+1},\partial_-U_{n+1},\rho).$$
Here the first isomorphism $a^{(1)}_n$ comes from the observation that the one-parameter semigroup $g^t: x\mapsto e^{-t\xi}(x), t\ge 0$ contracts $U-U_n$ to $\partial_-U-\partial_-U_n, n\ge 1$, while the second isomorphism $a_n^{(2)}$ is just the excision isomomorphism.

Next we would like to find ``good" representatives of the classes $\alpha_n, n\ge 0$ by finite linear combinations of smooth simplicial chains with values in $\rho$ . Such a representative can be encoded as a pair $(\gamma_n,\hat{\gamma}_n)$ where $\gamma_n\in C_\bullet(U_n, \rho), \hat{\gamma}_n\in C_\bullet(\partial_-U_n,\rho)$ subject to the conditions $\partial\hat{\gamma}_n=0,\partial \gamma_n=\hat{\gamma}_n$. Here  we abuse the notation and identify $\hat{\gamma}_n$ with its image under the natural embedding $i_n:\partial_-U_n\to {U}_n$.
In order to save the notation we will do the same thing below with other chains as well.

Let us explain in which sense the representatives are good. First we define the norm $|s|_\infty$ of a linear combination $s$ of simplicial chains with values in $\rho$ as the supremum over all simplicies of the norms of the corresponding vectors in $\rho$. Then one can analyze geometrically the isomorphisms $a_n^{(1)}$ and $a_n^{(2)}$ and ensure that norms $|\gamma_n|_\infty$ are uniformly bounded from above.
More precisely it goes such as follows. Using the fact that $U_n$ and $U_{n+1}$ are manifolds with corners and the flow $g^t$ defines a free action of the semigroup $\R_{\ge 0}$ on $U_n-U_{n+1}$ one can construct by induction chains $\gamma_n, \hat{\gamma}_n, n\ge 0$ as above as well as the chains $\delta_{n,n+1}\in C_\bullet(U_n,\rho), \hat{\delta}_{n,n+1}\in C_\bullet(\partial_-U_n,\rho)$ such that:

a) ${\gamma}_n-{\gamma}_{n+1}=\partial \hat{\delta}_{n,n+1}$;

b) $\partial\delta_{n,n+1}=\gamma_n-\gamma_{n+1}+\hat{\delta}_{n,n+1}$;

c) all norms $|\gamma_n|_\infty, |\hat{\gamma}_n|_\infty, |\delta_{n,n+1}|_\infty, |\hat{\delta}_{n,n+1}|_\infty$ are uniformly bounded.


 Using simplicial decomposition of $Y$ by smoothly embedded simplices we can present a closed  chain in 
 $C_\bullet(Y, \partial_+Y, D(\rho))$ as a finite sum $B:=\sum_j\tau_j\otimes b_j$, where $\tau_j$
 is a simplex and $b_j$ a flat section of the trivial local system $\rho^\ast_{|\tau_j}$. Then the  norm $|B|_\infty$ is equal to the maximum of the pointwise norms $|b_j(x)|, x\in \tau_j$.  We want to define the pairing of the homology class $[B]$ with $\alpha_0\in H_{loc,\bullet}(Y,\xi,\rho)=H_\bullet(U,\partial_-U,\rho)$.

Next we replace $Y$ by $Y^{ext}$ and will apply the flow $g^t$  to the above-defined chains extending them beyond the good neighborhood $U$. Notice that this flow moves the boundary of the chains towards $\partial_-Y^{ext}=\partial_-Y\times \{-\infty\}$. More precisely, we proceed such as follows.
 
 Let $W_n=U_n\cup e^{[-T_n,0]\xi}(\partial_-U_n):=U_n\cup(\cup_{0\le t\le T_n} g^{t\xi}(\partial_-U_n))$. Then  $\partial_-W_n=e^{-T_n\xi}(\partial_-U_n)$. To each quadruple $\sigma_{n,n+1}:=(\gamma_n,\hat{\gamma}_n, \delta_{n,n+1}, \hat{\delta}_{n,n+1})$ as above we associate a quadruple 
 $\sigma_{n,n+1}^{ext}:=(\gamma_n^{ext},\hat{\gamma}_n^{ext}, \delta_{n,n+1}^{ext}, \hat{\delta}_{n,n+1}^{ext})$ such that
 
 $$\gamma_n^{ext}=\gamma_n+e^{[-T_n,0]}(\hat{\gamma}_n)+g^{T_n}(\hat{\delta}_{n,n+1})\in C_\bullet(W_n,\rho)\subset C_\bullet(Y^{ext},\rho),$$
 $$\hat{\gamma}_n^{ext}=g^{T_n}(\hat{\gamma}_n), \hat{\delta}_{n,n+1}^{ext}=g^{T_n}(\hat{\delta}_{n,n+1}), \delta_{n,n+1}^{ext}=e^{[-T_n,0]}(\delta_{n,n+1}).$$
 
 The pair $(\gamma_n^{ext},\hat{\gamma}_n^{ext})$ is a closed relative chain in $C_\bullet(W_n, \partial_-W_n, \rho)$.


\vspace{3mm}

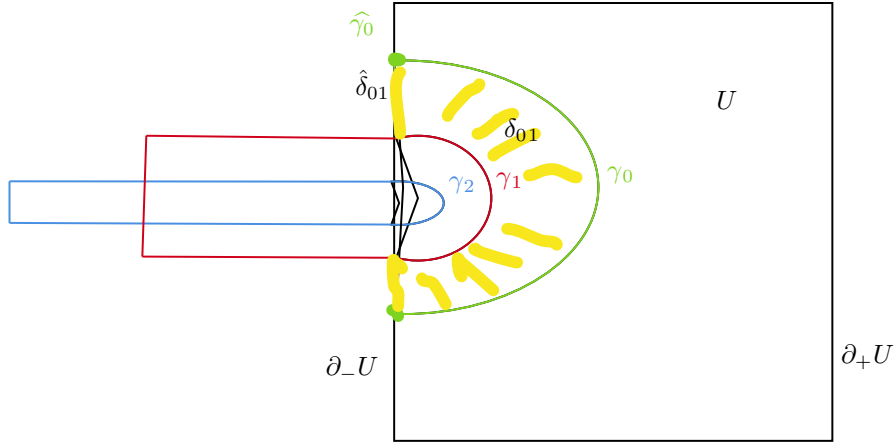
\tikzset{every picture/.style={line width=0.75pt}} 

\begin{tikzpicture}[x=0.75pt,y=0.75pt,yscale=-1,xscale=1]

\draw   (375,46) -- (596,46) -- (596,267) -- (375,267) -- cycle ;
\draw  [draw opacity=0] (376.05,174.6) .. controls (379.51,175.51) and (383.19,176) .. (387,176) .. controls (407.43,176) and (424,161.9) .. (424,144.5) .. controls (424,127.1) and (407.43,113) .. (387,113) .. controls (383.19,113) and (379.51,113.49) .. (376.05,114.4) -- (387,144.5) -- cycle ; \draw  [color={rgb, 255:red, 208; green, 2; blue, 27 }  ,draw opacity=1 ] (376.05,174.6) .. controls (379.51,175.51) and (383.19,176) .. (387,176) .. controls (407.43,176) and (424,161.9) .. (424,144.5) .. controls (424,127.1) and (407.43,113) .. (387,113) .. controls (383.19,113) and (379.51,113.49) .. (376.05,114.4) ;  
\draw [color={rgb, 255:red, 208; green, 2; blue, 27 }  ,draw opacity=1 ]   (376.05,114.4) -- (250,113) ;
\draw [color={rgb, 255:red, 208; green, 2; blue, 27 }  ,draw opacity=1 ]   (376.05,174.6) -- (248,174) ;
\draw  [draw opacity=0] (373.56,157.83) .. controls (374.84,157.94) and (376.16,158) .. (377.5,158) .. controls (389.93,158) and (400,153.08) .. (400,147) .. controls (400,140.92) and (389.93,136) .. (377.5,136) .. controls (376.16,136) and (374.84,136.06) .. (373.56,136.17) -- (377.5,147) -- cycle ; \draw  [color={rgb, 255:red, 74; green, 144; blue, 226 }  ,draw opacity=1 ] (373.56,157.83) .. controls (374.84,157.94) and (376.16,158) .. (377.5,158) .. controls (389.93,158) and (400,153.08) .. (400,147) .. controls (400,140.92) and (389.93,136) .. (377.5,136) .. controls (376.16,136) and (374.84,136.06) .. (373.56,136.17) ;  
\draw [color={rgb, 255:red, 74; green, 144; blue, 226 }  ,draw opacity=1 ]   (373.56,136.17) -- (181,136) ;
\draw [color={rgb, 255:red, 74; green, 144; blue, 226 }  ,draw opacity=1 ]   (373.56,157.83) -- (181,157) ;
\draw [color={rgb, 255:red, 74; green, 144; blue, 226 }  ,draw opacity=1 ]   (181,136) -- (181,157) ;
\draw [color={rgb, 255:red, 208; green, 2; blue, 27 }  ,draw opacity=1 ]   (250,113) -- (248,174) ;
\draw  [draw opacity=0] (375.98,202.96) .. controls (377.15,202.99) and (378.32,203) .. (379.5,203) .. controls (433.9,203) and (478,174.35) .. (478,139) .. controls (478,103.65) and (433.9,75) .. (379.5,75) .. controls (377.98,75) and (376.47,75.02) .. (374.97,75.07) -- (379.5,139) -- cycle ; \draw  [color={rgb, 255:red, 126; green, 211; blue, 33 }  ,draw opacity=1 ] (375.98,202.96) .. controls (377.15,202.99) and (378.32,203) .. (379.5,203) .. controls (433.9,203) and (478,174.35) .. (478,139) .. controls (478,103.65) and (433.9,75) .. (379.5,75) .. controls (377.98,75) and (376.47,75.02) .. (374.97,75.07) ;  
\draw [color={rgb, 255:red, 126; green, 211; blue, 33 }  ,draw opacity=1 ]   (374.97,75.07) -- (376.05,114.4) ;
\draw [color={rgb, 255:red, 126; green, 211; blue, 33 }  ,draw opacity=1 ]   (375.98,202.96) -- (376.05,174.6) ;
\draw  [color={rgb, 255:red, 126; green, 211; blue, 33 }  ,draw opacity=1 ][line width=4.5] [line join = round][line cap = round] (375,74) .. controls (376.42,74) and (381,75) .. (375,75) ;
\draw  [color={rgb, 255:red, 126; green, 211; blue, 33 }  ,draw opacity=1 ][line width=4.5] [line join = round][line cap = round] (377,204) .. controls (377,202.59) and (375.41,201) .. (374,201) ;
\draw  [color={rgb, 255:red, 248; green, 231; blue, 28 }  ,draw opacity=1 ][line width=4.5] [line join = round][line cap = round] (427,105) .. controls (427,105) and (427,105) .. (427,105) ;
\draw  [color={rgb, 255:red, 248; green, 231; blue, 28 }  ,draw opacity=1 ][line width=4.5] [line join = round][line cap = round] (435,102) .. controls (427.97,102) and (421.62,112) .. (417,112) ;
\draw  [color={rgb, 255:red, 248; green, 231; blue, 28 }  ,draw opacity=1 ][line width=4.5] [line join = round][line cap = round] (446,112) .. controls (444.05,112) and (426.9,121.1) .. (426,122) .. controls (425.67,122.33) and (425,122.53) .. (425,123) ;
\draw  [color={rgb, 255:red, 248; green, 231; blue, 28 }  ,draw opacity=1 ][line width=4.5] [line join = round][line cap = round] (418,87) .. controls (418,88.7) and (413.13,89.87) .. (412,91) .. controls (410.4,92.6) and (402,100.95) .. (402,103) ;
\draw  [color={rgb, 255:red, 248; green, 231; blue, 28 }  ,draw opacity=1 ][line width=4.5] [line join = round][line cap = round] (457,168) .. controls (457,165.98) and (451.62,166.41) .. (450,166) .. controls (444.46,164.62) and (437.86,160) .. (432,160) ;
\draw  [color={rgb, 255:red, 248; green, 231; blue, 28 }  ,draw opacity=1 ][line width=4.5] [line join = round][line cap = round] (437,177) .. controls (435.15,177) and (410.95,170) .. (416,170) ;
\draw  [color={rgb, 255:red, 248; green, 231; blue, 28 }  ,draw opacity=1 ][line width=4.5] [line join = round][line cap = round] (425,191) .. controls (419.67,186.33) and (414.47,181.51) .. (409,177) .. controls (408.27,176.4) and (407,174.06) .. (407,175) .. controls (407,177.42) and (407.24,182.24) .. (409,184) ;
\draw  [color={rgb, 255:red, 248; green, 231; blue, 28 }  ,draw opacity=1 ][line width=4.5] [line join = round][line cap = round] (401,198) .. controls (399.89,195.78) and (395.58,187.79) .. (394,187) .. controls (392.67,186.33) and (389,186.89) .. (389,185) ;
\draw  [color={rgb, 255:red, 248; green, 231; blue, 28 }  ,draw opacity=1 ][line width=4.5] [line join = round][line cap = round] (467,134) .. controls (464.07,134) and (460.06,130.58) .. (456,130) .. controls (449.82,129.12) and (446.47,133) .. (443,133) ;
\draw  [color={rgb, 255:red, 248; green, 231; blue, 28 }  ,draw opacity=1 ][line width=4.5] [line join = round][line cap = round] (378,81) .. controls (373.24,85.76) and (378,104.03) .. (378,112) ;
\draw  [color={rgb, 255:red, 248; green, 231; blue, 28 }  ,draw opacity=1 ][line width=4.5] [line join = round][line cap = round] (379,180) .. controls (377.67,177.35) and (376.09,178.18) .. (375,176) .. controls (373.72,173.45) and (374.22,185.09) .. (375,189) .. controls (375.24,190.18) and (376.83,190.81) .. (377,192) .. controls (377.33,194.31) and (377,196.67) .. (377,199) ;

\draw (536,89) node [anchor=north west][inner sep=0.75pt]   [align=left] {$\displaystyle U$};
\draw (599,217) node [anchor=north west][inner sep=0.75pt]   [align=left] {$\displaystyle \partial _{+} U$};
\draw (339,223) node [anchor=north west][inner sep=0.75pt]   [align=left] {$\displaystyle \partial _{-} U$};
\draw (429,103) node [anchor=north west][inner sep=0.75pt]   [align=left] {$\displaystyle \delta _{0}{}_{1}$};
\draw (354,78) node [anchor=north west][inner sep=0.75pt]   [align=left] {$\displaystyle \hat{\delta }_{01}$};
\draw (425,130) node [anchor=north west][inner sep=0.75pt]   [align=left] {$\displaystyle \textcolor[rgb]{0.82,0.01,0.11}{\gamma }\textcolor[rgb]{0.82,0.01,0.11}{_{1}}$};
\draw (401,132) node [anchor=north west][inner sep=0.75pt]   [align=left] {$\displaystyle \textcolor[rgb]{0.29,0.56,0.89}{\gamma }_{\textcolor[rgb]{0.29,0.56,0.89}{2}}$};
\draw (481,127) node [anchor=north west][inner sep=0.75pt]   [align=left] {$\displaystyle \textcolor[rgb]{0.49,0.83,0.13}{\gamma }_{\textcolor[rgb]{0.49,0.83,0.13}{0}}$};
\draw (351,49) node [anchor=north west][inner sep=0.75pt]   [align=left] {$\displaystyle \textcolor[rgb]{0.49,0.83,0.13}{\widehat{\textcolor[rgb]{0.49,0.83,0.13}{\gamma _{\textcolor[rgb]{0.49,0.83,0.13}{0}}}}}$};

\end{tikzpicture}

\vspace{3mm}

By Proposition \ref{integral estimate} one has $|\hat{\gamma}_n^{ext}|_\infty\le C_1e^{-n}$ for some $C_1>0$, since the integral of $\beta$ along a trajectory $e^{[-T_n,0]}(x), x\in \partial_-U_n$ is  bounded below by $n$ uniformly in $x$. It suffices to define the pairing of $B$ with $\gamma_n^{ext}$ and check that it has a well-defined limit as $n\to \infty$.

If we choose $\gamma_n, n\ge 1$ to be sufficiently general then the intersections of $\gamma_n^{ext}$ and  $B$ are transversal for all $n\ge 1$. Furthermore, the support of $\gamma_n^{ext}$ is separated from $\partial_+Y$. Thus we have a well-defined finite intersection index $\gamma_n^{ext}\cap B\in K$.
 


The intersection
$\gamma_n^{ext}\cap B$ has a limit as  $n\to +\infty$ for the following reasons. First we observe that
$$\gamma_n^{ext}-\gamma_{n+1}^{ext}=-e^{[-T_{n+1},-T_n]\xi}(\hat{\gamma}_{n+1})-e^{-T_n\xi}\hat{\delta}_{n,n+1}+\partial \delta_{n,n+1}.$$
Notice that $B\cap \partial \delta_{n,n+1}=0$ since the boundary of $B$ belongs to $\partial_+Y$.
Then we have the equality $\gamma_{n}^{ext}\cap B=\gamma_{n+1}^{ext} \cap B$ modulo an element which has small norm in $K$. More precisely, we observe that $|\gamma_n^{ext}\cap B-\gamma_{n+1}^{ext}\cap B|$ is bounded from above by the norm of the intersection of the RHS with $B$. But the latter is bounded from above by $const\cdot e^{-n}$. 

We conclude that $lim_{n\to \infty}\gamma_n^{ext}\cap B$ exists. One can show by similar considerations that the limit does not depend on a choice of representatives for the homology classes. This is the desired pairing $[\gamma_0]\cap [B]=\alpha_0\cap [B]$.



In order to prove that the pairing  is non-degenerate we use the ideas of [HL].  Let us first assume  that $\partial_{h}Y=\emptyset$. We will return to the general case later.

Then we have a well-defined homology class $[\Delta_{Y,\rho}]:=[\Delta_Y,\rho\boxtimes D(\rho)]\in H_\bullet(Y\times Y, (\partial_-Y\times Y)\sqcup (Y\times \partial_+Y),\rho\boxtimes D(\rho)$ where the class of the diagonal  represents inverse to the natural non-degenerate pairing $H_\bullet(Y,\partial_-Y,\rho)\otimes H_\bullet(Y,\partial_+ Y, D(\rho))\to K$.\footnote{The following observation will be useful below. Let $V$ and $W$ be objects of a symmetric monoidal category with the unit object ${\bf 1}$. The element $R\in V\otimes W$ defines an inverse to the pairing $S\in Hom(V\otimes W,{\bf 1})$ iff the composition $V\to V\otimes W\otimes V\to V$ is equal to  $id_V$, where the first arrow is $v\mapsto R\otimes v$ and the second one is $v_1\otimes w\otimes v_2\mapsto v_1\otimes S(v_2\otimes w)\in V\otimes {\bf 1}=V$, and the similar statement holds if we interchange $V$ and $W$. We will be interested in this observation when $R$ will be the cohomology class of the (shifted) diagonal and $S$ will be the $K$-valued intersection of cohomology classes.}

In order to proceed we need the following construction.  Let $U$ be a a good neighborhood of $\ZZ(\xi)$. There exists a smooth function $\phi: U\to \R_{\ge 0}$ such that $\xi(\phi)=0$ and $\phi=const>0$ near $\partial_hU$. We keep the notation for the natural extension of $\phi$ to $U^{ext}$.

We define the modified diagonal
$\Delta_U^{\phi}\in C_\bullet(U\times U,\rho\boxtimes D(\rho))$ such as follows. First we consider the map $U^{ext}\to U^{ext}\times U^{ext}$ defined by
$$x\mapsto (x,e^{-\phi(x)\xi}(x), x,e^{\phi(x)\xi}(x)).$$
Applying this map to $\Delta_{Y,\rho}$ we obtain a chain in $C_\bullet(U^{ext}\times U^{ext}, \rho\boxtimes D(\rho))$.

Second, let us choose the $const$ in the definition of $\phi$ sufficiently large. We will call such functions $\phi$ sufficiently large. Then applying to the above chain the map  $C_\bullet(U^{ext}\times U^{ext}, \rho\boxtimes D(\rho))\to  C_\bullet(U\times U, \rho\boxtimes D(\rho))$ coming from the natural retraction $U^{ext}\times U^{ext}\to U\times U$ we obtain the chain with the boundary in $\partial_-U\times U\cup U\times \partial_+U$. We denote this chain by $\Delta^{\phi}_{U,\rho}$. The cohomology class $[\Delta_{Y,\rho}^{\phi}]\in H_\bullet(U,\partial_-U,\rho)\otimes H_\bullet(U,\partial_+U,D(\rho))$ gives the inverse to the non-degenerate Poincar\'e pairing

$$H_\bullet(U,\partial_-U,\rho)\otimes H_\bullet(U,\partial_+U,D(\rho))\to K.$$

Let us apply this construction to the sequence of good neighborhoods  $U_n$ constructed above. We denote by $\phi_n$ the corresponding sufficiently large functions. We denote by $\hat{\phi}_n$ the functions on $W_n$ obtained by application to $\phi_n$ the flow of $\xi$. Then we  obtain a sequence of chains $\Delta_{W_n,\rho}^{\hat{\phi}_n}$. 

Let us  apply to the diagonal $\Delta_{Y,\rho}$ the flow $e^{t\xi}, t\in \R$. The diagonal shifted by this flow  consists of points $t:=(x,e^{t\xi}(x))$, where $x\in Y$ (technically speaking the point $e^{t\xi}(x)$  belongs to $Y^{ext}$).
Then we have $|(\Delta_{W_n,\rho}^{\hat{\phi}_n}-e^{-T_n\xi}(\Delta_{Y,\rho}))\cap B|\to 0$ as $n\to \infty$ for any closed chain $B\in C_\bullet(Y\times Y, (\partial_+Y\times Y)\cup (Y\times \partial_-Y), D(\rho)\boxtimes \rho)$. Hence  $\Delta_{W_n,\rho}^{\hat{\phi}_n}$ approximates a representative for the dual of the non-degenerate pairing. This proves the result in case when $\partial_{h}Y=\emptyset$.

In the case when $\partial_{h}Y\ne \emptyset$ we use the non-degenerate pairing $H_\bullet(Y,\partial_-Y,\rho)\otimes H_\bullet(Y,\partial_+Y\cup \partial_{h}Y,D(\rho))\to K$ together with the observation  that $H_\bullet(Y,\partial_+Y\cup \partial_{h}Y,D(\rho))\simeq H_\bullet(Y,\partial_+Y,D(\rho))$. This follows from the long exact sequence and the fact that $H_\bullet(\partial_{h}Y\cap \partial_+Y,D(\rho))=0$. In order to show this equality we use the argument similar to the above one applied to the shifted diagonal class $\Delta_{\partial_{h}Y,\rho}$ as well as the fact that $\xi$ nowhere vanishes on $\partial_{h}Y$.
This concludes our sketch of the proof of the Theorem. $\blacksquare$

\begin{rmk}\label{soft sheaf}
a) The following notion is useful if we cannot represent $Log(\rho)$ by a smooth $1$-form.
Let us define a (soft)  sheaf  $\FF_\xi$ with local sections being $C^0$ real-valued functions $g$ on $Y$ modulo constants such that $\xi(g)$ considered as a distribution is in fact a $C^0$  function.
 Then each element of $H^1(Y,\R)$ has a representative in $\Gamma(Y,\FF_\xi)$.  Furthermore, Theorem \ref{isomorphism for non-archimedean local systems} can be proved, if we replace $\beta$ by a global section $g$ of $\FF_\xi$ such that $\xi(g)>0$ outside of $\ZZ(\xi)$.
 
The advantage of the use of $\FF_\xi$ is based on the observation that the cohomology classes represented by $g\in \Gamma(Y,\FF_\xi)$ as above form an {\bf open} convex cone  $Cone_\xi\subset H^1(Y,\R)$. 

b) Cohomology $H^\bullet_{loc}(Y,\rho)$ and $H^\bullet(Y,\partial_-Y,\rho)$ as we vary $\rho$ give rise to $\Z$-graded algebraic coherent sheaves $\EE_{loc,(Y,\xi)}$
and $\EE_{glob,(Y,\xi)}$ respectively on the stack of ${\bf G}_m$-local systems on $Y$. Slightly modifying the proof of the Theorem \ref{isomorphism for non-archimedean local systems} one can show that $iso_{Y,\rho}$ can be extended to an isomorphism $iso_{Y}^{an}:\EE_{loc,Y}^{an}\simeq \EE_{glob,Y}^{an}$ of analytifications of these coherent sheaves over the tube domain $U_\xi$ corresponding to $Cone_\xi$. In particular  Theorem \ref{isomorphism for non-archimedean local systems} holds for all $\rho\in U_\xi$. Imposing some extra restrictions one can generalize it to higher rank local systems $\rho$. We will not do that in this paper.

Moreover, a particular choice of the non-archimedean field of $K$ is irrelevant. One can define $iso_{Y}^{an}$ at the level of adic spaces over 
$\Z$  in the sense of Huber, see [Hub] (in fact better to speak about adic stacks) . Even more generally, it can be extended to an isomorphism of the correspondent objects of triangulated dg-categories of coherent sheaves.

In the case of a pair $(X,\alpha)$ of a complex manifold $X$ endowed with a holomorphic $1$-form $\alpha$ we will use the notation $\EE_{loc,(X,\alpha)}$
and $\EE_{glob,(X,\alpha)}$, etc.  for the corresponding graded non-archimedean coherent sheaves.

\end{rmk}

\subsection{Comparison of real and holomorphic  cases}\label{real and complex cases}

Suppose $X,\alpha$ are the same as in Section \ref{de Rham local to global for 1-forms}.
Recall that in the holomorphic setting of Section \ref{Betti for 1-forms} we defined the global Betti cohomology $H^\bullet_{Betti,glob,t}(X,\alpha)=H^{\bullet}(\overline{X}_{cor}, D_{v,t}^{\R,-}\cup D_{log,t}^{\R,+}, E_{\alpha,t})$. In the RHS we can replace the local system $E_{\alpha,t}$ associated with the holomorphic $1$-form $\alpha/t$ by {\it any locally constant sheaf  $\rho$ of rank $1$} over an arbitrary field. In fact we will be interested in the non-archimedean local systems which belongs to the tube domain $U(C_{\alpha/t})$ defined in  Remark \ref{covering by two intersecting arcs}. We denote the corresponding cohomology groups by $H^\bullet_{Betti,glob,t}(X,\rho)$. Notice that  the exceptional divisor of the real blow-up depends on the parameter $t$, but the local system $\rho$ does not.
For  a general Hermitian metric on $X$ we denote by $v_{\alpha,t}$ the gradient vector field corresponding to the $1$-form $Re(\alpha/t)$. 
\vspace{2mm}

{\bf Assumption}

\vspace{2mm}
{\it 1) In the above notation there exist a compact smooth manifold with corners $Y=Y_{\alpha,t}\subset X$ and a Hermitian metric $h$ in its neighborhood such that the pair $(Y_{\alpha,t}, v_{\alpha,t})$ satisfies the properties 1)-4) from Section \ref{Betti cohomology for real one-forms}.  Moreover the pair $(Y_{\alpha,t}, \partial_-Y_{\alpha,t})$ is homotopy equivalent to the pair $(\overline{X}_{cor},  D_{v,t}^{\R,-}\cup D_{log,t}^{\R,+})$.  

2) More precisely for a decomposition of $\partial\overline{X}_{cor}=\overline{X}_{cor}-X$ into the union
of closed subsets $\cup_{\tau\in \{+,h,-\}}\partial_\tau  \overline{X}_{cor}$ defined below, there exist smooth non-negative functions $f_\tau, \tau\in  \{+,h,-\}$ in  a neighborhood of $\partial_\tau  \overline{X}_{cor}$ as well as sufficiently small positive numbers  $\epsilon_\tau, \tau\in \{+,h,-\}$ such that $f_\tau^{-1}(0)=\partial_\tau \overline{X}_{cor}, Y=\overline{X}_{cor}-\cup_{\tau\in \{+,h,-\}}\{x|f_\tau(x)<\epsilon_\tau\}$. Furthermore the decomposition $\partial Y=\cup_{\tau\in \{+,h,-\}}\partial_\tau Y$ is  obtained when we set $\partial_\tau Y=f_\tau^{-1}(\epsilon_\tau)$ with $\epsilon_\pm$  are much smaller than $\epsilon_h$.}

\begin{conj}\label{comparison of real and complex cases}
The above Assumption  always holds.

\end{conj}

Let us explain the definition of $\partial_\tau\overline{X}_{cor}$ where $\tau\in \{+,h,-\}$.\footnote{The reader can skip technical details below and simply accept the above Assumption when needed.}
For that we introduce four types of coordinates in a neighborhood of any point of $\overline{X}-X=D_v\cup D_h\cup D_{log}$. More precisely we denote by $z_i, z_i', z_i''$ the local coordinates at points of $D_v, D_{log}, D_h$ respectively, and by $z_i'''$ the remaining coordinates. We will not specify the sets of indices for each group of coordinates, if it does not lead to a confusion.

Now we define:

a) $\partial_-\overline{X}_{cor}=\partial_{v,-}\overline{X}_{cor}\cup \partial_{log,-}\overline{X}_{cor}$, where $\partial_{v,-}\overline{X}_{cor}$ is defined by the conditions $\#\{z_i\}\ge 1, \sum_ik_iArg(z_i)\in [-\pi/2,\pi/2]$, where $k_i\ge 1$ are the integers from the Proposition \ref{compactification for 1-forms} and $\partial_{log,-}\overline{X}_{cor}$ is defined by the conditions $\#\{z_i'\}\ge 1,$ and there exists $z_j'$ such that $Re(c_j)>0$ in the notation of the Proposition \ref{compactification for 1-forms}.

b) $\partial_+\overline{X}_{cor}=\partial_{v,+}\overline{X}_{cor}\cup \partial_{log,+}\overline{X}_{cor}$, where 
$\partial_{v,+}\overline{X}_{cor}$ is defined by the conditions $\#\{z_i\}\ge 1, -\sum_ik_iArg(z_i)\in [-\pi/2,\pi/2]$ with the same conventions as in a).
Similarly, $\partial_{log,+}\overline{X}_{cor}$ is defined by the conditions $\#\{z_i'\}\ge 1,$ and there exists $z_j'$ such that $Re(c_j)<0$.

c) The most complicated part is the definition of $\partial_h\overline{X}_{cor}$. It is covered by six closed subsets:

c1) The closure of the open stratum $\partial_{h,h}\overline{X}_{cor}$ is defined by the condition $\#\{z_i''\}\ge 1$.

c2) The stratum $\partial_{h,v}\overline{X}_{cor}$ corresponding to the intersection of horizontal and vertical divisors is defined by the condition   $\#\{z_i\}\ge 1, \sum_ik_iArg(z_i)\in \{-\pi/2,\pi/2\}$, i.e. $\prod_iz_i^{-k_i}$ is pure imaginary.

The intersection of the horizontal and logarithmic divisors gives rise to two strata:

c3) $\partial_{h,log}^{(1)}\overline{X}_{cor}$ consisting of points for which $\#\{z_i'\}\ge 1$ and there exists $z_j'$ such that $Re(c_j)=0$;

c4) $\partial_{h,log}^{(2)}\overline{X}_{cor}$ consisting of points for which $\#\{z_i'\}\ge 1$ and there exists $z_{j_1}'$ and $z_{j_2}'$ such that $Re(c_{j_1})<0$ and $Re(c_{j_2})>0$.

Intersection of horizontal, vertical and logarithmic divisors gives rise to the following two strata:

c5) $\partial_{h,v, log}^{(1)}\overline{X}_{cor}$ consisting of points for which  $\#\{z_i'\}\ge 1$ and  $\#\{z_i\}\ge 1$ and such that $\sum_ik_iArg(z_i)\in (-\pi/2,\pi/2)$ and there exists $z_j'$ such that $Re(c_j)<0$.

c6) $\partial_{h,v, log}^{(2)}\overline{X}_{cor}$ consisting of points for which  $\#\{z_i'\}\ge 1$ and  $\#\{z_i\}\ge 1$ and such that $-\sum_ik_iArg(z_i)\in (-\pi/2,\pi/2)$ and there exists $z_j'$ such that $Re(c_j)>0$.

Notice that $\partial\overline{X}_{cor}-\partial_h\overline{X}_{cor}=U_+\sqcup U_-$, where $U_\pm$ are open subsets of $\partial\overline{X}_{cor}$. Here $U_\pm=U_\pm^{(1)}\cup U_\pm^{(2)}$. The set $U_-^{(1)}$ consists of points such that  $\#\{z_i'\}=0, \#\{z_i''\}=0, \#\{z_i\}\ge 1$ satisfying the following conditions: $-\sum_ik_iArg(z_i)\in (-\pi/2,\pi/2)$ and for any $z_j'$ we have $Re(c_j)>0$.  The set $U_-^{(2)}$ consists of points such that  $\#\{z_i\}=0,  \#\{z_i''\}=0, \#\{z_i'\}\ge 1$ and for all $z_j'$ we have $Re(c_j)>0$.

The open set $U_+$ is defined similarly.  One can see that $\partial_\pm \overline{X}_{cor}=\overline{U}_\pm$. Moreover, the embeddings $\overline{U}_\pm\subset \overline{U}_\pm$ are homotopy equivalences.  In particular $\partial_\pm \overline{X}_{cor}$ are topological manifolds with corners.

The complement to the union of slightly smaller open subsets $U_+^0\subset U_+$ and $U_-^0\subset U_-$ is also a manifold with corners which is homotopy equivalent to $\partial_h\overline{X}_{cor}$, while the sets themselves will be homotopy equivalent to $\partial_- \overline{X}_{cor}$ or to $\partial_+ \overline{X}_{cor}$. 

Having the above-defined sets one {\it can hope} to define the distance functions $f_\tau, \tau\in \{+,h,-\}$ as well as the rest of the data from the above Assumption.

\begin{rmk}
Conjecture \ref{comparison of real and complex cases} holds in many special cases, e.g. in the case when $X$ is a complex curve, which we will consider below.
\end{rmk}

Since the discussion of all details will substantially increase the size of the paper, we will use the above Assumption leaving the proof of the Conjecture \ref{comparison of real and complex cases} to the interested reader. Then using the previously established  relation between local and global Betti cohomology in the non-archimedean setting  we deduce from the Conjecture \ref{comparison of real and complex cases}  that for any locally constant sheaf $\rho$ of rank $1$ over a non-archimedean field we have two isomorphisms of global and local Betti cohomology groups:
$$H_{Betti,glob,t}(X,\rho)\simeq H^\bullet(Y_{\alpha, t},\partial_-Y_{\alpha,t},\rho)$$
and

$$H^\bullet_{Betti,loc,t}(X,\rho)\simeq H^\bullet(\ZZ(\alpha),\phi_{W/t}(\rho))\simeq H^\bullet_{Betti,loc}(Y_{\alpha,t},\partial_-Y_{\alpha,t}, v_{\alpha,t},\rho),$$
where $dW=\alpha$ near $\ZZ(\alpha)$ and $W_{|\ZZ(\alpha)}=0$ (see Section \ref{Betti for 1-forms}).

\subsection{Betti global-to-local isomorphism in complex case}\label{Stokes rays for 1-forms}

Let us  discuss now the notion of Stokes rays in the case of holomorphic closed $1$-forms.
It is similar to the one in the case of functions.

\begin{defn}\label{definition of Stokes rays for 1-forms}
Stokes ray for the $1$-form $\alpha$ is a ray in $\R^2$ with the vertex at the origin (a.k.a. admissible ray) and such that the angle with the positive horizontal axis is equal to
 $Arg(Z(\gamma_{ij})):=Arg(\int_{\gamma_{ij}}\alpha)$, where $\gamma_{ij}$ is the isotopy class in $X$ of a path joining two points in the connected components $\ZZ_i(\alpha)$ and $\ZZ_j(\alpha)$ (we can have $i=j$).

\end{defn}

Let us call the angle $\theta$  {\it generic} if $\R_{\ge 0}\cdot e^{i\theta}=\R_{\ge 0}\cdot t$ is not a Stokes ray. Hence genericity of $\theta$ means that $\theta$ different from any of $Arg(Z(\gamma_{ij}))$. Thus the Betti cohomology and the corresponding $\C$-vector space  
$H^{\bullet}(U_{\varepsilon,j},U_{\varepsilon,j}\cap W_j^{-1}(\varepsilon\cdot e^{i\theta}),\Z)\otimes \C$ do not depend on $\theta$ as long as it is generic (see Section \ref{Betti for 1-forms}).
Now we can reformulate the non-archimedean Betti local-to global isomorphism for complex manifolds. 

\begin{prp-constr}\label{non-archimedean Betti local to global for holomorphic 1-forms}
Under the Assumption from Section \ref{real and complex cases}
let $\rho$ be a local system over a non-archimedean field $K$ such that $\rho$ belongs to the tube domain $U(Cone_{\alpha/t})$ from Remark \ref{covering by two intersecting arcs}. Then if $t$ does not belong to a Stokes ray, we have a well-defined isomorphism
$iso_t:=iso_{l_t}: H^\bullet_{Betti,loc,t}(X,\rho)\simeq H^\bullet_{Betti,glob,t}(X,\rho)$. Here $l_t=\R_{>0}\cdot t$.

\end{prp-constr}
{\it Proof.} The proof follows from the combination of the Theorem \ref{isomorphism for non-archimedean local systems}, comparison of local and global Betti cohomology in real and complex cases in Section \ref{real and complex cases}, and from the absence of saddle connections for $v_{\alpha,t}$ as long as $t$  does not belong to a Stokes ray. $\blacksquare$

We will  say  that $t\in \C^\ast$ is generic (or $\theta=Arg(t)$ is generic) if the corresponding $l_t$ is not a Stokes ray.

\begin{rmk}\label{extension to coherent sheaves}
The isomorphism $iso_t$ extends to the one of graded non-archimedean coherent sheaves on the tube domain $U(Cone_{\alpha/t})$ (see Remark \ref{soft sheaf}).

\end{rmk}

 
 For a Stokes ray $l_t=\R_{>0}\cdot t$ there are  two isomorphisms $iso_{t,\pm}$ obtained as limits in the non-archimedean sense of $iso_{t^\prime}$ for a non-Stokes ray $\R_{>0}\cdot t^\prime$ as $t^\prime$ is generic and approaches $l_t$ from the right for $+$ and  from the left for $-$. Then by the Remark \ref{extension to coherent sheaves} the automorphism $R_{l_t}:=iso_{t,-}^{-1}\circ iso_{t,+}$ gives rise to an automorphism of the graded non-archimedean analytic coherent sheaf $\EE^{an}_{loc, (X,\alpha/t)}$ over the tube domain $U(Cone_{\alpha/t})$ (see Remark \ref{soft sheaf}). 
 \begin{defn}\label{Rays for 1-forms}
 We call $R_{l_t}$ the {\it Stokes automorphism corresponding to the Stokes ray $l_t$}.
 \end{defn}
 Thus $R_{l_t}: H^\bullet_{Betti,loc,t}(X,\rho)\to H^\bullet_{Betti,loc,t}(X,\rho)$ is an isomorphism for any $\rho\in U(Cone_{\alpha/t})$.


\begin{rmk} \label{universal integer local system} Assume that the restriction of $\rho$ to the set of zeros $\ZZ(\alpha)$ is trivialized. Since rank $\rho$ is equal to $1$ we have $$H^\bullet_{Betti,loc,t}(X,\rho)\simeq H^\bullet_{Betti,loc,t}(X,\Z)\otimes K.$$
Then in the Proposition \ref{non-archimedean Betti local to global for holomorphic 1-forms} we can replace the local system over a non-archimedean field of arbitrary characteristic by the universal local system over the ring  $\Z[[C\cap \Gamma]]$ which consists of formal power series $\sum_{\gamma\in C\cap \Gamma}n_\gamma x^\gamma$. Here $\Gamma=H_1(X,\ZZ(\alpha), \Z)$ is a finitely generated abelian group, and $C\subset \Gamma\otimes \R$ is a rational strict convex cone, and $n_\gamma\in \Z$. Then we can speak about the universal Stokes automorphisms $R_{l_t}^{univ}\in Aut(H^\bullet_{Betti,loc,t}(X,\Z))[[\Gamma\cap \C]]$.
If we choose a $\Z$-basis in $H^\bullet_{Betti,loc,t}(X,\Z)/torsion$ then  $R_{l_t}^{univ}$ gives rise to an integer matrix consisting of integer block-diagonal matrices of the size $rk\,H^i_{Betti,loc,t}(X,\Z)\times rk\,H^j_{Betti,loc,t}(X,\Z), i,j\in I$.  Changing scalars to $\Q$ we can interpret $R_{l_t}$ as an element of $GL(k, \Q[[C\cap \Gamma]])$, where $k=rk\,H^\bullet_{Betti,loc,t}(X,\Z)$.

\end{rmk}




\subsection{Wall-crossing structure for $1$-forms}\label{WCS for 1-forms}

For simplicity of the exposition {\it we will restrict ourselves to the case of rank $1$ local systems trivialized at $\ZZ(\alpha)$}. Such local systems form a smooth scheme isomorphic to the torus $Hom(H_1(X, \ZZ(\alpha), \Z), {\bf G}_m)$ rather than a stack. Moreover the sheaf $\EE_{loc,(X,\alpha)}$ (see Remark \ref{soft sheaf}) becomes a vector bundle on this torus.

Let $\ZZ(\alpha)=\sqcup_{j\in J}\ZZ_j(\alpha)$ be a decomposition of the set of zeros into the union of connected components. In other words $J=\pi_0(\ZZ(\alpha))$.

Then the  WCS (we  refer to  [KoSo7] for the definition) on $S^1_\theta$ (equivalently on $\C^\ast_t$) is defined such as follows:

a) Constant local system of lattices $\Gamma_\theta=\Gamma=H_1(X, \ZZ(\alpha),\Z)/tors$;\footnote{One can allow torsion in the free abelian group which appears in the definition of the wall-crossing structure in [KoSo1] and [KoSo7]. In particular we can take above as $\Gamma$ the abelian group $H_1(X, \ZZ(\alpha),\Z)$. It will make the definition of our WCS slightly more complicated. For that reason we have factorized by the torsion.} The central charge $Z:\Gamma\to \C$ is given by $\gamma\mapsto \int_\gamma\alpha$.

b) Local system of $\Gamma_\theta$-graded Lie algebras $\g_\theta=\oplus_{(i,j)\in J\times J}V_{ij,\theta}\otimes W_{ij,\theta}$, where 
$$V_{ij,\theta}=V_{ij,\theta}^\Z\otimes \Q,$$ 
and $V_{ij,\theta}^\Z=Hom_{\Z-graded}(H^\bullet_{Betti, i,\theta}(X,{\Z}), H^\bullet_{Betti, j,\theta}(X,{\Z}))$. Similarly
$W_{ij,\theta}=W_{ij,\theta}^\Z\otimes \Q$, where $W_{ij,\theta}^Z=\Z^{\{\gamma\in \Gamma_\theta |\partial \gamma=[p_i]-[p_j], p_i\in \ZZ_i(\alpha), p_j\in \ZZ_j(\alpha)\}}$. 
Here 
$$\partial: H_1(X, \ZZ(\alpha), \Z)\to H_0(\ZZ(\alpha),\Z)\simeq \Z^{J}$$ 
is the boundary map. We denote by $\g_{\theta}^\Z\subset \g_\theta$ the lattice defined by $\g_\theta^\Z:=\oplus_{(i,j)\in \pi_0(Z(\alpha))\times \pi_0(Z(\alpha))}V_{ij,\theta}^\Z\otimes W_{ij,\theta}^\Z$.

The set of generators of the $\Q$-vector space $W_{ij,\theta}$ is a torsor over $H_1(X,\Z)$. Fixing $\gamma$ with $\partial \gamma=[p_i]-[p_j]$ we obtain the $\Gamma_\theta$-graded component $\g_{\theta,\gamma}\subset \g_\theta$, which is isomorphic over $\Q$ to $V_{ij,\theta}$. Notice that $\g_\theta$ is in fact an associative algebra, but we use the corresponding Lie algebra structure only. 

c)  The Stokes automorphisms $R_{l_t}^{univ}:=R_\theta^{univ}, \theta=Arg(t)$ make sense and are defined over $\Q$ (see Remark \ref{universal integer local system}). Then we have a WCS on $S^1_\theta$ (equivalently on $\C^\ast_t$). Let us denote it by $\sigma_{Betti}$.

\begin{rmk}\label{family of WCS}
Let us choose a sufficiently small tubular neighborhood $U_\epsilon$ of $Z(\alpha)$ and a holomorphic function $W:U_\epsilon\to \{|z|<\epsilon\}$ such that $\alpha=dW$ on $U_\epsilon$. Suppose that we have a continuous family of pairs $(X,\alpha)$ such that the homotopy type of the pair $(U_\epsilon, W^{-1}(\epsilon))$ does not change. Then we have a continuous family of the wall-crossing structures defined above. It follows that if we know the WCS for one value of the parameter space then we can compute it for any other value. Indeed, all what we use for that are the wall-crossing formulas.

\end{rmk}

We remark that if we consider a continuous family of pairs $(X,\alpha)$ such that $\pi_0(\ZZ(\alpha))$ is locally constant\footnote{In this case the cohomology with coefficients in the sheaf of vanishing cycles also stays locally constant.} one can show that the corresponding wall-crossing structures form a continuous family of WCS. The main idea is to consider the corresponding holomorphic function on the universal abelian cover of $X$.   In that case the continuity is clear from the results of the previous Section, where the Betti global-to-local isomorphism was given explicitly.

\subsection{Wall-crossing structures and their analyticity: revision}\label{analyticity of WCS}

Notice that for fixed $\theta$ the Lie algebra $\g_\theta$ is the Lie subalgebra of the Lie algebra of matrix-valued functions on the torus ${\bf T}_\Gamma=Hom(\Gamma, {\bf G}_m)$. In particular it can be embedded into the Lie algebra of vector fields on a bigger torus of dimension $rk\, H^\bullet(X, \ZZ(\alpha), \Z)+\sum_i rk\, H^\bullet_{Betti, i,\theta}(X,\Z)$. The embedding depends on the choice of a basis in local cohomology. The basis can have non-trivial monodromy, hence it cannot be chosen globally in a continuous way. Nevertheless one can slightly modify 
 the Definition 3.4.1 from our paper [KoSo12] and obtain a more general notion of analytic WCS which is good for the purposes of current paper. Main results of the loc. cit. hold in the more general framework including the open-closed property (Propositions 3.5.1, 3.5.2 from the loc.cit.). When we will use this more general set-up we will often give references to the corresponding less general results from [KoSo12].

\begin{thm}\label{analyticity of WCS for 1-forms} Under the Assumption from Section \ref{real and complex cases} and after extending the scalars from $\Q$ to $\C$ the above-constructed WCS $\sigma_{Betti}$ is analytic in the sense of [KoSo12].

\end{thm}

Let us discuss the idea of the proof of the Theorem \ref{analyticity of WCS for 1-forms} postponing the actual proof until Section \ref{rationality of WCS}. Although  full definition of the analytic WCS is a bit involved (see [KoSo12], Section 3.11), in our case the local system of lattices is constant and the dependence of stability data on $\g_\theta$ on the point of $S^1_\theta$ is such that it suffices to establish analyticity of the stability data on the graded Lie algebra $\g_\theta$, where $\theta$ is fixed. 

Recall that according to [KoSo12] if $Z(\Gamma)\subset \Z\oplus i\Z$ (such central charge is called {\it rational} in the loc.cit.) then for {\it some} finite cyclic cover $\R^2=\C=\cup_{i\in I}V_i$ by strict {\it rational} semiclosed sectors all transformations $A_{V_i}$ (see loc. cit.) are convergent in some {\it holomorphic} tube domains. We proved in the loc.cit. that then $A_V$ is analytic (i.e. given by a convergent series) for {\it any} strict {\it rational} sector $V$ and moreover analyticity in the rational case is in fact {\it equivalent} to the analyticity of all transformations $A_{V_i}$ for a {\it fixed} rational cyclic cover or for an {\it arbitrary} rational cyclic cover (see e.g.  Proposition 3.4.2 and Remark 3.4.3 in the loc.cit.).

We will show below that for a small deformation of the central charge to a rational one the transformation $A_V$ for a rational $V$ is given by a {\bf matrix-valued rational function}. This implies the Theorem \ref{analyticity of WCS for 1-forms}.

\begin{rmk}

On the other hand we do not know whether $A_V$ is analytic  in case if either $V$ is non-rational or $Z$ is non-rational. In this case the question of analyticity of $A_V$ seems to be a difficult problem related to the problem of approximation of a transcendental number by rational ones. We do not know if $A_V$ is analytic for all $V$ even in our case of matrix-valued rational functions.
\end{rmk} 

We will illustrate the Theorem \ref{analyticity of WCS for 1-forms} in a couple of non-trivial examples in Section \ref{WCS for 1-form on curve}.  
First we will need to revisit the notion of stability data and analytic stability data. Before doing that in the following subsections let us make the a comment which will be used later on.

\begin{rmk}\label{remark about rationality of WCS}
Let us fix two sets of complex numbers: $z_i, 1\le i\le k$ and $s_j, 1\le j\le l$ as well as the collection of finite-dimensional complex vector spaces $V_i, 1\le i\le k$.
Define the lattice $\Gamma=\Gamma_0\oplus \Z^l$, where $\Gamma_0\subset \Z^k$ is the kernel of the map $(n_1,...,n_k)\mapsto n_1+...+n_k$. Consider the Lie algebra $\g=End(V, \C[q_1^{\pm 1},...,q_l^{\pm 1}])$, where $V=\oplus_{1\le i\le k}V_i$.
The grading on $\g=\oplus_{\gamma\in \Gamma}\g_\gamma$ is defined in such a way that the component $V_{i_1}^\ast\otimes V_{i_2}\otimes \C\cdot q_1^{m_1}...q_l^{m_l}$ has the grading $\gamma=(a_i,b_j)$, where $a_i=\delta_{ii_1}-\delta_{ii_2}, b_j=k_j$. One can define the central charge by the formula $Z((a_1,...,a_k), (b_1,...,b_l))=\sum_ia_iz_i+\sum_jb_js_j$, but we will not use it.

Let us use the notation and terminology from [KoSo12]. With an admissible cone $C\subset \Gamma\otimes \R$ we associate a pronilpotent Lie algebra $\g_C$ as well as the corresponding pronilpotent group $G_C=exp(\g_C)$. If a rational hyperplane defines a disjoint decomposition $C=C_-\sqcup C_0\sqcup C_+$ then we have the factorization $G_C=G_{C_-}G_{C_0}G_{C_+}$.  It is shown in the loc. cit. that analytic elements in $G_C$ are exactly those which are products of analytic elements in the above factorization formula. It is exactly this property allowed us to prove that the set of analytic WCS is open and closed in the topological space of all WCS. Similarly to the loc.cit. one can prove the same result for our Lie algebra $\g$ provided we replace analyticity by {\bf rationality} with respect to the variables $q_j, 1\le j\le l$. Of course this {\bf rationality property} of WCS implies analyticity.

\end{rmk}

\subsubsection{Stability data and hyperbolic hyperplanes}

To simplify the exposition we will discuss the simplest case of the wall-crossing structure, namely the case of stability data on a fixed graded Lie algebra. Generalization to the case when we have a {\it trivial} local system of lattices and {\it non-trivial} local system of graded Lie algebras is straightforward. This will be sufficient for our main case of holomorphic $1$-forms.

We will assume that our stability data has rank $2$, i.e. $Z_\R$ is a surjective map $\Gamma_\R\to \C=\R^2$. The rank $1$ case is easy (see [KoSo12]).
By the Support Property there exists a
quadratic form $Q$  on $\Gamma_\R\simeq \R^n$ of signature $(2,n-2)$ such that $Q$ is negative on $Ker\, Z_\R$ and positive on the support of the wall-crossing structure.

Having all this in mind, let us start with some linear algebra data, which we get for free from stability data of rank $2$ on a graded Lie algebra. Namely, let $\Gamma$ be a free abelian group of rank $n$ and $\Gamma_\R=\Gamma\otimes \R\simeq \R^n$. Assume that $\Gamma_\R$ is endowed with a non-degenerate quadratic form of signature $(2, n-2)$. We will also assume a choice of locally-constant orientation on all $2$-dimensional subspaces which have signature $(2,0)$ with respect to the quadratic form $Q$.

\begin{defn}\label{hyperbolic hyperplanes} A hyperbolic (or $Q$-hyperbolic if we want to stress the dependence on $Q$)  hyperplane is a cooriented hyperplane $H\subset \Gamma_\R$ such that the restriction $Q_{|H}$  has signature $(1,n-2)$.
\end{defn}

The reader should remember that coorientation is a part of the structure, although this is not reflected in the notation.

\begin{rmk}\label{inductive limit over quadratic forms}
We can get rid of the dependence of the quadratic form by taking the inductive limit over all of those $Q$ which satisfy the above property, as we did in [KoSo1]. For that reason we will often skip from the notation  the quadratic form $Q$ and  speak about hyperbolic hyperplanes rather than $Q$-hyperbolic. 
\end{rmk}

We  denote by $\UU$ the set of pairs hyperbolic hyperplanes. Equivalently $\UU$ can be described as a set of such closed half-spaces $P=P_H\subset \Gamma_\R=\R^n$ that the restriction of $Q$ on the boundary of $H=\partial P$ has signature $(1,n-2)$. The set $\UU$ carries a natural Hausdorff topology.

 Notice that the point $u\in \UU$ gives rise to a ray $l=l_u$ in the open domain $\{Q>0\}:=\{x\in \Gamma_\R|Q(x)>0\}$. The ray $l$ can be thought of as the only ray in the half-space orthogonal to its boundary and pointing out in the direction where $Q>0$.
This allows us to define a binary relation on the points of $\UU$ which is described below. Informally it can be thought of as a partial clockwise order,
similar to the relation on $S^1=\R/2\pi \Z$ given by $\{(\theta_1, \theta_2)\in \R^2|\theta_1-\theta_2\in [0,\pi)\}$.

Before giving the definition let us explain the underlying geometry.
Let $H$ be a hyperbolic hyperplane and $W\subset H$ an $(n-2)$-dimensional subspace such that $Q_{|W}$ has signature $(0,n-2)$ (i.e. restriction of $Q$ on $W$ is negative and $W$ is generic hyperplane in $H$ with respect to this property). We have $\Gamma_\R=W\oplus W^\perp$ and  $H=W\oplus c_W$, where $c_W\subset W^\perp\simeq \R^2$ is a cooriented line in the oriented $Q$-orthogonal complement $W^\perp\simeq \R^2$.

\begin{defn}\label{admissible pairs} A pair $(H_1,H_2), H_i\in \UU, i=1,2$ is called admissible if either $H_1=H_2$ or the restriction of $Q$ to $W:=H_1\cap H_2$ has signature $(0,n-2)$, and $H_2$ is obtained from $H_1$ by the clockwise rotation in the plane $W^\perp\simeq \R^2$ by an angle $0<\theta<\pi$.
We write $H_1\le H_2$ for an arbitrary admissible pair and $H_1<H_2$ in case if $H_1\ne H_2$.

\end{defn}

If $(H_1,H_2)$ is an admissible pair then for the corresponding half-spaces $P_{H_i}, i=1,2$ we construct a strict convex cone $C_{H_1,H_2}$ defined as  the convex hall of the points of the nonconvex cone $C_{H_1,H_2}^{ncv}:=(P_{H_2}-P_{H_1})\cap \{Q>0\}$.

\vspace{3mm}

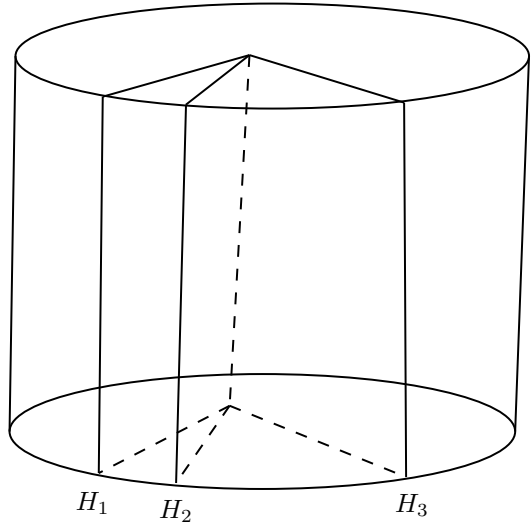
\tikzset{every picture/.style={line width=0.75pt}} 

\begin{tikzpicture}[x=0.75pt,y=0.75pt,yscale=-1,xscale=1]

\draw   (167,49.5) .. controls (167,34.86) and (224.98,23) .. (296.5,23) .. controls (368.02,23) and (426,34.86) .. (426,49.5) .. controls (426,64.14) and (368.02,76) .. (296.5,76) .. controls (224.98,76) and (167,64.14) .. (167,49.5) -- cycle ;
\draw   (164,239) .. controls (164,222.98) and (221.08,210) .. (291.5,210) .. controls (361.92,210) and (419,222.98) .. (419,239) .. controls (419,255.02) and (361.92,268) .. (291.5,268) .. controls (221.08,268) and (164,255.02) .. (164,239) -- cycle ;
\draw    (167,49.5) -- (164,239) ;
\draw    (426,49.5) -- (419,239) ;
\draw    (211,70) -- (209,260) ;
\draw    (363,73) -- (364,262) ;
\draw    (285,49) -- (211,70) ;
\draw    (285,49) -- (363,73) ;
\draw  [dash pattern={on 4.5pt off 4.5pt}]  (275,226) -- (209,260) ;
\draw  [dash pattern={on 4.5pt off 4.5pt}]  (285,49) -- (275,226) ;
\draw  [dash pattern={on 4.5pt off 4.5pt}]  (275,226) -- (364,262) ;
\draw    (285,49) -- (253,74) ;
\draw    (253,74) -- (248,265) ;
\draw  [dash pattern={on 4.5pt off 4.5pt}]  (275,226) -- (248,265) ;

\draw (187,267) node [anchor=north west][inner sep=0.75pt]   [align=left] {$ $};
\draw (238,272) node [anchor=north west][inner sep=0.75pt]   [align=left] {$\displaystyle H_{2}$};
\draw (196,268) node [anchor=north west][inner sep=0.75pt]   [align=left] {$\displaystyle H_{1}$};
\draw (357,269) node [anchor=north west][inner sep=0.75pt]   [align=left] {$\displaystyle H_{3}$};

\end{tikzpicture}

\vspace{3mm}

We can extend the notion of admissibility of a pair of cooriented hyperbolic hyperplanes to a structure of  simplicial set. Namely, $0$-simplices are elements of $\UU:=\UU^{(0)}$, the set of $1$-simplices $\UU^{(1)}\subset \UU\times \UU$ are admissible pairs $(H_1,H_2), H_1\le H_2$.  Furthermore $2$-simplices are  {\it admissible triples} $(H_1,H_2,H_3)\in \UU^{(2)}\subset \UU^3$ i.e. such triples  that $(H_1,H_2), (H_2,H_3), (H_1,H_3)$ are admissible pairs and $C^{ncv}_{H_1,H_3}=C^{ncv}_{H_1,H_2}\sqcup C^{ncv}_{H_2,H_3}$. Similarly one defines  $k$-simplices. We will use only the case when $k\le 2$.

Suppose that we are given a $\Gamma$-graded Lie algebra $\g=\oplus_{\gamma\in \Gamma}\g_\gamma$ over $\Q$. Then to each  cone $C_{H_1,H_2}$ we can assign a pronilpotent group $G_{C_{H_1,H_2}}$ (see [KoSo1], [KoSo12] for more details). This assignment satisfies the property that for an admissible triple
$(H_1,H_2,H_3)$ one has $G_{C_{H_1,H_2}}\subset G_{C_{H_1,H_3}}\supset G_{C_{H_2,H_3}}$. 

Suppose that to each cone $C_{H_1,H_2}$ we have assigned an element $A_{H_1,H_2}\in G_{C_{H_1,H_2}}$, and for an admissible triple 
$(H_1,H_2,H_3)$ this assignment satisfies the cocycle property 
$$A_{H_1,H_2}A_{H_2,H_3}=A_{H_1,H_3}.$$

Summarizing,  to the  graded Lie algebra $\g$ and the quadratic form $Q$ we associate the data consisting of the collections of group elements $A_{H_1,H_2}\in G_{C_{H_1,H_2}}, H_i\in \UU, i=1,2$ satisfying the above cocycle conditions. Let us call these data a {\it hyperbolic stability data on $\g$}. 

Hyperbolic stability data would give an ordinary stability data on $\g$ in the sense of [KoSo1], [KoSo12] in case when we also had  a central charge. In fact there is an open set of  central charges compatible with $Q$. More precisely, we will consider central charges $Z: \Gamma\to \C$ such that $Q_{|Ker\, Z_\R}$ has signature $(0,n-2)$ and the above-mentioned orientation on the $Q$-orthogonal complement $(Ker\, Z_\R)^\perp\simeq \R^2$ coincides with the standard orientation on $C\simeq \R^2$.
Such central charge $Z$ is called {\it compatible with $Q$}.  
Let $l\subset \R^2$ be a ray with the vertex at the point $(0,0)$ (such rays are called admissible in [KoSo12]). Then $H(l):=Z^{-1}_\R(l\cup -l)$ is a hyperbolic hyperplane. The coorientation is inherited from the clockwise orientation of $\R^2$. We call such hyperplanes {\it vertical} because they contain $Ker\,Z_\R$.

Now we can say that hyperbolic stability data give rise to stability data  on $\g$ with any central charge $Z$ which is compatible with $Q$. Moreover the support of this stability data belongs to the set $\{Q\ge 0\}$. 

Thus we have a map which to every admissible in the sense of [KoSo12] semiclosed sector $V\subset \R^2$ bounded by the rays $l_1<l_2$ (clockwise order)  associates a strict convex cone $C(V)$ which is bounded by the intersection $H({l_1})\cap H(l_2)\cap \{Q>0\}$. The transformation $A_V$ from [KoSo1] is defined as the one  associated with $C(V)$. In the above notation it coincides with    $A_{H(l_1), H(l_2)}$.

\vspace{3mm}

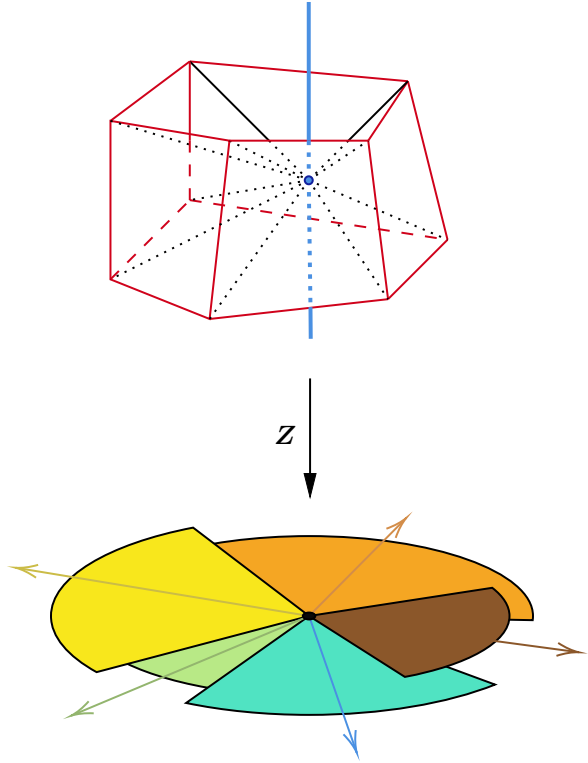
\tikzset{every picture/.style={line width=0.75pt}} 

\begin{tikzpicture}[x=0.75pt,y=0.75pt,yscale=-1,xscale=1]

\draw  [fill={rgb, 255:red, 184; green, 233; blue, 134 }  ,fill opacity=0.39 ] (239.34,369.74) .. controls (196.74,366.95) and (162.42,355.78) .. (150.51,341.22) -- (260.16,329.91) -- cycle ;
\draw  [fill={rgb, 255:red, 245; green, 166; blue, 35 }  ,fill opacity=0.23 ] (187.33,299.09) .. controls (206.99,293.16) and (232.41,289.58) .. (260.16,289.58) .. controls (322.53,289.58) and (373.1,307.64) .. (373.1,329.91) .. controls (373.1,330.62) and (373.04,331.32) .. (372.94,332.02) -- (260.16,329.91) -- cycle ;
\draw  [fill={rgb, 255:red, 248; green, 231; blue, 28 }  ,fill opacity=0.14 ] (152.99,358.25) .. controls (138.49,350.2) and (130,340.44) .. (130,329.91) .. controls (130,310.41) and (159.16,293.52) .. (201.68,285.31) -- (260.16,329.91) -- cycle ;
\draw  [color={rgb, 255:red, 208; green, 2; blue, 27 }  ,draw opacity=1 ] (160,80) -- (220,90) -- (210,180) -- (160,160) -- cycle ;
\draw [color={rgb, 255:red, 208; green, 2; blue, 27 }  ,draw opacity=1 ]   (210,180) -- (300,170) ;
\draw [color={rgb, 255:red, 208; green, 2; blue, 27 }  ,draw opacity=1 ]   (290,90) -- (300,170) ;
\draw [color={rgb, 255:red, 208; green, 2; blue, 27 }  ,draw opacity=1 ]   (290,90) -- (310,60) ;
\draw [color={rgb, 255:red, 208; green, 2; blue, 27 }  ,draw opacity=1 ]   (330,140) -- (300,170) ;
\draw [color={rgb, 255:red, 208; green, 2; blue, 27 }  ,draw opacity=1 ]   (310,60) -- (330,140) ;
\draw [color={rgb, 255:red, 208; green, 2; blue, 27 }  ,draw opacity=1 ]   (160,80) -- (200,50) ;
\draw [color={rgb, 255:red, 208; green, 2; blue, 27 }  ,draw opacity=1 ]   (200,50) -- (310,60) ;
\draw    (200,50) -- (240,90) ;
\draw    (310,60) -- (280,90) ;
\draw  [dash pattern={on 0.84pt off 2.51pt}]  (240,90) -- (260,110) ;
\draw  [dash pattern={on 0.84pt off 2.51pt}]  (280,90) -- (260,110) ;
\draw  [dash pattern={on 0.84pt off 2.51pt}]  (260,110) -- (160,160) ;
\draw  [dash pattern={on 0.84pt off 2.51pt}]  (160,80) -- (260,110) ;
\draw  [dash pattern={on 0.84pt off 2.51pt}]  (260,110) -- (330,140) ;
\draw  [dash pattern={on 0.84pt off 2.51pt}]  (260,110) -- (300,170) ;
\draw [color={rgb, 255:red, 208; green, 2; blue, 27 }  ,draw opacity=1 ]   (200,50) -- (200,90) ;
\draw [color={rgb, 255:red, 208; green, 2; blue, 27 }  ,draw opacity=1 ] [dash pattern={on 4.5pt off 4.5pt}]  (200,90) -- (200,120) ;
\draw [color={rgb, 255:red, 208; green, 2; blue, 27 }  ,draw opacity=1 ] [dash pattern={on 4.5pt off 4.5pt}]  (160,160) -- (200,120) ;
\draw [color={rgb, 255:red, 208; green, 2; blue, 27 }  ,draw opacity=1 ] [dash pattern={on 4.5pt off 4.5pt}]  (200,120) -- (330,140) ;
\draw  [dash pattern={on 0.84pt off 2.51pt}]  (260,110) -- (290,90) ;
\draw  [dash pattern={on 0.84pt off 2.51pt}]  (260,110) -- (200,120) ;
\draw  [dash pattern={on 0.84pt off 2.51pt}]  (260,110) -- (210,180) ;
\draw  [dash pattern={on 0.84pt off 2.51pt}]  (220,90) -- (260,110) ;
\draw  [fill={rgb, 255:red, 80; green, 227; blue, 194 }  ,fill opacity=0.09 ] (353.94,364.53) .. controls (330.26,373.96) and (297,379.82) .. (260.16,379.82) .. controls (237.71,379.82) and (216.58,377.64) .. (198.14,373.81) -- (260.16,329.91) -- cycle ;
\draw  [fill={rgb, 255:red, 139; green, 87; blue, 42 }  ,fill opacity=0.22 ] (352.54,315.68) .. controls (358.14,320.03) and (361.25,324.84) .. (361.25,329.91) .. controls (361.25,343.18) and (339.92,354.73) .. (308.48,360.66) -- (260.16,329.91) -- cycle ;
\draw    (260.67,210) -- (260.67,268) ;
\draw [shift={(260.67,270)}, rotate = 270] [fill={rgb, 255:red, 0; green, 0; blue, 0 }  ][line width=0.08]  [draw opacity=0] (12,-3) -- (0,0) -- (12,3) -- cycle    ;
\draw [color={rgb, 255:red, 208; green, 2; blue, 27 }  ,draw opacity=1 ]   (220,90) -- (290,90) ;
\draw [color={rgb, 255:red, 74; green, 144; blue, 226 }  ,draw opacity=1 ][line width=1.5]    (260,20) -- (260,90) ;
\draw [color={rgb, 255:red, 74; green, 144; blue, 226 }  ,draw opacity=1 ][line width=1.5]    (260.87,174.01) -- (261,190) ;
\draw [color={rgb, 255:red, 74; green, 144; blue, 226 }  ,draw opacity=1 ][line width=1.5]  [dash pattern={on 1.69pt off 2.76pt}]  (260,90) -- (260,110) ;
\draw [color={rgb, 255:red, 74; green, 144; blue, 226 }  ,draw opacity=1 ][line width=1.5]  [dash pattern={on 1.69pt off 2.76pt}]  (260,110) -- (260.87,174.01) ;
\draw  [color={rgb, 255:red, 17; green, 31; blue, 157 }  ,draw opacity=1 ][fill={rgb, 255:red, 74; green, 144; blue, 226 }  ,fill opacity=1 ] (257.99,110) .. controls (257.99,108.89) and (258.89,107.99) .. (260,107.99) .. controls (261.11,107.99) and (262.01,108.89) .. (262.01,110) .. controls (262.01,111.11) and (261.11,112.01) .. (260,112.01) .. controls (258.89,112.01) and (257.99,111.11) .. (257.99,110) -- cycle ;
\draw [color={rgb, 255:red, 74; green, 144; blue, 226 }  ,draw opacity=0.66 ]   (260.16,329.91) -- (283.54,397) ;
\draw [shift={(284.2,398.89)}, rotate = 250.78] [color={rgb, 255:red, 74; green, 144; blue, 226 }  ,draw opacity=0.66 ][line width=0.75]    (10.93,-3.29) .. controls (6.95,-1.4) and (3.31,-0.3) .. (0,0) .. controls (3.31,0.3) and (6.95,1.4) .. (10.93,3.29)   ;
\draw [color={rgb, 255:red, 139; green, 87; blue, 42 }  ,draw opacity=0.75 ]   (260.16,329.91) -- (394.22,347.96) ;
\draw [shift={(396.2,348.23)}, rotate = 187.67] [color={rgb, 255:red, 139; green, 87; blue, 42 }  ,draw opacity=0.75 ][line width=0.75]    (10.93,-3.29) .. controls (6.95,-1.4) and (3.31,-0.3) .. (0,0) .. controls (3.31,0.3) and (6.95,1.4) .. (10.93,3.29)   ;
\draw [color={rgb, 255:red, 147; green, 181; blue, 110 }  ,draw opacity=0.96 ]   (260.16,329.91) -- (142.05,379.45) ;
\draw [shift={(140.2,380.23)}, rotate = 337.25] [color={rgb, 255:red, 147; green, 181; blue, 110 }  ,draw opacity=0.96 ][line width=0.75]    (10.93,-3.29) .. controls (6.95,-1.4) and (3.31,-0.3) .. (0,0) .. controls (3.31,0.3) and (6.95,1.4) .. (10.93,3.29)   ;
\draw [color={rgb, 255:red, 211; green, 141; blue, 73 }  ,draw opacity=1 ]   (260.16,329.91) -- (307.8,281.83) ;
\draw [shift={(309.2,280.4)}, rotate = 134.73] [color={rgb, 255:red, 211; green, 141; blue, 73 }  ,draw opacity=1 ][line width=0.75]    (10.93,-3.29) .. controls (6.95,-1.4) and (3.31,-0.3) .. (0,0) .. controls (3.31,0.3) and (6.95,1.4) .. (10.93,3.29)   ;
\draw [color={rgb, 255:red, 203; green, 188; blue, 71 }  ,draw opacity=1 ]   (260.16,329.91) -- (114.18,305.88) ;
\draw [shift={(112.2,305.56)}, rotate = 9.35] [color={rgb, 255:red, 203; green, 188; blue, 71 }  ,draw opacity=1 ][line width=0.75]    (10.93,-3.29) .. controls (6.95,-1.4) and (3.31,-0.3) .. (0,0) .. controls (3.31,0.3) and (6.95,1.4) .. (10.93,3.29)   ;
\draw  [fill={rgb, 255:red, 10; green, 9; blue, 9 }  ,fill opacity=1 ] (257.01,329.91) .. controls (257.01,328.95) and (258.42,328.18) .. (260.16,328.18) .. controls (261.9,328.18) and (263.31,328.95) .. (263.31,329.91) .. controls (263.31,330.87) and (261.9,331.65) .. (260.16,331.65) .. controls (258.42,331.65) and (257.01,330.87) .. (257.01,329.91) -- cycle ;

\draw (241,232.4) node [anchor=north west][inner sep=0.75pt]    {$\boldsymbol{Z}$};

\end{tikzpicture}

\vspace{3mm}

The above discussion together with the results of Sections 2 of [KoSo1], [KoSo12] implies the following result.

\begin{prp}\label{stability data via hyperplanes} 
The above map gives rise to a bijection between the set of hyperbolic stability data on $\g$ with the fixed central charge $Z$ which compatible with the quadratic form $Q$ and the set of  stability data on $\g$ with the central charge $Z$ and support belonging to $\{Q\ge 0\}$.

\end{prp}

The Proposition \ref{stability data via hyperplanes} implies that fixing a collection of group elements $(A_{H_1,H_2})_{(H_1,H_2)\in \UU^{(2)}}$ we obtain a map from an open subset in $Hom(\Gamma,\C)$ to $Stab(\g)$. The topology on $Stab(\g)$ can be characterized by the property that this map is a local homeomorphism. This allows us to endow the set of hyperbolic stability data with a Hausdorff topology, so that the map in the Proposition \ref{stability data via hyperplanes} becomes a homeomorphism.

We remark that the above-mentioned group elements $A_{H_1,H_2}$ can be extended to a larger class of pairs $(H_1,H_2)$.

\begin{defn}\label{bounded collection} We say that an ordered collection $(H_i)_{i\in i}, H_i\in \UU$ is bounded if there exists $H_1^\prime, H_2^\prime\in \UU$ such that the triple $(H_1^\prime, H_i, H_2^\prime)\in \UU^{(2)}$ for any $i\in I$.

\end{defn}

Then

$$A_{H_1,H_2}=A_{H_1^\prime, H_1}^{-1}A_{H_1^\prime, H_2}=A_{H_1,H_2^\prime}A_{H_2,H_2^\prime}^{-1}\in G_{C_{H_1^\prime, H_2^\prime}}.$$ 
We see that the cocycle condition now holds for bounded triples.

Let $\g=Vect_\Gamma$ be the Lie algebra of vector fields on the algebraic torus $Hom(\Gamma, \C^\ast)$. The following proposition is a  translation of the notion of analytic stability data to the language of hyperbolic stability data. 

\begin{prp}\label{analyticity for hyperbolic stability}
1) Analyticity of stability data on $Vect_\Gamma$ with the central charge compatible with the given quadratic form $Q$ is equivalent to analyticity of transformations $A_{H_1,H_2}, H_i\in \UU, i=1,2$ for rational $H_1$ and $H_2$.

2) Moreover, let $(H_i)_{i\in I\simeq \Z/nZ}, H_i\in \UU$ be a cyclically ordered collection such that any pair $(H_i,H_{i+1}), i\in I$ is admissible,
$\sqcup_{i\in I}C^{ncv}_{H_i, H_{i+1}}=\{x\in \Gamma_\R|Q(x)>0\}$, and $\cap_{i\in I}H_i=W$ has codimension $2$.

Then analyticity of all transformations $A_{H_i, H_{i+1}}$ is equivalent to the analyticity of $A_{H,H^\prime}$ for any pair $(H,H^\prime)\in \UU^{(2)}$, where $H$ and $H^\prime$ are rational.

3) Finally, for any central charge $Z$ compatible with $Q$ analyticity of the stability data with the central charge $Z$ is equivalent to the analyticity of all $A_{H_i,H_{i+1}}$ from the cyclic cover as in 2).

\end{prp}

\begin{rmk}\label{new set up}
The above Proposition holds in the set up discussed at the beginning of Section \ref{analyticity of WCS} (which is a generalization of the one from [KoSo12]).

\end{rmk}

\subsubsection{Rational hyperplanes and analyticity of $\sigma_{Betti}$}

The above discussion holds for an arbitrary wall-crossing structure. Let us return to our  particular WCS, namely to  $\sigma_{Betti}$. Recall that for any admissible ray $l$ a wall-crossing structure assigns the automorphism $A_l$.  In the case of $\sigma_{Betti}$ we have a formula for $A_l$, namely $A_l=iso_{l,-}^{-1}iso_{l,+}$, where $iso_{l,\pm}$ are isomorphisms of the non-archimedean Betti local cohomology with Betti global cohomology for the rays ${l,-}$ (resp. ${l,+}$) which are sufficiently (in fact infinitesimally) close to $l$ and are on the left (resp. on the right) from $l$ for the clockwise orientation of $\R^2$. Indeed, if $l$ is generic then $A_l=id$ because the Proposition \ref{non-archimedean Betti local to global for holomorphic 1-forms} implies that we can identify canonically fibers of $\mathcal{H}^\bullet_{Betti,loc}(X,\rho)$ on the left and on the right of $l$. On the other hand, if $l$ is a Stokes ray then $A_l=R_l$ is the corresponding Stokes automorphism, which is given  by the above formula. This definition implies that the cocycle $A_{H_1,H_2}$ is the coboundary on the vertical hyperplanes.  More precisely the following proposition holds.

\begin{prp}\label{trivialization for sectors}
If $V$ is a semiclosed admissible sector bounded by $l_1$ and $l_2$  in the clockwise order, with $l_1$ included in $V$ and $l_2$ excluded from $V$ then $A_V=iso_{l_2,+}^{-1}iso_{l_1,+}$.

\end{prp}
{\it Proof.} Follows from the factorization formula $A_V=\overrightarrow{\prod}_{l\subset V}A_l$, where the product is taken in the clockwise order over all admissible rays in $V$ combined with the definition of $A_l$ given above. Notice that since $l_2$ is not included in the semiclosed sector $V$, the last term in the clockwise product can be identified with $iso_{l_2,+}^{-1}$.        $\blacksquare$

For any cooriented hyperbolic  hyperplane $H^\prime<H(l)$, we define the transformation $iso_{H^\prime,+}$ by the formula
$$iso_{H^\prime,+}=iso_{l,+}A_{H^\prime, H(l)}.$$ 
The RHS a priori depends on $l$, but the cocycle property implies that it does not.
Indeed, let $l_1>l$ be another  admissible ray which is clockwise on the left from $l$ then 
$$iso_{l_1,+}A_{H^\prime, H(l_1)}=iso_{l_1,+}A_{H(l), H(l_1)}A_{H^\prime, H(l)}=$$
$$iso_{l_1,+}iso_{l_1,+}^{-1}iso_{l,+}A_{H^\prime, H(l)}=iso_{l,+}A_{H^\prime, H(l)}.$$

By interchanging $l$ and $l_1$ we see that the computation gives the same result for $l_1<l$ as well. 

\begin{cor}\label{trivialization for all hyperbolic hyperplanes}
If the pair $(H_1, H_2)$  is bounded then
$$A_{H_1,H_2}=iso_{H_2,+}^{-1}iso_{H_1,+}.$$

\end{cor}
{\it Proof.} Follows from the cocycle condition and Proposition \ref{trivialization for sectors}. $\blacksquare$

Thus we see that that the cocycle $A_{H_1,H_2}$ is always a coboundary. Notice that by our conventions $iso_{H(l),+}=iso_{l,+}$.

Let us now recall that  for generic ray $l$ the isomorphism $iso_l=iso_{l,+}=iso_{l,-}$ was interpreted in Remark \ref{universal integer local system} as an isomorphism of the local and global Betti cohomology groups with coefficients in the universal local system over the ring $\Z[[C\cap \Gamma]]=\{\sum_{\gamma\in C\cap \Gamma}n_\gamma x^\gamma, n_\gamma\in \Z\}$, where $C\subset \Gamma_\R$ is a strict convex cone depending on the ray $l$. In this case we can speak about holomorphic analyticity of $iso_l$ and $iso_{l, \pm}  $ provided they are defined over the subring of convergent series. Having all that in mind we can formulate the following result.

\begin{thm}\label{semirational analyticity}
Under the  Assumption of Section \ref{real and complex cases} we have the following:
for any rational $H^{\prime}\in \UU$ the map $iso_{H^\prime,+}$ is a complex analytic isomorphism.

\end{thm}

{\it Proof.} Theorem follows from a stronger result about rationality, which we will discuss in Section \ref{de Rham for pair of 1-forms}. $\blacksquare$

It follows from Theorem \ref{semirational analyticity} that for any  $(H_1, H_2)\in \UU^{(1)}$ the transformation $A_{H_1,H_2}$ is analytic.

Notice that if $Z$ is rational, we can take $H^\prime=Z^{-1}(l)$, where $l$ is a rational admissible ray. Then the Theorem \ref{semirational analyticity} together with Proposition \ref{trivialization for sectors} imply that $A_V$ is analytic for any admissible open or semiclosed rational sector $V$. Hence 
$\sigma_{Betti}$ is analytic  in agreement with Proposition 3.4.2 from [KoSo12].
Furthermore analyticity property in  Theorem \ref{semirational analyticity} is an open and closed condition in the topology on the space of central charges (or stability data). Same is true for the original analyticity property from [KoSo12].

\begin{cor}\label{semirational analyticity implies analyticity}
Under the Assumption from Section \ref{real and complex cases}  the WCS $\sigma_{Betti}$ is analytic. In particular the 
Theorem \ref{analyticity of WCS for 1-forms} holds.
\end{cor}

{\it Proof.} We can choose a rational central charge $Z^\prime$ which belongs to the same connected component as $Z$ and such that $Q$ is negative on $Ker\, Z^\prime_\R$. Choose a cyclically ordered collection of rational admissible rays $l_i, i\in I$ such the semiclosed sectors $V_i, i\in I$ bounded by $l_i<l_{i+1}$ form a cyclic covering of $\R^2$. We have the transformations $A_{H^\prime(l_i), H^\prime(l_{i+1})}$, where $H^\prime(l_i)=(Z^\prime)^{-1}(l_i)$. The Corollary \ref{trivialization for all hyperbolic hyperplanes} implies that these transformations are analytic.
Hence the WCS associated with $Z^\prime$ is analytic by  Proposition 3.4.2 from [KoSo12]. The results of Section 3.5 from the loc.cit imply  that the whole connected component containing this rational WCS consists of analytic WCS.\footnote{The reader should not forget that although we refer to the results from [KoSo12] we actually mean their generalizations mentioned at the beginning of Section \ref{analyticity of WCS}.} In  particular $\sigma_{Betti}$ is analytic. $\blacksquare$

\begin{prp}\label{analyticity implies semirational analyticity}
If  $iso_{H_1,+}$ is analytic for one rational $H_1\in \UU$ then it is analytic for any rational $H_2\in \UU$.

\end{prp}
{\it Proof.} We need to prove that $iso_{H_2}$ is analytic for any  $H_2\in \UU$. Choose rational central charge $Z^\prime$ such that $Ker\,Z^\prime_\R\subset H_1$. Construct a chain of rational hyperplanes $H_i\in \UU, i\in I=\{i_1,..., i_m\}$ such that $H_{i_1}=H_1, H_{i_m}=H_2$ and such that all consecutive pairs $H_{i_k}, H_{i_{k+1}}$ are bounded. Since $iso_{H_1}=iso_{H_{i_1}}$ is analytic and $A_{H_{i_1}, H_{i_2}}$ is analytic we conclude that $iso_{H_{i_2}}$ is analytic. Then we continue by induction. This concludes the proof. $\blacksquare$

Under the Assumption from Section \ref{real and complex cases} we can say that $iso_{Betti}$ gives rise to an isomorphism of {\it holomorphic bundles} on a small punctured disc in $\C^\ast_t$, i.e. $\mathcal{H}^\bullet_{Betti,glob}(X,\alpha)\simeq \mathcal{H}^\bullet_{Betti,loc}(X,\alpha)$.

\begin{defn}\label{meromorphic structure}
For a holomorphic bundle $E$ over a germ of a punctured disc a meromorphic structure at $t=0$ is an equivalence class of extensions of $E$ at $t=0$, where the equivalence means that  extensions belong to the same orbit of the group of germs at $t=0$ of meromorphic gauge transformations.

\end{defn}

Meromorphic structure gives rise to a finitely generated $\C\{t\}[t^{-1}]$-module of germs of global sections of $E$.

Then  the holomorphic vector bundle  $\mathcal{H}_{Betti,glob}(X,\alpha)$ carries a meromorphic structure at $t=0$. 
More precisely, using the Stokes automorphisms we can glue a {\it new} holomorphic graded vector bundle $\mathcal{H}_{Betti,loc,WCS}^\bullet(X,\alpha)$ over a germ of punctured disc. It  carries a meromorphic structure at $t=0$ since the gluing automorphisms have exponential decay at $t=0$. We denote by $H^\bullet_{Betti,loc,WCS}(X,\alpha)$  the corresponding  graded $\C\{t\}[t^{-1}]$-module of sections of $\mathcal{H}_{Betti,loc,WCS}^\bullet(X,\alpha)$.

\subsection{An approach to analyticity in the complex case}\label{non-archimedean approach to analyticity}

One can try to approach analyticity naively, by attempting to extend the non-archimedean global-to-local isomorphism to the complex analytic one. For that one needs estimates on the coefficients of arising series. 

\begin{defn} \label{thimbles in complex case}
Suppose that $\alpha$ has only simple zeros. Then unstable sets for $Re(\alpha)$ (see Definition \ref{thimbles in real case}) are called thimbles. 
\end{defn}
Each thimble is contractible and moreover it is homeomorphic to $\R^m$.

\begin{que}\label{analyticity of Betti local to global}
Is it true that isomorphisms $iso_t, iso_{t,+}, iso_{t,-}$ defined above give complex analytic isomorphisms of coherent sheaves on holomorphic tube domains associated with open convex subsets $V_\C\subset Cone_{\alpha/t}\subset H^1(X, \ZZ(\alpha),\R)$ such that $V_\C$ contains $[1,+\infty]\cdot V_\C$ and furthermore $(0,+\infty)\cdot V_\C=Cone_{\alpha/t}$?

\end{que}

If the answer to the Question \ref{analyticity of Betti local to global} is positive then transformations $A_V$ are complex analytic for any strict sector $V$ (not necessarily a rational one). On the other hand a positive answer to the Question \ref{analyticity of Betti local to global} follows from  exponential bounds on certain volumes which we will discuss below.

Let us assume that all zeros of $\alpha$ are isolated and simple. Then   thimbles  are just stable subsets for $\ZZ(\alpha)$. They are contractible, hence the restriction of the $1$-form $Re(\alpha/t)$ to a thimble $Th_{i,t}$ emanated from  $x_i\in \ZZ(\alpha)$ is an exact real-valued $1$-form. Let $f_i/t$ be the smooth function on $Th_{i,t}$ such that $d(f_i/t)=Re(\alpha/t)_{|Th_{i,t}}$. We say that $Th_{i,t}$ {\it agrees with $\alpha$ in the direction $\theta=Arg(t)$} if $f_i(x)/t\to -\infty$ provided $x\to \infty, x\in Th_{i,t}$. We may assume that $f_i/t<0$ on $Th_{i,t}$.
Then the pairing of a de Rham cohomology class with the Betti homology class defined by $Th_{i,t}$ is well-defined as long as $\int_{Th_{i,t}}e^{f_i/t}vol<\infty$, where $vol$ is the restriction of a chosen volume form on $X$ to $Th_{i,t}$. The integral can be written as $\int_{0}^{+\infty}e^{-s}V_i(s)ds$, where $V_i(s)$ is the volume of the ``level set'' $f_i/t=s$.  The function $V_i(s)$ increases as $s\to +\infty$, and we need to ensure that  $V_i(s)$ grows no more than exponentially.

The desired exponential bounds hold trivially if $dim_\C X=1$, and can be shown with more efforts in the case $dim_\C X=2$. On the other hand they are not known (and maybe not true) if $dim_\C X\ge 3$. The reason for that is the complicated dynamics of the trajectories of the gradient vector field $v_{\alpha,t}$ corresponding to $\alpha/t$ near the closed subset $\ZZ(\alpha)$.

Despite of the positive answer in small dimensions we have doubts about the positive answer to the Question \ref{analyticity of Betti local to global} in general.
Indeed in  [KoSo12] we asked a similar question about analyticity of factors $A_{V_1}, A_{V_2}$ in the formula $A_{V_1\sqcup V_2}=A_{V_1}A_{V_2}$ in the case of irrational central charges but could not answer that question. In any case the answer to the Question \ref{analyticity of Betti local to global} is not easy to obtain.

\section{Rationality of WCS for $1$-forms and generalizations}\label{applications}

We start with two  examples which illustrate the general theory. These examples are interesting by themselves.

\subsection{ Motivating examples: wall-crossing structures for  $1$-forms on complex curves} \label{WCS for 1-form on curve}

Let $(C, \alpha)$ be a pair consisting of a smooth possibly non-compact complex algebraic curve $C$ of genus $g$ and a holomorphic $1$-form $\alpha$ with  the set of  zeros $\ZZ(\alpha)=\{x_1,...,x_k\}$, having multiplicities $n_1,...,n_k$ respectively, $n_i\ge 1, 1\le i\le k$. We assume that $\alpha$ has poles at all punctures. Let $\Gamma=H_1(C, \ZZ(\alpha), \Z)$ denote the integer relative homology group and $Z: \Gamma\to \C$ the central charge given by $\gamma\mapsto \int_\gamma \alpha$.  Then analyticity of the wall-crossing structure  can be checked directly for sufficiently small $|t|$, since there are  exponential bounds on the volumes of thimbles. 
 
 Let $\overline{C}$ be a smooth compact curve containing $C$ and  $\overline{C}-C=D_v\sqcup D_{log}$, where $D_v=\{y_1,...,y_l\}, D_{log}=\{w_1,...,w_m\}$ are different points. Notice that $D_h=\emptyset$ for curves.
 
 In other words $\overline{C}-C$ consists of poles of $\alpha$, which we split into the union of logarithmic and higher order poles. Equivalently, we can think of $\alpha$ as of meromorphic $1$-form on $\overline{C}$.
 Then near each point $w_i$ we have $\alpha=c_id log(z-w_i)+...$, where $c_i\ne 0$ and dots denote terms holomorphic in the local coordinate $z-w_i$.

\subsubsection{Wall-crossing structure related to $\Gamma$-function}\label{gamma-function}
Let $\alpha=({1\over{x}}-1)dx$ considered as $1$-form on $X=\C^\ast\subset \overline{X}=\C{\bf P}^1$. Then $\ZZ(\alpha)=\{1\}, D_v=\{\infty\}, D_{log}=\{0\}$.

Notice that $\alpha$ can be integrated  along the thimble $L_+=(0,+\infty)$ with the volume form on $\C^\ast$ equal to $dx/x$.  This thimble gives rise to the chain which is a difference of two defined by the rays $[1,+\infty)$ and $(0,1]$ emanating from the zero $x=1$ of the form $\alpha$.

When restricted  to $L_+$ we can write $\alpha=df, f=log(x)+1-x$. Then the corresponding version of the modified exponential integral for $Re(t)>0$ becomes 
$I^{mod}(t)={1\over{\sqrt{2\pi t}}}\int_{L_+}e^{{1\over{t}}f(x)}{dx\over{x}}={\Gamma(\lambda)\over{\sqrt{2\pi}e^{-\lambda}\lambda^{\lambda-1/2}}}$, where $\lambda=1/t$. This expression has an asymptotic expansion which belongs to $\C[[t]]$ and gives rise to a resurgent series.
Similar considerations hold for $Re(t)<0$ and $L_-=(-\infty, 0)$.  Then $I^{mod}(t)=1/I^{mod}(-t)$.
Let us set $I_R^{mod}(t)={1\over{\sqrt{2\pi t}}}e^{1/t}t^{1/t}\Gamma(1/t)$ for $Re(t)>0$ and $I_L^{mod}(t)=1/I_R^{mod}(-t)$ for $Re(t)<0$. Both functions can be extended analyticially to small open neighborhoods of rays $i\R_{>0}$ and $i\R_{<0}$. Then the WCS is equivalent to the following Riemann-Hilbert problem which connects these two functions:
$$I_L^{mod}(t)=I_R^{mod}(t)(1-exp(-2\pi i/t)), t\in i\R_{>0},$$
$$I_R^{mod}(t)=I_L^{mod}(t)(1+exp(-2\pi i/t))^{-1}, t\in i\R_{<0}.$$

There are infinitely many critical points of the antiderivative of  $\alpha$ on the universal covering of $\C^\ast$ given by $log(x)=2\pi i k, k\in \Z$.
The corresponding ``intersection indices"  of different thimbles are given by the following formulas:
$n_{ab}=-1, b=a-1$, $n_{ab}=+1, b>a$ and $n_{ab}=0$ otherwise. Here $a,b\in \Z$.
Notice that in this case there is no problem with convergence of integrals, which is present in general. We have a sequence of morphisms for $t\notin i\R$:
$$H_{Betti,loc,t}^1(X, \alpha)\to H^1_{Betti,glob,t}(X, \alpha)\to H_{DR,t}^1(X,\alpha)\to H^1(\overline{X}_{cor},D_v, E_{\alpha,t}),$$ 
where the first arrow  is the Betti global-to-local isomorphism, the second one is the isomorphism between global Betti and de Rham cohomology, while the last arrow corresponds to the integration over  non-compact cycles. It is an isomorphism only if ${1\over{t}}\notin \{-1,-2,-3,...\}$. More precisely, if $1/t$ does not belong to this countable set then the form $dx/x$ gives rise to a class in the global de Rham cohomology, and it can be paired with the thimble $L_+$. Looking at the above explicit formula for $I^{mod}(t)$ we see that the ``prohibited" integer numbers are   poles of the $\Gamma$-function.

\subsubsection{Square-tiled surfaces}\label{square-tiled surfaces}

 Let us now discuss the relationship of the Betti global-to-local isomorphism and  wall-crossing structures in a special case. For simplicity we  will {\it assume that all zeros of $\alpha$ are simple} (considerations in the general case are similar), and $C=\overline{C}$ (i.e. $\alpha$ has no singularities).
 
 \begin{defn} The pair $(C,\alpha)$ is called square-tiled if the cohomology class $[\alpha]\in H^1(C, \ZZ(\alpha), \C)$ in fact belongs to $H^1(C, \ZZ(\alpha), \Z\oplus i\Z)$
 
 \end{defn}
  
 The form $\alpha$ gives rise to flat metric $|\alpha|^2=\alpha\otimes\overline{\alpha}$ on the Riemann surface  $C$ with singularities at the points $x_i, 1\le i\le k$. The Riemann surface $C$ can be decomposed into the union of finitely many ``geodesic squares'' with the sides being geodesic segments of the metric of length $1$.  This explains the terminology, which agrees with the one in the theory of  dynamical systems.

Let us discuss the construction of the wall-crossing structure in the non-archimedean framework in the square-tiled case.
Recall the character torus ${\bf T}_\Gamma=Hom(\Gamma, {\bf G}_m)$. Its $\C$-points parametrize $\C^\ast$ local systems on $C$ trivialized at all zeros $x_i\in \ZZ(\alpha)$.  For a non-archimedean field $K=\overline{\C((s))}$ and a local system
$\rho\in {\bf T}_\Gamma^{an}(K)$ we have a point $Log(|\rho|)\in \Gamma^\vee\otimes \R$ associated with the ``tropical map'' $Log|\bullet|: {\bf T}_\Gamma^{an}(K)\to \Gamma^\vee\otimes \R\simeq \R^{N}, N=rk\, \Gamma$, where $|\bullet|$ is induced by the non-archimedean norm on $K$ and the ${\bf T}_\Gamma^{an}$ refers to the non-archimedean analytic space in the sense of Berkovich associated with ${\bf T}_\Gamma$.

For $t\in \C^\ast$ we denote as before by $Cone_{\alpha/t}\subset \Gamma_\R$ a sufficiently narrow open strict convex cone containing the ray $\R_{>0}\cdot [Re(\alpha/t)]\in H^1(C,\ZZ(\alpha), \R)=\Gamma_\R^\ast$.

 Consider the set of local systems $\rho\in {\bf T}_\Gamma^{an}(K)$ for which $Log(|\rho|)\in Cone_{\alpha/t}$.
 Such local systems give rise to an open subset (``tube domain'') $U(Cone_{\alpha/t})\subset {\bf T}_\Gamma^{an}(K)$.
 The Proposition \ref{non-archimedean Betti local to global for holomorphic 1-forms} immediately gives the following.
  
  \begin{prp}\label{isomorphisms of sheaves}
 For any  square-tiled surface $(C,\alpha)$ we have the following: 
 
 The isomorphisms $iso_{t,\pm}, t\in \{\pm 1, \pm i\}$ extend to isomorphisms of analytic coherent sheaves $\EE^{an}_{loc,(C,\alpha)}\simeq \EE^{an}_{loc, (C,\alpha)}$ over $U(Cone_{\alpha/t})$.

 \end{prp}
 
 \begin{thm}\label{rationality of isomorphisms of coherent sheaves}
The above isomorphisms of analytic coherent sheaves are induced by the corresponding isomorphisms of algebraic coherent sheaves 
$$\EE_{loc,(C,\alpha)}\simeq \EE_{loc, (C,\alpha)}$$
on a Zariski open subset of ${\bf T}_\Gamma$.

\end{thm}


Let us illustrate the Theorem \ref{rationality of isomorphisms of coherent sheaves}  in the case when $C=\overline{C}$ is a genus $2$ smooth connected projective curve and $\alpha$ is a regular $1$-form with two simple zeros marked on the figure below by ${\bf v_1},{\bf v_2}$. In the case $\Gamma=H_1(C, \{{\bf v_1}, {\bf v_2}\},\Z)\simeq \Z^5$.

 \

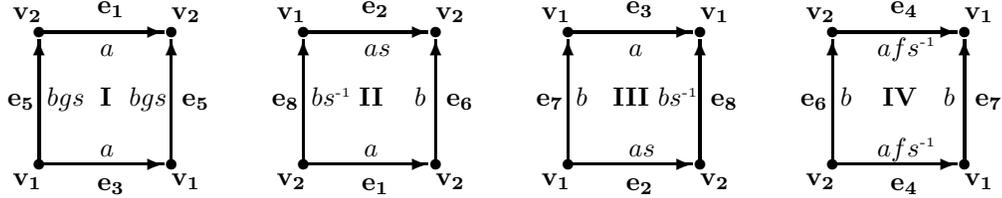
 \[
\centering
\begingroup
\setlength{\unitlength}{0.85pt} 
\makebox[\textwidth][c]{%
  \begin{picture}(370,70)(-12,0)
    \put(0,0){\circle*{4}}
    \put(0,50){\circle*{4}}
    \put(50,0){\circle*{4}}
    \put(50,50){\circle*{4}}
    
    \put(100,0){\circle*{4}}
    \put(100,50){\circle*{4}}
    \put(150,0){\circle*{4}}
    \put(150,50){\circle*{4}}
    
    \put(200,0){\circle*{4}}
    \put(200,50){\circle*{4}}
    \put(250,0){\circle*{4}}
    \put(250,50){\circle*{4}}
    
    \put(300,0){\circle*{4}}
    \put(300,50){\circle*{4}}
    \put(350,0){\circle*{4}}
    \put(350,50){\circle*{4}}
    
    \thicklines
    \put(0,0){\vector(1,0){47}}
    \put(0,0){\vector(0,1){47}}
    \put(0,50){\vector(1,0){47}}
    \put(50,0){\vector(0,1){47}}
    
    \put(100,0){\vector(1,0){47}}
    \put(100,0){\vector(0,1){47}}
    \put(100,50){\vector(1,0){47}}
    \put(150,0){\vector(0,1){47}}
    
    \put(200,0){\vector(1,0){47}}
    \put(200,0){\vector(0,1){47}}
    \put(200,50){\vector(1,0){47}}
    \put(250,0){\vector(0,1){47}}
    
    \put(300,0){\vector(1,0){47}}
    \put(300,0){\vector(0,1){47}}
    \put(300,50){\vector(1,0){47}}
    \put(350,0){\vector(0,1){47}}
    
    \put(-10,-8){$ \bf v_1$}
    \put(50,-8){$\bf v_1$}
    \put(-10,55){$ \bf v_2$}
    \put(50,55){$\bf v_2$}
    
    \put(90,-8){$\bf v_2$}
    \put(150,-8){$\bf v_2$}
    \put(90,55){$\bf v_1$}
    \put(150,55){$\bf v_2$}
    
    \put(190,-8){$\bf v_1$}
    \put(250,-8){$\bf v_2$}
    \put(190,55){$\bf v_1$}
    \put(250,55){$\bf v_1$}
    
    \put(290,-8){$\bf v_2$}
    \put(350,-8){$\bf v_1$}
    \put(290,55){$\bf v_2$}
    \put(350,55){$\bf v_1$}
    
    \put(22,-10){$ \bf e_3$}
    \put(22,57){$\bf e_1$}
    \put(-12,22){$\bf e_5$}
    \put(54,22){$\bf e_5$}
    
    \put(122,-10){$\bf e_1$}
    \put(122,57){$\bf e_2$}
    \put(88,22){$\bf e_8$}
    \put(154,22){$\bf e_6$}
    
    \put(222,-10){$\bf e_2$}
    \put(222,57){$\bf e_3$}
    \put(188,22){$\bf e_7$}
    \put(254,22){$\bf e_8$}
    
    \put(322,-10){$\bf e_4$}
    \put(322,57){$\bf e_4$}
    \put(288,22){$\bf e_6$}
    \put(354,22){$\bf e_7$}
    
    \put(23,22){$\bf I$}
    \put(121,22){$\bf II$}
    \put(217,22){$\bf III$}
    \put(319,22){$\bf IV$}
    
    \put(3,22){$bgs$}
    \put(34,22){$bgs$}
    \put(23,3){$a$}
    \put(23,41){$a$}
    
    \put(103,22){$bs^{\mbox{\tiny -1}}$}
    \put(142,22){$b$}
    \put(123,3){$a$}
    \put(123,41){$as$}
    
    \put(203,22){$b$}
    \put(234,22){$bs^{\mbox{\tiny -1}}$}
    \put(223,3){$as$}
    \put(223,41){$a$}
    
    \put(303,22){$b$}
    \put(342,22){$b$}
    \put(317,3){$afs^{\mbox{\tiny -1}}$}
    \put(317,41){$afs^{\mbox{\tiny -1}}$}
   \end{picture}%
}
\endgroup
\]

 \
 
 Here we have 2 vertices $\{{\bf v_1},{\bf v_2}\}$, 8 edges $\{{\bf e_1},\dots,{\bf e_8}\}$\\
  and 4 squares $\{{\bf I},{\bf II},{\bf III},{\bf IV}\}$. The symbols written near edges {\it inside} squares denote holonolomy of the local system trivialized at vertices (zeroes of 1-form $\alpha$), defined universally over the ring $\Z[a^{\pm 1},b^{\pm 1},f^{\pm 1},g^{\pm 1},s^{\pm 1}]$, which is the group ring of the abelian group $\Gamma$. The central charge $Z: \Gamma\to \C$ in this notation is given by $Z(a)=1, Z(b)=i=\sqrt{-1}, Z(f)=Z(g)=Z(s)=0$.
 
 We  explain below the construction of the isomorphism 
 $$(iso_{1,+}^\ast)^{-1}=\lim_{\epsilon\to 0^+}(iso_{t_\epsilon^{hor}}^\ast)^{-1}$$
  between homology groups $H_{1,loc, t_\epsilon^{hor}}(C, \rho)$ and $H_{1, glob, t_\epsilon^{hor}}(C, \rho)$ where $t=t_\epsilon^{hor}=e^{i\epsilon}$. We assume that $tan(t_\epsilon^{hor})\notin \Q$ and $\epsilon$ is very small positive number.   More precisely, this isomorphism should be understood as isomorphism of fibers of vector bundles over $(K^\times)^5=Hom(\Z^5, K^\times)$ with the fibers which are the corresponding local and global Betti homology.

First we explain the global Betti homology groups. They can be calculated by means of the chain complex $C_0\to C_1\to C_2$, where $C_0=C_0(C,\rho)=K\langle {\bf v_1,v_2}\rangle, C_1=C_1(C,\rho)=K\langle {\bf e_1,...,e_8} \rangle, C_2=C_2(C,\rho)=K\langle{\bf I,II,III,IV}\rangle$, i.e. they are $K$-vector spaces spanned by the vertices, edges and faces of the decomposition of $C$ into squares.
 
  The differential in the chain complex  can be directly read from the picture. E.g. we have
 $$\partial ({\bf e_1}) ={\bf v_2}-a\cdot {\bf  v_2}, \quad \partial ({\bf e_2}) ={\bf v_1}-as\cdot {\bf  v_2},\quad \dots, \quad \partial  ({\bf e_8}) ={\bf v_2}-bs^{-1}\cdot {\bf  v_1}$$ 
and 
 $$ \partial({\bf I}) =({\bf e_5}+bgs\cdot {\bf e_1})-({\bf e_3}+a \cdot{\bf e_5}),\quad \partial({\bf II}) =({\bf e_8}+bs^{-1}\cdot {\bf e_2})-({\bf e_1}+a \cdot{\bf e_6}),\quad \dots $$

 Let us explain local Betti homology.
 The abelian group $H_{1,loc, t_\epsilon^{hor}}(C, \alpha)$ has basis consisting of two thimbles emanating from the zeros ${\bf v_1}$ and ${\bf v_2}$.
This is explained in more detail  below, where the horizontal rays emanating from ${\bf v_1}$ are denoted by $\gamma_{1,\pm}^{\to}$ and the corresponding thimble is denoted by $\gamma_1^{\to}$.
 
 Namely, the thimble $\gamma_1^{\to}$ with the small positive slope emanating from vertex $\bf v_1$ consists of two rays $\gamma_{1,+}^{\to}$ and $\gamma_{1,-}^{\to}$. As a chain it is equal to the {\it difference} $[\gamma_{1,+}^{\to}-\gamma_{1,-}^{\to}]$ of rays endowed with the natural orientations. 
 
 \
 
 The ray $\gamma_{1,+}^{\to}$ starts in the south-west corner of square $\bf I$, and the ray  $\gamma_{1,-}^{\to}$ starts in the south-west corner of square $\bf III$. If the slope is positive and very small then for a long time $\gamma_{1,+}^{\to}$  stays in the horizontal cyclinder obtaines by gluing two vertical sides of square $\bf I$ with itself:
 
 
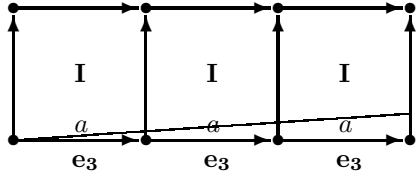
\begin{picture}(200,70)(-70,0)
 \put(0,0){\circle*{4}}
  \put(0,50){\circle*{4}}
  \put(50,0){\circle*{4}}
  \put(50,50){\circle*{4}}

  \put(100,0){\circle*{4}}
  \put(100,50){\circle*{4}}
  \put(150,0){\circle*{4}}
  \put(150,50){\circle*{4}}
  
   \thicklines
  \put(0,0){\vector(1,0){47}}
  \put(0,0){\vector(0,1){47}}
  \put(0,50){\vector(1,0){47}}
  \put(50,0){\vector(0,1){47}}
  
  \put(50,0){\vector(1,0){47}}
  \put(50,50){\vector(1,0){47}}

  \put(100,0){\vector(1,0){47}}
  \put(100,0){\vector(0,1){47}}
  \put(100,50){\vector(1,0){47}}
  \put(150,0){\vector(0,1){47}}
  
  \put(23,22){$\bf I$}
  \put(73,22){$\bf I$}
  \put(123,22){$\bf I$}
  
\put(22,-10){$ \bf e_3$}
\put(72,-10){$ \bf e_3$}
\put(122,-10){$ \bf e_3$}

\put(23,3){$a$}
\put(73,3){$a$}
\put(123,3){$a$}
  
  \thinlines
  \qbezier(0,0)(75,5)(150,10)

 \end{picture}
 
 \
 
 \
 
and the contribution of $\gamma_{1,+}^{\to}$  approximates the sum of the geometric progression
$$ {\bf e_3}+a\cdot {\bf e_3}+a^2\cdot {\bf e_3}+\dots=\frac{1}{1-a}\cdot {\bf e_3} $$

\

Similarly, ray $\gamma_{1,-}^{\to}$ stays for a long time in the horizontal cyclinder obtained by gluing  vertical sides of squares ${\bf III},{\bf II}$ and $\bf IV$:

 
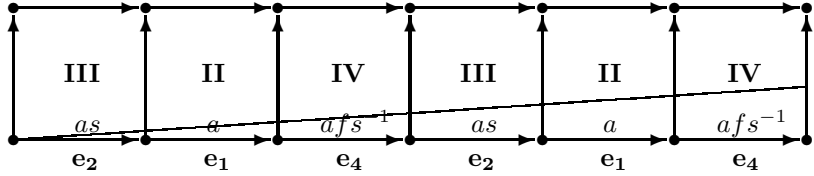
\begin{picture}(300,70)(0,0)
 \put(0,0){\circle*{4}}
  \put(0,50){\circle*{4}}
  \put(50,0){\circle*{4}}
  \put(50,50){\circle*{4}}

  \put(100,0){\circle*{4}}
  \put(100,50){\circle*{4}}
  \put(150,0){\circle*{4}}
  \put(150,50){\circle*{4}}
  
 \put(200,0){\circle*{4}}
  \put(200,50){\circle*{4}}
  \put(250,0){\circle*{4}}
  \put(250,50){\circle*{4}}
  
  \put(300,0){\circle*{4}}
  \put(300,50){\circle*{4}}
  
   \thicklines
  \put(0,0){\vector(1,0){47}}
  \put(0,0){\vector(0,1){47}}
  \put(0,50){\vector(1,0){47}}
  \put(50,0){\vector(0,1){47}}
  
  \put(50,0){\vector(1,0){47}}
  \put(50,50){\vector(1,0){47}}

  \put(100,0){\vector(1,0){47}}
  \put(100,0){\vector(0,1){47}}
  \put(100,50){\vector(1,0){47}}
  \put(150,0){\vector(0,1){47}}
  
  \put(150,0){\vector(1,0){47}}
  \put(150,50){\vector(1,0){47}}
  
  \put(200,0){\vector(1,0){47}}
  \put(200,0){\vector(0,1){47}}
  \put(200,50){\vector(1,0){47}}
  \put(250,0){\vector(0,1){47}} 
  
 \put(250,0){\vector(1,0){47}}
  \put(250,50){\vector(1,0){47}}
  \put(300,0){\vector(0,1){47}} 
  
  \put(19,22){$\bf III$}
  \put(71,22){$\bf II$}
  \put(120,22){$\bf IV$}
  \put(169,22){$\bf III$}
  \put(221,22){$\bf II$}
  \put(270,22){$\bf IV$}
  
\put(22,-10){$ \bf e_2$}
\put(72,-10){$ \bf e_1$}
\put(122,-10){$ \bf e_4$}

\put(172,-10){$ \bf e_2$}
\put(222,-10){$ \bf e_1$}
\put(272,-10){$ \bf e_4$}

\put(23,3){$as$}
\put(73,3){$a$}
\put(116,3){$afs^{-1}$}

\put(173,3){$as$}
\put(223,3){$a$}
\put(266,3){$afs^{-1}$}

  \thinlines
  \qbezier(0,0)(75,5)(150,10)
  \qbezier(150,10)(225,15)(300,20)
  
 \end{picture}
 
 \
 
 \
 
 In this way we obtain the sum consisting of 3 geometric series
  \begin{multline*} {\bf e_2}+as\cdot {\bf e_1}+a^2s\cdot {\bf e_4}+a^3f\cdot {\bf e_2}+a^4fs\cdot{\bf e_1}+\dots=\\
 = {1\over 1-a^3 f}\cdot{\bf e_2}+{as\over a-a^3 f}\cdot {\bf e_1}+{a^2 s\over 1-a^3 f}\cdot {\bf e_4}\end{multline*}

 The limit of chain $\gamma_1^{\to}$ as the slope tends to zero, is
 \begin{multline*}(iso_{1,+}^\ast)^{-1}(\gamma_1^{\to})=lim_{\epsilon\to 0}(iso_{t_\epsilon^{hor}}^\ast)^{-1}(\gamma_{1}^{\to})= [\gamma_{1,+}^{\to}-\gamma_{1,-}^{\to} ]=\\ = [\frac{1}{1-a}\cdot {\bf e_3} -{1\over 1-a^3 f}\cdot{\bf e_2}-{as\over a-a^3 f}\cdot {\bf e_1}
  -{a^2 s\over 1-a^3 f}\cdot {\bf e_4}].\end{multline*}

One can check by a direct calculation that the last expression in the square bracket is closed under $\partial$.

Similarly we can deal with the horizontal thimbles emanating from the vertex $\rm\bf{v}_2$ and calculate
$(iso_{1,+}^\ast)^{-1}(\gamma_2^\to)$.

Then looking at the above explicit formulas as expressions of parameters $a,b,f,g,s$ one can show that  $(iso_{1,+}^\ast)^{-1}(\gamma_1^{\to})$ and $(iso_{1,+}^\ast)^{-1}(\gamma_2^\to)$ form a basis in $H_1(C,\rho)$ on a Zariski open subset of ${\bf T}_\Gamma$.

Similarly one can consider thimbles close to vertical direction. Then in the end we obtain formulas for all $iso_{t,-}, t\in \{\pm 1,\pm i\}$.


Let $V_1$ be the open first quadrant in $\R^2$, and $A_{V_1}$ the transition matrix between ``vertical'' and ``horizontal'' bases. Then $A_{V_1}$ is given by the following $2\times 2$ matrix with coefficients which are rational functions in $a,b, f, g, s$.

 


The matrix can be expressed as 
\[
\begin{pmatrix}
  M_{11} & M_{12} \\
  M_{21} & M_{22}
\end{pmatrix},
\]
where the individual entries are given by:
\begin{align*}
  M_{11} &= 1 + \frac{(a^3 f-a^2 f) (b^2 g-b^3 g) + (a-a^3 f) (b^3 g-b)}{(1-a^3 f) (1-b^3 g)}, \\[1ex]
  M_{12} &= \frac{s(a-a^3 f)b^2 + s(a^3 f-a^2 f)b}{(1-a^3 f) (1-b^3 g)} - \frac{s(a^2 f-a f)b}{(1-a^3 f)(1-b)}, \\[1ex]
  M_{21} &= \frac{s^{-1}a^2 (b^3 g-b) + a (b^2 g-b^3 g)}{(1-a^3 f) (1-b^3 g)} - \frac{s^{-1}a (b g-b^2 g)}{(1-a)(1-b^3 g)}, \\[1ex]
  M_{22} &= 1 + \frac{a^2 b^2+a b}{(1-a^3 f) (1-b^3 g)} - \frac{a^3 fb}{(1-a^3 f)(1-b)} - \frac{a b^3 g}{(1-a) (1-b^3 g)}.
\end{align*}

Analogously one computes other  matrices $A_{V_j}, 2\le j\le 8$, where $V_j$ are either other quadrants in $\R^2$ or coordinate half-axes. The collection $A_{V_j}, 1\le j\le 8$ determines the wall-crossing structure. Clearly it is analytic because all entries of the matrices $A_{V_j}$ are rational functions on ${\bf T}_\Gamma$.  This result is an illustration of the general  Proposition \ref{rationality of isomorphisms} discussed in the Section \ref{rationality of WCS}.

\subsection{Rationality and analyticity of WCS}\label{rationality of WCS}
As we explained above,  we can generalize the theory of local and global Betti and de Rham cohomology from the case of holomorphic closed $1$-forms to the case of complex-valued closed $1$-forms which are holomorphic near the set of zeros.  In particular one can hope to have the corresponding version of the analytic WCS. Now the set of zeros $\ZZ$ does not have to be the set of isolated points. Furthermore, we may vary the cohomology class $[\alpha]\in H^1(X,\ZZ,\C)$ in an open subset in the vector space $H^1(X,\ZZ,\C)$.  In this new framework we have two types of parameters: the topological one, which is the cohomology class $[\alpha]$ and the non-topological, which is a choice of the complex structure near $\ZZ$.
Then we can deform the given pair to the one for which the cohomology class $[\alpha]\in H^1(X,\ZZ, \Z\oplus i\Z)$, i.e. it is integer.  Such a closed $1$-form gives rise to a map $f: X\to T^2=\R/\Z\times \R/\Z$ such that $f(\ZZ)=(0,0)$ and $(0,0)$ is the only critical value. The generic fiber is a smooth real $(2n-2)$-dimensional manifold. 

In this setting to a generic  admissible ray $l$ we can assign the corresponding isomorphism $iso_l$ defined similarly to the one in Proposition \ref{non-archimedean Betti local to global for holomorphic 1-forms}. For a rational admissible ray $l$ we have analogs of the isomorphisms $iso_{l,\pm}$ defined as limits from left and right of the corresponding generic isomorphisms as well as the isomorphism $A_l=iso_{l,-}^{-1}iso_{l,+}$. It enjoys certain integrality property (see Remark \ref{universal integer local system}). As in the case of  holomorphic $1$-forms those isomorphisms can be thought of as isomorphisms of analytic coherent sheaves given by the cohomology groups with coefficients in non-archimedean local systems  over the non-archimedean field $K$ (see  Section \ref{non-archimedean Betti local to global}) trivialized at $\ZZ$ (see Remark \ref{soft sheaf} b)). The corresponding  non-archimedean tube domains $U(Cone_{\alpha/t})$ are associated with  sufficiently narrow cones containing cohomology class $[\alpha/t]$.

\begin{prp}\label{rationality of isomorphisms} Let $l=l_t=\R_{>0}\cdot t, t\in \C^\ast$ be a rational ray in $\R^2$.
Assume that  $X$ is compact and $[\alpha]\in H^1(X,\ZZ, \Z\oplus i\Z)$. Then  the isomorphisms $iso_{l,\pm}$ are analytic isomorphisms of the algebraic coherent sheaves over the tube domain $U(Cone_{\alpha/t})\subset {\bf T}_\Gamma:=Hom(\Gamma, {\bf G}_m), \Gamma=H_1(X, \ZZ,\Z)$ which map the space of rational sections into the space of rational sections. In particular the corresponding WCS enjoys the rationality property  and hence it is analytic.

\end{prp}

{\it Proof.} Let $\rho$ be a  rank $1$ local system over a non-archimedean field $K$ over $X$, trivialized at the set of zeros $\ZZ=\ZZ(\alpha)$. Then by Leray spectral sequence one has $H^\bullet(X,\rho)\simeq H^\bullet(T^2, \R f_\ast(\rho))$. Notice that the direct image $F_\rho:=\R f_\ast(\rho)$ is a complex of constructible sheaves  on $T^2$ which  is a local system outside of the point $p=(0,0)$. 
As such it can be described in terms of some linear algebra data. More precisely, let $V^\bullet$ be a generic stalk of the direct image sheaf and $V_0^\bullet$ be the  stalk at $p$. They are complexes of finite-dimensional $K$-vector spaces. Then the above-mentioned linear algebra data consist of a pair of automorphisms (in the derived sense) $A, B: V^\bullet\to V^\bullet$, a morphism of complexes $i:V_0^\bullet\to V^\bullet$ and a homotopy $h$ between $ABi$ and $BAi$. \footnote{In the case when the complexes are concentrated in degree $0$ this amounts to the condition that the image of $i$ belongs
to the subspace $V^{inv}$ consisting of the invariants with respect to the commutant $\langle ABA^{-1}B^{-1}\rangle$ of the free group generated by $A$ and $B$.}  We may assume that $A$ corresponds to the meridian of $T^2$ while $B$ corresponds to the equator of $T^2$.
Equivalently, these data can be described in terms of the representation of a dg-quiver with two vertices and four arrows corresponding to $A,B,i,h$, subject to the above restrictions on $A,B,i,h$. Thinking of the quiver as of a small dg-category $\CC$, we get a dg-functor $\CC\to Perf_{dg}(Vect_K)$ to the dg-category of perfect $K$-modules.

Let us now consider $\rho$ as a variable point of ${\bf T}_\Gamma$. Then we have a universal local system $\rho_{univ}$ on $X$ with the fiber which is a rank $1$ free module over $\OO({\bf T}_\Gamma)=\Z[\Gamma]$. The corresponding sheaf $F_{\rho_{univ}}$ is a complex of constructible sheaves of perfect $\Z[\Gamma]$-modules on the torus $T^2$. 

Repeating the above considerations with vector spaces replaced by  $\Z[\Gamma]$-modules we will get a functor $\CC\to Perf(\Z[\Gamma]-mod)$. The fact that $A$ and $B$ preserve the category of perfect $\Z[\Gamma]$-modules relies on the properness of the map $f$.

In what follows in order to simplify the exposition we will restrict considerations to closed $K$-points of $Spec(\Z[\Gamma])$.

Since $K$ is a field the category $Perf_{dg}(Vect_K)$ is quasiequivalent to the dg-category $Vect_K^\Z$ of $\Z$-graded finite-dimensional $K$-vector spaces. The latter is understood as the category of complexes of finite-dimensional $K$-vector spaces endowed with the trivial differential. 
Then we can replace the above functor $\CC\to Perf_{dg}(Vect_K)$ by its ``minimal model\rq\rq{}, i.e. by a functor from $\CC$ to the above category of complexes with trivial differential.
Then  we can write $V^\bullet=\oplus_{j\in \Z}V^{j}, V_0^\bullet=\oplus_{j\in \Z}V_0^{j}$ which are  direct sums of finite-dimensional vector spaces. The above-mentioned linear algebra data are encoded in $A,B\in Aut_{Vect_K^\Z}(V^\bullet), i\in Hom_{Vect_K^\Z}(V_0^\bullet, V^\bullet), h\in Hom_{Vect_K^\Z}(V_0^\bullet, V^\bullet[-1])$.

Then cohomology $H^\bullet(X,\rho_{univ})$ are isomorphic to the cohomology of the complex
$$C^\bullet=V_0^\bullet\oplus V^\bullet[-1]\oplus V^\bullet[-1]\oplus V^\bullet[-2]$$
endowed with the differential
$$d_{C^\bullet}=
\begin{pmatrix}
0&0\\
(1-A)i&0&0&0&\\
(1-B)i&0&0&0\\
h&-(1-B)&(1-A)&0\\
\end{pmatrix}
.
$$

Consider also the complex $C^\bullet_{can}=V_0^\bullet\oplus V^\bullet[-1]$ endowed with the differential

$$d_{C^\bullet_{can}}=
\begin{pmatrix}
0&0\\
i&0\\
\end{pmatrix}
.
$$

Having in mind that in the proof of the  Theorem \ref{isomorphism for non-archimedean local systems} we worked with chains rather than cochains it is more convenient to consider dual complexes $(C^\bullet)^\ast$ and $(C_{can}^\bullet)^\ast$ endowed with dual differential $d_{C^\bullet}^\ast$ and $d_{C_{can}^\bullet}^\ast$ respectively.

Assume that $(1-A)$ and $(1-B)$ are invertible (this is true e.g. if both $A$ and $B$ have small non-archimedean norm).

There are two quasi-isomorphisms of complexes
$\phi_A^\ast: (C_{can}^\bullet)^\ast\to (C^\bullet)^\ast, \phi_B^\ast: (C_{can}^\bullet)^\ast\to (C^\bullet)^\ast$
given explicitly by the formulas
$$\phi_A^\ast(u_0,u)=
\begin{pmatrix}
1&0\\
0&(1-A^\ast)^{-1}u\\
0&0\\
0&0\\
\end{pmatrix}
$$

and

$$\phi_B^\ast(u_0,u)=
\begin{pmatrix}
1&0\\
0&0\\
0&(1-B^\ast)^{-1}u\\
0&0\\
\end{pmatrix}
.
$$
Here the notation like $A^\ast$ etc. means the dual map.
One checks directly that the following matrix is a morphism of complexes and it is a homotopy inverse to $\phi_A^\ast$
$$\psi_A^\ast=
\begin{pmatrix}
1&0&h^\ast(1-A^\ast)^\ast&0\\
0&1&(1-B^\ast)(1-A^\ast)&0\\
\end{pmatrix}
$$
and similarly one has the homotopy inverse $\psi_B^\ast$.

Writing $(1-A^\ast)^{-1}$ as a sum $\sum_{m\ge 0}(A^\ast)^m$ one can see that on the cohomology $H^\bullet(X,\rho_{univ})$ the  dual to the map  $\phi_A^\ast$ can be identified with $iso_{l_{vert},+}$, where $l_{vert}$ (resp. $l_{hor}$) are the rays in $\R^2$ corresponding to the meridian (resp. equator) directions on $T^2$.
This follows from the  identification of the powers $(A^\ast)^m u$ with the terms of the extension of the chain $u$ in the Theorem \ref{isomorphism for non-archimedean local systems}\footnote{ It is also helpful  to look at the explicit formulas in Section \ref{square-tiled surfaces}.}
Similarly, the dual $\phi_B$ to $\phi_B^\ast$ induces on the cohomology $H^\bullet(X,\rho_{univ})$ the map $iso_{l_{hor},+}$.

Consider the composition
$$\psi_A^\ast\circ \phi_B^\ast=
\begin{pmatrix}
1&h^\ast(1-A^\ast)^{-1}(1-B^\ast)^{-1}\\
0&(1-A^\ast)(1-B^\ast)(1-A^\ast)^{-1}(1-B^\ast)^{-1}\\
\end{pmatrix}
$$

Notice that the  matrix rationally depends on $A$ and $B$. Since $A$ and $B$ depend rationally on the point of 
${\bf T}_\Gamma$ the Proposition follows. $\blacksquare$

\begin{cor}\label{rationality and analyticity}

1) In the previous notation the isomorphisms $A_{H(l_1), H(l_2)}$ are series on the torus given by matrices with entries which are rational functions with coefficients in $\Q$.  

2) Theorem \ref{analyticity of WCS for 1-forms} holds.

\end{cor}


\begin{rmk}\label{interpretation via wheel of lines}

One can interpret the rationality property of WCS in a different way, using the geometric picture with the wheel of projective lines. Namely for the torus ${\bf T}_\Gamma\simeq (\C^\ast)^n$ one can consider the algebraic coherent sheaves $\EE_i=H^i({\bf T}_\Gamma, \rho_{univ}), i\ge 0$. A choice of central charge and a compatible wheel of rational cones gives a partial compactification $\overline{\bf T}_\Gamma$ of ${\bf T}_\Gamma$ which contains a wheel of projective lines (see [KoSo7], [KoSo12]). Then each $\EE_i$ has a canonical extension to this partial compactification, and moreover it is an algebraic vector bundle in a Zariski open neighborhood of the wheel of projective lines. 

\end{rmk}

\begin{rmk}\label{rationality in non-compact case}
Proposition \ref{rationality of isomorphisms} holds in a larger generality, when $X$ is non-compact but has a ``good compactification" similar to the one we considered in Section \ref{case of 1-forms}. In this case we  do not have global algebraic coherent sheaves $\EE_i$ as above. Instead we have only a $1$-parameter family of algebraic coherent sheaves depending on the angle $\theta\in \R/2\pi \Z$. For a sufficiently small arc these coherent sheaves form a  constant family thus giving rise to a single coherent sheaf which is canonically extended to the corresponding part of the wheel of projective lines. 
More pedantically this can be spelled out such as follows. On the wheel of projective lines $\cup_i {\bf P}^1_i$ endowed with the standard analytic topology we define a sheaf of commutative algebras $\OO^{new}$ with the stalk over a point $x$ consisting of such germs of the stalk $\OO_{\overline{\bf T}_\Gamma, Zar,x}$ that their restrictions to the wheel of projective lines (considered as a complex analytic space) is locally constant. Then the above-mentioned extension will be a sheaf of locally-free finitely generated modules over $\OO^{new}$. In particular it gives a local system on the wheel $\cup_i{\bf P}^{1,an}$ which is obtained via the pull-back from the one on the standard circle $S^1$ (see [KoSo7], [KoSo12]).
\end{rmk}

\subsection{Wall-crossing formulas and non-commutative identities}\label{WCF and nc identities}
Here we present computations which explain how the wall-crossing formulas can be derived in the framework of the Proposition \ref{rationality of isomorphisms} making a further simplifying assumption $V_0=0$. It turns out that the geometric details of that Proposition are not very important, and the wall-crossing formulas appear as corollaries of some non-commutative identities.

\begin{prp}
1) (first 6-term identity)  If $(1-B)$ and $(1-AB)$ are both invertible then
$$(1-B)^{-1}(1-A)=(1-AB)^{-1}(1-A)\cdot (1-B)^{-1}(1-BA).$$

2) (second 6-term identity) If
$(1-AB)$ and $(1-BA)$ are both invertible  then
$$(1-A)(1-BA)^{-1}(1-B)=(1-B)(1-AB)^{-1}(1-A)$$

\end{prp}

{\it Proof.} We will prove here the second $6$-terms identity:
$$(1-A)(1-BA)^{-1}(1-B)=(1-B)(1-AB)^{-1}(1-A)$$

We remark that the meaning of this identity is that the noncommutative expression in the LHS is {\it symmetric} in $A,B$.

Let us take the inverse of both sides and prove the corresponding identity:
$$(1-B)^{-1}(1-BA)(1-A)^{-1}\stackrel{?}{=}(1-A)^{-1}(1-AB)(1-B)^{-1}.$$
If we introduce the notation $X:=1-A,\quad Y:=1-B$ then we want to prove the identity
$$Y^{-1}(1-(1-Y)(1-X))X^{-1}\stackrel{?}{=}X^{-1}(1-(1-X)(1-Y))Y^{-1}$$
$$Y^{-1}(Y+X-YX)X^{-1}\stackrel{?}{=}X^{-1}(X+Y-XY)Y^{-1}.$$
But 
$$X^{-1}+Y^{-1}-1=Y^{-1}+X^{-1}-1. \quad \blacksquare$$

Let us rewrite the first $6$-term identity such as follows:

$$h(B,A)=h(AB,A)\cdot h(B,BA),$$
where $$\quad h(B,A):=(1-B)^{-1}(1-A).$$

Let us iterate the above identity:
$$ h(B,A)=\\
h(AB,A)h(B,BA)=$$
$$(AAB,A)h(AB,ABA)h(BAB,BA)h(B,BBA)=\\
...$$

$$  
h(\r01,\b11)=\\~\\
h(\r12,\b11)h(\r01,\b12)=\\~\\
h(\r23,\b11)h(\r12,\b23)h(\r13,\b12)h(\r01,\b13)=\\~\\$$
$$h(\r34,\b11)h(\r23,\b34)h(\r35,\b23)h(\r12,\b35)h(\r25,\b12)h(\r13,\b25)h(\r14,\b13)h(\r01,\b14)=\\\dots,$$

where 
$$ 
\r01:=B, \quad \b11:=A,\,\,\, \\~\\
\r {p+p'} {q+q'}:=\b {p'}{q'}\cdot \r{p}{q},\qquad \b {p+p'} {q+q'}:=\r{p}{q}\cdot \b {p'}{q'}$$

Let us replace $h(B,A)$ by $\tilde h(B,A):=(1-A)(1-B)^{-1}$, i.e.
$h(B,A)\rightsquigarrow \tilde h(B,A)$. Both expressions coincide in a certain limit. We get:

$$ h(\r34,\b11)h(\r23,\b34)h(\r35,\b23)h(\r12,\b35)h(\r25,\b12)h(\r13,\b25)h(\r14,\b13)h(\r01,\b14)\rightsquigarrow \\~\\$$
$$\hh(\r34,\b11)\hh(\r23,\b34)\hh(\r35,\b23)\hh(\r12,\b35)\hh(\r25,\b12)\hh(\r13,\b25)\hh(\r14,\b13)\hh(\r01,\b14)=\\~\\$$
$$(1-\b11)(1-\r34)^{-1}(1-\b34)(1-\r23)^{-1}(1-\b23)\dots(1-\r01)^{-1}=\\~\\$$
$$(1-\b11)\gg34\gg23\gg35\gg12\gg25\gg13\gg14(1-\r01)^{-1}$$
where
$$ \gg{p}{q}:=(1-\r{p}{q})^{-1}(1-\b{p}{q}).$$

In the limit
$$
h(\r01,\b11)=(1-\b11)\cdot \prod_{{p\over q}\in (0,1)\searrow} \gg{p}{q}\cdot(1-\r01)^{-1}.$$

This is the desired wall-crossing formula. The non-commutative rational expressions in the ordered product can be interpreted as factors $A_l$
corresponding to different rays in the clockwise factorization formula $A_V=\prod^{\to}_lA_l$, where $V$ is the first quadrant in $\R^2$.

\section{Generalizations}\label{generalizations}

\subsection{Morse-Novikov theory of a pair of closed $1$-forms}\label{Morse-Novikov for pair of 1-forms}

Recall the framework of Morse-Novikov theory.
Let $(Y,\beta)$ be a pair consisting of a compact closed manifold $Y$ and a closed $1$-form $\beta$ which has only isolated Morse zeros. 
The deformation theory of the pair $(Y,\beta)$ is unobstructed, and hence the moduli space of such pairs can be locally embedded into $H^1(Y, \ZZ(\beta), \R)$ by the map $(Y,\beta)\mapsto [\beta]$. Then Novikov version of the Morse theory assigns to such a pair and a generic Riemannian metric a cochain complex over the Novikov ring $Nov$. This can be generalized to the setting of the Section \ref{Betti cohomology for real one-forms}. Namely, $Y$ can be not necessarily closed, and the form $\beta$ can have not necessarily Morse and even not necessarily isolated zeros.
Moreover, the local system $E_\beta$ corresponding to $\beta$ can be replaced by the one which is close to $E_\beta$ in non-archimedean sense.
Instead of the  differential in the Morse-Novikov complex  we now have a spectral sequence converging from local Betti cohomology to the global one.

In the holomorphic setting we have  a complex manifold endowed with a holomorphic closed $1$-form $\alpha$. Let $\beta_\theta=Re(e^{-i\theta}\alpha), \theta\in \R/2\pi \Z$. Then for a generic $\theta$ and arbitrary hermitian metric there are no saddle connections tangent to the vector field corresponding to $\beta_\theta$.  Hence as we have already discussed,  we have canonical isomorphism between local and global Betti cohomology. In the special case when $\alpha$ has isolated Morse (i.e. simple) zeros, the Morse-Novikov complex for $\beta_\theta$ has trivial differential. Hence in the holomorphic case the main question is not about the differential but about the arising wall-crossing structure.

From the $C^\infty$ point of view we have {\it two} real closed $1$-forms: $\alpha_1=Re(\alpha)$ and $\alpha_2=Im(\alpha)$. It is interesting to generalize the contents of Section \ref{case of 1-forms} to the case of the pair of real closed $1$-forms, maybe  satisfying some extra conditions. In this subsection we discuss this possibility. 

Thus we assume that we are given a triple $(Y,\alpha_1,\alpha_2)$ such that $Y$ is a $C^\infty$ compact closed manifold and $\alpha_i, i=1,2$ are smooth closed $1$-forms on $Y$. We impose the following restrictions on the pair $\alpha_1, \alpha_2$ which are automatically satisfied for the real and imaginary parts of a holomorphic closed $1$-form:

a) the closed sets of zeros $\ZZ(\alpha_i), i=1, 2$ satisfy the property $\ZZ(\alpha_1)=\ZZ(\alpha_2):=\ZZ$;

b) outside of the common set of zeroes $\ZZ$ the forms $\alpha_1$ and $\alpha_2$ are pointwise independent (i.e. they span at each point is a $2$-dimensional
real vector space);

c) there exist functions $f_i, i=1,2$ in a neighborhood of $\ZZ$ such that each $f_i$ is equal to zero on $\ZZ$ and $df_i=\alpha_i, i=1,2$.

If a) and b) are satisfied we  will call a smooth vector field $\xi$ on $Y$  {\it gradient-like} for our triple if pointwise we have $\alpha_1(\xi)>0$ outside of $\ZZ$ and $\alpha_2(\xi)=0$. 

It is easy to see that gradient-like vector fields always exist. Indeed, we can always construct such a vector field outside of $\ZZ$. Then multiplying it by a smooth function which has a sufficiently fast decay on $\ZZ$ we extend the vector field to $\ZZ$.
In fact the set of gradient-like vector fields is a non-empty open convex cone, and hence it is contractible.


\begin{rmk}
We expect a generalization of the above set-up as well as the following one to the case of manifolds with corners.
In particular there is a compactification $\overline{Y}$ of $Y$ to a manifold with corners.  We assume that $\overline{Y}-Y=\partial_\pm Y\cup \partial_hY$, and we assume that there is a manifold with corners containing $\ZZ$ which has a similar splitting of its boundary, and all these data satisfy the same properties as in  Section \ref{Betti cohomology for real one-forms}
\end{rmk}

Now we can repeat considerations of  Section \ref{Betti cohomology for real one-forms}. More precisely, let us assume that the image of the map $\int(\alpha_1+i\alpha_2): H_1(Y, \ZZ,\Z)\to \C$ does not contain $\R_{>0}$. Then there are no saddle connections, i.e. 
there is no trajectory $\phi: \R_\tau\to Y$  such that as $\tau\to +\infty$ and $\tau\to -\infty$ the family of points $\phi(\tau)$ has a limit which is a point of $\ZZ$. It follows that the results of Section \ref{non-archimedean Betti local to global} hold. In particular we have the corresponding canonical isomorphism $H^\bullet_{Betti,glob}(Y,\rho)\simeq H^\bullet_{Betti,loc}(Y,\rho)$. Here $\rho$ is a non-archimedean local system trivialized at $\ZZ$ which belongs to an appropriate non-archimedean tube domain.

Therefore one has the wall-crossing structure on the moduli space of non-archimedean rank one local systems on $Y$ trivialized at $\ZZ$, similarly to the one in Section \ref{WCS for 1-forms}. Furthermore, let us consider for a smooth manifold $S$ a smooth family $(Y_s,\alpha_{1,s}, \alpha_{2,s}), s\in S$ of triples as above endowed with a non-linear connection which identifies $\ZZ_s$ for different $s\in S$. Then we have a family of  closed $1$-forms $\beta_s$ which are obtained as generic linear combinations $c_1(s)\alpha_{1,s}+c_2(s)\alpha_{2,s}$ (this means that at each point the rank of the vector space they span is either $0$ or $2$).
Assume that:

i) For each $s\in S$ the closed $1$-forms $\alpha_{i,s}, i=1,2$ satisfy the properties a)-c) above;

ii) Using c) consider a unique smooth in $s\in S$ family of smooth functions $f_{i,s}, i=1,2$ on a sufficiently small neighborhood $U_\epsilon(\ZZ)$ of $\ZZ$ such that the restriction of $f_{i,s}$ to each connected component of $\ZZ$ is equal to zero, $df_{i,s}=\alpha_{i,s}, i=1,2$. Then the condition b) above implies that all elements of the family of maps $f_s=(f_{1,s}, f_{2,s}): U_\epsilon(\ZZ)\to \R^2=\C$ have the only critical value $0\in \C$.

Moreover we have a smooth in $s\in S$ family of local systems of abelian groups on $S^1\subset \C$ given by $\phi_{f_s}(\underline{\Z}_Y)$, where $\phi_{f_s}$ refers to the sheaf of vanishing cycles of $f_s$.  This sheaf is a direct sum of the corresponding sheaves associated with each connected component of $\ZZ$. Then we obtain a smooth over $S$ family of WCS.

We can construct a continuous family as above deforming  $\alpha_1$ and $\alpha_2$  in such a way that the common set of zeros and above-described local systems of sheaves of vanishing cycles do not change, but  the cohomology class $[\alpha_1+\sqrt{-1}\alpha_2]\in H^1(Y, \ZZ, \C)$ can change.\footnote{For that we perturb $\alpha_1$ and $\alpha_2$ by a closed relative $1$-form which is equal to zero in a neighborhood of $\ZZ$ but has non-trivial cohomology class.}
Then the assumptions i)-ii) imply that the elements of the corresponding family of WCS depend only on the cohomology class  $[\alpha_1+\sqrt{-1}\alpha_2]$. In particular, these WCS are isomorphic as long as we require in addition that this class does not change. Similarly to the example of square-tiled surfaces we get a convex cone in $H^1(Y, \ZZ,\R)$. It is acted by the group $SL_2(\R)$, so that everything looks  as a higher-dimensional analog of the theory of abelian differentials on a curve (see [KoZo]). It is tempting to ask in this framework about higher-dimensional analogs of the results from the theory of abelian differentials, e.g. those about Lyapunov exponents (see e.g. [EsKoMoeZo]).\footnote{In the higher-dimensional case one might need slightly different restrictions on the deformations. In particular, only the homotopy type of the common set of zeros $\ZZ$ does not change, and furthermore for the continuous family of closed sets $\ZZ_s$ over the deformation base $s\in S$ enjoys the continuity property rather than semicontinuity.}

\subsection{De Rham cohomology for a pair of $1$-forms}\label{de Rham for pair of 1-forms}

Finally we would like to make a similar to the above generalization  of the de Rham cohomology. 

We assume that we have a pair of closed $1$-forms as in the Section \ref{Morse-Novikov for pair of 1-forms} which satisfy the conditions a) and b) from the loc.cit. together with the following condition:

c) there exists a germ of a complex structure at the common set of zeros $\ZZ$ such that the form $\alpha:=\alpha_1+i\alpha_2$ is holomorphic.

Moreover we will assume that we are given a choice of such a complex structure.

For simplicity let us suppose  that our $C^\infty$ manifold $X$ is compact, so we do not have problems with a choice of compactification.

We define a sheaf of graded vector spaces $\underline{\Omega}^\bullet_{X,\alpha_1,\alpha_2}$ in such a way that its germ at $x\in X-\ZZ$ coincides with the germ of the sheaf $\underline{\Omega}^\bullet_{X,C^\infty}\otimes \C$ of complexified smooth differential forms, and its germ at $z\in \ZZ$ consists of holomorphic differential forms with the respect to the given complex structure.
We endow $\underline{\Omega}^\bullet_{X,\alpha_1,\alpha_2}$ with the family of differentials  $td+\alpha\wedge(\bullet), t\in \C$.

Then for each $t\in \C$ we define the de Rham cohomology
$$H^\bullet_{DR,t}(X,\alpha_1,\alpha_2)=\mathbb{H}^\bullet(X, (\underline{\Omega}^\bullet_{X,\alpha_1,\alpha_2}, td+\alpha\wedge(\bullet))).$$

In this way we obtain a coherent sheaf on $\C_t$. Moreover all the results relating de Rham global and local cohomology can be generalized to this new setting.

\subsection{Generalization to the real case symplectic case: allowable pairs of symplectic forms}\label{allowable symplectic forms}

If $X$ is a complex manifold then the standard complex symplectic structure $\omega^{2,0}_{T^\ast X}$ on the cotangent bundle $T^\ast X$ gives rise to the {\it pair} of $C^\infty$ symplectic forms $(Re(\omega^{2,0}_{T^\ast X}), Im(\omega^{2,0}_{T^\ast X}))$ which are generically $\R$-linear independent.

\begin{defn}\label{allowable pair of symplectic forms}
A pair of real symplectic forms $(\omega_1, \omega_2)$ on a smooth real manifold $M$ is called allowable if for any $(c_1,c_2)\in \R^2-\{(0,0)\}$ the form $\omega_{c_1,c_2}=c_1\omega_1+c_2\omega_2$ is symplectic.
\end{defn}

One can prove  the following fact of linear algebra (it can be also derived directly from considerations in a chosen basis of $V$).

\begin{prp}\label{pencil of 2-forms} If $V$ is a finite-dimensional real vector space endowed with a pencil $\omega_{c_1,c_2}=c_1\omega_1+c_2\omega_2$ of skew-symmetric real bilinear forms such that $\omega_{c_1,c_2}$ is non-degenerate for $(c_1, c_2)\in \R^2-\{(0,0)\}$ then:

1) dimension of $V$ is divisible by $4$;

2) There is a pair of vector subspaces $L_1,L_2$ of $V$ such that $V=L_1\oplus L_2$, each $L_i, i=1,2$ is Lagrangian with respect to both symplectic forms $\omega_1$ and $\omega_2$ and such that $\omega_i, i=1,2$ gives rise to an isomorphism $\nu_i: L_1\simeq L_2^\ast$;

3) the linear map $A=\nu_2^{-1}\nu_1:L_1\to L_1$ has the spectrum which belongs $\C-\R$.

\end{prp}

\begin{cor}\label{almost complex structure}
Under the assumptions of the above proposition there is a linear map $J: V\to V$ such that $J^2=-id_V$ and furthermore $\omega_1(Jv, v)>0$ for $v\in V-\{0\}$ and $\omega_2(Jv, v)=0$ for any $v\in V$.

\end{cor}

It looks plausible that
the set of maps $J$ from the  Corollary \ref{almost complex structure} is a {\it contractible} open subset of a closed set of such maps $J$  that $J^2=-id_V, \omega_2(Jv, v)=0$ for any $v\in V$. We did not try to prove this result.

It follows that if a smooth real manifold $M$ carries a pair of allowable symplectic structures $(\omega_1,\omega_2)$ then $dim_\R M=4n$, and there is an almost complex structure $J$ on $M$ such that $\omega_1(Jv, v)>0, v\ne 0$ and $\omega_2(Jv, v)=0$ on each tangent space.

Assume that $M$ carries a pair of allowable of symplectic forms $(\omega_1, \omega_2)$ and we are given a closed subset $L\subset M$  such that $L$ is a singular Lagrangian with respect to both $\omega_i, i=1,2$ and moreover for any $\gamma\in H_2(M,L, \Z)$ one has $\int_\gamma(\omega_1+i\omega_2)\notin \R_{>0}$. Then $L$ is {\it unobstructed} in the sense that there are no $J$-holomorphic discs with boundary on $L$ for a contractible set of almost complex structures $J$ on $M$ such that $\omega_1(Jv, v)>0, v\ne 0$ and $\omega_2(Jv, v)=0$. 

\begin{rmk}\label{embedding of local Fukaya category}
In the language of Fukaya categories which we will discuss without details in the next section, 
the observation that $L$ is unobstructed means that there exists a fully faithful embedding of the {\it local} Fukaya category $\FF_{L,loc}$  associated with $L$ into the {\it global} Fukaya category $\FF(M, \omega+iB)$, where $\omega=\omega_1, B=\omega_2$.  Replacing the form $\omega_1+i\omega_2$ by the rotated form $e^{i\theta}(\omega_1+i\omega_2)$ we obtain a family of fully faithful embeddings $i_\theta: \FF_{L,loc}\to \FF(M, e^{i\theta}(\omega_1+i\omega_2))$ as long as $\theta\in \R/2\pi \Z$ does not belong to the countable set of Stokes directions (those for which  $\int_\gamma(\omega_1+i\omega_2)\in \R_{>0}$ for some $\gamma\in H_2(M,L, \Z)$. 

Let us define Stokes isomorphisms by the formula $T_{\theta_{st}}=i_{\theta_+}^{-1}\circ i_{\theta_-}$, where $\theta_+>\theta_{st}$ (resp. $\theta_-<\theta_{st}$) is the slope of a ray which is close to $\theta_{st}$. Using the collection of Stokes isomorphisms $T_{\theta_{st}}$ one can define the wall-crossing structure similarly to the holomorphic case.
One can hope that  previously discussed results and conjectures about analyticity of this WCS in holomorphic case  can be generalized to the case of allowable pairs of real symplectic forms. 
\end{rmk}

\begin{rmk}\label{analyticity of WCS for allowable pairs}
In the notion of analytic WCS the crucial role is played by rational central charges. In the case of $C^\infty$ complex-valued closed $1$-form $\alpha_\C=\alpha_1+i\alpha_2$ the rationality constraint (up to a rational rescaling) means that both $\alpha_1$ and $\alpha_2$ have integer periods. In the case of a pair of allowable $2$-forms and a simultaneously Lagrangian subset $L$ the corresponding property means that the integrality of the cohomology classes $[\omega_i]\in H^2(M,L,\R), i=1,2$. In particular we have two prequantization bundles with connections, both trivialized over $L$. One can hope for an analog of the map to the torus $T^2$ as in Section \ref{rationality of WCS}. But instead of rationality one will probably have analyticity only.

\end{rmk}

\begin{rmk}\label{analogy with allowable metrics}
We mention here a peculiar analogy of the story with two symplectic structures as above with allowable complex metrics from [KoSe]. Indeed, assume that $J$ is integrable. Then instead of considering a pair of real symplectic forms on a $4n$-dimensional real manifold $M$ one can consider a holomorphic symplectic form $\omega_\C=\omega_1+i\omega_2$ on a complex manifold $M_\C$ (for that we need to assume that $M$ is real analytic). Then simultaneous  Lagrangian subvarieties for $\omega_1,\omega_2$ are simply $\omega_\C$-Lagrangian subvariaties in $M_\C$ (again we should consider real analytic Lagrangians in order to be able to speak about their complexification). Then allowable pairs of symplectic forms  are analogous to allowable complex metrics in the loc.cit.

Deformation of the pair $(\omega_1,\omega_2)$ in such a way that the  cohomology class $[\omega_1+i\omega_2]\in H^2(M,\C)$ is preserved is equivalent to the deformation of $M$ as a totally real submanifold of $M_\C$. Furthermore there is an analog of this observation if we replace $M$ by a pair $(M,L)$ where $L$ is Lagrangian with respect to both $\omega_1$ and $\omega_2$.

This real analog of the holomorphic story will be discussed elsewhere.
\end{rmk}


\section{Comparison isomorphisms  and 
Generalized Riemann-Hilbert Correspondence}\label{Fukaya categories}

In this section we are going to discuss conjectural categorifications of  the Betti and de Rham cohomology as well as the categorification of the comparison isomorphisms between them.
We will explain that the comparison isomorphisms between Betti and de Rham cohomology  are  incarnations of the local and global generalized Riemann-Hilbert correspondences (global GRHC for short). The exposition in this section  will be sketchy, since we are going to discuss our generalized Riemann-Hilbert correspondence  elsewhere.

\subsection{Summary of comparison isomorphisms}\label{comparison isomorphisms for 1-forms}

In this subsection we will summarize previously discussed comparison isomorphisms between de Rham and Betti cohomology. Later we will explain how they fit into a more general categorical picture.

Let  us start with a summary of  the cohomology theories associated with the pair $(X,\alpha)$, where $X$ is a complex manifold and $\alpha$ is holomorphic closed $1$-form, and then discuss the relationships between them. 
So far we have defined the following:

i) The {\it global} de Rham cohomology as a graded algebraic coherent sheaf $\mathcal{H}_{DR,glob}^\bullet(X,\alpha)$ on $\C_t$. The space of global sections ${H}_{DR,glob}^\bullet(X,\alpha)$ is a finitely generated $\C[t]$-module.

\

ii) The local de Rham  cohomology space $H_{DR,loc}^\bullet(X,\alpha)$ as graded finitely generated $\C[[t]]$-module.

\

iii) The global  Betti cohomology  $\mathcal{H}^\bullet_{Betti,glob}(X,\alpha)$ as a graded analytic coherent sheaf on $\C^\ast_t$, which is a graded holomorphic vector bundle over a germ of punctured disc.
The space of its global sections gives rise to a graded free module $H^\bullet_{Betti,glob}(X,\alpha)$ over the algebra $\mathcal{A}:=\lim_{\epsilon\to 0}\OO^{an}(D_\epsilon-\{0\})$, where $D_\epsilon$ is a disc of radius $\epsilon$ with center at $t=0$.

\

iv) The locally constant sheaf of graded finitely generated abelian groups $\mathcal{H}_{Betti,loc,\Z}^\bullet(X,\alpha)$ of local Betti cohomology as well as the corresponding graded algebraic vector bundle  $\mathcal{H}_{Betti,loc}^\bullet(X,\alpha)$ on $\C^\ast_t$ endowed with a connection with regular singularieties at $t=0$ and $t=\infty$. The  space of global sections of the latter is denoted by ${H}_{Betti,loc}^\bullet(X,\alpha)$ (we forget about the connection and consider it just as a $\C[t,t^{-1}]$-module). By changing  scalars we obtain a $\C((t))$-module $H^\bullet_{Betti,loc}(X,\alpha)\otimes_{\C[t,t^{-1}]}\C((t))$. 

\

v) Using the analytic WCS associated with $\alpha$ we constructed a graded holomorphic vector bundle $\mathcal{H}_{Betti,loc,WCS}^\bullet(X,\alpha)$ over a germ of punctured disc as well with a meromorphic structure at $t=0$ as well as the corresponding graded $\C\{t\}[t^{-1}]$-module ${H}_{Betti,loc,WCS}^\bullet(X,\alpha)$ of its global sections. We recall that $\mathcal{H}_{Betti,loc,WCS}^\bullet(X,\alpha)$ does not carry a natural connection.

\

We  know that
the $\C[t]$-module ${H}_{DR,glob}^\bullet(X,\alpha)$ of  global sections of 
$$\mathcal{H}_{DR,glob}^\bullet(X,\alpha)$$ satisfies the property that the natural morphism 
$$\phi_{DR}: {H}_{DR,glob}^\bullet(X,\alpha)\otimes_{\C[t]} \C[[t]]\to{H}_{DR,loc}^\bullet(X,\alpha)$$
 is an isomorphism. We keep the same notation for the isomorphism obtained from this one by extension of scalars to $\C((t))$.

We also have  isomorphisms: 
$$\phi_{Betti,loc}:{H}_{Betti,loc,WCS}^\bullet(X,\alpha)\otimes_{\C\{t\}[t^{-1}]} \C((t))$$
$$\simeq {H}_{Betti, loc}^\bullet(X,\alpha)\otimes_{\C[t,t^{-1}]} \C((t)),$$
and 
$$\phi_{WCS}: {H}_{Betti,loc,WCS}^\bullet(X,\alpha)\otimes_{\C\{t\}[t^{-1}]}{\mathcal A}\simeq H_{Betti,glob}^\bullet(X,\alpha).$$

The isomorphism $\phi_{WCS}$ induces a meromorphic structure at $t=0$ on the holomorphic vector bundle $\mathcal{H}^\bullet_{Betti,glob}(X,\alpha)$. We denote the corresponding $\C\{t\}[t^{-1}]$-module of germs of meromorphic sections by  $H_{Betti,glob}^\bullet(X,\alpha)^{mer}$. Then we have an isomorphism of $\C\{t\}[t^{-1}]$-modules

$$\phi_{WCS}^{mer}: {H}_{Betti,loc,WCS}^\bullet(X,\alpha)\simeq H_{Betti,glob}^\bullet(X,\alpha)^{mer}.$$

 We denote by $\phi_{Betti}$ the   isomorphism obtained from $\phi_{Betti,loc}\circ (\phi_{WCS}^{mer})^{-1}$ by extension of scalars:
 $$\phi_{Betti}: {H}_{Betti,glob}^\bullet(X,\alpha)^{mer}\otimes_{\C\{t\}[t^{-1}]} \C((t))\simeq {H}_{Betti, loc}^\bullet(X,\alpha)\otimes_{\C[t,t^{-1}]} \C((t)).$$

Furthermore, we have the global Riemann-Hilbert correspondence (in the terminology of this paper  it is the global de Rham-to-Betti isomorphism) which amounts to an isomorphism of analytic coherent sheaves over $\C^\ast_t$
$$\mathcal{RH}_{glob}: \mathcal{H}_{DR,glob}^\bullet(X,\alpha)\simeq  \mathcal{H}_{Betti,glob}^\bullet(X,\alpha).$$
Hence we have an isomorphism of the corresponding modules of sections:
$$RH_{glob}: {H}_{DR,glob}^\bullet(X,\alpha)\otimes_{\C[t]}{\mathcal A}\simeq  {H}_{Betti,glob}^\bullet(X,\alpha).$$

Also  the local Riemann-Hilbert correspondence $RH_{loc}$ for the {local}  de Rham and Betti cohomology gives rise to an isomorphism: 
$$RH_{loc}: {H}_{DR,loc}^\bullet(X,\alpha)\otimes_{\C[[t]]} \C((t))\simeq {H}_{Betti,loc}^\bullet(X,\alpha)\otimes_{\C[t,t^{-1}]} \C((t)).$$ 

Let us combine all the above information and formulate the following {\it comparison of isomorphisms conjecture}.

\begin{conj} \label{comparison of isomorphisms diagram} 

1) The isomorphism $\mathcal{RH}_{glob}$ respects the meromorphic structures at $t=0$, thus inducing the isomorphism of $\C\{t\}[t^{-1}]$-modules of sections
$$RH_{glob}^{mer}: {H}_{DR,glob}^\bullet(X,\alpha)\otimes_{\C[t]}\C\{t\}[t^{-1}]\simeq  {H}_{Betti,glob}^\bullet(X,\alpha)^{mer}.$$

2) Let us denote by ${RH}_{glob}^{form}$  the extension of $RH_{glob}^{mer}$ obtained by changing scalars to $\C((t))$. 
Then the following diagram is commutative
\[ \begin{tikzcd}
{H}_{DR,glob}^\bullet(X,\alpha)\otimes_{\C[t]} \C((t))\arrow{r}{{RH}_{glob}^{form}} \arrow[swap]{d}{\phi_{DR}}&{H}_{Betti, glob}^\bullet(X,\alpha)^{mer}\otimes_{\C\{t\}[t^{-1}]} \C((t))  \arrow{d}{\phi_{Betti}}  \\%
{H}_{DR, loc}^\bullet(X,\alpha)\otimes_{\C[[t]]} \C((t)) \arrow{r}{RH_{loc}}& H_{Betti,loc}^\bullet(X,\alpha)\otimes_{\C[t,t^{-1}]} \C((t))
\end{tikzcd}
\]

\end{conj}


\begin{rmk}
The Conjecture \ref{comparison of isomorphisms diagram} can be thought of as a version of  the Deligne-Malgrange Riemann-Hilbert correspondence in the case of irregular $D$-modules associated with closed $1$-forms. 
\end{rmk}

We notice that the Conjecture \ref{comparison of isomorphisms diagram} is based on the highly non-trivial property of analyticity of the WCS associated with a holomorphic closed $1$-form. It is not clear for us at this time how to approach the Conjecture \ref{comparison of isomorphisms diagram}.

\begin{rmk}\label{generalization to non-trivial connections}
The above-discussed story about comparison isomorphisms admits the following immediate generalization. Instead of the family of connections $d+{\alpha \over t}$ let us consider a local system of free rank $1$ modules over $\C\{t\}$ and twist the corresponding family of connections by adding  ${\alpha \over t}$. Then the Conjecture \ref{comparison of isomorphisms diagram} has a generalization with the following modification of two opposite vertices of the commutative diagram: for global de Rham cohomology is now a module over $\C\{t\}$ rather than over $\C[t]$, and local Betti cohomology is now a module over $\C\{t\}[t^{-1}]$ rather than over $\C[t,t^{-1}]$.

\end{rmk}

\subsection{Fukaya categories of  complex symplectic manifolds}\label{Fukaya categories of complex manifolds}

 Let $M$ be a complex manifold of dimension $2n$ endowed with a holomorphic symplectic form $\omega^{2,0}$. We split $\omega^{2,0}$ as a sum $\omega+iB$. Here $\omega=Re(\omega^{2,0})$ endows $M$ with a {\it real} symplectic manifold and the closed $2$-form $B=Im(\omega^{2,0})\in H^2(M,\R)$ gives rise to the $B$-field. The latter is the cohomology class  $[B]_{U(1)}\in H^2(M,U(1))=H^2(M,\R/2\pi \Z)$ obtained via the exponential map. 
 
Fukaya category  of  a general non-compact complex symplectic manifold $M$  depends on  additional choices, which include a {\it partial} compactification which we call {\it log extension} $M_{log}$. \footnote{The notion of log extension generalizes the notion of log symplectic manifold.} For the sake of brevity we will omit the dependence of the Fukaya category on $M_{log}$ from the notation. In general a log extension is non-unique. By definition it  is a Poisson manifold containing $M$ such that $M\subset M_{log}$ is an open symplectic leaf. Furthermore it is required that  $M_\infty:=M_{log}-M$ is a normal crossing divisor such that at any point of $M_\infty$ there are local coordinates on $M$ such that 
$$M_\infty=\{(x_1,...,x_{2n})|\prod_{1\le i\le k}x_i=0\}, \quad 1\le k\le n,$$

$$\omega^{2,0}=\sum_{1\le i\le k}{dx_i\over {x_i}}\wedge dx_{i+n}+\sum_{k+1\le i\le n}dx_i\wedge dx_{i+n}.$$

The action of $\C^\ast$ on the holomorphic symplectic form given by $\omega^{2,0}\mapsto \omega^{2,0}_t:=\omega^{2,0}/t, t\in \C^\ast$ gives rise to a family of {\it global} Fukaya categories with the fibers $\FF_t(M)=\FF(\omega_t,[B_t]_{U(1)})$, where $\omega_t=Re(\omega^{2,0}/t)$ and $B_t=Im(\omega^{2,0}/t)$. \footnote{Sometimes we will abuse the notation and denote this category by $\FF(\omega_t+i B_t)$. Moreover the initial notation $\FF(\omega_t,[B_t]_{U(1)})$ is also slightly misleading, since the properly defined category should depend on a {\it representative} of the $B$-field.} In fact this is an {\it analytic} family of $\A$-categories.\footnote{The notion of analytic family of $\A$-categories will be also discussed elsewhere. Roughly such a family over an analytic stack $S$ is a sheaf of $\A$-categories enriched over the sheaf of symmetric monoidal categories ${Perf}_S$ of  perfect complexes of $\OO_S$-modules.}

\begin{rmk}\label{gerbe}
There is a useful generalization of the above considerations. Namely let us choose $[B_{init}]\in H^2(M,U(1))$ which we treat as a class of a $U(1)$-gerbe on $M$. Then we have the corresponding family of global Fukaya categories $\FF(\omega_t,[B_t]_{U(1)}+[B_{init}])$.

\end{rmk}

Some (but not all) objects of $\FF_t(M)$ are pairs $(L,\rho)$, where $L$ is a $C^\infty$  unobstructed oriented Lagrangian submanifold of $(M,\omega_t)$ with vanishing Maslov class,  such that the closure of $L$ in $M_{log}$ is compact (and satisfying some other properties near $M_\infty$ which we skip here), and $\rho$ is a smooth complex vector bundle on $L$ endowed with a unitary connection whose curvature is given by the $B$-field. More precisely $\rho$ should belong to a category of local systems on $L$ twisted by the $\C^\ast$-gerbe on $L$ with the class 
$$([B_t]_{U(1)}+[B_{init}])_{|L}+\overline{St}_2(L),$$ 
where $\overline{St}_2(L)$ is the image of the second Stieffel-Whitney class $St_2(L)\in H^2(L, \Z/2\Z)$ under the natural homomorphism $H^2(L, \Z/2\Z)\to H^2(M,U(1))$.

\begin{rmk}\label{unobstructedness}
Unobstructedness of $L$ means that  there are no non-trivial pseudo-holomorphic discs with boundary on $L$ with respect to a compatible with $\omega_t$ almost complex structure. Examples of unobrstructed Lagrangians include exact Lagrangians in exact symplectic manifolds as well as graphs of closed $1$-forms in cotangent bundles. In the first case exactness of $L$ implies that for any non-trivial pseudo-holomorphic disc $\phi: (D,\partial D)\mapsto (M,L)$ the area $\int_D\phi^\ast(\omega)=0$, hence $\phi$ must be a constant map. In the second case $L=graph(\alpha)$ the vanishing of area follows from the following calculation:
$$\int_D\phi^\ast(\omega)=\int_{D}\phi^\ast(d\lambda_{Liov})=\int_{\partial D}(\pi\circ\phi)^\ast (\alpha)=\int_D(\pi\circ\phi)^\ast (d\alpha)=0,$$
where $\lambda_{Liov}$ is the Liouville form, and $\pi :T^\ast X\to X$ is the natural projection.

\end{rmk}

Morphisms in $\FF_t(M)$ are defined using the partial wrapping near the divisor $M_\infty$. More precise definition of an appropriate Fukaya category as well as the  inclusion to the story the case of singular Lagrangian subvarieties $L$ will be discussed elsewhere.

\begin{rmk} \label{local Fukaya category in general}
In what follows we will also need the notion of local Fukaya category. In this paper we will need its simplest version. Similarly to the notion of the global Fukaya category it can be formulated in $C^\infty$ setting. 
Having a possibly singular Lagrangian subvariety $L$ of a real symplectic manifold we can define  the partially wrapped  Fukaya category with appropriate stops (see [Syl]) of a sufficiently small Liouville (better Weinstein) neighborhood of $L$.  It is expected that such a neighborhood always exists. In the holomorphic case assuming that $L$ itself is a holomorphic Lagrangian subvariety, one can take a Stein neighborhood of $L$. Then the local Fukaya category $\FF_{L,loc}$  is a  category of modules of finite rank over a $\Z$-graded $\A$-algebra over $\Z$ which is the endomorphism algebra of a generator of the partially wrapped Fukaya category.

\end{rmk}

\subsection{The case of cotangent bundle}\label{Fukaya and Betti}

For $M=T^\ast X$ where $X$ is a smooth complex algebraic variety endowed with a regular function $f$ we propose the following choice of $M_{log}$ having in mind applications to exponential integrals. Recall the set up of Section \ref{data for exp int}. For simplicity we will assume that $D_0=\emptyset$. First we choose a smooth compactification $\overline{X}$ of $X$, such that $D=\overline{X}-X$ is a simple normal crossing divisor. We assume that $D=D_v\sqcup D_h$ where at the components of the ``vertical divisor" $D_v$ the function $f$ has poles of order at least $1$, and it does not have poles at the components of the ``horizontal" divisor $D_h$. Moreover, we assume that the restriction of  $f$ to all strata of $D_h-D_v$ has no critical points. In what follows we will choose $M_{log}$ in such a way that the the closures of $L_0=X$ and $L_1=graph(df)$
 in $M_{log}$ are compact, and their intersection set can be identified with $Crit(f)\subset T^\ast X$.\footnote{Here we use the assumption that $f$ has no critical points on $D_h$.} 
 

Let $M_1=T^\ast_{log\,D}\overline{X}$ be the logarithmic cotangent bundle of $X$
(i.e. its sections are logarithmic $1$-forms with respect to $D$).

Let now $g: T^\ast _{log\, D_h}(\overline{X}-D_v)\to T^\ast_{log\, D_h} (\overline{X}-D_v)$ be the  automorphism given by the shift by $df$, i.e. $g(x,p)=(x, p+df(x))$  on $T^\ast X$ and the natural extension of this formula to $D_h$. Let us define $M_{log}$ by gluing two copies of $M_1$ along the common open part $ T^\ast_{log\, D_h} (\overline{X}-D_v)$ to which we apply the automorphism $g$, i.e. $M_{log}=M_1\cup_g M_1$. Using valuative criterion one can show that $M_{log}$ is a smooth separated scheme, hence it is a smooth algebraic variety.  It is easy to see that indeed the closures $\overline{L}_i\subset M_{log}, i=0,1$ are compact.

It is natural for several reasons  to modify the $B$-field by adding to it another $B$-field corresponding to the natural pull-back to $M=T^\ast X$ of the Stieffel class $St_2(X)\in H^2(X,\Z_2)$. The corresponding $B$-field $[B_{init}]\in H^2(M, U(1))$ does not depend on $t$. The modified $B$-field is $[B_t]_{U(1)}+[B_{init}]\in H^2(M, U(1))$.

Then we have the corresponding analytic family of Fukaya categories $\FF_t(M)=\FF(\omega_t, [B_t]_{U(1)}+[B_{init}])$ associated with the above data.

Notice that  $M=T^\ast X$ is an {\it exact symplectic manifold}: $\omega^{2,0}=d\lambda$, where $\lambda$ is the holomorphic Liouville form. Consequently it is exact as a real symplectic manifold. Lagrangian submanifolds $L_0$ and $L_1$ are exact Lagrangians and hence they are unobstructed. The class of the $U(1)$-gerbe on $L_i, i=0,1$ given by $([B_t]_{U(1)}+[B_{init}])_{|L}+\overline{St}_2(L)$ is trivial. Indeed the restriction of $B_t$ to either of these submanifolds (considered as real Lagrangian submanifolds of $(M,\omega_t)$) is zero, and $\overline{St}_2(L)=[B_{init}]_{|L}$ and $2St_2(L)=0$. Therefore any local system on $L_i, i=0,1$ gives an object of $\FF_t(M)$. In particular trivial rank one local systems give such objects which we denote by abusing notation by $L_0$ and $L_1$.

Since  $L_0\cap L_1\subset M$ in order to compute $Hom_{\FF_t(M)}(L_0,L_1)$ we do not need the partial wrapping. Instead it suffices to consider small Hamiltonian perturbations with Hamiltonian functions which are compactly supported on $X$, and the support contains $L_0\cap L_1=Crit(f)$. We may assume that the corresponding symplectomorphisms are ``very close" to the identity.
It follows that $Hom(L_0, L_1)$ in the Fukaya category $\FF_t(M)$ can be defined over $\Z$. 
In the special case when  $f$ has only Morse critical points we obtain that $Hom(L_0, L_1)\simeq\Z^{L_0\cap L_1}[-dim_\C\,X]$ (the isomorphism is not canonical for $t$ belonging to Stokes rays).

\begin{prp} \label{Fukaya and Betti comparison} In the above notation for any $t\in \C^\ast$ there is a natural isomorphism 
$$Ext_{\FF_t(M)}^\bullet(L_0, L_1)\simeq H^\bullet(X, f^{-1}(z),\Z)$$
 as long as $z$ belongs to the ray $\R_{<0}\cdot t$ and $|z|$ is sufficiently large.

\end{prp}
{\it Sketch of the proof.}  For simplicity let us assume that $t\in \R_{>0}$.

We start with some general considerations in the $C^\infty$ case. For that reason we change the notation  and speak about a compact $C^\infty$ manifold $Y$ and its cotangent bundle $T^\ast Y$. In the definition of the Fukaya category of $T^\ast Y$ which we briefly outlined above, we use the $B$-field $[B_{init}]$ which is the pullback to $T^\ast Y$ of $St_2(Y)$. Then the zero section $Y$ gives rise to an object of the Fukaya category $\FF(T^\ast Y)$.
Recall that  for  a smooth compact manifold $Y$ the cohomology $H^\bullet(Y, \Z)$ is isomorphic to $Ext^\bullet_{\FF(T^\ast Y)}(Y,Y)$. We can relax the condition that $Y$ is compact such as follows. Let now $Y$ be the interior of a smooth manifold $\overline{Y}$ with corners $\partial\overline{Y}$, and $g: \overline{Y}\to \R$ be a  smooth {\it proper}  Morse function such that the restriction of $g$ on all boundary strata has no critical points, and such that the cardinality $|Crit(g)|$ of the set of critical points of $g$ is finite. Then the Morse complex $Morse^\bullet(g)$ computes the cohomology $H^\bullet(Y, g^{-1}(-\infty)):=H^\bullet(Y, g^{-1}(-z)),\Z)$, where $z>0$  is sufficiently large.
On the other hand $Morse^\bullet(g)$ is quasi-isomorphic to  $Hom_{\FF(T^\ast Y)}(Y, graph(dg))$.
This can  be shown similarly to [KoSo6], where the corresponding conjecture of Fukaya and Oh was proved. Therefore we see that $Ext^\bullet_{\FF(T^\ast Y)}(Y, graph(dg))\simeq H^\bullet(Y, g^{-1}(-\infty)))$.

Let us return to our case when $X$ is a smooth  algebraic variety over $\C$. Let us first assume that
$f: X\to \C$ is a regular function with finitely many Morse critical points. 

Let $Y:=X_\R$ denote the  $C^\infty$ manifold underlying $X$. Then the real symplectic manifold $T^\ast X_\R$ is isomorphic to $(M,Re(\omega^{2,0}))$. Under this isomorphism $graph(df)\subset M$ is identified with $graph(d(Re(f)))\subset T^\ast X_\R$.
Moreover this induces the isomorphism for any $t\in \R_{>0}$
$$Ext^\bullet_{\FF(T^\ast X_\R)}(X_\R,graph\, d(Re(f)))\simeq Ext^\bullet_{\FF_t(M)}(L_0,L_1).$$

Similarly we have an isomorphism
$$H^\bullet(X, f^{-1}(z),\Z)\simeq H^\bullet(X_\R, Re(f)^{-1}(z),\Z).$$

In order to apply the above considerations   we need an appropriate partial compactification $\overline{X}_\R$ such that $Re(f)$ extends to a proper function without critical points on the boundary.  Namely we define $\overline{X}_\R$  as the closure of $graph(f)\subset X\times \C$ in the space $\overline{X}_{cor}\times \R\times i\overline{\R}$, where $\overline{\R}:=[-\infty,+\infty]\supset \R$ and $\overline{X}_{cor}$ is the oriented real blow-up of $\overline{X}$ at $\overline{X}-X$ (see Section \ref{compactification for 1-forms}).
Then the extended function $Re(f)$ does not have critical points on the boundary strata since the holomorphic function $f_{|D_h}$ has no critical points. 


Finally, if  $f$ is not Morse we can replace $Re(f)$ by its small perturbation which is Morse, i.e.  define $g_1=Re(f)+h$, such that $h$ has compact support, and $g_1$ is Morse. The $graph(dRe(f))$ which is Hamiltonian isotopic to $L_1$ should be replaced by the $graph(dg_1)$ which is also Hamiltonian isotopic to $L_1$. Function $Re(f)$ is replaced by $g_1$ and the rest of the above argument is changed in the obvious way. $\blacksquare$

\begin{rmk}\label{local version}
There is a similar isomorphism for the local Fukaya category $\FF_{L_0\cup L_1,loc,t}:=\FF(U_\epsilon(L_0, L_1), \omega_t+i(B_t+B_X))$ associated with a small Liouville neighborhood $U_\epsilon(L_0, L_1)$ of the union $L_0\cup L_1$ with $\A$-structure defined by means of the infinitesimal wrapping. The corresponding $Ext$-groups are canonically isomorphic to the cohomology groups of the critical locus of $f$ (the critical locus coincides with $L_0\cap L_1$), with coefficients in the sheaf of vanishing cycles of $f$.

\end{rmk}

\subsection{Betti global-to-local isomorphism and wheels of projective lines}\label{comparison isomorphisms and wheels of lines}

In this subsection we are going to utilize  the description of stability data  and wall-crossing structures in terms of wheels of $\C{\bf P}^1$'s in a toric or toric-like varieties proposed for the first time in [KoSo7] and later revisited in [KoSo12]. In [KoSo12]  we also introduced the notion of  analytic wall-crossing structure and explained how it can be described in terms of the wheels of projective lines. In order to simplify the exposition we are going to discuss a special case of stability data and analytic stability data of rank 2 only, instead of more general wall-crossing structures.

In [KoSo12], Section 2.8 we considered stability data on the graded Lie algebra   of algebraic vector fields on the torus $(\C^\ast)^n$. This Lie algebra contains the graded Lie subalgebra of vector fields which are linear in the first $m\le n$ variables, which is the same as the Lie algebra of $\mathfrak{gl}(m)$-valued functions on $(\C^\ast)^{n-m}$. The  graded Lie algebras in the present paper are of this type (see Section \ref{WCS for 1-forms}).

We start with a wheel of strict rational cones in the vector space $\Gamma\otimes \R,  \,  \Gamma=H^1(X,\ZZ(\alpha),\Z)$ such that if we apply to each of the cones the central charge $Z:\gamma\mapsto \int_\gamma \alpha$ we obtain a strict sector in $\R^2=\C$. The support of our stability data is contained in the union of these cones. Then  the corresponding group elements $A_V$ can be used (cf. [KoSo12]) to glue a vector bundle
on the formal neighborhood $\XX^{form}_\C$ of a wheel of projective lines $\cup_i{\bf P}^1_i$ in a toric variety containing an open orbit isomorphic to the torus $Hom(\Gamma, {\bf G}_m)$. 
Analyticity of the stability data (see Sections \ref{WCS for 1-forms}, \ref{analyticity of WCS}) implies that this bundle extends to the germ $\XX^{an}$ of a $\C$-analytic neighborhood of the wheel of projective lines. 

Integrality of the group elements (see Remark \ref{extension to coherent sheaves}) can be reformulated as a statement that this bundle is the pullback of the one on the formal scheme $\XX^{form}_\Z$ over $\Z$ which is the formal neighborhood of the wheel of projective lines $\cup_i{\bf P}^1_{\Z,i}$ in the toric variety understood as a scheme over $\Z$.  Then we have  a surjective morphism $\XX^{form}_\C\twoheadrightarrow \XX^{form}_\Z$ which is equal to $id$ on the wheel.

Therefore our  stability data can be equivalently described in the following way:

1) Finite rank local system of lattices $\EE_{S^1}$ on the circle $S^1$.

2) Analytic vector bundle $\EE^{an}$ on $\XX^{an}$.

3) Bundle of locally-free $\Z$-modules $\EE_\Z^{form}$ on $\XX_\Z^{form}$ endowed with:

3a) an isomorphism
$$(\XX^{form}_\C\hookrightarrow  \XX^{an})^\ast(\EE^{an})\simeq (\XX^{form}_\C\twoheadrightarrow \XX^{form}_\Z)^\ast(\EE_\Z^{form}),$$

3b) an isomorphism of the bundle of locally-free $\Z$-modules
$(\cup_i {\bf P}^1_{\Z,i}\to \XX_\Z^{form})^\ast(\EE_\Z^{form})$ with the one on the wheel $\cup_i {\bf P}^1_{\Z,i}$ which is trivial on each ${\bf P}^1_{\Z,i}$ and is naturally identified with the local system $\EE_{S^1}$.\footnote{Identification comes from the fact that a bundle of locally-free $\Z$-modules on the wheel $\cup_i {\bf P}^1_{\Z,i}$ which is trivial on each ${\bf P}^1_{\Z,i}$ gives rise to a collection of free $\Z$-modules on the cyclically ordered set of double points of the wheel as well as isomorphisms of these modules for any adjacent pair of such points.}

In our  case these data can be derived from those in Section \ref{WCS for 1-forms}. In particular $\EE_{S^1}=\mathcal{H}^\bullet_{Betti,loc,\Z}(X,\alpha)$. {\it We denote by $\mathcal{H}_{Betti, loc, WCS}^\bullet(X,\alpha)$ the corresponding analytic vector bundle on $\XX^{an}$}. This notation agrees with the similar one in the previous subsection.

With our central charge $Z: \gamma\mapsto \int_\gamma \alpha$ we associate a map $exp_Z: \C^\ast_t\to \XX^{an}$ which on the open subtorus $(\C^\ast)^n\subset \XX^{an}$ can be described as the composition of the map $t\mapsto {1\over t}\int \alpha: H_1(X,\ZZ(\alpha),\Z)\to \C$ with the exponential map $exp: \C\to \C^\ast$. 

Then Betti global-to-local isomorphism means that
$$exp_Z^\ast(\mathcal{H}_{Betti,loc,WCS}^\bullet)=\mathcal{H}_{Betti,glob}^\bullet(X,\alpha).$$

Furthermore the graded vector bundle $\mathcal{H}_{Betti,glob}^\bullet(X,\alpha)$ is endowed with a meromorphic structure at $t=0$ induced from the natural meromorphic structure on the LHS (i.e. we can speak about analytic sections with no more than polynomial growth at $t=0$). In particular the corresponding space of sections is a module over $\C\{t\}[t^{-1}]$.

Via the map $exp_Z$  one assigns to an analytic bundle $\EE$ on $\XX^{an}$ satisfying the property  that the restriction $\EE_{|\C{\bf P}^1_i}$ is trivial for any projective line $\C{\bf P}^1_i$,  the following data:

(a) Local system of $\C$-vector spaces $\EE_{S^1}$ on $S^1$.

(b) Free $\C\{t\}[t^{-1}]$-module $E^{mer}$ together with an isomorphism of $\C((t))$-modules
$$E^{mer}\otimes_{\C\{t\}[t^{-1}]}\C((t))\simeq RH_{loc}^{-1}(\EE_{S^1}),$$
where $RH_{loc}^{-1}(\EE_{S^1})$ is the vector bundle with the regular singular connection on $Spec(\C((t)))$ (equivalently free $\C((t))$-module) corresponding to the local system $\EE_{S^1}$ under the local Riemann-Hilbert correspondence.
Then the pull-back $exp_Z^\ast(\EE)$ is naturally identified with $E^{mer}\otimes_{\C\{t\}[t^{-1}]}\mathcal{A}$ where the ring $\mathcal{A}$ was defined in the previous subsection (see condition (iii) there).

\begin{rmk}\label{limits of diagrams}

More abstractly the data (a), (b) are obtained via the  pull-backs associated with the natural functors
$$\OO^{an}_{\C^\ast}-mod \leftarrow \C\{t\}[t^{-1}]-mod\to \C((t))-mod\leftarrow Conn(Spec(\C((t)))\simeq Loc(S^1),$$
where the very right functor is the local Riemann-Hilbert correspondence on the formal punctured disc.

\end{rmk}

In the following section we will ``categorify" this description so that vector bundles will be replaced by categories.

\subsection{Families of Fukaya categories and wheels of projective lines}\label{analytic families of Fukaya categories}

In this subsection we are going to discuss a conjectural categorification of the Betti global-to-local isomorphism in the framework of Section \ref{comparison isomorphisms and wheels of lines}. We will work in a bigger generality than just the one of cotangent bundles. This should explain to the reader the relationship between local and global categories in general Holomorphic Floer theory. The detailed discussion of this topic will be presented elsewhere. For that reason we will omit many details and precise definitions. Nevertheless we think that the discussion below will be helpful.

Let $(M,\omega^{2,0})$ be a complex symplectic manifold and $L\subset M$ a holomorphic Lagrangian subvariety (possibly singular).  As before we use the notation $\omega_t=Re(\omega^{2,0}/t), B_t=Im(\omega^{2,0}/t),t \in \C^\ast$. Then one has the {\it analytic} family of $\A$-categories $\FF_{glob}^{an}:=\FF_M$ of global Fukaya categories on $\C^\ast_t$ with the fiber $\FF_{glob,t}^{an}=\FF_{M,t}:=\FF_t(M)=\FF(M,Re(\omega^{2,0}/t)+iIm(\omega^{2,0}/t))$. We also have the {\it formal} family $\FF^{form}_{glob}$ over the Novikov ring $Nov_\Z=\Z[[T^\C]]$ (see Definition \ref{Novikov rings}). 

On the other hand, with the Lagrangian subvariety $L$ one can associate a local system on $S^1$ of $\A$-categories $\FF_{loc,\Z}:=\FF_{L,loc,\Z}$ over $\Z$.
Moreover, if $\theta$ is not the direction of a Stokes ray, i.e. does not belong to the set $\{\int_\gamma \omega^{2,0}|\gamma\in H_2(M,L,\Z)\}$ then there is a fully faithful embedding of $\A$-categories which are fibers of our families over $t=e^{i\theta}$:

$$i_\theta: \FF_{loc,\Z,\theta}\otimes_\Z \Z[[T^\C]]\to \FF_{glob,\theta}^{form}.$$

Setting $\Gamma=H_2(M,L,\Z)$ and taking as the central charge $Z:\Gamma\to \C$ given by $Z(\gamma)=\int_\gamma \omega^{2,0}$ we see that the above functor can be thought of as a categorification the Betti global-to-local isomorphism.
One can then ask about categorification of the wall-crossing structure. In order to do that one should ensure the Support Property.
The latter means that there exists a conical closed subset $Supp\subset \Gamma\otimes \R$ such that the intersection of $Supp$ with the kernel $Ker\,Z_\R$ of the $\R$-linear extension of $Z$ to $\Gamma\otimes \R$ is equal to zero. We should assume that the homology classes of pseudo-holomorphic discs with respect to a compatible with $\omega^{2,0}$ pseudo-hyperk\"ahler structure belong to $Supp$.

Then the categorification of the Stokes isomorphisms are the elements
$g_\theta=i_{\theta_-}^{-1}\circ i_{\theta_+}$ for $\theta\in Z(\Gamma)$ and $\theta_+$ (resp. $\theta_-$) being the slope of the ray which is on the right (resp. on the left) of the Stokes ray with the slope $\theta$.

As we will explain in one of the subsequent papers on the project in fact $g_\theta$ comes via a pullback of an element $\hat{g}_\theta\in Aut(\FF_{loc,\Z,\theta}\otimes_\Z \Z[[T^N]])$, where $N$ denote the monoid $Z^{-1}(e^{i\theta}\cdot \R_{\ge 0})\cap Supp\cap \Gamma.$

The data $\hat{g}_\theta$ are equivalent to the family of $\A$-categories $\FF_{loc,WCS}:=\FF_{L,loc,WCS}$ on $\XX_\Z^{form}$ such that its restriction to the wheel $\cup_i{\bf P}^1_i$ is obtained in the natural way from the local system $\FF_{loc,\Z}$ on $S^1$ (cf. a similar discussion in Section \ref{comparison isomorphisms and wheels of lines}). 

Now we can formulate the conjecture  which can be thought of as  a categorification of the analytic WCS revisited in  Section \ref{comparison isomorphisms and wheels of lines}:
\begin{conj}\label{equivalence of pull-backs}
There exists canonical analytic family of $\A$-categories $\FF_{loc}^{an}:=\FF_{L,loc}^{an}$ on $\XX^{an}$ and an equivalence of sheaves of categories
$$(\XX^{form}_\C\to \XX^{an})^\ast(\FF_{loc}^{an})\simeq (\XX^{form}_\C\to \XX_\Z^{form})^\ast(\FF_{loc,WCS}).$$
\end{conj}

Next conjecture can be thought of as a  categorical analog of the Betti global-to-local isomorphism in the set-up of 
Section \ref{comparison isomorphisms and wheels of lines}:

\begin{conj}\label{equivalence of pull-backs on the line}
There is a fully faithful embedding of families of categories
$$exp_Z^\ast(\FF_{loc}^{an})\hookrightarrow \FF_{glob}.$$

\end{conj}

This embedding induces a fully faithful functor on the categories which are  global sections of our families. 

Notice that an embedding $L\to L^\prime$ induces a fully faithful embedding of local Fukaya categories endowed with autoequivalences $\hat{g}_\theta$, i.e. $\FF_{L,loc,WCS}\hookrightarrow \FF_{L^\prime,loc, WCS}$.

\begin{rmk}\label{relation to 1-forms}
In the case of closed $1$-form $\alpha$ we take $L=L_0\cup L_1$, where $L_0=X, L_1=graph(\alpha)$. There are no pseudo-holomorphic discs with boundary on $L_i, i=0,1$, hence they define canonical objects $E_0^{Betti},E_1^{Betti}$ and $E_0^{Betti,loc},E_1^{Betti,loc}$\footnote{In the previous subsection we sloppily denoted by $L_i, i=0,1$ the objects which we now denote by  $E_i^{Betti}, i=0,1$.} in the  categories of global sections of the corresponding families on $\C^\ast$. E.g. $E_0^{Betti},E_1^{Betti}\in Ob(\Gamma(\C^\ast,\FF_{glob}))$. Then Stokes isomorphisms $g_\theta$ give rise to the automorphisms of $Ext^\bullet (E_0^{Betti,loc}, E_1^{Betti,loc})$. The latter appear in the definition of the WCS in Section \ref{WCS for 1-forms}.

\end{rmk}

Summarizing, we have  four different types of Fukaya categories:

i) $\FF_{loc,\Z}$, which is a local system of $\Z$-linear $\A$-categories on $S^1$ (equivalently, we can speak about a local system on $\C^\ast$  equivariant with respect to dilations).

ii) $\FF_{loc,WCS}$ which is a family of $\Z[[T^N]]$-linear $\A$-categories on $\XX_\Z^{form}$.

iii) $\FF_{loc}^{an}$ which is an analytic family of $\A$-categories over $\XX^{an}$.

iv) $\FF_{glob}$ which is an analytic family of $\A$-categories over $\C^\ast_t$. They are subject to (conjectural) relations described above.

The categories i)-iii) depend on a choice of a complex Lagrangian subvariety $L\subset M$, while the category iv) does not depend  on $L$ (it depends only on $(M,\omega^{2,0})$ and some additional global choices). 

Consider  the ringed space $(\C^\ast, \OO_{\C^\ast}^{an})$. Let us define the ringed space $(\C,\OO_{\C,0}^{mer})$ such as follows.
Outside of $t=0$ it coincides with $(\C^\ast, \OO_{\C^\ast}^{an})$, and the stalk of $\OO_{\C,0}^{mer}$ at $t=0$ consists of germs of functions meromorphic at $t=0$.  Then the choice of Lagrangian subvariety $L$ gives rise to an extension of the sheaf $\FF_{glob}$  on $(\C^\ast, \OO_{\C^\ast}^{an})$ to the one on $(\C,\OO_{\C,0}^{mer})$.  Hence we have another sheaf of categories

v) $\FF_{L,glob}^{mer}$ which is the above-mentioned extension. Its stalk $\FF_{L,glob,0}^{mer}$ at $t=0$ is an $\A$-category over the field $\C\{t\}[t^{-1}]$. 

At the level of $Ext$-groups this extension gives rise to a meromorphic structure at $t=0$ on the holomorphic vector bundles. In particular this explains conceptually the meromorphic structure on the global Betti cohomology from Section \ref{comparison isomorphisms for 1-forms}.

\subsection{Global and local holonomic $DQ$-modules}\label{holonomic DQ-modules}

Now let us discuss the de Rham side. 
Unlike the Betti/Fukaya case we will have essentially only two categories: the local one depending on $(M,L)$ and the global one depending on $M$ only (and some additional global choices). There will be no wall-crossings on the de Rham side.
We start by recalling some general facts on deformation quantization.

For an arbitrary complex symplectic manifold $(M,\omega^{2,0})$ there is a canonically defined sheaf of triangulated $\A$-categories  linear over $\C[[t]]$. It is a sheaf of triangulated $\A$-categories  consisting of complexes of modules over a sheaf of algebroids.  The latter is locally isomorphic to the category of free rank one modules  over a formal deformation quantization $\OO_{M,t}$ of the sheaf of algebras ${\mathcal O}_M$ of analytic functions on $M$. 
Similarly to the Fukaya category case this sheaf can be twisted with a class $[B_{init}]\in H^2(M,\C^\ast)$. 
Furthermore this sheaf of triangulated $\A$-categories is endowed with a sheaf of $t$-structures, so we can speak about the corresponding sheaf of abelian categories.\footnote{In fact there is a notion of abelian $\A$-category which is appropriate in this context.}
 We denote by $\CC=\CC(M, [B_{init}])$ the abelian category of  global sections of this sheaf of abelian categories. The category of global sections can be trivial. 
Objects of $\CC$ are called {\it $DQ$-modules} (see [KasSch1]).\footnote{In what follows we will skip $[B_{init}]$ from notation.}

\

Next we will give a definition of the local category of holonomic $DQ$-modules. For that let us fix a closed analytic Lagrangian subset $L\subset M$, i.e. a possibly singular analytic subvariety which is Lagrangian on the smooth locus. 
We denote by $Hol_{L,loc}$ the full subcategory of $\CC$ consisting of locally finitely generated $DQ$-modules which are topologically supported on $L$. 
\footnote{More precisely the $DQ$-modules should be also  regular along $L$  in a sense which generalizes  the one in [AgKash]. 
We will not discuss this notion in current paper.} We call objects of $Hol_{L, loc}$ {\it holonomic $DQ$-modules supported  on $L$}. 
The inclusion $L\hookrightarrow L^\prime$ of closed Lagrangian analytic subsets induces a fully faithful embedding of the corresponding local categories.

If  $L$ is smooth and a trivialization of the $U(1)$-gerbe with the class $[B_{init}]_{|L}+\overline{St}_2(L)$ is chosen\footnote{In case $[B_{init}]=0$ the choice of trivialization is equivalent to the choice of $K_L^{1/2}$.} then $Hol_{L,loc}$ contains a canonical object $E_L$, which is locally isomorphic to $\OO_L[[t]]$.

\

In order to define the category $Hol_{glob}$ of global holonomic $DQ$-modules we will assume for simplicity that $(M,\omega^{2,0})$ is an affine algebraic symplectic manifold endowed with a Poisson compactification $\overline{M}$ such that $\overline{M}-M$ is a Poisson ample divisor. We define  $M_{log}$ as a Zariski open subset of $\overline{M}$ which satisfies the properties from the beginning of Section \ref{Fukaya categories of complex manifolds}.

In [Ko2] under the assumptions $H^i(\overline{M},\OO)=0,i=1,2$ the deformation quantization algebra $\OO_t(M)$ was defined as a filtered algebra, flat over $\C[[t]]$. The filtration modulo $t$ is given by the order of pole of a regular function at $\overline{M}-M$.
E.g. in the case $M=T^\ast \C=\C^{2}$ and $\overline{M}=\C{\bf P}^{2}$ the $k$-th term of filtration of the algebra $\OO_t(M)$ consists of expressions
$$\sum_{i,j\ge 0, i+j\le k}a_{ij}x^i(t\partial_x)^j, \quad a_{ij}\in \C[[t]].$$

Notice that in this example the algebra $\OO_t(M)$ contains a subalgebra over $\C[t]$ consisting of the above expressions in which all $a_{ij}\in \C[t]$. Based on this and similar examples it was conjectured (see [Ko2], Section 2.3, Conjecture 1) that the algebra $\OO_t(M)$ contains the natural filtered subalgebra $\OO_t^{an}(M)$ defined over the algebra $\OO^{an}(D)$ of analytic functions on an open disc $D\subset \C_t$ containing the point $t=0$. \footnote{More optimistic conjecture from [Ko2] says that one can replace $D$ by the whole line $\C$.} {\it In what follows we will suppose that the above assumptions from [Ko2] as well as the above-mentioned conjecture from loc. cit. are satisfied.}

For a  finitely-generated  $\OO_t^{an}(M)$-module $E$ one can define the notion of ``good filtration" $E_{\le m}$ starting with a finite set of generators and applying the canonical filtration on $\OO_t^{an}(M)$. It is explained in [Ko2] that such a choice allows us to extend $E$ to an object $\overline{E}$ on a properly defined deformation quantization category of $\overline{M}$. One can spell out the property that the ``support of $\overline{E}$ at infinity" does not intersect $\overline{M}-M_{log}$.  Also one can define the notion of $E$ being holonomic using the growth condition  $rk\, E_{\le m}\le C\cdot m^{dim_\C M/2}$.
Finally we define the category of global holonomic $DQ$-modules $Hol_{glob}:=Hol(M)$ as the category consisting of holonomic modules satisfying the above ``support at infinity" condition. For a given partial compactification $M_{log}$ the category $Hol_{glob}$ does not depend on a choice of the full compactification $\overline{M}$.

\

Next we discuss the comparison between $Hol_{L,loc}$ and $Hol_{glob}$. We will treat $Hol_{glob}$ as the category of global sections of a sheaf of categories on the ringed space $(D,\OO_D^{an})$. By taking the pull-back to the formal disc $\widehat{D}$ we obtain the category $Hol_{glob}^{form}$ which is linear over $\C[[t]]$. This category is well-defined even without the assumptions from [Ko2]. Objects from $Hol_{glob}^{form}$ are supported on {\it algebraic} Lagrangian subvarieties such that their closures in $\overline{M}$ are contained in $M_{log}$. Let us choose such a subvariety $L$. Then we have two categories over $\C[[t]]$:

a) The category $Hol_{L,loc}$ which depends on the formal symplectic neighborhood in $M$ of {\it analytic} Lagrangian subvariety $L$.

b) The category $Hol^{form}_{L,glob}$ of objects of $Hol^{form}_{glob}$ supported on $L$.

There is a natural  functor $Hol^{form}_{L,glob}\to  Hol_{L,loc}$. Roughly speaking we forget the condition of polynomial growth of
$L$ at infinity.

{\it We claim that if $\overline{L}$ is smooth at infinity and intersects strata of $M_{log}-M$ transversally then this functor gives an equivalence of categories. This functor categorifies the  local-to-global de Rham isomorphism.}

\

Let us illustrate this discussion in the case when $M=T^\ast X,  L=L_0\cup L_1$ where $L_0=X, L_1=graph(df)$,  and $[B_{init}]$ is the pull-back of $St_2(X)$. Here $f$ satisfies the assumptions of Section \ref{data for exp int}. We choose $M_{log}$ as in the Section \ref {Fukaya and Betti}. Notice that the $U(1)$-gerbes with classes $[B_{init}]_{|L_i}+\overline{St}_2(L_i), i=0,1$ are canonically trivialized. Therefore we have canonical objects in $ Hol_{L_i,loc}, i=0,1$, and hence the corresponding objects $E_i^{DR,loc}\in Hol_{L,loc}$. Similar considerations hold for an arbitrary complex symplectic manifold.

One can check directly that 
$$Ext^\bullet(E_0^{DR,loc}, E_1^{DR,loc})\simeq \mathbb{H}^\bullet(X, (\Omega_X[[t]], td+df\wedge (\bullet)))=H_{DR,loc}^\bullet(X,f).$$
Similarly we have objects $E_i^{DR,glob}\in Hol_{glob}, i=0,1$.  In the more traditional language they are holonomic algebraic $t$-$D$-modules given by $\OO_X$ and $e^{f\over t}\OO_X$ respectively. Analogously to the local case we have an isomorphism of $\C\{t\}$-modules

$$Ext^\bullet(E_0^{DR,glob}, E_1^{DR,glob})\simeq \mathbb{H}^\bullet(X_{Zar}, (\Omega_X\{t\}, td+df\wedge (\bullet)))=$$
$$=H_{DR,glob}^\bullet(X,f)\otimes_{\C[t]}\C\{t\}.$$

The global-to-local de Rham isomorphism
$$Ext^\bullet(E_0^{DR,glob}, E_1^{DR,glob})\otimes_{\C\{t\}}\C[[t]]\simeq Ext^\bullet(E_0^{DR,loc}, E_1^{DR,loc}).$$

It gives the same comparison isomorphism 
$$H^\bullet_{DR,glob}(X,f)\otimes_{\C\{t\}}\C[[t]]\simeq H^\bullet_{DR,loc}(X,f)$$
as in Proposition 2.3.4 d).

\subsection{Categorical analogs of comparison isomorphisms}\label{categorical analogs}

In this subsection we will explain how the cohomological isomorphisms from Sections \ref{comparison isomorphisms for 1-forms}, \ref{comparison isomorphisms and wheels of lines} are related to certain categorical statements concerning four different categories and four functors between them.  For  general complex symplectic manifolds the appropriate framework is the one  of  the 
{\it generalized Riemann-Hilbert correspondence} (RH-correspondence for short) 
which we will discuss in detail in  subsequent papers on the project. 
Generalized  RH-correspondence is formulated as an equivalence of  categories associated with the Betti and de Rham sides  in global and local cases. In this subsection we 
will discuss the simplest version, postponing a more general case to the subsequent papers.

Let us start with the local RH-correspondence. Let $L$ be a complex Lagrangian subvariety of $(M,\omega^{2,0})$. Therefore it is a Lagrangian subvariety for the symplectic form $Re(\omega^{2,0})$. Recall (see Remark \ref{local Fukaya category in general}) that the category $\FF_{L,loc}$ can be described as the category of finite rank modules over a finite type $\A$-algebra $A_{L,\Z}$ over $\Z$. Varying real symplectic structure $Re(e^{i\theta}\omega^{2,0}), \theta\in \R/2\pi \Z$ we obtain an automorphism $\psi\in Aut(A_{L,\Z})$.  It induces an autoequivalence  of the category of finite rank $A_{L,\Z}$-modules. It can be identified with the monodromy of the corresponding  local system of categories   on $S^1$ (or equivalently on $\C^\ast$). We are going to associate with $\psi$ a new $\A$-category over $\C((t))$.

Let ${\bf B}=\overline{\C((t))}[\log\,t]$. Let $\psi_{\bf B}: {\bf B}\to {\bf B}$ be the automorphism given by 
$$t^{1/k}\mapsto e^{2\pi i/k}t^{1/k}, \quad \log\,t\mapsto \log\, t+2\pi i.$$
We define the $\A$-category $\widehat{\FF}_{L,loc}^{form}$ as the category of modules over the algebra $(A_{L,\Z}\otimes_\Z{\bf B})[T,T^{-1}]_{tw}$ of twisted Laurent polynomials, where $Ad(T)$ acts as the automorphism  $\psi\otimes \psi_{\bf B}$, and 
 which are perfect over ${\bf B}$. 
 This category is linear over the fixed points subalgebra ${\bf B}^{\psi_{\bf B}}\simeq \C((t))$. 
 Define a full subcategory ${\FF}_{L,loc}^{form}\subset \widehat{\FF}_{L,loc}^{form}$ by considering only 
 modules obtained by the extension of scalars from modules over the algebra $(A_{L,\Z}\otimes_\Z{\bf B}_+)[T,T^{-1}]_{tw}$, 
 where ${\bf B}_+\subset {\bf B}$ is the subring of series containing only non-negative powers of $t$ and $\log\,t$.

 Then the conjectural {\it local Riemann-Hilbert correspondence} says that there is an equivalence functor
$$RH_{loc}=RH_{L,loc}: Hol_{L,loc}\otimes_{\C[[t]]}\C((t))\stackrel{\sim}{\to} \FF_{L,loc}^{form}.$$


\

 \begin{exa}\label{local RH for smooth L}
 Let us assume that $L$ is smooth, connected with chosen based point $l_0$  and $[B_{init}]=\overline{St}_2(L)$. Then the dg-algebra $A_{L,\Z}$ is isomorphic to $Chains(\Omega(L,l_0))$ and in this case $\psi=id_{A_{L,\Z}}$. Recall that in this set up we have a canonical object $E_L^{DR}\in Hol_{L,loc}$. Then it gives a functor from the category $Hol_{L,loc}$ to the category of sheaves of $\C[[t]]$-modules on $L$. Namely to a sheaf $\underline{V}$ of $\C[[t]]$-modules on $L$ corresponding to an object of $Hol_{L, loc}$ we associate a sheaf $\underline{Hom}_{Hol_{L, loc}}(E_L^{DR},\underline{V})$ of $\C[[t]]$-modules. After tensoring with $\C((t))$ we obtain a $\C((t))$-local system on $L$ which is obtained from a local system of $\C[[t]]$-modules by the change of scalars. Such a local system gives an object of $\FF_{L,loc}^{form}$. In this way we obtain an example of $RH_{loc}$.

 \end{exa}

The global RH-correspondence is much more mysterious. We will discuss it in a separate paper. Roughly it claims an equivalence of categories
$$RH_{glob}: Hol_{glob}\otimes_{\C\{t\}}{\mathcal A}\stackrel{\sim}{\to} \FF_{glob},$$
where the algebra ${\mathcal A}$ of germs at $t=0$ of analytic functions on $\C^\ast_t$ with arbitrary growth as $t\to 0$ was defined in Section \ref{comparison isomorphisms for 1-forms}.

\begin{exa}\label{relation to the usual RH}
Let $X$ be a smooth affine algebraic variety. Let us choose a compactification $\overline{X}\supset X$ by a simple normal crossing divisor $D$. Then we define a log-extension of $M=T^\ast X$ by $M_{log}=T^\ast_{log\, D}\overline{X}$. Then $Hol_{glob}$ is defined over $\C[t]$ its restriction to any point $t=t_0\in \C^\ast$ can be identified with the category of holonomic $D$-modules on $X$ with regular singularities and microlocal support in $X$ (i.e. algebraic vector bundles on $X$ with regular singular flat connections). The category $\FF_{glob}$ restricted to $t=t_0$ can be identified with the category of $\C$-local systems on $X$. Then we expect that $RH_{glob}$ becomes the usual Riemann-Hilbert correspondence.
\end{exa}

Next we propose a conjecture which is a categorical extension of  Conjecture \ref{comparison of isomorphisms diagram} part 1).

\begin{conj}\label{RH for categories with meromorphic structures}
For a given $L$ the functor $RH_{glob}$ gives rise to the functor $RH_{L,glob}^{mer}$ which identifies the subcategory 
$$Hol_{L,glob}\otimes_{\C\{t\}}\otimes \C\{t\}[t^{-1}]\subset Hol_{glob,L}\otimes_{\C\{t\}}{\mathcal A}$$ 
with the the category 
$$\FF_{L,glob,0}^{mer}\subset \FF_{glob}.$$

Here $Hol_{L,glob}$ is the full subcategory of $Hol_{glob}$ consisting of objects whose restriction to $t=0$ are coherent sheaves supported on $L$, and the category $\FF_{L,glob,0}^{mer}$ was defined as item v) at the end of Section \ref{analytic families of Fukaya categories}.

\end{conj}

Assuming the above conjecture one can formulate a categorical extension of Conjecture \ref{comparison of isomorphisms diagram} part 2).

\begin{conj}\label{diagram of categories}

The following diagram is commutative

\[ \begin{tikzcd}
{Hol}_{L,glob}\otimes_{\C\{t\}} \C((t))\arrow{r}{{RH}_{glob}^{mer}} \arrow[swap]{d}{\phi_{DR}}&{\FF}_{L,glob,0}^{mer}\otimes_{\C\{t\}[t^{-1}]} \C((t))  \arrow{d}{\phi_{Betti}}  \\%
Hol_{L,loc}\otimes_{\C[[t]]} \C((t)) \arrow{r}{RH_{loc}}& \FF_{L,loc}^{form}
\end{tikzcd}
\]
Here we keep the notation for the extension of the functor $RH_{glob}^{mer}$ to the $\C((t))$-linear categories.

\end{conj}

Let us assume the Conjecture \ref{diagram of categories} for $M=T^\ast X$, and let $L_0=X, L_1=graph(\alpha)$, where $\alpha$ is a closed algebraic $1$-form as before.
Then $Ext$-groups between two canonical objects associated with $L_0$ and $L_1$ in the corresponding categories in   we obtain the commutative diagram of cohomology groups from the Conjecture \ref{comparison of isomorphisms diagram}. Thus the Conjecture 
\ref{diagram of categories} implies the latter.

\part{Infinite-dimensional exponential integrals}

The reader should be aware that except of the rigorously developed theory of quantum wave functions the contents of this part of the paper is mostly conjectural and sometimes speculative. Our aim is to apply the ideas from the first part of the paper in the infinite-dimensional framework.

\section{Quantum wave functions}\label{infinite-dimensional case}

Infinite-dimensional exponential integrals (a.k.a. Feynman integrals or functional integrals) are not well-defined as objects of mathematics e.g. because there is no top degree volume form (Feynman measure). Also the integration cycle is infinite-dimensional, and we do not have an appropriate integration theory in this case. 

Trying to make sense of the infinite-dimensional exponential integral mathematician can use various approaches.
For example one can consider the formal expansion of the functional integral with respect to the coupling constant typically present in the problem (``Planck constant''). This formal expansion is given in terms of the formal Feynman rules and hence is well-defined mathematically. These formal expansions are typically divergent, but their Borel resummations often give analytic functions in the dual variable corresponding to the coupling constant, at least outside of a countable set.  In other words, these series are {\it resurgent}. 

Another approach utilizes the replacement of the initial ``integration cycle" to another one for which  the infinite-dimensional integral can be given some meaning. After that  one can try to make  sense of the analytic continuation of the integral over this new cycle with respect to coupling constants. 

The above-mentioned approaches (as well as  several others) can be considered as an attempt to give a meaning e.g. to the Feynman path integral as an infinite-dimensional exponential period by giving a precise meaning to its de Rham and Betti sides.

In this section we will propose an approach to analytic continuation of the Feynman integrals in terms of deformation quantization.
It can be thought of as a holomorphic version of the equivalence of Lagrangian and Hamiltonian formalisms in quantum mechanics.
We are going to reformulate the complexified infinite-dimensional exponential integral as a collection of data which can be computed non-perturbatively in terms of the {\it finite-dimensional data} consisting of quantum wave functions and their pairings. We propose to use the  pairing of quantum wave functions as a definition of the {\it analytically continued} infinite-dimensional exponential integral. Quantum wave functions are associated with holomorphic Lagrangian subvarieties (and some additional data). From the perspective of the Betti side of the RH-correspondence the possibility to continue analytically the exponential integral is based on the  analyticity of the corresponding wall-crossing structure (see [KoSo12]). The main conjecture of the loc.cit. claims  that analyticity of the WCS implies resurgence of the related formal series.

The underlying geometry in this section will be again the one of a pair of finite-dimensional complex Lagrangian submanifolds in a complex symplectic manifold. If the Lagrangian submanifolds are unobstructed then the corresponding WCS is determined by the vector spaces of vanishing cycles associated with intersection points as well as by the virtual number of pseudo-holomorphic discs  with boundary on our Lagrangian submanifolds. In case when we have  exponential bounds on the numbers of such discs, the corresponding WCS is analytic. Comparison of Betti and de Rham pictures predicts that the pairings at the intersection points of the quantum wave functions associated with our Lagrangian submanifolds give rise to formal series which should be resurgent.
As we will see the proper definition of the pairing requires a certain normalization.

{\it Remark about the notation.}
Since our exposition in this section is motivated by physics, we will change the notation for  the parameter from $t$ to the ``Planck constant" notation $\hbar$, although in practice $\hbar$ can be any parameter of the theory (``a  coupling constant").

\subsection{Path integral with holomorphic Lagrangian boundary conditions}\label{path integrals}

Let us briefly explain what we would like to achieve in the case of the standard complex symplectic vector space. 

Let ${\bf C}^{2n}$ be the standard complex vector space with coordinates $({\bf q},{\bf p})=(q_1,\dots,q_n,p_1,\dots,p_n)$ endowed with the standard holomorphic symplectic form $\omega^{2,0}=\sum_i dp_i\wedge dq_i$. We remark that in some cases it is useful to think of this symplectic manifold as  the  universal cover of $(\C^\ast)^{2n}$. This is useful e.g. in the complexified Chern-Simons theory in the $1$-dimensional space-time formulation, since the action functional is well-defined on the universal cover only.

The ``space of fields" by definition  is the space of $C^\infty$ maps
$$\varphi:[0,1]\to \C^{2n},\quad \varphi(t)=({\bf q}(t),{\bf p}(t))=(q_1(t),\dots,p_n(t))$$

The boundary conditions for the path integral with the action $S(\phi)$ are specified by a choice of two holomorphic Lagrangian submanifolds $L_0,L_1\subset \C^{2n}$. This means that we consider smooth maps $\varphi:[0,1]\to \C^{2n}$ such that 
$$\varphi(0)\in L_0,\quad \phi(1)\in L_1.$$

The action functional (in the first order formalism) is 
$$S(\varphi)=\int_0^1\sum_i p_i(t){d q_i(t)\over dt}+\int_0^1H({\bf q}(t),{\bf  p}(t),t)dt+f_0(\phi(0))-f_1(\phi(1)),$$
where 
$$H:\C^{2n}\times [0,1]\to \C$$
is holomorphic in complex coordinates $(\bf q,\bf p)$ and $C^\infty$ in real coordinate $t\in[0,1]$ and $f_j,j=0,1$ are primitives of the closed $1$-forms $(\sum_i p_idq_i)_{|L_j}$.

The case $H\equiv 0$ is already very interesting. The arising theory is topological. Examples include WKB theory and the complexified  Chern-Simons theory. If the matrix of second derivatives $({\partial^2 H\over \partial p_i\partial p_j})_{1\le i,j\le n}$ is non-degenerate, then one can exclude ${\bf p}(t)$ using the equations of motions, and get the Lagrangian $\mathcal L({\bf q},\dot {\bf q},t)$ as the Legendre transform of $H$ depending on first derivatives (hence the term ``first order formalism").

The action functional on the space of fields is a holomorphic function on an infinite-dimensional complex manifold. Its critical points are solutions of  the Euler-Lagrange equation
$${\delta S(\varphi) \over \delta \phi}=0 \iff \dot p_i={\partial H \over\partial q_i}, \dot q_i=-{\partial H \over\partial p_i}$$
Thus it defines a family of partially defined holomorphic symplectomorphisms 
$$g_t:\{\text{an open dense domain in }\C^{2n}\}\to \C^{2n},\quad t\in [0,1],\,g_0=\text{identity map}.$$
Critical points are identified with intersection points
$$\text{Crit}(S)\simeq g_1(L_0)\cap L_1.$$

Typically all intersections  are transversal, and there are countably many of them. Critical points   $x_j$ of $S$ are naturally labeled by the set of intersection points of Lagrangians $L_0\cap L_1$.

Symbolically, the Feynman path integral is written as 
$$\int e^{S(\varphi)\over \hbar}\mathcal D\varphi$$

The infinite-dimensional cycle of integration is not well-defined as well as the ``Feynman measure" $\mathcal D\varphi$. By analogy with  finite-dimensional exponential integrals we can expect that regardless of the definition of the Feynman path integral (or more general functional integrals) the answer depends  on the ``relative homology class of the integration cycle" (whatever this means). In several examples including e.g. the complexified Chern-Simons theory, the integration cycle  can be chosen as an integer linear combination of the 
 {\it infinite-dimensional} Lefschetz thimbles. The notion of Lefschetz thimble is well-defined in the infinite-dimensional case, and one can give some  meaning to the corresponding exponential integral.

In what follows  we will offer a formalism which assigns to each critical point $x_j$ of $S$ at which our holomorphic Lagrangian submanifolds intersect transversally,  a well-defined asymptotic series in $\hbar$. This series should be thought of as the infinite-dimensional analog of the RHS in  the properly normalized stationary phase expansion formula and will have the form: \footnote{The factor $(2\pi \hbar)^{-{n\over 2}}$ is conventional in  quantum mechanics. In the finite stationary phase expansion one would have the factor $(2\pi \hbar)^{{n\over 2}}$. We do not understand the geometric origin of this discrepancy.}
 
$$\intop_{\overset{\text{local Lefschetz thimble}}{\text{\tiny outcoming from } x_j} }e^{S(\phi)\over \hbar}\mathcal D\phi  \,\underset{\hbar\to 0}{\sim}\,e^{S(x_j)\over \hbar}(2\pi \hbar) ^{-{n\over 2}}\cdot(c_{0,j}+c_{1,j}\hbar+c_{2,j}\hbar^2+\dots).$$
Furthermore, we will explain that the formal series which appear in the RHS are expected to be resurgent for very general reasons, and moreover they should give rise to  multivalued analytic functions, which are related to each other for different critical points.

We will also discuss how  to go beyond the case when the Lagrangian subvarieties intersect transversally, i.e. when the approach with Lefschetz thimbles does not work.  

\subsection{Space of paths and the potential}\label{space of paths}

Part of the data from the Section \ref{path integrals} can be spelled out in a coordinate-free way.

Namely, consider a complex symplectic manifold $(M,\omega^{2,0}), dim_\C M=2n$ and a pair of complex Lagrangian submanifolds $L_0, L_1$. We {\it do not} require that $L_0$ intersects $L_1$ transversally. In general  $L_0\cap L_1$ is an analytic subset of $M$.

We denote by $P(L_0,L_1)$ the set of real smooth paths $\varphi:[0,1]\to M$ such that $\varphi(0)\in L_0, \varphi(1)\in L_1$. 
This is an infinite-dimensional complex manifold. 

The manifold $P(L_0,L_1)$ carries a holomorphic closed $1$-form $\eta$ given by the integration of $\omega^{2,0}$. Namely for any path in $P(L_0,L_1)$ the integral of $\eta$ over a path is defined as the integral of $\omega^{2,0}$ over the real two-dimensional membrane in $M$ corresponding to this path.
Closedness of $\eta$ follows from the one of $\omega^{2,0}$ and the vanishing of $\omega^{2,0}$ restricted to $L_i, i=0,1$.

Zeros of $\eta$ are constant paths, i.e. maps of the interval $[0,1]$ to the intersection $L_0\cap L_1$.
 In case if the Hamiltonian $H\ne 0$ we need to add to $\eta$ the holomorphic exact $1$-form given by the differential of the holomorphic function on $P(L_0,L_1)$ equal to $\int_0^1Hdt$. One can check that in the case $M=\C^{2n}$ considered in Section \ref{path integrals} the form $\eta$ is equal to $dS$ where $S=S(\varphi)$ was defined in loc.cit.

\subsection{Sheaf of vanishing cycles on the space of paths}\label{infinite-dimensional vanishing cycles}

Let $S$ be an analytic function in a neighborhood of the locus of zeros $\ZZ(\eta)$ such that $\eta=dS$. We may assume that $S$ vanishes on $\ZZ(\eta)$.

After a choice of square roots of the canonical bundles  $K_{L_i}^{1/2}, i=0,1$
one can define the sheaf of vanishing cycles $\phi_{S}(\Z)$ on $\ZZ(\eta)$. Rigorous definition requires a proof of the existence of the so-called orientation data a.k.a. spin structure data (see [KoSo1]). Detailed discussion of the sheaf of vanishing cycles for the above action $S$ will be given elsewhere.

In this subsection we will discuss the sheaf of vanishing cycles in the infinite-dimensional case in abstract terms.
Namely, let  us  assume that $X$ is a complex  manifold, which can be infinite-dimensional, and $f:X\to \C$ an analytic function. 
We do not want to make 
these notions more precise at this time, but they are clear in the finite-dimensional case. The reader can keep in mind two main examples: the space of paths considered above and  the complexification of the Chern-Simons functional which we will discuss later.

Let us assume that the second derivative $d^2(f):T_X\to T^{\ast}_X$ is a Fredholm operator at each critical point. 
Then locally the  set $Crit(f)$ of critical points of $f$ is a finite-dimensional analytic space.

Under this assumption we can replace the pair $(X,f)$ near each point $x\in Crit(f)$ by a finite-dimensional model.  

Here is the construction.
Let  us pick a vector subspace $T^{vert}\subset T_{X,x}$ of finite codimension such that
the restriction of $d^2(f)$ to $T^{vert}$ is non-degenerate.
It follows  that we can extend $T^{vert}$ to a  ``vertical'' analytic subbundle ${\mathcal T}^{vert}$ which gives rise to a fibration $p: U\to B$ of an open neighborhood
 $x\in U$ over a finite-dimensional manifold $B$. Then on each fiber $U_b:=p^{-1}(b)\subset U$ of the fibration which 
 is ``sufficiently close'' to $x$ there is a unique Morse critical point of the function $f_b:=f_{|U_b}$.
 We define the function $\overline{f}:B\to \C$ by the condition that $\overline{f}(b)$ is the critical value of $f_b$. The projection $p$ identifies locally $Crit(f)$ with $Crit(\overline{f})$, We declare that the stack of the sheaf of vanishing cycles $\phi_f({\Z})$ at $x$ is equal to  $\phi_{\overline{f}}({\Z})_{p(x)}$.
 
 \begin{rmk} The sheaf $\phi_{\overline{f}}({\Z})$ depends on a choice of $T^{vert}$. For two choices of $T^{vert}$ the corresponding sheaves of vanishing cycles are isomorphic but the isomorphism is  defined up to a sign. In general in order to have unambiguous definition of  $\phi_{\overline{f}}({\Z})$ one has to choose the above-mentioned ``orientation data"  (see e.g.[KoSo1], [KoSo5]).

 \end{rmk}


For $X=P(L_0, L_1)$ let us take  $f=S$ as above and assume that $H=0$. Then  the set $Crit(f)$ is  the set of constant paths, hence is identified with $L_0\cap L_1$. The tangent space at the constant path corresponding to $m\in L_0\cap L_1$ is identified with the vector space of $C^\infty$ maps $\varphi: [0,1]\to T_mM$ such that $\varphi(0)\in T_mL_0$ and $\varphi(1)\in T_mL_1$. The subspace $T^{vert}$ can be chosen as the intersection of the kernels of finitely many linear maps of the type
$$\varphi\mapsto \int_0^1\langle\varphi(t),\psi(t)\rangle dt,$$
where $\psi\in C^\infty[0,1]\otimes T_m^\ast M$ and $\langle\bullet,\bullet \rangle$ denote the canonical pairing of the tangent and cotangent spaces.\footnote{In order to be on the safe side  we have chosen  $\psi$ to be  {\it smooth} functions. Calculations show that a choice of more singular functions or distributions e.g. the delta-functions at $t=0$ and $t=1$ does not guarantee existense of Morse function on the fibers of $p$.}

One can show that in this case  the sheaf of vanishing cycles $\phi_{S}({\Z})$ is canonically isomorphic to the {\it finite-dimensional} sheaf of vanishing cycles $\phi_{(L_0,L_1)}$ associated with the pair $(L_0, L_1)$, and it is supported on $L_0\cap L_1$. It  is locally equal to the sheaf of the vanishing cycle $\phi_F({\Z})$ of the function $F$ which is determined by the property that if  $L_0$ is identified with the zero section of a cotangent bundle then $L_1=graph(dF)$. Alternatively one can interpret the sheaf $\phi_{(L_0,L_1)}$ in terms of (sheaves of) local Fukaya categories. Namely the stalk of  $\phi_{(L_0,L_1)}$ at any $x\in L_0\cap L_1$ is equal to $Ext^\bullet_{\FF_{(L_0\cup L_1)\cap B_\epsilon(x),loc}}(L_0\cap B_\epsilon(x), L_1\cap B_\epsilon(x))$, where is $B_\epsilon(x)$ is a small open ball with the center at $x$. Also we have an isomorphism 
$$\R\Gamma(L_0\cap L_1, \phi_{(L_0,L_1)})\simeq Ext^\bullet_{\FF_{L_0\cup L_1,loc}}(L_0, L_1).$$
Here we identify Lagrangian subvarieties with the corresponding objects of the local Fukaya category.

In the more general case when $H\ne 0$ the sheaf of vanishing cycles of the functional $S$ can be identified with the one corresponding to the intersection $g_1(L_0)\cap L_1$.

In the case when the Lagrangian submanifolds intersect transversally one can define a canonical up to a sign local integration cycle (``infinite-dimensional thimble").

In the remaining part of this section we are going to propose a rigorous finite-dimensional formalism which gives in the end formal series in $\hbar$ which are ``morally" equal to the asymptotic expansions of the path integrals over the local integration cycles.

\subsection{Reminder on deformation quantization and $\ast$-products}{\label{star-product}

If one can get rid of ambiguities in the path integral at the level of formal expansions in $\hbar$, one should have a QFT on the $1$-dimensional space-time manifold with the coordinate $t\in [0,1]$. After that  one should be able to add observable at any intermediate moment $t\in (0,1)$.
From the classical limit point of view  the observables should correspond to functions on $\C^{2n}$ or more generally on a complex symplectic manifold $M$. 

We will assume from now on that the Hamiltonian $H$  vanishes.\footnote{Generalization to the case $H\ne 0$ will be discussed later.} Then  the corresponding theory is {\it topological}, i.e. invariant under orientation-preserving diffeomorphisms of $[0,1]$. This implies that observables form an {\it associative algebra}. In the  formal expansion in the parameter $\hbar$ the associative product should be given by  the $\ast$-product on the space of observables.  
Thus we can use the language of deformation quantization of complex symplectic manifolds. It has already appeared in this paper when we spoke about the de Rham side of various comparison isomorphisms and categorical equivalences.

Let $(M,\omega^{2,0})$ be a complex symplectic manifold as before. We will often restrict ourselves to the case when $M$ is a smooth complex algebraic variety.

Recall that the deformation quantization of $(M,\omega^{2,0})$ can be defined
at the different levels:

 1)  (Most concrete) a $\ast$-product on $\mathcal O_M[[\hbar]]$, i.e. a $\C[[\hbar]]$-linear associative product $\ast$  given by
 $$f\ast g=fg+\hbar {\{f,g\}\over 2} +\dots=\sum_{m\ge 0}\hbar^m B_m(f,g),$$
 where $B_m$ is a bi-differential operator of order $(m,m)$ for each $m$. Here $B_0(f,g)=fg$ is the usual product and $B_1(f,g)={1\over 2}\{f,g\}$ is given by the Poisson bracket associated with $\omega^{2,0}$. 
 
 Recall that any two $\ast$-products are locally equivalent via transformations of the form $f\mapsto f+\sum_{k\ge 1}\hbar^kA_k(f)$, where $A_k$ is a differential operator of order less or equal than $k$.

 2) (Less abstract) a  sheaf of algebras $\mathcal O_{M,\hbar}$ over $\mathbb C[[\hbar]]$, together with an isomorphism
 $$\mathcal O_{M,\hbar}/\hbar\mathcal O_{M,\hbar}\simeq \mathcal O_M,$$
 and such that $\mathcal O_{M,\hbar}$ is locally isomorphic to the sheaf  $\mathcal O_M[[\hbar]]$ endowed with a $\ast$-product as in 1).

 3)  (Most abstract) a  sheaf of $\C[[\hbar]]$-linear categories $\CC_M$ (see e.g. [KasSch1], [Ko2]) which are locally equivalent to the sheaf of categories of $\OO_{M,\hbar}$-modules from 2) together with a choice of equivalence $\CC_M/\hbar \CC_M\simeq \OO_M-mod$.

 Equivalence classes of such sheaves of categories bijectively correspond to elements of ${1\over{\hbar}}H^2(M,\C[[\hbar]])$ of the form ${[\omega^{2,0}]\over\hbar}+\dots$. These elements are called Deligne classes of the corresponding sheaves of categories. Moreover there exists a {\it canonical} sheaf of categories $\CC_{M}^{can}$ (not just the equivalence class) with the Deligne class ${[\omega^{2,0}]\over\hbar}$.

 \begin{exa}\label{examples of Deligne classes}
  Let $X$ be a smooth complex manifold, $M=T^\ast X$ endowed with the standard symplectic structure. In this case $[\omega^{2,0}]=0$.
   
 a) The sheaf of $\hbar$-microdifferential operators $\EE_{M,\hbar}$ is the natural extension of the notion of $\hbar$-differential operators \footnote{It has different names in different papers see e.g. [KasSch1].} is an example of the sheaf $\OO_{M,\hbar}$ with the Deligne class equal to 
 $$-{1\over 2}c_1(T^\ast X)\in H^2(X,\C)\subset H^2(M,{1\over \hbar}\C[[\hbar]]).$$
 
 b) Given an element $\delta\in H^1(X, \Omega^{1,cl}_X)\simeq H^1(X,\OO_X^\times/\C^\times)$ one can twist the sheaf $\EE_{M,\hbar}$ and get a new sheaf of algebras $\EE_{M,\hbar}^\delta$. Namely if $\delta$ is locally represented by an invertible function $f$ (modulo scalars) we can define the corresponding automorphism of $\EE_{M,\hbar}$  by $Ad(f): P\mapsto fPf^{-1}$. The Deligne class of $\EE_{M,\hbar}^\delta$ is equal to $-{1\over 2}c_1(T^\ast X)+\partial(\delta)$ where $\partial$ is the coboundary map from $H^1(X,\Omega^{1,cl}_X)\simeq H^1(X,\OO_X/\C)$ to $H^2 (X,\C)=H^2(T^\ast X, \C)$. Notice that for a given element of  $Pic(X)\otimes_\Z \C$ one can speak about the corresponding sheaf of twisted $\hbar$-microdifferential operators via the map $Pic(X)\simeq H^1(X,\OO_X^\ast)\stackrel{dlog}{\to}H^1(X,\Omega^{1,cl}_X)$.  The  group $Pic(X)\otimes_\Z \C$ can be thought of as the group of equivalence classes of tensor products of complex powers of line bundles on $X$. 
  
 c)  According to b) one can speak about the sheaf of $\hbar$-microdifferential operators  on a complex power of a given line bundle on $X$. Then  the above-mentioned sheaf of categories $\CC_M^{can}=\CC_{T^\ast X}^{can}$ has a distinguished generator $G$ for which the sheaf of endomorphisms $\underline{End}(G)$ is canonically identified with the sheaf of $\hbar$-microdifferential operators acting on the line bundle $K_X^{1/2}$.

 d) Furthermore one can replace $Pic(X)\otimes_\Z \C$ by $Pic(X)\otimes_\Z (\C\oplus \C\cdot {1\over \hbar})$.  The corresponding $\hbar$-microdifferential operators form a sheaf of $\C[[\hbar]]$-algebras on the twisted cotangent bundle $T^\ast_\lambda X$. Here $\lambda\in H^1(X,\Omega^{1,cl}_X)$ is the defined by the composition
 
$$ Pic(X)\otimes_\Z (\C\oplus \C\cdot {1\over \hbar})\twoheadrightarrow Pic(X)\otimes_\Z\C\cdot{1\over \hbar}\stackrel{\hbar\cdot}{\to}Pic(X)\otimes_\Z\C\to H^1(X,\Omega^{1,cl}_X).$$

  \end{exa}
  
  \begin{rmk}\label{sheaves of categories with gerbes}
  
One can slightly generalize the most abstract case 3) by replacing the sheaf of categories $\OO_M-mod$ by an $\OO_M^\times$-gerbe. This means that we omit in the definition of the sheaf of categories the choice of equivalence $\CC_M/\hbar \CC_M\simeq \OO_M-mod$. Such a generalization will be useful later in relation to the notion of quantum wave function structure.

\end{rmk}

 Let us comment on 1)-3). 
 
 It is known (see e.g. [Ko2]) that given a Poisson structure $\pi$ on a, say,  smooth complex  algebraic variety $P$ one can define a canonical map $Spec(\C[[\hbar]])\to MC(R\Gamma(P, \wedge^\bullet T_P))$, where $MC(R\Gamma(P, \wedge^\bullet T_P))$ is the space of Mauer-Cartan elements in the DGLA  $R\Gamma(P, \wedge^\bullet T_P)$ corresponding to the sheaf of polyvector fields. The  map is given by the canonical path $\hbar\pi$. By the formality theorem (see e.g. [Ko2]) the space $MC(R\Gamma(P, \wedge^\bullet T_P))$ can be identified with the space of Maurer-Cartan elements in the corresponding Hochschild cochain complex, controlling the deformation theory of the category of perfect $\OO_P$-modules. 
 
 Let us return to our symplectic case, so that $P=M, \pi^{-1}=\omega^{2,0}$.
The  path $\hbar\mapsto \hbar\pi$ gives rise to the above-mentioned canonical sheaf of categories $\CC_M^{can}$, flat over $\C[[\hbar]]$ and identified modulo $\hbar$ with $\OO_M-mod$. It is not true in general that $\CC_M^{can}$ is equivalent to the category $\OO_{M,\hbar}-mod$ for some sheaf of quantized algebras $\OO_{M,\hbar}$ (the obstruction is discussed in [NeTsy]).

But even if we know that there is such an equivalence, the equivalence functor $\Phi: \CC_M^{can}\simeq \OO_{M,\hbar}-mod$ is {\it not} defined canonically. A choice of $\OO_{M,\hbar}$ and the functor $\Phi$ means that we consider the deformation theory of a  generator of $\CC_M^{can}$. This generator corresponds to what physicists call {\it canonical coisotropic brane} $\mathcal{B}_{cc}$. 

In the next subsection we will axiomatize the functor $\Phi$. The corresponding structure on $M$ is named {\it quantum wave function structure}.

For the applications to path integrals the case 2) will be sufficient, although in practice the sheaf of algebras $\OO_{M,\hbar}$ often arises from the $\ast$-product  on $\OO_M$, i.e. as in 3).

In particular, in the framework of  Section \ref{path integrals} (i.e. when the symplectic manifold is $\C^{2n}$) we are going to use the  {\it Moyal $\ast$-product}:
$$f \ast_{Moyal} g=\left[\exp\Bigl({\hbar\over 2}\sum_i \bigl(\partial_{p_i}\otimes \partial_{q_i}-\partial_{q_i}\otimes \partial_{p_i}\bigr)\Bigr) (f\otimes g)\right]_{|Diagonal\,\,\C^{2n}\subset \C^{2n}\times \C^{2n}}. $$

The Moyal $\ast$-product  is  invariant under affine symplectic transformations, i.e. elements of the group $Sp(2n,\C)\ltimes \C^{2n}$:
$$\begin{pmatrix}{\bf q} \\ {\bf p}\end{pmatrix}  \mapsto T\cdot \begin{pmatrix}{\bf q} \\ {\bf p}\end{pmatrix}+\begin{pmatrix}{\bf q}_0 \\ {\bf p}_0\end{pmatrix},\qquad T\in Sp(2n,\C),  \begin{pmatrix}{\bf q}_0 \\ {\bf p}_0\end{pmatrix}\in \C^{2n}.$$
Let $\mathcal O_\hbar:=\OO_{\C^{2n},\hbar}$ be the corresponding sheaf of algebras.

\begin{rmk} There is another canonical $\ast$-product which is defined on the cotangent bundle of any complex manifold, and consists of micro-differential operators on functions (or half-densities), i.e functions in coordinates $(q_1,\dots,q_n)$, derivatives $(\hbar\partial_{q_1},\dots,\hbar\partial_{q_n})$ and formal series in $\hbar$. The polynomial part of $\Gamma(\C^{2n},\mathcal O_{\C^{2n},\hbar})$ is 
 $$ \C[q_1,\dots,q_n][\hbar\partial_{q_1},\dots,\hbar\partial_{q_n}][[\hbar]]$$
In the case of the cotangent bundle to $\C^n$ this $\ast$-product  is given by
$$ f \ast_{micro} g=\left[\exp\Bigl({\hbar}\sum_i \partial_{p_i}\otimes \partial_{q_i}\Bigr) (f\otimes g)\right]_{|Diag\,\,\C^{2n}\subset \C^{2n}\times \C^{2n}} $$
It is covariant with respect to   changes of coordinates in $\C^n_{\bf q}$ lifted uniquely as symplectomorphisms of  $\C^{2n}_{\bf q,p}=T^*\C^n_{\bf q}$.
\end{rmk}

\subsection{Quantized boundary conditions as modules}\label{left and right modules}

In the case of  general complex symplectic manifold $(M,\omega^{2,0})$ one can think that the path integral with the action $S$  and observables $\psi_i, 0\le i\le n$ can be symbolically written  as
$$\int_{\varphi\in P(L_0, L_1)}e^{S(\varphi)/\hbar}\psi_0(\varphi(t_0))\psi_1(\varphi(t_1))...\psi_n(\varphi(t_n)){\mathcal D}\varphi,$$
where $0=t_0<t_1<t_2<...<t_n=1$ are marked points on the interval $[0,1]$, $\psi_i$ are ``observables" and $\psi_0\in \OO(L_0), \psi_1,....,\psi_{n-1}\in \OO(M),  \psi_n\in \OO(L_1)$.
Here ${\mathcal D}\varphi$ is the ill-defined ``Feynman measure" on the space of maps $P(L_0, L_1)$. 
The path integral discussed previously corresponds to the case when there are no observables.

Heuristically we should also specify the infinite-dimensional integration cycle. According to the discussion in Section \ref{infinite-dimensional vanishing cycles} this choice is equivalent to a choice of an element 
$\gamma_{loc}\in \R Hom_{\Z-mod}(Ext^\bullet_{\FF_{L_0\cup L_1,loc}}(L_0,L_1),\Z)$. \footnote{Unfortunately the tradition calls an element of the cohomology group ``the vanishing cycle", although it should be better called ``the vanishing cocyle".}

We would like to think of the above path integral in the following heuristic way. We have an algebra $A$ over $\C[[\hbar]]$,  a right $A$-module $E_0^{right}$ and left $A$-module $E_1^{left}$ as well as the linear map $\langle \bullet, \bullet \rangle_{\gamma_{loc}}: E_0^{right}\otimes_AE_1^{left}\to {\bf B}$, where the ring ${\bf B}=\overline{\C((\hbar))}[\log\,\hbar]$ was introduced in Section \ref {Fukaya categories}. Then for $\psi_0\in E_0^{right}, \psi_1,...,\psi_{n-1}\in A, \psi_n\in E_1^{left}$ we have the ``topological correlator"
$$\langle\psi_0,\psi_1...\psi_{n-1}(\psi_n)\rangle_{\gamma_{loc}}\in {\bf B}.$$

E.g. in the case 1) from the previous subsection we have $A=(\OO(M)[[\hbar]], \ast)$. The path integral interpretation implies that $1_M\in \OO(M)\subset A$ is the unit for the $\ast$-product.
Furthermore the above integral can be understood as a way to endow $E_0^{right}=\OO(L_0)[[\hbar]], E_1^{left}=\OO(L_1)[[\hbar]]$ with the  structures of right and left $A$-modules, such that  modulo $\hbar$ they coincide with the usual $\OO(M)$-module structures on $\OO(L_i),i=0,1$. In particular one has canonical generators $1_{L_0}^{right}\in E_0^{right}, 1_{L_1}^{left}\in E_1^{left}$ corresponding to the functions identically equal to $1$.
This example can be considered as a motivation for the future discussion of the notion of quantum wave function. 

Then the above path integral without observables can be interpreted as the pairing 
$\langle 1_{L_0}^{right},1_{L_1}^{left}\rangle_{\gamma_{loc}}\in {\bf B}$. This pairing is covariant with respect to the natural action of the group $\Z$ acting as the monodromy on $\gamma_{loc}$ and ${\bf B}$. Getting rid of the dependence on $\gamma_{loc}$ we can interpret the pairing as an element $\langle 1_{L_0}^{right},1_{L_1}^{left}\rangle\in (Ext^\bullet_{\FF_{L_0\cup L_1,loc}}(L_0,L_1)\otimes {\bf B})^\Z$.

Recall that we have two canonical objects $E_{L_i}^{DR,loc}, i=0,1$ in the local category $Hol_{L_0\cup L_1,loc}$.
By the local Riemann-Hilbert correspondence discussed in Section  \ref{Fukaya categories}  the pairing $\langle 1_{L_0}^{right},1_{L_1}^{left}\rangle$ can be identified with an element 

$$\mu_{L_0,L_1}\in Ext^{\bullet}_{Hol_{L_0\cup L_1,loc}}(E_{L_0}^{DR,loc}, E_{L_1}^{DR,loc}).$$

Recall that the category $Hol_{L_0\cup L_1,loc}$ is fully faithfully embedded in the category $D^b(A-mod)$, where $A=(\OO(M)[[\hbar]], \ast)$. Under this embedding the object $E_{L_1}^{DR,loc}$ is identified with $E_1^{left}$ whereas $E_0^{DR,loc}$ is identified with the dual object ${\bf D}(E_0^{right})=Ext^n_{mod-A}(E_0^{right},A)$.

Notice that we have an   isomorphism 
$$E_0^{right}\otimes_A E_1^{left}\simeq Ext^\bullet_{A-mod}({\bf D}(E_0^{right}), E_1^{left})[n].$$

We suggest that under this isomorphism the element $\mu_{L_0,L_1}$ corresponds to $1_{L_0}^{right}\otimes 1_{L_1}^{left}$.

Recall that  the quantum mechanics interprets our  path integral with $H=0$ as the pairing $\langle \psi_{L_0}|\psi_{L_1}\rangle$ of the quantum wave functions corresponding to the boundary conditions $L_0$ and $L_1$ respectively. We are going to develop a theory of quantum wave functions which will identify this pairing with the one discussed above.

One can also ask about the global analogs of the above data, i.e. if there exist:

1) an element $\mu\in Ext^n_{Hol_{glob}}(E_{L_0}^{DR,glob}, E_{L_1}^{DR, glob})$ corresponding to the volume form $vol_X$, where $n={dim_\C M\over{2}}$;

2) a class $\gamma\in Ext^n_{\FF_{glob,\hbar}}(E_{L_0}^{Betti, glob,\hbar}, E_{L_1}^{Betti, loc, \hbar})$ corresponding to the ill-defined global integration cycle in the infinite-dimensional path integral. 

In order to make sense of 1) by analogy with Section \ref{Fukaya categories} one should assume that the quantized algebra and the modules are defined over the ring $\C\{\hbar\}$ of analytic germs.\footnote{The same assumption is presumably necessary for the resurgence of the arising formal series in $\hbar$.}

Finally, in the case $H\ne 0$ there is a modification of the above proposal based on the notion of transport of quantum wave functions which we will discuss in Section \ref{transport of wave functions}.

\subsection{Reminder on Harish-Chandra pairs}\label{Harish-Chandra pairs}

 
 In this as well as the next subsection we are going to discuss the notion of {\it quantum wave function structure} which underlies the concept of quantum wave function mentioned previously.
 A (left) quantum wave function $\psi_L:=\psi_L^{left}$ will be  a cyclic vector  in the (left) holonomic $DQ$-module of WKB type associated with a complex Lagrangian submanifold $L\subset M$  (see e.g. [AgSch], [KasSch1] ). In order to define the holonomic $DQ$-module corresponding to  $L$ one has to fix a square root $K_L^{1/2}$. But a choice of  cyclic vector is not automatic. It depends on the additional structure on $M$ which we will define below and call it the quantum wave function structure.

We are going to utilize the approach which is based on the ideas of formal differential geometry of Gelfand-Kazhdan (see [GeKazh]). This subsection is devoted to the reminder of the related notion of Harish-Chandra pair. Before giving formal definitions we explain our main example of Harish-Chandra pair.

Let us endow $\C^{2n}$  with the standard holomorphic symplectic form, and let $\OO_{\C^{2n},\hbar}$ be the corresponding sheaf of quantized analytic functions endowed with the Moyal product. We are going to consider the associated formal Weyl algebra $W_q=\C[[q_1,...,q_n]][[\hbar\partial_{q_1},...,\hbar\partial_{q_n}]][[\hbar]]$ as a Lie algebra over $\C$. We endow $W_q$ with a grading such that $deg\,\hbar=2, deg\, q_i=1$. Then $deg\,\hbar\partial_{q_i}=1$. 

Consider the naturally graded Lie algebra 
$$\g_{max}=\prod_{k\ge -1}\g_{max}^k:={1\over{\hbar}}W_q/\C\cdot{1\over{\hbar}}.$$ 

Then we have a surjective  map  of vector spaces (not of Lie algebras) $\g_{max}\to \C^{2n}$, which is the projection to the vector space of elements of degree $-1$.  One can check by a direct computation that the kernel  of this map is a Lie subalgebra of $\g_{max}$. The degree zero component $\g_{max}^0\subset \g_{max}$ is naturally isomorphic to the  Lie algebra $Lie(Sp(2n, \C))$.
This Lie algebra   is in fact a  Lie algebra of a proalgebraic  group $G_{max,+}=\C^\ast\times Sp(2n,\C) \ltimes G_1$, where $G_1$ is a pronilpotent group.   The factor $\C^\ast$ corresponds to the central element $1=\hbar/\hbar\in \g_{max}$ of degree $0$. Notice that $G$ acts on $\g$ via the adjoint representation and on $\C^{2n}$ via $Sp(2n,\C)$.  
We have a short exact sequence of $\g$-modules
$$0\to Lie(G_{max,+})\to \g_{max}\to \C^{2n}\to 0.$$
 
 We would like to think of the pair $(\g_{max}, G_{max,+})$ as of the Harish-Chandra pair analogous to the pair $(Der(\C[[x_1,...,x_k]]), Aut_0(\C[[x_1,...,x_k]])$, i.e. the pair consisting of the Lie algebra of formal vector fields on $\C^k$ and the group of formal automorphisms preserving the origin (see [GeKazh]).  In the framework of deformation quantization of symplectic manifolds it was probably V. Drinfeld who first suggested to use the formalism of Harish-Chandra pairs. Let us recall the corresponding terminology.
 
 \begin{defn}\label{Harish-Chandra pair} Let $G_+$ be an affine proalgebraic group and $\g$ be a Lie algebra which is dual to a countably-dimensional Lie coalgebra. We will say that the pair $(\g, G_+)$ is a Harish-Chandra pair if the following data a), b) satisfying the condition c) are given:
 
 a) A continuous action of $G_+$ on $\g$.
 
 b) An embedding of $Lie(G_+)$ to $\g$ such that $\g/Lie(G_+)$ is a finite-dimensional vector space.

 c) The infinitesimal action of $Lie(G_+)$ on $\g$ given by a) coincides with the adjoint $Lie(G_+)$-action given by b).

 \end{defn}
 
 \begin{rmk}\label{dual to enveloping algebra} If $\g=Lie(G)$ is a finite-dimensional Lie algebra of a Lie group then the formal spectrum of the dual $U(\g)^\ast$ to its universal enveloping algebra can be identified with the formal neighborhood of $1\in G$. This geometry corresponds to the Harish-Chandra pair $(\g,\{1\})$. In general the Harish-Chandra pair $(\g, G_+)$ could be understood as a formal thickening of the proalgebraic group $G_+$ (i.e. it is a formal group in some directions).
 
 \end{rmk}

Harish-Chandra pairs naturally form a category.  

Let now $X$ be  a smooth proalgebraic variety  $X$ and $(\g,G_+)$ be a Harish-Chandra pair.
\begin{defn}\label{structure on manifold}
 Suppose we are given a $G_+$-torsor $\widehat{X}\to X$ together with a homomorphism $\g\to Vect(\widehat{X})$ of Lie algebras  of which agrees with the given actions of $G_+$ on $\g$ and $\widehat{X}$. In this case we will speak about  $(\g,G_+)$-action on $\widehat{X}$.

If in addition the induced linear map $\g\to T_x\widehat{X}$ is an isomorphism for each point $\widehat{x}\in \widehat{X}$ we will speak about
 free transitive action of $(\g, G_+)$ on $\widehat{X}$. In the latter case we will also call $\widehat{X}$ a $(\g,G_+)$-torsor on $X$, and say that $X$ is endowed with a  {\it $(\g,G_+)$-structure}.
 \end{defn}

In order to display in the notation the Harish-Chandra pair acting on $\widehat{X}$ we will sometimes write $\widehat{X}_{(\g,G_+)}$ instead of just $\widehat{X}$. We will omit the subscript in case if the Harish-Chandra pair is clear from the context.

\begin{rmk}\label{Harish-Chandra for non-algebraic manifolds}
The  above definition can be naturally reformulated in the cases when $X$ is  a  $C^\infty$ manifold or a complex analytic one. We will be using these versions without further comments.

\end{rmk}

\begin{defn}\label{strict morphisms}
A morphism $f: (\g_1,G_{1,+})\to (\g_2,G_{2,+})$  of Harish-Chandra pairs is called strict is the induced map  $\g_1/Lie(G_{1,+})\to \g_2/Lie(G_{2,+})$ is an isomorphism.
\end{defn}

Let $f$ be a  strict morphism. If $X$ carries a $(\g_1, G_{1,+})$-structure, then it carries the induced $(\g_2, G_{2,+})$-structure.
\begin{defn}\label{lift of structure}
 In this case we will say that the $(\g_1, G_{1,+})$-structure on $X$ is a lift of the $(\g_2, G_{2,+})$-structure on $X$.
\end{defn}

Notice that the lift (it is sometimes called in the literature  the reduction or  descend of Harish-Chandra pairs) of a structure is non-canonical. 

One can naturally define the notion of morphism of torsors  associated with different Harish-Chandra pairs.
Then the above-defined lift  of structures gives rise to  a morphism $\widehat{X}_{(\g_1, G_{1,+})}\to \widehat{X}_{(\g_2, G_{2,+})}$ of torsors on $X$.

One can form a category of $(\g, G_+)$-structures on $X$  with objects being different Harish-Chandra pairs $(\g, G_+)$.
The natural notion of a morphism uses strict morphisms of Harish-Chandra pairs.

\begin{defn}\label{central extensions}
A strict morphism $f$ is called central extension if the induced homomorphisms $\g_1\to \g_2$ and $G_{1,+}\to G_{2,+}$ are epimorphisms whose kernels are central and moreover the action of the ``kernel group" $G_f:=Ker(G_{1,+}\to G_{2,+})$ on $\g_1$ is trivial.

\end{defn}

\begin{rmk}\label{Deligne class}
If $X$ carries a $(\g_2,G_{2,+})$-structure and $f: (\g_1,G_{1,+})\to (\g_2,G_{2,+})$ is a central extension then we have a naturally defined  ``obstruction class" $\eta:=\eta_f\in H^2(X, G_f)$.
The lift of $(\g_2,G_{2,+})$-structure to $(\g_1,G_{1,+})$-structure exists if and only if $\eta_f=0$. 
In this case the set of classes of isomorphisms of lifts is a torsor over $H^1(X, G_f)$. The group of automorphisms of any lift is naturally isomorphic to $H^0(X, G_f)$.

\end{rmk}

Let  us illustrate the above discussion in the example of symplectic structure on a complex analytic manifold of the fixed dimension $2n$. 

 Let $\g_{cl}$ be the Lie algebra of Hamiltonian vector fields on the formal completion of $\C^{2n}$ at $0\in \C^{2n}$, and $G_{cl,+}$ be the proalgebraic group of formal symplectomorphisms of $\C^{2n}$ preserving the origin. For a complex symplectic manifold $(M,\omega^{2,0})$
let $\widehat{M}_{cl}$ be the $G_{cl,+}$-torsor over $M$ whose points are pairs consisting of a point of $M$ and  formal Darboux coordinates at the point. 
It is easy to see that a $(\g_{cl}, G_{cl,+})$-structure on  $M$  is equivalent to a symplectic structure on $M$.

The Harish-Chandra pair $(\g_{cl}, G_{cl,+})$ has a natural central extension $(\g_{cl}^{\prime\prime}, G_{cl,+}^{\prime\prime})$ (the reason for the notation will become clear later). Here $\g_{cl}^{\prime\prime}=\C[[p_1,...,p_n,q_1,...,q_n]]$ is endowed with the standard Poisson bracket and $G_{cl,+}^{\prime\prime}=\C\times G_{cl,+}$. Then the obstruction class $\eta$ is equal to $[\omega^{2,0}]\in H^2(M,\C)$.

\begin{rmk}\label{special cases of lifts}
a) If $M$ carries $(\g_{cl}^{\prime\prime}, G_{cl,+}^{\prime\prime})$-structure then the induced symplectic structure is exact.

b) Notice that there is a natural embedding of the group $2\pi i\Z$ into the first factor of  $G_{cl,+}^{\prime\prime}$. If $M$ carries a $(\g_{cl}^{\prime\prime}, G_{cl,+}^{\prime\prime}/2\pi i\Z)$ structure then the induced symplectic structure is the curvature form of the natural prequantization line bundle on $M$. Hence the class $[\omega^{2,0}]$ belongs to the image of $H^2(M,2\pi i\Z)$ in $H^2(M,\C)$.

\end{rmk}

 \subsection{Quantum wave function structures}\label{QWFS}

In this subsection we will apply the ideas of the previous one in the case of $(\g_{max}, G_+)$-structures.
Notice that have a natural epimorphism of Harish-Chandra pairs $(\g_{max},G_{max,+})\to (\g_{cl}, G_{cl,+})$.   Let $(M,\omega^{2,0})$ be a complex symplectic manifold. As we explained above this means that $M$ carries a  $(\g_{cl}, G_{cl,+})$-structure.

 \begin{defn}\label{manifold with wave functions}
A quantum wave function structure (QWFS for short) on $M$  is a lift of the $(\g_{cl}, G_{cl,+})$-structure on $M$ to a $(\g_{max}, G_{max,+})$-structure on $M$.
  
 \end{defn}
 
 Given such a lift we denote by $\widehat{M}:=\widehat{M}_{(\g_{max},G_{max,+})}$ the corresponding proalgebraic complex manifold. By definition it is  endowed with a free action of $G_{max,+}$ such that
  $\widehat{M}/G_{max,+}\simeq M$.
 Notice that by definition we have a $G_{max,+}$-equivariant morphism of Lie algebras $\g_{max}\to T_{\widehat{M}}$ which identifies $\g_{max}$ with the tangent space at any point, and which is compatible with the induced action of $Lie(G_{max,+})$.

 Thus the quantum wave function structure can be thought of as a refinement of the notion of symplectic structure. 
 
 \begin{rmk}\label{smooth star-product}
 One can give a similar definition of the quantum wave function structure for any real $C^\infty$ or complex algebraic symplectic manifold.
 \end{rmk}
 
There are two more Harish-Chandra pairs closely related to $(\g_{max}, G_{max,+})$ which are of independent interest.
 
 Consider a pair $(\g_{max}^\prime, G^\prime_{max,+})=({1\over{\hbar}}W_q/{1\over{\hbar}}\C[[\hbar]], G_{max,+}/exp(\C[[\hbar]]))$, where we take the quotient of the group $G_{max,+}$ by the central subgroup. 
Notice that $\g_{max}^\prime\simeq Der(W_q)$.

Let $\g_{max}^{\prime\prime}={1\over{\hbar}}W_q$ and $G^{\prime\prime}_{max,+}$ be the connected simply-connected proalgebraic group such that $Lie(G^{\prime\prime}_{max,+})=\g_{max}^{\prime\prime}\times_{\g_{max}^{\prime}}Lie(G_{max,+}^{\prime})$.
In plane words elements of $Lie(G^{\prime\prime}_{max,+})$ are such elements of ${1\over{\hbar}}W_q$ that the induced modulo $\hbar$ formal Hamiltonian vector fields vanish at $0$.

Then we have a sequence of strict morphisms of Harish-Chandra pairs:

$$(\g_{max}^{\prime\prime},G_{max,+}^{\prime\prime})\xrightarrow{f_1}(\g_{max},G_{max,+})\xrightarrow{f_2}(\g_{max}^{\prime},G_{max,+}^{\prime})\xrightarrow{f_3}(\g_{cl},G_{cl,+}).$$

Moreover $f_1,f_2,f_2\circ f_1$ are central extensions and the corresponding kernel groups are such as follows:
$$G_{f_1}={\C\over\hbar}\oplus2\pi i\Z,$$
$$G_{f_2}=\C^\ast\times {\hbar}\C[[\hbar]],$$
$$G_{f_2\circ f_1}={1\over \hbar}\C[[\hbar]].$$

 
Notice that a lift  of the $(\g_{cl}, G_{cl,+})$ structure on $M$ to a $(\g_{max}^\prime, G_{max,+}^\prime)$-structure on $M$ is equivalent to a choice of the   sheaf of quantized algebras $\OO_{M,\hbar}$. 

Recall (see Section \ref{star-product} and e.g. [BezKal])  that a sheaf of quantized algebras  gives rise to a class in $H^2(M,{1\over \hbar}\C[[\hbar]])$ with the initial term ${[\omega^{2,0}]\over \hbar}$ called the Deligne class of the quantization. In our language the Deligne class admits an alternative description, namely it is equal to   $\eta_{f_2\circ f_1}$ (we omit the proof). The fact that the initial term of the class $\eta_{f_2\circ f_1}$ is equal to ${[\omega^{2,0}]\over \hbar}$ can be seen from the following commutative diagram or equivalently, the morphism of central extensions \footnote{It is clear from the diagram  that the Harish-Chandra pair  $(\g_{cl}^{\prime\prime},G_{cl,+}^{\prime\prime})$ corresponds to the Harish-Chandra pair $(\g_{max}^{\prime\prime}, G_{max,+}^{\prime\prime})$, hence the similarity in the notation.}:

\[ \begin{tikzcd}
(\g_{max}^{\prime\prime}, G_{max,+}^{\prime\prime})\arrow{r} \arrow[swap]{d}&(\g_{max}^\prime,G_{max,+}^\prime) \arrow{d}  \\%
(\g_{cl}^{\prime\prime},G_{cl}^{\prime\prime}) \arrow{r}& (\g_{cl}, G_{cl,+})
\end{tikzcd}
\]

Since we have the natural epimorphism of Harish-Chandra pairs $(\g_{max}, G_{max,+})\to (\g_{max}^\prime, G^\prime_{max,+})$ we conclude that QWFS on $M$ gives rise to the sheaf $\OO_{M,\hbar}$. 

For a given QWFS we have the obstruction class $\eta_{f_1}\in H^2(M, {\C\over \hbar}\oplus 2\pi i\Z)$. The Deligne class of the sheaf $\OO_{M,\hbar}$ associated with our QWFS is equal to the image of $\eta_{f_1}\in H^2(M, {1\over \hbar}\C[[\hbar]])$ induced by the natural inclusion ${\C\over\hbar}\oplus 2\pi i\Z\hookrightarrow {1\over \hbar}\C[[\hbar]]$ which sends $2\pi i\in 2\pi i \Z$ to $2\pi i\in \C\subset {1\over \hbar}\C[[\hbar]]$. Therefore the Deligne class can be written as
$${[\omega^{2,0}]\over \hbar}+\hbar^0\cdot j(\beta),$$
where $\beta\in H^2(M,2\pi i\Z)$ and $j: H^2(M,2\pi i\Z)\to  H^2(M,\C))$ is the natural homomorphism.

Furthermore if we have a $(\g_{max}^{\prime\prime}, G^{\prime\prime}_{max,+})$-structure on $M$ then we have the corresponding QWFS and the sheaf $\OO_{M,\hbar}$. The Deligne class in this case is equal to zero. In particular $[\omega^{2,0}]=0$ and the symplectic structure on $M$ is exact. For that reason a QWFS which can be lifted to $(\g_{max}^{\prime\prime}, G^{\prime\prime}_{max,+})$-structure will be called {\it exact}.

\begin{rmk}\label{Deligne class of categories}

Recall (see Section \ref{star-product}) that one can consider sheaves $\CC_M$ of $\C[[\hbar]]$-linear categories on $M$  locally on $M$ equivalent to the standard category  $\OO_{\C^{2n},\hbar}=(\OO_{\C^{2n}}, \ast_{Moyal})$-mod, {\bf without} the choice of an identification  modulo $\hbar$  with the sheaf of categories $\OO_M$-mod.  Autoequivalences of such sheaves of categories and automorphisms of autoequivalences are described  locally by the cross-module
$$\OO_{\C^{2n},\hbar}^\times\to Aut^0(\OO_{\C^{2n},\hbar}),$$
where $\OO_{\C^{2n},\hbar}^\times$ is the sheaf of invertible quantized functions, and $Aut^0(\OO_{\C^{2n},\hbar})$ is the group of such automorphisms of the quantum algebra which are equal to $id$ modulo $\hbar$. This cross-module is equivalent to another one, namely to
$$\C[[\hbar]]^\times\to \{1\},$$
where $\C[[\hbar]]^\times$ denote the group of invertible elements of the algebra $\C[[\hbar]]$.
Hence the equivalence classes of the above sheaves $\CC_M$  are parametrized by elements of 
$H^2(M, \C[[\hbar]]/2\pi i\Z).$ \footnote{This implies that the sheaf $\CC_M/\hbar \CC_M$ is a sheaf $\OO_M-mod$ twisted by a $\OO_M^\times$-gerbe obtained from a $\C^\ast$-gerbe via the natural embedding $\C^\ast\to \OO_M^\times$.}
For the canonical sheaf of categories $\CC_M^{can}$ from Section \ref{star-product} the corresponding class is equal to zero. 

In the case when we have an identification $\CC_M$ modulo $\hbar$ with the sheaf  $\OO_{M}$-mod,  the corresponding deformation class of the sheaf $\CC_M$ is equal to the image of the difference of the Deligne class of $\CC_M$ and ${[\omega^{2,0}]\over \hbar}$ under the natural map $H^2(M,\C[[\hbar]])\to H^2(M, \C[[\hbar]]/2\pi i\Z)$. Hence  in the case when we have a QWFS, the corresponding class is zero and the sheaf of categories $\CC_M$ is equivalent to $\CC_M^{can}$. One can show that the equivalence of $\CC_M$ with $\CC_M^{can}$ is also canonical.

\end{rmk}

\begin{rmk}\label{QWFS and generators}

Let us explain  why the notion QWFS can be reformulated as a choice of generator of $\CC_M^{can}$. First we remark that having a QWFS is equivalent to having a triple $T_1$=(QWFS, sheaf of quantized algebras $\OO_{M,\hbar}$, isomorphism of $\OO_{M,\hbar}$ with the sheaf of quantized algebras coming from the QWFS). Indeed, the last two elements of the triple form a non-empty contractible groupoid.

Similarly a choice of generator $G$ of the sheaf of categories $\CC_M^{can}$ is equivalent to a choice of a triple $T_2$=(generator $G$, sheaf of quantized algebras $\OO_{M,\hbar}$, an isomorphism of $\OO_{M,\hbar}$ with the sheaf $\underline{End}(G)$). 

The groupoids of triples $T_1$ and $T_2$ are naturally fibered over the groupoid of sheaves of quantized algebras $\OO_{M,\hbar}$ via the forgetful morphism. Hence we have to show that the fibers are equivalent.

In the case of triples $T_1$ the fiber over a fixed sheaf $\OO_{M,\hbar}$ is non-empty iff the image of the Deligne class of $\OO_{M,\hbar}$ in $H^2(M,\C^\ast\times \hbar\C[[\hbar]])$ vanishes. In this case the set of isomorphism classes of objects of the fiber is a torsor over $H^1(M, \C^\ast\times \hbar\C[[\hbar]])$. The  automorphism group of each object is $H^0(M, \C^\ast\times \hbar\C[[\hbar]])$. All this follows from the obstruction theory of the above-defined central extension $f_2$.

In the case of triples $T_2$ the fiber over a fixed sheaf $\OO_{M,\hbar}$ is equivalent to the groupoid of equivalences between the category $\OO_{M,\hbar}-mod$ and $\CC_M^{can}$.
Again it is non-empty iff the image of the Deligne class of $\OO_{M,\hbar}$ in $H^2(M,\C^\ast\times \hbar\C[[\hbar]])$ vanishes.
The set of isomorphism classes of objects of the fiber is a torsor over $H^1(M, \C^\ast\times \hbar\C[[\hbar]])$, and the automorphism group of each object is $H^0(M, \C^\ast\times \hbar\C[[\hbar]])$. This follows from the Remark \ref{Deligne class of categories}.
Furthermore one can show that the fiber in the case of $T_1$-triples is {\bf canonically} isomorphic to the corresponding fiber  in the case of $T_2$-triples. Therefore the notion of QWFS is indeed equivalent to the notion of generator of $\CC_M^{can}$.
\end{rmk}

\begin{rmk}\label{twisting by gerbes}

Given a $\C^\ast$-gerbe ${\mathcal G}$ on $M$ with the class $cl({\mathcal G})\in H^2(M,\C^\ast)$ one can twist the sheaf of categories $\CC_M$ by ${\mathcal G}$ using the natural morphism $\C^\ast\to Aut(\CC_M)$. Similarly one can speak about a QWFS twisted by ${\mathcal G}$.\footnote{One can replace here $\C^\ast$-gerbes by $\C[[\hbar]]^\times$-gerbes.}

Then we can generalize the Remark \ref{QWFS and generators} such as follows: the notion of QWFS twisted by ${\mathcal G}$ is equivalent to a choice of  generator of the sheaf of categories $\CC_M^{can}$ twisted by ${\mathcal G}$. In the notation of  Section \ref{holonomic DQ-modules} a choice of the gerbe $\mathcal{G}$ corresponds to the choice of $[B_{init}]$. 
The existence of  the sheaf of quantized algebras $\OO_{M,\hbar}$  imposes the following restriction: $[B_{init}]$ belongs to the image of $H^2(M,\C)$ under the exponential map $H^2(M,\C)\to H^2(M,\C^\ast)$. Equivalently, the image of $[B_{init}]$ under the coboundary map $H^2(M,\C^\ast)\to H^3_{tors}(M,\Z)$ is trivial.

In particular we can apply the procedure of twisting by gerbe to the case of  the sheaf of $\hbar$-microdifferential operators on $M=T^\ast X$. The latter has the Deligne class equal to $-{1\over 2}c_1(T^\ast X)$ (see Example \ref{examples of Deligne classes}  a)).

\end{rmk}

The notion of QWFS (possibly twisted by a gerbe) will be used later in order to spell out resurgent properties of perturbative expansions of functional integrals and their hypothetical interpretation via the Riemann-Hilbert correspondence.

Notice that in the  Remark \ref{twisting by gerbes} the Deligne class has  at most two terms: the one with $\hbar^{-1}$ and the one with $\hbar^0$. This property holds beyond the case of cotangent bundles, e.g. it holds in the case of quantum tori or twisted cotangent bundles (see Example \ref{examples of Deligne classes} d)). There is another class of examples for which the Deligne class is more complicated. Those include filtered quantizations of smooth affine symplectic varieties {\it analytically depending on $\hbar$}. For this class of examples which we do not describe explicitly, the Deligne class is expected to belong to
$H^2(M, {1\over \hbar}\C\{\hbar\})\subset H^2(M, {1\over \hbar}\C[[\hbar]])$, and moreover can be extracted from the first non-trivial term of the Hodge filtration on the periodic cyclic homology of the quantized algebra (see [Ko6], Theorem 1.36.1, where the cohomology class of the symplectic form should be replaced by the Deligne class). For such examples we expect that  the classes of  gerbes defining twisted QWFS belong to $H^2(M,\C\{\hbar\}^\times)\subset H^2(M,\C[[\hbar]]^\times)$.


\

Let us give few examples of QWFS.

1) Suppose that the symplectic manifold $(M,\omega^{2,0})$ is endowed with an affine structure for which the symplectic form $\omega^{2,0}$ is covariantly constant. In this case we say that $M$ carries an {\it affine symplectic structure} (examples are $\C^{2n}, (\C^{\ast})^{2n}$ or a compact complex symplectic torus). The corresponding Harish-Chandra pair is $(sp(2n)\ltimes \C^{2n}, Sp(2n))$. Then we have a natural morphism of Harish-Chandra pairs 
$$(sp(2n)\ltimes \C^{2n}, Sp(2n))\to (\g_{max}, G_{max,+}).$$
Therefore any manifold with affine symplectic structure carries a natural QWFS. The associated sheaf of quantized algebras $\OO_{M,\hbar}$ is naturally isomorphic to  $\OO_M[[\hbar]]$ endowed with the Moyal $\ast$-product. 

1a) Let us assume that the affine symplectic structure on $M$ is conical. This means that $M$ is endowed with a vector field $Eu$ (Euler vector field) such that $Eu$ preserves the affine structure, its linear part is equal to the identity operator (hence $Lie_{Eu}(\omega^{2,0})=2\omega^{2,0}$). In other words in local affine coordinates 
$$Eu=\sum_{1\le i\le 2n}(x_i+c_i)\partial_{x_i},\quad \omega^{2,0}=\sum_{1\le i<j\le 2n}\omega_{ij}dx_i\wedge dx_j,\qquad c_i,\omega_{ij}\in \C.$$
Then the QWFS can be canonically lifted to $(\g_{max}^{\prime\prime},G_{max,+}^{\prime\prime})$-structure and hence it is exact.

2) Assume that $(M,\omega^{2,0})$ is a symplectic manifold endowed with a Lagrangian foliation, i.e. an integrable Lagrangian subbundle of the tangent bundle $T_M$. This structure can be described in terms of the following Harish-Chandra pair $(\g,K)$ where 
$$\g=\C[[q_1,...,q_n]]/\C\rtimes \bigoplus_{1\le i\le n}{1\over \hbar}\C[[q_1,...,q_n]]\hbar\partial_{q_i}\subset \g_{max}.$$
The algebra Lie $\g$ is naturally isomorphic to the Lie algebra of Hamiltonian vector fields on the formal neighborhood of $0\in \C^{2n}$ preserving the coordinate foliation spanned by $\partial_{p_i}, 1\le i\le n$.
The group $K$ is defined as the group of  automorphisms of this formal neighborhood preserving  symplectic structure, foliation and $0\in \C^{2n}$. Then $(\g,K)$ is naturally embedded in $(\g_{max}, G_{max,+})$. Hence $M$ carries a QWFS.

2a) Let $M=T^\ast X$ with the standard symplectic structure. Then $M$ carries a natural Lagrangian foliation by cotangent fibers.
One can show that the corresponding  QWFS is exact (this is due to the existence of the  Euler vector field which expands $\omega^{2,0}_{T^\ast X}$ and preserves the  foliation). The corresponding sheaf of quantized algebras $\OO_{M,\hbar}$ is naturally isomorphic to the sheaf of $\hbar$-microdifferential operators acting on $K_X^{1/2}$ (see Example \ref{examples of Deligne classes} c)).

On the contrary in the general case of a twisted cotangent bundle (see Example \ref{examples of Deligne classes} d) ) the corresponding QWFS is not necessarily exact.


3) More generally, a  QWFS on $M$ arises naturally from a  {\it semi-affine symplectic structure} on $M$ (see [KoSo11], Sect. 8.1). This follows from the observation that there exists a homomorphism $\g_{m,n}\to W_q$, where $\g_{m,n}$ is the Lie algebra from the loc. cit., Sect. 8.1. This construction includes those from the previous two examples as special cases.

 4) In all above examples class $\beta\in H^2(M,2\pi i\Z)$ defined above the Remark \ref{Deligne class of categories} was equal to zero. There are examples when $\beta\ne 0$. E.g. let $M=T^\ast X$ and $\beta$ is the pull-back of $c_1(\LL)$, where $\LL$ is a holomorphic line bundle on $X$. The corresponding sheaf of quantized algebras is isomorphic to the sheaf of $\hbar$-microdifferential operators on $\mathcal{L}\otimes K_X^{1/2}$.


 \begin{rmk}\label{Hilbert-90}
 
a)   A $(\g_{max}^{\prime\prime}, G^{\prime\prime}_{max,+})$-structure on $M$ can be interpreted as symmetries of a family over formal punctured disc 
 of infinite-dimensional vector spaces. In a similar way a $(\g_{max}^{\prime}, G^{\prime}_{max,+})$-structure on $M$ can be interpreted as the symmetries modulo rescaling by $\C[[\hbar]]^\times$, while a $(\g_{max}, G_{max,+})$-structure on $M$ corresponds to the symmetries modulo an overall rescaling.

b)  Informally speaking a $(\g_{max}^{\prime\prime}, G^{\prime\prime}_{max,+})$-structure can be thought of as an analog of the Hilbert-90 theorem since we claim that the first cohomology with coefficients in some infinite-dimensional general linear group is trivial.

c) In the case of $(\g_{max}, G_{max,+})$-structure the  vector spaces over $\C((\hbar))$ discussed in a) and b) are defined only up to rescaling by $\C^\ast\subset \C((\hbar))^\times$. \footnote{There is a clash of notation throughout the paper: the set of non-zero complex numbers (or the corresponding algebraic variety) is denoted by $\C^\ast$, while in most of the other cases the set of invertible elements of an algebra $A$ is denoted by $A^\times$. The notation $V^\ast$ is reserved for the dual to $V$ vector spaces, with the only exception $V=\C$. We hope this will not lead to a confusion, since we never use the notation $\C^\ast$ for the dual vector space to $\C$.}
  
   \end{rmk}

 There is a version of the above constructions which will be useful later, when we will discuss the appearance of half-forms in the story. Namely, let $\sigma: W_q\to W_q^{opp}$ be an isomorphism of $\C$-algebras defined by the involution $\sigma: \hbar \mapsto -\hbar, \hat{q}_i=q_i\mapsto \hat{q}_i, \hat{p}_i=\hbar\partial_{q_i}\mapsto \hat{p}_i, 1\le i\le n$. Then we have a graded Lie subalgebra $\g_{max}^{odd}=\{a\in {1\over{\hbar}}W_q|\sigma(a)=-a\}/\C\cdot{1\over{\hbar}}\subset \g_{max}$ and  the corresponding group $G_{max,+}^{odd}\subset G_{max,+}$. The Harish-Chandra pair $(\g_{max}^{odd}, G_{max,+}^{odd})$  maps surjectively to the Harish-Chandra pair $(\g_{cl}, G_{cl,+})$. 
 
 The direct sum of graded components of $\g_{max}^{odd}$ of degrees $-1$ and $0$ has a basis
 $${\hat{p}_i\over \hbar}, {\hat{q}_i\over\hbar},{\hat{p}_i\hat{p}_j\over\hbar},  {\hat{q}_i\hat{q}_j\over \hbar}, {\hat{q}_i\hat{p}_j+\hat{p}_j\hat{q}_i\over {2\hbar}}, \quad 1\le i, j\le n,$$
 and is naturally isomorphic to the Lie algebra of affine symplectic vector fields on $\C^{2n}$. 
 The factor $1/2$ in these formulas is responsible for the appearance of half-densities in the next subsection. 
 A large class of examples of the Harish-Chandra pairs  $(\g_{max}^{odd}, G_{max,+}^{odd})$ arises from the semi-affine symplectic structures.

Let us illustrate the relevance of the new class of Harish-Chandra pairs   $(\g_{max}^{odd}, G_{max,+}^{odd})$   by considering a family of homomorphisms $\varphi_\lambda$ of the Lie algebra $Vect_n\ltimes \Omega^{1,cl}$ (this Lie algebra consists of symplectic vector fields on $\C^{2n}$) into $\g_{max}$ given by $\varphi_\lambda(\xi,\alpha)={1\over{\hbar}}\xi+{\lambda \over{\hbar}}\xi div(\xi)+ {f_\alpha\over{\hbar}}$, where $df_\alpha=\alpha$. Using the natural action of $W_q$ on $\C[[q_1,...,q_n]]((\hbar))$ we obtain a $\lambda$-family of such actions of $Vect_n\ltimes \Omega^{1,cl}$. But only for $\lambda={1\over{2}}$ the image of the corresponding homomorphism into $End(\C[[q_1,...,q_n]]((\hbar))$ belongs to the image of $\g_{max}^{odd}$.
 
 As follows from the discussion in the Section \ref{star-product} a choice of QWFS gives rise  to a choice of the functor $\Phi$ which gives an equivalence of the abstractly defined sheaf of categories $\CC_M$ and the category of $\OO_{M,\hbar}$-modules. Since a choice of the QWFS is equivalent to the choice of a lift of the $(\g_{cl}, G_{cl,+})$-structure on $M$ to the $(\g_{max}, G_{max,+})$-structure on $M$, the latter structure gives also a choice of a generator of $\CC_M$, namely the sheaf $\OO_{M,\hbar}$. 
 
 One can ask what in these terms means a lift of the $(\g_{cl}, G_{cl,+})$-structure on $M$ to the $(\g_{max}^{odd}, G_{max,+}^{odd})$-structure on $M$? 
 Notice that the canonical sheaf of categories $\CC_M$ is  automatically endowed with an  equivalence functor $\Psi: \sigma^\ast(\CC_M)\to \CC_M^{opp}$ induced by the homomorphism $\sigma: \C[[\hbar]]\to\C[[\hbar]], \hbar\mapsto -\hbar$. There is an ``opposite'' functor $\Psi^{opp}: \sigma^\ast(\CC_M^{opp})\to \CC_M$, and we have $\Psi^{opp}\circ\Psi=id$. Let $\EE\in Ob(\CC_M)$ be a generator, i.e. a generator of each category $\CC_M(U)$, where $U$ is an open subset, satisfying natural compatibility conditions for different open subsets as well as the condition $\EE_M=\OO_M$ mod $\hbar$. Then the lift from a  $(\g_{cl}, G_{cl,+})$-structure on $M$ to a $(\g_{max}^{odd}, G_{max,+}^{odd})$-structure on $M$ means that there is a isomorphism $\phi_{\EE}:\EE\to \Psi(\sigma^\ast(\EE))$ which is equal to $id$ modulo $\hbar$ and such that $\phi_{\EE}^{opp}=\phi_{\EE}$. In this case we will say that our QWFS is {\it symmetric}.

 \subsection{Remarks about the relation to quantized integrable systems}
In this subsection we will make few rather informal comments on the relation of QWFS with quantized integrable systems.

We start with an observation that in the framework of Example 2) from Section \ref{QWFS}, i.e. for a symplectic manifold endowed with Lagrangian foliation there exists a canonical subsheaf of  commutative $\C[[\hbar]]$-algebras of the sheaf $\OO_{M,\hbar}$ of quantized algebra of analytic functions. This subsheaf is canonically isomorphic to the sheaf of $\C[[\hbar]]$-valued functions on $M$ constant along the foliation. In the case when the foliation is a Lagrangian {\it fibration} over an $n$-dimensional complex manifold $B$ the above structure is essentially a structure of quantized integrable system (cf. e.g. [GV]). In general the groupoid of quantized integrable systems such that their classical limits are  Lagrangian fibrations over the same base $B$ is equivalent to the groupoid of classical integrable systems over $B$ depending formally on $\hbar$. \footnote{One can generalize this correspondence allowing $B$ to vary formally with respect to $\hbar$, or more generally to allow Lagrangian foliations instead of fibrations.} Nevertheless the above construction of the canonical quantization of a classical integrable system is often not adequate. E.g. in the case when the Lagrangian foliation has singularities,  the above canonically quantized sheaf of algebras in general does not admit an extension to singularities. All that will be discussed in the forthcoming joint paper of the first author and Alexander Soibelman.

\

Let us consider a quantum integrable system such that its classical limit is a possibly singular Lagrangian fibration $\pi: M\to B$.  Assume that the sheaf of commutative subalgebras of $\OO_{M,\hbar}$ is identified (as a sheaf of $\C[[\hbar]]$-algebras) with $\pi^\ast(\OO_B)[[\hbar]]$. Then we have a canonical holomorphic family of objects of the category $\OO_{M,\hbar}-mod$ parametrized by $B$, such that the object corresponding to $b\in B$ is $E_b=\OO_{M,\hbar}\otimes_{\pi^\ast(\OO_B)}\pi^\ast(\C_b[[\hbar]])$, where $\C_b$ is the sky-scraper sheaf sitting at $b$. 

Conversely, suppose we have a quantized sheaf of categories $\CC_M$ and a holomorphic family of objects $(E_b)_{b\in B}$ such that in the classical limit this family becomes the family $(\OO_{\pi^{-1}(b)})_{b\in B}$ for a singular Lagrangian fibration $\pi:M\to B$. We claim that in this case one gets a canonical sheaf of quantized algebras $\OO_{M,\hbar}$ and a subsheaf of maximal commutative $\C[[\hbar]]$-subalgebras of $\OO_{M,\hbar}$ which is canonically isomorphic to $\pi^\ast(\OO_B)[[\hbar]]$.
Indeed, locally on $M$ after a choice of generator $G$ of $\CC_M$ (so $\CC_M$ is canonically equivalent to the sheaf of categories  $\underline{End}(G)^{op}-mod$) we can describe a holomorphic family of objects over $B$ as a module $V$ over $\underline{End}(G)^{op}\otimes_{\C[[\hbar]]} \pi^\ast(\OO_B[[\hbar]])$ which is locally free of rank $1$ over $\underline{End}(G)^{op}$.  One can see that $V$ (considered locally as an object of $\CC_M$ which is endowed with a morphism $\pi^\ast(\OO_B)\to V$) does not depend on the choice of $G$ up to a canonical isomorphism. Thus $V$ gives rise to a canonical global family of objects of $\CC_M$ which are locally generators. We define the sheaf of $\C[[\hbar]]$-algebras as $\OO_{M,\hbar}=\underline{End}(V)^{op}$. The morphism  $\pi^\ast(\OO_B)\to V$ gives an embedding of sheaves of algebras $\pi^{\ast}(\OO_B)\to \OO_{M,\hbar}$. In physics terminology the family of objects $(E_b)_{b\in B}$ is called a {\it canonical coisotropic brane} ${\mathcal B}_{cc}$. In the ``real life" examples of quantized integrable systems of Hitchin type  with  the gauge group $GL(m)$, the variety $B$  appears as the locus of {\it opers} in $M_{DR}$, which is a component of the moduli space of holonomic $\hbar-D$-modules on the base curve.

We see that the sheaf of quantized algebras on $M$  arises naturally when one has just an ${1\over 2}dim_\C M$-dimensional holomorphic  family of ``holonomic objects"  with an appropriate behavior in the classical limit.\footnote{More generally if the classical limit of the family of objects gives a singular Lagrangian fibration on a $k$-fold ramified covering of $M$ then instead of the sheaf of algebras $\OO_{M,\hbar}$ it is plausible to expect to have a sheaf of algebras on $M$ locally isomorphic to the sheaf $Mat(k\times k,\OO_{M,\hbar})$ containing a subsheaf of commutative subalgebras.} Nevertheless there is no natural QWFS in this situation, since e.g. the Deligne class has in general more than two terms as a series in $\hbar$.

 \subsection{Quantizable compactifications and QWFS}\label{quantizable compactifications and QWFS}
 
Let us summarize part of the previous discussion from the point of view of different versions of deformation quantization discussed in Section \ref{star-product}. Recall that a QWFS gives rise to a sheaf $\OO_{M,\hbar}$ of algebras over $\C[[\hbar]]$ quantizing $\OO_M$ together with an equivalence of the category of $\OO_{M,\hbar}$-modules with the canonical sheaf of $\C[[\hbar]]$-linear categories $\CC_M$ associated to $(M,\omega^{2,0})$. The corresponding Deligne class of such quantizations is ${[\omega^{2,0}]\over \hbar}$, while in general we can get any formal path ${[\omega^{2,0}]\over \hbar}+\sum_{i\ge 0}[\alpha_i]\hbar^i, [\alpha_i]\in H^2(M,\C)$. In order to reduce from general deformation quantizations to those which admit QWFS one has to choose a generator of $\CC_M$. We describe below a class of examples where such a generator can be described more explicitly.
 
 Recall that according to [Ko2] (see also Section \ref{holonomic DQ-modules}) one can construct a filtered algebra $\OO_\hbar(M)$ over $\C[[\hbar]]$  provided $M$ is an algebraic affine complex symplectic manifold which admits a {\it quantizable compactification}. The latter means that there is a smooth projective variety $\overline{M}\supset M$ such that:
 
 a) $H^1(\overline{M}, \OO_{\overline{M}})=H^2(\overline{M}, \OO_{\overline{M}})=0$;
 
 b) the divisor $D=\overline{M}-M$ is ample and with normal crossings;
 
 c)  the Poisson bivector field corresponding 
 to the symplectic structure on $M$ is tangent to $D$ in a sense defined in loc.cit. and extends to $\overline{M}$ making it into a Poisson manifold with the open symplectic leaf $M$.
 
 \begin{conj}\label{QWFS from semiflat quantization} Under the above assumptions there is a canonical QWFS on $M$.

 \end{conj}

 Recall that under the above assumptions a)-c)  factors of  the filtered $\C[[t]]$-algebra $\OO_\hbar(M)$  are free $\C[[\hbar]]$-modules of finite rank. There is an involution which maps $\OO_\hbar(M)$ into the opposite algebra in such a way that $\hbar\mapsto -\hbar$, i.e.   we have $\OO_{-\hbar}(M)\simeq\OO_\hbar(M)^{opp}$ in the obvious notation. In terms of formal differential geometry this means that one has the corresponding Harish-Chandra pair $(\g^{odd}_{cl}, G_{cl,+}^{odd})$ acting on $\widehat{M}$. This Harish-Chandra pair is a quotient of  $(\g^{odd}_{max}, G_{max,+}^{odd})$.
 
 Then the Conjecture \ref{QWFS from semiflat quantization} can be spelled out more precisely.
 
\begin{conj}\label{existing of lift} Under the above assumptions a)-c) there is a lift of  $(\g^{odd}_{cl}, G_{cl,+^{odd}})$-structure on $M$ to a  $(\g^{odd}_{max}, G_{max,+}^{odd})$-structure on $M$.

\end{conj}
 
 Another problem which should be  related to the resurgence of arising series and  which will be discussed later on, is whether the $\C[[\hbar]]$-algebra $\OO_\hbar(M)$ can be defined over the smaller ring $\C\{\hbar\}\subset \C[[\hbar]]$ of germs of convergent series at $\hbar=0$. In such a case the above-mentioned formal series  ${[\omega^{2,0}]\over \hbar}+\sum_{i\ge 0}[\alpha_i]\hbar^i, [\alpha_i]\in H^2(M,\C)$ should be convergent in $H^2(M,\C)$. 
 
 \begin{conj}\label{analyticity of algebra}
 Under the assumptions a)-c) above the algebra $\OO_\hbar(M)$ over $\C[[\hbar]]$ is obtained by the extension of scalars from an algebra $\OO_\hbar(M)^{an}$ over $\C\{\hbar\}$ and the isomorphism $\OO_{-\hbar}(M)\simeq \OO_\hbar(M)^{opp}$ is obtained by the extension of scalars from the one for $\OO_\hbar(M)^{an}$.
 
 \end{conj}

 \subsection{Smooth Lagrangian submanifolds and quantum wave functions}\label{quantum wave functions}

 Let $L\subset M$ be a smooth Lagrangian subvariety. Let us define $\widehat{L}_{cl}\subset \widehat{M}_{cl}$ as the set of pairs $(m, (p_1,p_2,...,p_n, q_1,..., q_n))$, where $m\in M$ and $(p,q)=(p_1,p_2,...,p_n, q_1,..., q_n)$ is a formal Darboux coordinate system at $m$ such that the completion of $L$ at $m$ is given by the equations $p_1=...=p_n=0$.

Recall the Harish-Chandra pair $(\g_{cl},G_{cl,+})$ describing symplectic manifolds. There is an embedding of Harish-Chandra pairs $({\p}_{cl}, P_{cl,+})\hookrightarrow (\g_{cl},G_{cl,+})$, where we single out $\p_{cl}$ and $P_{cl,+}$ by the additional condition that vector fields and symplectomorphisms stabilize the standard Lagrangian subspace $L_{st}=\{p_1=...=p_n=0\}$. We have $dim_\C({\p}_{cl}/Lie(P_{cl,+}))=n$.
It follows from the definitions that $\widehat{L}_{cl}$ carries the natural $({\p}_{cl}, P_{cl,+})$-action. 

We have the following diagram of morphisms of Harish-Chandra pairs 
$$(\g_{max}, G_{max,+})\to (\g_{cl}, G_{cl, +})\leftarrow ({\p}_{cl}, P_{cl,+}).$$

 We denote by $({\p}_{max}, P_{max,+})$ the Harish-Chandra pair which is the fiber product of these morphisms. One can see that there is an embedding of Harish-Chandra pairs $(\p_{max}, P_{max, +})\subset (\g_{max}^{\prime\prime}, G_{max, +}^{\prime \prime})$.

 Then we can summarize the above constructions as the following statement.
 
 \begin{prp}\label{QWFS on thickening of Lagrangian}
 A $({\p}_{max}, P_{max,+})$-structure on $L$ is equivalent to the following data:
 
 a) A symplectic structure on the formal neighborhood (thickening) of $L$.
 
 b) A quantum wave function structure on the symplectic thickening from a).

 \end{prp}
 
 We can give a slightly different description of $\p_{cl}$. Notice that as the Lie algebra of Hamiltonian vector fields preserving $L_{st}$ the Lie algebra $\p_{cl}$ can identified with the Poisson ideal $I_{Lag}=\oplus_{1\le i\le n}\C[[p_1,...,p_n,q_1,..., q_n]]p_i\subset \C[[p_1,...,p_n, q_1,..., q_n]]$, since under the natural map from the space of (formal) functions to Hamiltonian vector fields the ideal $I_{Lag}$ goes to $\p_{cl}$.

 Let $P_{Lag}\subset Sp(2n, \C)$ be the subgroup of linear symplectic transformations preserving the Lagrangian subspace $L_{st}$ (i.e. it is the stabilizer of $L_{st}$ in $Sp(2n,\C))$. There is a natural short exact sequence
 $$1\to P_1\to P_{cl, +}\to P_{Lag}\to 1,$$
 where $P_1$ is a pronilpotent group.
 It follows that there is a surjective homomorphism $P_{max,+}\to P_{Lag}$ with the pronilpotent kernel. We have also the following short exact sequence of groups
 
$$1\to P_{Lag,1}\to P_{Lag}\to GL(n,\C)\to 1,$$
where $P_{Lag,1}$ is a nilpotent group, and $GL(n,\C)$ should be thought of as $Aut(L_{st})$, the group of linear automorphisms of the vector space $L_{st}$.
 
Next we would like to define a $2:1$ covering of the group $P_{max,+}$ and discuss its algebraic and geometric meaning. We start with a remark that the above short exact sequence implies that the following fact for the fundamental group: $\pi_1(P_{max,+})\simeq \Z$ for any $n\ge 1$. Next we notice that $\p_{max}$ is naturally embedded into the Lie algebra $\g_{max}^{\prime\prime}={1\over{\hbar}}W_q$. The Weyl algebra $W_q$ can be interpreted as the algebra of formal pseudo-differential
operators over $\C[[\hbar]]$ generated by $q_i, \hbar\partial_{q_i}, 1\le i\le n$. Then the $\C[[\hbar]]$-module $E:=\C[[q_1,...,q_n]][[\hbar]]$ is naturally a $W_q$-module.

\begin{prp}\label{action of Lie algebra}
For any $f\in \p_{max}$ we have:

a) $[f,\bullet]\in Der(W_q)$, where $[\bullet,\bullet]$ denote the commutator in $W_q$.

b) Multiplication by $f$ gives rise to an element of $End_{\C[[\hbar]]}(E)$.

\end{prp}
{\it Proof.} Part a) is obvious since $[{1\over{\hbar}}W_q, W_q]\subset W_q$. In order to prove b) consider a short exact sequence
$$0\to W_q\to \p_{max}\to {1\over{\hbar}}I_{Lag}\to 0.$$
Notice that ${1\over{\hbar}}W_q/W_q={1\over{\hbar}}\C[[p_1,...,p_n,q_1,...,q_n]][[\hbar]]\supset I_{Lag}$. Then it suffices to check that $\partial_{q_i}={1\over{\hbar}}(\hbar\partial_{q_i}), 1\le i\le n$ acts on $E$, which is obvious. 

In particular we see that although a priori multiplication by $f$  defines a map $E\to {1\over{\hbar}}E$, in fact it maps $E$ into $E$. This concludes the proof. $\blacksquare$

\begin{cor}\label{action of Lie algebra on pair}
There is a natural homomorphism from $\p_{max}$ to the Lie algebra of derivations of the pair $(W_q, E)$.

\end{cor}

There is an embedding of Lie algebras $Lie(Sp(2n, \C))\to \g_{max}^{\prime\prime}={1\over{\hbar}}W_q$ such that the image is the Lie subalgebra
generated by the elements ${1\over{\hbar}}q_iq_j, {1\over{\hbar}}(\hbar\partial_{q_i}\cdot \hbar\partial_{q_j}), 1\le i\ne j\le n, {1\over{\hbar}}(\hbar\partial_{q_i}+{\hbar\over{2}}), 1\le i\le n$. Notice that degree zero part of the graded Lie algebra ${1\over{\hbar}}W_q$ can be canonically identified with $Lie(Sp(2n, \C))\oplus \C\cdot 1$. This ensures that the above embedding is essentially unique. Then we see that we have an embedding $Lie(GL(n,\C))\to Lie(P_{max,+})$ such that the image is identified with the Lie algebra generated by $q_j\partial_{q_j}, 1\le i\ne j\le n, q_i\partial_{q_i}+1/2, 1\le i\le n$. Presence of the summand $1/2$ leads to the following conclusion: the group $P_{max,+}$ does not acts on the pair $(W_q,E)$, but this pair is acted by the  group $P_{max,+}^{(2)}$, which is $2:1$ covering of $P_{max,+}$.

Let us introduce the Harish-Chandra pairs $(\p_{cl}, P_{cl, +}^{(2)}), (\p_{max}, P_{max,+}^{(2)})$ where the upper index refers to the $2:1$ covering of the corresponding group. We have the natural morphism of Harish-Chandra pairs $(\p_{max}, P_{max,+}^{(2)})\to (\p_{cl}, P_{cl, +}^{(2)})$. In concrete terms an element of $P_{max,+}^{(2)}$ is an element $g\in P_{max,+}$ together with the choice of the square root $\sqrt{Jac(g_{|T_0L_{st}})}$, where
$Jac(g_{|T_0L_{st}})$ is the determinant of the Jacobi matrix of the restriction of $g$ to the tangent space at $0$ of the Lagrangian subspace $L_{st}$.

Notice that the  $(\p_{cl}, P_{cl, +}^{(2)})$-structure on $L$ is equivalent to the choice of Lagrangian structure on $L$ together with a choice of the square root $K_L^{1/2}$, where $K_L$ is the line bundle of the top degree forms (canonical line bundle). 

We can now summarize the above discussion in the following way.

\begin{prp}\label{bundles of algebras and modules}

Assume now that $M$ is endowed with a QWFS, and $L$ carries a $(\p_{cl}, P_{cl, +}^{(2)})$-structure.The lift of this structure to $(\p_{max}, P_{max,+}^{(2)})$-structure obtained via the above morphism of Harish-Chandra pairs endows the formal neighborhood $U_L$ of $L$ with an infinite-dimensional  bundle of quantized algebras $\OO_{U_L,\hbar}\to U_L$ endowed with a flat connection  as well as a bundle $\EE_L\to L$ of $\OO_{U_L,\hbar}$-modules, also endowed with a flat connection. 

\end{prp}

Notice that the bundle of algebras over $U_L$ gives rise to the bundle of germs of algebras over $L$. Therefore the above structure depends only on the germ of $L$ in an ambient symplectic manifold (maybe formal) as long as the latter is endowed with a QWFS.

\begin{defn}\label{definition of quantum wave functions}
Suppose that $M$ is endowed with a QWFS.
If $L$ satisfies the assumption of the Proposition \ref{bundles of algebras and modules} we will say that $L$ carries a  bundle of quantum wave functions. Same terminology will be applied to a germ of Lagrangian submanifold.
Quantum wave function supported on $L$ is an element of the vector space of flat global sections $\Gamma_{flat}(L, \EE_L)$ of the above bundle of modules.\footnote{We  also can (and will) speak about a local quantum wave function $\psi_U$ which is associated  with an open subset  $U\subset L$. }

\end{defn}

If we want to stress the support of a quantum wave function we will denote it by $\psi_L$.



\begin{rmk}\label{Borel-Weil and wave functions}
 Let $G$ be a complex Lie group and $G_+\subset G$ a complex Lie subgroup. Given a character $\chi: G_+\to \C^\ast$ one can construct a line bundle $\LL_\chi\to G/G_+$. The space of global sections $E_\chi:=\Gamma(G/G_+,\LL_\chi)$ is naturally a $G$-module. Explicitly it is given by the space of functions $f$ on $G$ such that $f(g\cdot g_+)=\chi(g_+)f(g), g\in G, g_+\in G_+$.
Suppose that we just have a Harish-Chandra pair $(\g,G_+)$ (this means that the  group $G$ corresponding to $\g$ may not exist). We still  can speak about the formal thickening of $G_+$ in the (non-existing) group $G$. It corresponds to the formal completion along the fibers of the normal bundle. The formal thickening is a partially formal manifold (cf. Remark \ref{dual to enveloping algebra}). Then the corresponding ``functions" are formal series in the variables $x_i\in \g/Lie(G_+)$ which are semi-invariant in the direction $G_+$ with respect to the character $\chi$. 


\end{rmk}

Now we can discuss the ``odd'' version of the above construction. Namely, we define $\p_{max}^{odd}\subset \p_{max}\subset {1\over{\hbar}}W_q$ as the intersection $\g_{max}^{odd}\cap \p_{max}$.  Similarly we have a subgroup $P_{max,+}^{odd}\subset P_{max,+}$.

The Harish-Chandra pair $(\p_{max}^{odd}, P_{max,+}^{odd})$ naturally acts on the $\C[[\hbar]]$-module $\C[[q_1,...,q_n]][[\hbar]]$. The action is the restriction of the natural action of the Harish-Chandra pair $(\p_{max}, P_{max}^{(2)})$. Hence the Harish-Chandra pair   $(\p_{max}^{odd}, P_{max,+}^{odd})$ acts on the quotient $\C[[q_1,...,q_n]][[\hbar]]/\hbar\C[[q_1,...,q_n]][[\hbar]]$. We can describe this action explicitly. The action of $\p_{max}^{odd}$ comes from the action of $Der(\C[[q_1,...,q_n]])$ on $\C[[q_1,...,q_n]]$ and induced by the composition of the Lie algebra homomorphisms
$$\p_{max}^{odd}\to \p_{cl}\to  Der(\C[[q_1,...,q_n]])\to End(\C[[q_1,...,q_n]]),$$
where the first two arrows are natural epimorphisms, while the last one is given by the map $\xi=\sum_{1\le i\le n}\xi_i\partial_{q_i}\mapsto \sum_{1\le i\le n}\xi_i\partial_{q_i}+{1\over{2}}\partial_{q_i}(\xi_i)=\xi+{1\over{2}}div(\xi)$.

Even more explicitly, notice that 
$$\g_{max}^{odd}=\prod_{k=-1,1,3,5,..., m=0,1,2,..., (k,m)\ne (-1,0)}\hbar^k Sym^m(\C^{2n}),$$ 
and we have the surjective map $\g_{max}^{odd}\to \prod_{k\ge 1}\hbar^k Sym^m(\C^{2n})$. It follows that $\p_{max}^{odd}$ maps surjectively onto ${1\over{\hbar}}I_{Lag}=\p_{cl}$. Thus $\p_{max}^{odd}$ acts on $\C[[q_1,...,q_n]]$ via its projection to $\p_{cl}$. This action is given explicitly on the monomial ${1\over{\hbar}}Sym(q_{i_1}...q_{i_m}(\hbar\partial_{q_{i_m}})):=\sum_{\sigma\in S_m}{1\over{m!}}{1\over{\hbar}}\sigma(q_{i_1}...q_{i_m}(\hbar\partial_{q_{i_m}})=\xi+{1\over{2}}div(\xi)$,
where $\xi=q_{i_1}...q_{i_m}\partial_{q_{i_m}}$.

From the point of view of Remark \ref{Borel-Weil and wave functions} this can be illustrated such as follows. Consider the Lie subalgebra of $Lie(P_{max,+})$ isomorphic to $\mathfrak{gl}(n)\oplus \C$, where the Lie algebra $\mathfrak{gl}(n)$ is spanned by ${1\over{2\hbar}}(q_i\hbar\partial_{q_j}+\hbar\partial_{q_i}q_j)_{1\le i,j,\le n}$ and $\C$ is spanned by $1={1\over{\hbar}}\cdot \hbar$. The algebra $\mathfrak{gl}(n)\oplus \C$ has a character $\chi$ given by $\chi(a,\lambda)={1\over{2}}tr(a)+\lambda$. The lift of this character to the Harish-Chandra pair $(\p_{max}^{odd}, P_{max}^{odd})$ does not contain the summand $\lambda$ since $1\notin \p_{max}^{odd}$.


Thus we obtain the following result.

\begin{prp}\label{half-forms} The lift of the natural Harish-Chandra module structure on $Spec(\C[[q_1,...,q_n]])$ from $(\p_{cl}, P_{cl,+}^{(2)})$ to $(\p_{max}^{odd}, P_{max,+}^{odd})$ identifies the former canonically with the module of formal half-forms on the formal completion $\widehat{\C}_0^n$.

\end{prp}

If $M$ is endowed with a QWFS one can generalize the above  considerations for an arbitrary smooth Lagrangian submanifold $M$ or for a germ of a smooth Lagrangian submanifold.

\subsection{Explicit formulas for quantum wave functions}\label{explicit formulas for wave functions}

Let us explain the notion of  quantum wave function in ``concrete'' terms. Suppose we are given a holomorphic Lagrangian submanifold $L$ together with a choice of $K_L^{1/2}$. Then
for any Darboux coordinate system $(p,q):=(p_1,...,p_n, q_1,...,q_n)$ we can choose a germ of holomorphic function $\mathcal F_0=\mathcal F_0(q)$ such that $L=\{(q,\partial \mathcal F_0/\partial q)\}$.  In a bit more invariant way it is a choice of a germ of holomorphic function $\mathcal F_0$ on $L$ such that $d\mathcal F_0=\lambda_{|L}$, where $\lambda=pdq$ is the Liouville $1$-form.

Then one can think informally that  a quantum wave function supported on $L$ can be written as $\psi_L=e^{\mathcal F_0/\hbar}\sum_{m\ge 0}\mu_m\hbar^m:=e^{\mathcal F_0/\hbar}\mu(\hbar)$, where $\mu_m, m\ge 0$  are holomorphic half-densities in coordinates $q_1,...,q_n$ (i.e. on $\C^n_{q_1,...,q_n}$). More precisely, under the change 
of Darboux coordinates $\mu_0$ transforms as a half-density on $L$ and each $\mu_m, m\ge 1$ is transformed as a half-density on $L$ plus a correction term depending only on the terms $\mu_i, i<n$.

\begin{rmk}\label{identification of wave function modules}
The choice of $\FF_0$ is not canonical, since we can add a constant: $\FF_0\mapsto \FF_0+c, c\in \C$. On the other hand the module of quantum wave functions supported on $L$ should not depend on this ambiguity. In order to achieve that we identify modules corresponding to $\FF_0$ and $\FF_0+c$ by formal multiplication by $e^{c\over \hbar}$:
$$e^{\mathcal F_0/ \hbar}\sum_{m\ge 0}\mu_m\hbar^m \stackrel{e^{c/\hbar}\cdot}{\mapsto} e^{(\mathcal F_0+c)/\hbar}\sum_{m\ge 0}\mu_m\hbar^m.$$

\end{rmk}



If $\mu_0$ does not vanish (i.e. if $\psi_L$ is a generator of the corresponding holonomic $DQ$-module) then we will write the quantum wave function as
$$\psi_L=\exp\bigl(\sum_{g\ge 0}\hbar^{g-1} \mathcal F_g\bigr)\cdot (dq_1\wedge dq_2\wedge\dots \wedge dq_n)^{1/2} $$
where $\mathcal F_0,\mathcal F_1,\mathcal F_2,\dots$ are functions of $\bf q$. We may assume that $\mathcal F_1$ is defined modulo $2\pi i \Z$
as long as we have fixed a section of $K_L^{1/2}$ which we think of as $exp(\mathcal F_1)$.

For a {\it given} choice of $\mathcal F_0$ and of $K_L^{1/2}$, the sheaf of quantum wave functions (given locally by an arbitrary sequence functions $G_0,G_1,\dots$ in  coordinates $\bf q$) is a sheaf of modules over $\mathcal O_{\C^{2n},\hbar}$. This is a  {\it holonomic} $\mathcal O_{\C^{2n},\hbar}$-module.

\begin{rmk}\label{three levels of abstraction}
One can describe  quantum wave functions at the three levels which correspond to the three levels of the deformation quantization from  Section \ref{star-product}:

 1) (Most concrete, assuming the odd version of the QWFS) a choice of a wave function $\psi$ with $G_0$ non-vanishing (i.e. given by the exponential form parameterized by $(\mathcal F_0,\mathcal F_1,\dots)$) defines a left  ideal $I=I_\psi$, as the annihilator of $\psi$ by the $\mathcal O_{M,\hbar}$-action.

2) (Intermediate) if we choose a sheaf of algebras $\mathcal O_{M,\hbar}$ over $\C[[\hbar]]$, then this object can be {\it non-canonically} presented as $\mathcal O_{M,\hbar}/I$ where $I\subset \mathcal O_{M,\hbar}$ is a sheaf of left ideals, s.t. $I\,mod\, \hbar$ defines $L$.

 3) (Most abstract) for arbitrary complex symplectic manifold $(M,\omega^{2,0})$ and a holomorphic Lagrangian submanifold $L\subset M$ {\it together} with a choice of $K_L^{1/2}$, there is certain {\it canonical} holonomic object 
 in the canonical sheaf of categories $\CC_M^{can}$ (cf.  [KasSch1]).

Notice that left ideal $I$ {\it does not} define $\psi$ uniquely: in terms of the sequence of functions $\mathcal F_0,\mathcal F_1,\dots$ one has freedom
$$(\mathcal F_0,\mathcal F_1,\dots )\rightsquigarrow (\mathcal F_0+c_0,\mathcal F_1+c_1,\dots ),\qquad c_i\in \C\,\, \forall i=0,1,\dots$$

\end{rmk}

Next we are going to make more precise  the informal discussion of the local pairing between quantum wave functions from Section \ref{left and right modules}.
In the de Rham language it corresponds to the composition of morphisms of sheaves (and spaces of their global sections):
$$\underline{Ext}^n(\EE_L, \OO_{M,\hbar})\otimes_{\OO_{M,\hbar}} \underline{Hom}(\OO_{M,\hbar}, \EE_{L^{\prime}})\to \underline{Ext}^n(\EE_L,\EE_{L^\prime}).$$

Let $x\in L\cap L^\prime$ be a transversal intersection point. Fix a left quantum wave function $\psi_L^{left}$ corresponding to $L$ and a right quantum wave function $\psi_{L^\prime}^{right}$ corresponding to $L^\prime$ \footnote{Right quantum wave functions can be defined  replacing $\OO_{M,\hbar}$ by $\OO_{M,\hbar}^{op}$. In the case when  $\psi_L^{left}$ and $\psi_{L^\prime}^{right}$ give nowhere vanishing half-densities in the classical limit, then they correspond to what we denoted in Section \ref{left and right modules} by $1_{L^{left}}$ and $1_{L\prime}^{right}$. }. They define the corresponding local sections of the above tensor factors.
The stalk of the sheaf $\underline{Ext}^n(\EE_L,\EE_{L^\prime})$ at the point $x$ is canonically (up to a sign) isomorphic to $ \C[[\hbar]]$. {\it After a choice of this isomorphisms we have  a well-defined  local pairing $\langle\psi_L^{left}, \psi_{L^\prime}^{right}\rangle_x\in \C[[\hbar]]$}.

In the case when a QWFS is symmetric one can identify right quantum wave functions with the left ones via the change of the quantization parameter $\hbar\mapsto -\hbar$. {\it Unless we say otherwise will assume from now on that our QWFS is symmetric.}

In coordinates the pairing of quantum wave functions can be described such as follows. Consider two Lagrangian submanifolds $L_1=\{p_1=...=p_n=0\}$ and $L_2=\{q_1=...=q_n=0\}$ in the formal completion at $0\in \C^{2n}$. The corresponding canonical modules over the formal Weyl algebra $W_q$ are denoted by $E_1$ (right  module) and $E_2$ (left module) respectively. Then $$E_1\simeq 1\cdot \C[[q_1...,q_n, \hbar\partial_{q_1},...,\hbar\partial_{q_n}]][[\hbar]], E_2=\C[[q_1...,q_n, \hbar\partial_{q_1},...,\hbar\partial_{q_n}]][[\hbar]]\cdot \delta_q,$$ 
where $\delta_q$ denotes the delta-function in variables $q_1,.., q_n$. There is a natural pairing $\langle \bullet,\bullet\rangle: E_1\otimes_{W_q} E_2\to \C[[\hbar]]$ which can be identified with a unique pairing $I_{E_1}\backslash W_q\otimes_{W_q}W_q/I_{E_2}\to \C[[\hbar]]$ normalized by the condition $\langle 1,\delta_q\rangle=1\in \C[[\hbar]]$. Here $I_{E_1}$ and $I_{E_2}$ denote left and right ideals corresponding to the cyclic $W_q$-modules $E_1$ and $E_2$ respectively, and the above-mentioned generators $1,\delta_q$ can be identified with the images of $1\in W_q$ under the natural quotient maps. Since we assumed that our QWFS was symmetric we may convert left $W_q$-module $E_1$ into the right $W_q$ module, so $E_1\otimes_{W_q}E_2$ is well-defined and can be naturally identified with the above tensor product of quotient modules.
Explicitly the pairing can be written down such as follows
$$\langle \sum_{I,k}a_{I,k}q^I\hbar^k, \sum_{J,k}b_{J,k}p^J\hbar^k\rangle=\sum_{I}a_{I,k_1}b_{I,k_2}|I|! \hbar^{k_1+k_2}.$$

Here $I,J$ are multi-indices and summation runs over $k\ge 0, |I|\ge 0, |J|\ge 0$.

Let us now illustrate the relationship between the pairing of quantum wave functions and local exponential integrals assuming  that $M=\C^{2n}:=\C^{2n}_{{\bf q}, {\bf p}}$, and  both Lagrangians $L_0, L_1$ have   one-to-one projections to the coordinate space $\C^n_{\bf q}$.

  Notice that a  choice of   quantum wave function supported on $L\subset \C^{2n}$ is the same as a choice of a  holomorphic function $\mathcal F_0'$ on $L$\footnote{Such a choice is local if $L$ does not have a one-to-one projection to $ \C^n_{\bf q}$} such that 
   $$d{\mathcal F}_0'=\left(-\sum_i p_i dq_i\right)_{|L}.$$
  
If $L_0,L_1\subset \C^{2n}$ are two Lagrangian submanifolds as above intersecting {\it transversally} at some point $({\bf q}_\alpha,{\bf p}_\alpha)\in \C^{2n}$, we define  $\langle \psi_{L_0}|$ and $|\psi_{L_1}\rangle $ which are resp. right and left quantum wave functions supported on $L_0$, resp. $ L_1$\footnote{Here and below we follow the bra and ket notation of physicists which are very convenient when one wants to distinguish between left and right wave functions.}. We defined above the {\it de Rham local pairing} (determined up to a sign $\pm$) 
$$\pm\langle \psi_{L_0},\psi_{L_1}\rangle_{({\bf q}_\alpha,{\bf p_\alpha})}:=\pm\langle \psi_{L_0}|\psi_{L_1}\rangle_{({\bf q}_\alpha,{\bf p_\alpha})}\in \exp\bigl({1\over \hbar}(\mathcal F_0+\mathcal F_0')({\bf q}_\alpha,{\bf p}_\alpha)\bigr)\C[[\hbar]].$$
Since  both $L_0$ and $L_1$ project  one-to-one to $\bf q$-coordinates near the intersection point $({\bf q}_\alpha,{\bf p}_\alpha)\in \C^{2n}$ we have
$$ L_0=\text{ graph of } d\mathcal F_0,\qquad L_1=\text{ graph of } -d\mathcal F_0' $$
 
 Then both  $|\psi\rangle$ and $\langle \psi'|$ can be viewed as half-densities:
 $$\langle\psi_{L_0}|=\exp\bigl(\sum_{g\ge 0}\hbar^{g-1} \mathcal F_g\bigr)\cdot (dq_1\wedge dq_2\wedge\dots \wedge dq_n)^{1/2}, $$
 $$|\psi_{L_1}\rangle=\exp\bigl(\sum_{g\ge 0}\hbar^{g-1} \mathcal F_g'\bigr)\cdot (dq_1\wedge dq_2\wedge\dots \wedge dq_n)^{1/2},$$
 and the function $\mathcal F_0+\mathcal F_0'$ in $\bf q$ has Morse critical point at ${\bf q}_\alpha\in \C^n$. Then we define the {\it Betti local pairing} at $({\bf q}_\alpha,{\bf p}_\alpha)\in L_0\cap L_1\subset \C^{2n}$  as a formal integral over the local Lefschetz thimble, i.e.
 $$ \langle \psi_{L_0}|\psi_{L_1}\rangle_{({\bf q}_\alpha,{\bf p}_\alpha)}:=(2\pi \hbar)^{-{n\over 2}}\intop_{\overset{\text{Lefschetz thimble}}{\text{\tiny near } {\tiny \bf q}_\alpha} } \exp\bigl(\sum_{g\ge 0}\hbar^{g-1} (\mathcal F_g+\mathcal F_g')\bigr)\cdot (dq_1\wedge dq_2\wedge\dots \wedge dq_n)$$
We know that this pairing is covariant with respect to affine symplectic transformations, hence it can be defined {\it without} the assumption that the projection of both $L_0,L_1$ to $\C^n_{\bf q}$ near  $({\bf q}_\alpha,{\bf p}_\alpha)$ is locally one-to-one.
{\it Betti and de Rham local pairings agree with each other}.

The prefactor $(2\pi \hbar)^{-n/2}$ ensures that the paring takes value in $\C[[\hbar]]$. In this way we obtain the formula which agrees with the expected formula in Section \ref{path integrals}.
More precisely in the case  $H=H({\bf q},{\bf p},t)\equiv 0$, i.e. when
 $$S(\phi)=\int_0^1\sum_i p_i(t){d q_i(t)\over dt} $$
 we  {\it define} in the above notation the formal path integral over the local Lefschetz thimble to be equal to the de Rham or Betti  local pairing
  $$\int e^{S(\phi)\over \hbar}\mathcal D\phi=\langle \psi_{L_0}|\psi_{L_1}\rangle $$
(here we omit the Lefschetz thimbles from the notation).

\begin{rmk} \label{non-smooth Lagrangians}
 a) In [KoSo11] we constructed a quantum wave function only as the sheaf of ideals, without fixing constants, which are $(\omega_{g,0})_{g=0,1,\dots}$ in the language of topological recursion.

b) Local parings of quantum wave functions are expected to give resurgent  perturbative series. This means resurgence of the perturbative expansions of the  corresponding  exponential integrals. Therefore not every choice of constants is good. Here are a couple of proposals based on experiments:

 b1) First, if we ignore constants $(c_0,c_1,\dots)$ as above (or the generating series $\sum_{i\ge 0} c_i \hbar^i\in \C[[\hbar]]$). In other words we consider the ideal $I$ instead of the quantum wave function, this ideal should come from an analytic in $\hbar,  |\hbar|\ll 1$ family of {\it algebraic} differential/difference equations in $\C^n$ or $(\C^\ast)^n$. 
 
 b2) Second, the normalization of constants $(c_0,c_1,\dots)$ should come from a certain  normalization of the quantum wave function {\it at infinity}. We will illustrate this point in the example of the WKB expansion.

\end{rmk}

\begin{rmk}\label{non-transversal intersections}
 Assume that $L_0$ and $L_1$ intersect non-transversally at some point $({\bf q}_\alpha,{\bf p}_\alpha)\in \C^{2n}$ but the intersection point is {\bf isolated}. Let  $\mu=\mu_\alpha\ge 2$ be the Milnor number of the corresponding singularity (e.g. $\mu=k-1$ if $L_1=\text{graph} (d x^k), L_0=\text{graph} (d 0), \,k\ge 3$). Then the Lefschetz thimble  should be replaced by a vanishing cycle (there are $\mu$ linearly independent  cycles). The resulting integral has asymptotic expansion in {\bf rational} powers of $\hbar$ and integer non-negative powers of $\log \hbar$.
This can be generalized to a {\bf non-isolated} intersection, where one should replace vanishing cycles by the global sections of an appropriate constructible sheaf of local vanishing cycles on $L_0\cap L_1$.

One can hope that in the case of {\bf singular} Lagrangian varieties $L_0,L_1$  the corresponding quantum wave functions still exist, and they are formal exponents of more general (possibly multivalued) expressions in $(q_1,\dots,q_n)$ and $\hbar$.

\end{rmk}

\begin{rmk}\label{pairing of wave functions and Fukaya}
Here is yet another categorical interpretation of the pairing of quantum wave functions. Recall that Holomorphic Floer Theory predicts that to a singular complex Lagrangian submanifold $L$ one can construct two  objects $E_L^{Betti,loc}\in \FF_{L,loc}$ and $E_L^{DR,loc}\in Hol_{L,loc}$. 

Let us take $L=L_0\cup L_1$, where $L_i, i=0,1$ are two transversally intersecting smooth complex Lagrangian submanifolds. Choose a point $x\in L_0\cap L_1$. Then we have  embeddings $\FF_{L_0,loc,x}\to \FF_{L_0\cup L_1,loc,x}\leftarrow \FF_{L_1, loc,x}$, where the notation means that we consider as an ambient symplectic manifold a small neighborhood of $x$ in $M$. The local pairing $\langle\psi_{L_0},\psi_{L_1}\rangle_x$ takes values in $Ext^n(E_{L_0}^{Betti,loc}, E_{L_1}^{Betti,loc})$ which is isomorphic after changing of scalars to  $Ext^n(E_{L_0}^{DR,loc}, E_{L_1}^{DR,loc})$. Outside of Stokes rays local categories can be identified with the global ones, so we have the global pairing of quantum wave functions.

These considerations can be applied to the non-transversal case as well. In particular, taking $L=L_0=L_1$ we obtain that $\langle \psi_L|\psi_L\rangle$ corresponds to the integration of top degree volume form $\psi_L\cdot \psi_L\in H^n(L)$ over the fundamental class. In other words our pairing of quantum wave functions contains the theory of periods of volume forms.

\end{rmk}

\subsection{Example:WKB asymptotics}\label{example of wkb}

Consider the following family of differential equations depending analytically on $\hbar$:
$$\Bigl[-\left(\hbar {d\over dx}\right)^2+ (x^4+1)\Bigr] \psi(x)=0.$$
We normalize the WKB solution by the condition at $+\infty$:
$$\lim_{x\to +\infty} \psi(x)\cdot \exp\Bigl({x^3\over 3 \hbar}\Bigr) =1\text{ as a series in } \hbar .$$
In other words, 
$$\mathcal F_0(x)=-x^3/3\text{ for } x\to +\infty\\
\lim_{x\to +\infty} {\mathcal F}_g(x)=0\text{ for }g\ge1.$$
This gives a multi-valued  analytic wave function on the punctured elliptic curve
$$L_0:=\{(q,p)\in \C^2|p^2=q^4+1, p\ne 0\},\,\,q=x$$
with  branches which differ from each other by the value of quantum periods:
$$H_1(L_0,\Z)\to \C[[\hbar]].$$

Quantum periods  are equal modulo $\hbar$ to the classical periods of the 1-form $p \,dq_{|L_0}$. 
The quantum periods for small loops about ramification points are equal to $\hbar/4$.

If we {\it fix} $x=x_0\in \C$ then the formal expression in $\hbar$ (which is the ``value" of WKB solution at $x=x_0$):
$$ \exp\Bigl(\sum_{g\ge 0} \mathcal F_g(x_0) \,\hbar^{g-1}\Bigr)$$
can be identified with $\langle \psi'|\psi\rangle$ where $\psi'$ corresponds to  $``\delta(x-x_0)"$, 
i.e. it is an analytic right quantum wave function associated with the Lagrangian submanifold
$$ L_1:=\{(q,p)\in \C^2|q=x_0\}$$
More precisely, it is the quantum wave function obtained 
by applying to the most basic function (constant) $$\psi(x):=1=\exp(0),\quad \mathcal F_0=\mathcal F_1=\dots=0$$ 
an appropriate element of $SL(2,\C)\ltimes\C^2$ which moves the line $q=0$ to the line $q=x_0$.

\subsection{Transport of quantum wave functions}\label{transport of wave functions}

So far we considered the case corresponding to the action functional $\int pdq$. We would like to have a generalization to the case when the action functional is equal to $\int pdq +Hdt$, where $H=H(p,q,t)$ is a time dependent Hamiltonian.   If we have a time-independent Hamiltonian $H=H(p,q)$ the Hamiltonian formalism of the quantum mechanics says that the path integral for the path $[0,t]\to M$ is given by the expression $\langle \psi_0|e^{-t\widehat{H}}|\psi_1\rangle$, where 
$\widehat{H}$ is the corresponding quantum Hamiltonian. This agrees with  the  ``topological" case $H=0$  discussed previously, which corresponds to the pairing $\langle \psi_0|\psi_1\rangle$ in the bra-ket notation. Mathematically this leads us to the notion of the transport of quantum wave functions which we are going to discuss below. In this subsection we will be using the terminology and results from Section \ref{Harish-Chandra pairs}, so the reader should consult it if necessary.

Suppose that a complex symplectic manifold $(M, \omega^{2,0})$ is endowed with a QWFS. Recall that it gives the following
structures on $M$:

1) A sheaf of non-commutative $\C[[\hbar]]$-algebras $\OO_{M,\hbar}$ which quantizes the
sheaf of holomorphic functions $\OO_M$. 

2) To a complex Lagrangian submanifold $L$ endowed with a choice of $K_L^{1/2}$
a canonical sheaf of $\OO_{M,\hbar}$-modules $\EE_L$. By definition elements of $E_L:=\Gamma(L, \EE_L)$  are 
quantum wave functions supported on $L$. We use the notation  $\psi_L$  for
such a quantum wave function.

Given a smooth family of quantum Hamiltonians 
$$(\widehat{H}_t)_{t\in [0,1]}\in \Gamma(M, \OO_{M,\hbar})\widehat{\otimes} C^\infty([0,1])$$
we have the corresponding smooth  family of classical Hamiltonians obtained from the
quantum family by taking their images modulo $\hbar$:

$$(H_t)_{t\in [0,1]}\in \Gamma(M, \OO_{M})\widehat{\otimes} C^\infty([0,1]).$$

The classical family gives raise to a  vector field  $v={d\over {dt}}-\{H_t, \bullet\}$ on $M\times [0,1]$.

Let $L\subset M$ be a complex Lagrangian submanifold.  Let us denote by $U\subset L\times [0,1]$ the subset consisting of pairs $(x,t)$ such that the orbit $t^\prime\mapsto \phi_{t^\prime}(x)$ of the vector field $v$ is well-defined for $t^\prime\in [0,t]$. The subset $U$ is open and contains a neighborhood of $L\times \{0\}$.\footnote{Notice that $U$ is a manifold of ``mixed" type: it is complex analytic along $L$ and $C^\infty$ along $t$. Notice also that if a point $m\in M$ does not flow to infinity under the flow of $v$ then all points of the fiber $\widehat{M}_{(\g_{max}, G_{max,+})} \to M$ do not flow to infinity. Indeed $G_{max,+}$ is a projective limit of finite-dimensional algebraic groups, hence the above claim follows from the corresponding one in the finite-dimensional case.
}
Then we have a smooth family of Lagrangian
submanifolds $L^t=\{\phi_t(x)| (x,t)\in U\}\subset M$ depending on $t\in [0,1]$. Alternatively $L^t=\phi(L\times \{t\})\cap U$, where $\phi(x,t):=\phi_t(x), x\in L, t\in [0,1]$.

Suppose we are given a quantum wave function $\psi_L$ supported on $L=L^{t=0}$. Our goal is to define a family of quantum wave functions $\psi_{L^t}$ supported on $L^t$ obtained by a ``parallel transport" of $\psi_L$ defined by the family $(\widehat{H}_t)_{t\in [0,1]}$. Informally we would like to solve the Schr\"odinger equation
$$\hbar {d\over {dt}}\psi_{L^t}=\widehat{H}_t\psi_{L^t}.$$
The change $\widehat{H}_t\mapsto \widehat{H}_t+{c(t)}$ with $c(t)\in C^{\infty}([0,1])\otimes \C$ leads to the multiplication of the quantum wave function by $e^{\int_0^tc(s)ds\over \hbar}\in e^{\C\over \hbar}$ which does not change the quantum wave function  as an element of the corresponding module (see Remark \ref{identification of wave function modules}).

In the remaining part of this subsection we will give a formal treatment of the parallel transport.  First let us explain what do we mean by a parallel transport of quantum wave functions.

Let $WF$ denote an infinite-dimensional vector bundle on $U$ such that the fiber $WF_{(x,t)}$ is the vector space of quantum wave functions on the formal completion of $L^t$ at the point $\phi_t(x)\in L^t$. As any bundle of infinite jets the bundle $WF$ carries a natural flat connection along $L$. We would like to extend it in the direction of $t$ in such a way that the total connection $\nabla$ will be flat. The quantum wave function $\psi_L$ gives a flat section of $WF_{|L\times \{0\}}$. Then there will be a unique flat section of $WF$ extending $\psi_L$. We will define $\psi_{L^t}$ as its restriction to $\phi(L\times \{t\})\cap U=L^t$  and call it {\it the parallel transport} of $\psi_L$.

In order to define the flat connection  we will use the following two general constructions. 

{\it Construction 1:  Lie algebra of symmetries of a $(\g,G_+)$-structure.}

Let $(\g, G_+)$ be a Harish-Chandra pair. Suppose that a smooth proalgebraic (or complex, $C^\infty$, etc.) variety $X$ is endowed with a $(\g, G_+)$-structure, so we have a $(\g, G_+)$-torsor $\widehat{X}\to X$. 
 Consider the  Lie algebra $Vect_{(\g, G_+)}(\widehat{X})$ of global  vector fields on $\widehat{X}$ commuting with the action of $(\g, G_+)$, so the Lie algebra $Vect_{(\g, G_+)}(\widehat{X})$  preserves our $(\g, G_+)$-torsor in the natural sense. This Lie algebra admits a different description.
 
Since the Lie algebra $\g$ is a  $(\g, G_+)$-module with respect to the adjoint action,  it gives rise to a bundle of infinite-dimensional Lie algebras $\widehat{X}_{ad}\to X$ endowed with a  flat connection.  The corresponding Lie algebra of  flat sections is naturally isomorphic to $Vect_{(\g, G_+)}(\widehat{X})$. For example, if $\g$ is the Lie algebra of formal vector fields on a coordinate vector space then $Vect_{(\g, G_+)}(\widehat{X})$ coincides with the Lie algebra  $Vect(X)$ of the vector fields on the manifold $X$. Yet another description of $Vect_{(\g, G_+)}(\widehat{X})$ is given in the following Proposition.

\begin{prp}\label{adjoint torsor and vector fields}
The set $Maps(\widehat{X}, \widehat{X}_{ad})$ of morpshisms of proalgebraic varieties endowed with  $(\g, G_+)$-action is naturally identified with the set $Vect_{ (\g, G_+)}(\widehat{X})$.
\end{prp}
{\it Proof.} The vector bundle $\widehat{X}_{ad}\to X$ carries a formal flat connection. The corresponding vector space of flat sections is naturally identified with the space of vector fields $Vect_{ (\g, G_+)}(\widehat{X})$. The result follows. $\blacksquare$

{\it Construction 2: Flat connection associated with a family of vector fields commuting with a $(\g,G_+)$-structure.}

Basic idea can be explained in the following elementary example. Let $G$ be an algebraic  group, $\pi: E\to B$  be a $G$-torsor and $\rho: G\to Aut(V)$ be a representation of $G$. We denote by ${\mathcal V}$ the associated vector bundle on $B$. Assume that we are given a family of vector fields $\xi_t, t\in [0,1]$ on $E$ commuting with the action of $G$. It induces  a family of vector fields $\eta_t, t\in [0,1]$ on $B$, so we have the corresponding non-autonomous differential equation on $B$:
$${dx(t)\over dt}=\eta_t(x(t)).$$

Consider a  $G$-torsor over  $B\times [0,1]$ with the fiber of $(b,t)$ equal to $\pi^{-1}(x(t))$, where $x(t)$ is the trajectory of $\eta_t$ such that $x(0)=b$. Then the family of vector fields $\xi_t$ gives rise to a  $G$-connection on this torsor along the second factor $[0,1]$. It induces the corresponding  connection along $t$ on the induced vector bundle on $B\times [0,1]$.

There is a version of the above construction. Suppose that $(\g,G_+)$ is a Harish-Chandra pair acting on a $G_+$-torsor $\pi: E\to B$ and $V$ be a $(\g,G_+)$-Harish-Chandra module. Given a family of vector fields $\xi_t$ as above which commute with the $(\g,G_+)$-action, we see that the corresponding vector fields $\eta_t$ are tangent to the leaves of the formal foliation  $\mathcal{F}$ of $B$ defined by $\g$. In this way we obtain a flat connection on the vector bundle on $B\times [0,1]$ along the foliation $\FF^\prime$ spanned by $\mathcal{F}$ and $\partial_t$. Summarizing,  on each leaf of $\FF^\prime$ (or on an open subset of the leaf) we obtain a vector bundle endowed with a flat connection.

\

Now we can formulate the main result.

\begin{prp-constr}\label{parallel transport} Suppose that $M$ is endowed with QWFS.
Let $\psi_L$ be a quantum wave function supported on a Lagrangian submanifold $L$. Then for each $t\in [0,1]$ there is a unique quantum wave function $\psi_{L^t}$ obtained from $\psi_L$ by the parallel transport in the sense explained above by utilizing the flat connection constructed in the proof below.
\end{prp-constr}

{\it Proof.} We apply the above two constructions in the following situation. We take $B$ to be the space of formal Lagrangian germs in $M$. In other words a point of $B$ is a pair consisting of a point $m\in M$ and a formal germ of a  Lagrangian submanifold at $m$. We have a natural morphism $E:=\widehat{M}_{(\g_{max}, G_{max, +})}\to B$ which is a $G_{max,+}$-torsor. The pro-algebraic manifold $B$ is endowed with a formal foliation coming from the $(\p_{max}, P_{max,+})$-action on $E$ (see Section \ref{quantum wave functions}). Global leaves of this formal foliation correspond to (non-formal) complex Lagrangian submanifolds  $L\subset M$.

By  Proposition \ref{adjoint torsor and vector fields} the set  $Vect_{(\g_{max}, G_{max,+})}(\widehat{M}_{(\g_{max}, G_{max, +})})$  contains $\OO_\hbar(M)=\Gamma(M,\OO_{M,\hbar})$, i.e. the set of quantum Hamiltonians. Indeed the  $(\g_{max}, G_{max,+})$-module structure on $\g_{max}$ associated with the adjoint action comes from the one for the Weyl algebra $W_q$.  
It follows that the set of flat  sections of $\widehat{M}_{ad}$ can be identified with the set of elements of the form ${1\over{\hbar}}\widehat{H}$ modulo ${\C\over{\hbar}}$, where  $\widehat{H}\in \OO_\hbar(M)$.

Then we have a family of vector fields $\xi_t, t\in [0,1]$ on $E$ induced by the family of quantum Hamiltonians $(\widehat{H}_t/\hbar)_{t\in [0,1]}$. We can understand this family as a single vector field $\xi_{tot}$ on $E\times [0,1]$ where $\xi_{tot}(e,t)=\xi_t(e)\boxplus \partial_t$. Let us denote by $E_U$ the $P_{max,+}$-principal bundle over $U$ with the fiber over $(x,t)$ given by the set of points $\widehat{m}\in \widehat{M}_{(\g_{max},G_{max,+})}$ such that the corresponding point $\widehat{m}_{cl}\in  \widehat{M}_{(\g_{cl},G_{cl,+})}$ gives the formal symplectic coordinate system on $M$ at the point $\phi_t(x)$ such that the formal germ of $L^t$ at $\phi_t(x)$ is given by the equation $p_1=...=p_n=0$.
We have a tautological  embedding $E_U\hookrightarrow E\times [0,1]$ such that $\xi_{tot}$ is tangent to the image. Therefore we have a free transitive action of the Harish-Chandra pair $(\p_{max}\oplus \R\cdot \xi_{tot}, P_{max,+})$ on $E_U$.

A choice of $K_L^{1/2}$ gives the $2:1$-cover of the restriction of $E_U$ to $L\times \{0\}\subset U$. It can be naturally extended to the $2:1$ cover $E_U^{(2)}\to E_U$. We have a free transitive action of the Harish-Chandra pair $(\p_{max}\oplus \R\cdot \xi_{tot}, P_{max,+}^{(2)})$ on $E_U^{(2)}$.

Consider the  $(\p_{max}\oplus \R\cdot \xi_{tot}, P_{max,+}^{(2)})$-Harish-Chandra module ${V}$  modeled  in the symplectic coordinate $(p_1,...,p_n, q_1,...,q_n)$ on $V=\C[[q_1,...,q_n]][[\hbar]]$ (see Section \ref{quantum wave functions}) where $\xi_{tot}$  acts trivially. 
Then we obtain a flat connection on the vector bundle $WF\to U$ associated with the principal $P_{max,+}^{(2)}$-bundle $E_U^{(2)}$ and the above module $V$ consider a $P_{max,+}^{(2)}$-module.
$\blacksquare$

\begin{rmk}\label{symmetric parallel transport}
If our QWFS is symmetric then the parallel transport can be defined similarly to the above for the ``odd" versions of  Harish-Chandra pairs, in particular for $(\g_{max}^{odd}, G_{max,+}^{odd})$ and  $(\p_{max}^{odd}, P_{max,+}^{odd})$.

\end{rmk}

One can restate the parallel transport of quantum wave functions slightly differently. Namely, the above-constructed flat connection is uniquely defined by the following Schr\"odinger equation which should hold for any Lagrangian submanifold $L^\prime$ transversal to $L^t$ at $x_t\in L^t\cap L^\prime,  t\in [0,1]$:
$$\hbar{d\over{dt}}(\langle\psi_{L^t},\psi_{L^\prime}\rangle_{x_t})=\langle(\widehat{H}_t-H_t(x_t)\cdot 1)(\psi_{L^t}),\psi_{L^\prime}\rangle_{x_t}.$$

\begin{rmk}\label{coisotropic transport}
One can generalize the above considerations to the case of $DQ$-modules with {\bf coisotropic} support. 

\end{rmk}

\section{Quantum wave functions and resurgence}\label{resurgence}

\subsection{Wheels of Lagrangians and resurgent series}\label{wheels of Lagrangians}
In this subsection we will deal with categorical considerations which we will combine later with the  theory of quantum wave functions.

Let $M$ be a smooth complex  affine variety of  dimension $2n$ endowed with an algebraic symplectic form $\omega^{2,0}$.
 
We will formulate below a conjecture about resurgence of certain series associated to the following data:

1) $L_0, L_1,..., L_k=L_0$ be a cyclically ordered set of complex Lagrangian submanifolds of  $(M,\omega^{2,0})$ together with a choice of square roots $K_{L_i}^{1/2}$ of the canonical line bundles $K_{L_i}$. For each $L_i$ let us choose a point $x_{i,i+1}\in L_i\cap L_{i+1}$ assuming that the intersection is transversal at $x_{i,i+1}$, and denote by ${\bf x}$ the cyclically ordered collection $(x_{i,i+1})_{i\in \Z/k}$.

2) $\EE:=(E_{i}:=E_{L_i})_{i\in \Z/k}$  a collection of modules of finite type over the $\C[[\hbar]]$-algebra $\OO_\hbar(M)$ such that for each corresponding sheaf $\EE_i$ the sheaf $\EE_{i}/\hbar \EE_i$ is a pure coherent sheaf set-theoretically supported on $L_i$. 
We will write in this case $L_i=Supp(\EE_{i})$. We will assume  that generically all $\EE_{i}/\hbar \EE_i$ have rank $1$.

3) A collection of elements $\mu=(\mu_{i,i+1})_{i\in \Z/k}$ such that $\mu_{i,i+1}\in Ext^n(E_{i}, E_{i+1})$. 

4) A collection of homotopy classes of paths $\gamma=
(\gamma_i)_{i\in \Z/k}$ on $L_i$ such that $\gamma_i$ joins the intersection points $x_{i-1,i}$ and $x_{i,i+1}$ for each $i\in \Z/k$.

Recall again that with a complex symplectic manifold $(M,
\omega^{2,0})$ one can associate canonically a sheaf of $\C[[\hbar]]$-linear categories $\CC_M^{can}$ of (left) $DQ$-modules.. 
With a complex Lagrangian submanifold $L\subset M$ endowed with a choice of $K_L^{1/2}$ one can associate 
a sheaf of categories of holonomic $DQ$-modules supported on $L$.
Assuming that $M$  is endowed with a QWFS the above sheaf of categories $\CC_M^{can}$ is canonically equivalent to the sheaf of categories of $\OO_{M,\hbar}$-modules over some sheaf of quantized algebras $\OO_{M,\hbar}$. Furthermore there is 
a sheaf $\EE_{L}^{can}$ 
of holonomic $DQ$-modules over the sheaf of $\C[[\hbar]]$-algebras ${\OO}_{M,\hbar}$ such that $\EE_L^{can}$  is supported on $L$ and 
$rk(\EE_{L}^{can}/\hbar \EE_L^{can})=1$.

Then {\it locally in analytic topology} on each $L_i, i\in \Z/k$ there is an isomorphism of sheaves $\EE_{i}\simeq \EE_{L_i}^{can}$.  Such  isomorphisms form a sheaf of  torsors $Iso(\EE_i, \EE_{L_i}^{can})$ over $\C[[\hbar]]^{\times}$ which is isomorphic to the sheaf of torsors of rank $1$ locally constant sheaves (local systems) over $\C[[\hbar]]$ on $L_i$. Probably the best way to see this is  to use the conjectural generalized local RH-correspondence over $\C[[\hbar]]$, which is a theorem in this particular case. Namely, each sheaf $\EE_{L_i}^{can}$ corresponds under the local RH-corrsepondence to the rank $1$ local system on $L_i$.
Such sheaves of torsors are classified by $H^1(L_i, \C[[\hbar]]^\times)\simeq Rep(\pi_1(L_i), \C[[\hbar]]^\times)$. The corresponding local systems are determined by their monodromies.

In concrete terms the torsor is described in the following way. For any $i\in \Z/k$ and an element $\delta_i^{cl}\in H_1(L_i,\Z)$ we have an element $V_{\delta_i^{cl}}(\EE_{i})\in \C[[\hbar]]^\times$. By analogy with the conventional terminology of resurgent functions, (see e.g. [DelPh]) we call it the {\it  Voros coefficient} (a.k.a {\it Voros symbol}) corresponding to $\EE_{i}$ and $\delta_i^{cl}$. Clearly $V_{\gamma_i^{cl}+\delta_i^{cl}}(\EE_{i})=V_{\gamma_i^{cl}}(\EE_{i})\cdot V_{\delta_i^{cl}}(\EE_{i})$ for any two $\gamma_i^{cl}, \delta_i^{cl}$ as above.

With the data 1)-4) we would like to associate an element $Tr_{{\bf x},\gamma}(\EE, \mu)\in \C[[\hbar]]$ which agrees with Voros symbols in the following sense. For a collection 
$\delta^{cl}=(\delta_i^{cl})_{i\in \Z/k}$ one has (for the natural action $(\gamma_i)_{i\in \Z/k}\mapsto (\gamma_i+\delta_i^{cl})_{i\in \Z/k}$ of closed $1$-homology on relative homology classes ) we have:

$$Tr_{{\bf x},\gamma+\delta^{cl}}(\EE, \mu)=Tr_{{\bf x},\gamma}(\EE, \mu)\cdot \prod_{i\in \Z/k}V_{\delta_i^{cl}}(\EE_i).$$

The construction goes such as follows.

 
The collection $\EE=(E_i)_{i\in \Z/k}$ determines Voros symbols $V_{\delta_i^{cl}}(\EE_{i})\in \C[[\hbar]]^\times, i\in \Z/k$.
 Let $\rho_i$ be the rank one local system on $L_i$ corresponding to $E_{i}$ under the local RH-correspondence. Then the $\C[[\hbar]]$-module $Ext^n(E_{i}, E_{i+1})$ is naturally decomposed into the direct sum over the connected components of $L_i\cap L_{i+1}$. The summand corresponding to the intersection point $x_{i,i+1}$ is  identified with the free $\C[[\hbar]]$-module $Hom({\rho_i}_{|x_{i,i+1}}, {\rho_{i+1}}_{|x_{i,i+1}})$. Here  $Hom$ is taken in the category of free $\C[[\hbar]]$-modules, and the notation $\rho_{|x}$ means the fiber of the local system $\rho$ at the point $x$.  Let $pr_{i,i+1}(\mu)$ denote the projection of $\mu_{i,i+1}$ onto the direct summand $\C[[\hbar]]\cdot [x_{i,i+1}]$ of $Ext^n(E_{i}, E_{i+1})$ corresponding to $x_{i,i+1}$. By the local RH-correspondence we obtain an element $pr_{i,i+1}^{Betti} \in Hom({\rho_i}_{|x_{i,i+1}}, {\rho_{i+1}}_{|x_{i,i+1}})$. 
 
Finally $Tr_{{\bf x},\gamma}(\EE,\mu)\in \C[[\hbar]]$ is defined as the trace of the cyclic compositions of operators $pr_{i,i+1}^{Betti}, i\in \Z/k$ combined with the  consecutive parallel transports (holonomies) along the paths $\gamma_i, i\in \Z/k$.  This definition is manifestly cyclically invariant.

We remark that by definition the series $Tr_{{\bf x},\gamma}(\EE,\mu)$ agrees with Voros symbols in the above sense.

\begin{conj}\label{resurgence of monodromy for wheels}
Assume that the global quantized algebra of functions $\OO_\hbar(M)$ as well as all elements $\mu_{i,i+1}, i\in \Z/k$ are defined over the ring of analytic germs $\C\{\hbar\}$. Then

a) Voros coefficients $V_{\delta_i^{cl}}(\EE_{i})$  are resurgent series in $\hbar$.

b) Each series  $Tr_\gamma(\EE,\mu)$ is resurgent in $\hbar$.
\end{conj}

\begin{rmk}\label{higher rank DQ-modules}
The above considerations as well the Conjecture can be generalized to the case when $\EE_{i}/\hbar \EE_i$ are pure coherent sheaves of higher rank supported ion $L_i, i\in \Z/k$.

\end{rmk}

Let us comment on the Conjecture \ref{resurgence of monodromy for wheels}.  One should  think of  $Tr_\gamma(\EE,\mu)$ as formal expansion at $\hbar=0$ of the holomorphic section of a non-linear bundle on $\C_\hbar$ with the fiber,
which is isomorphic to the complex toric variety $\mathcal{X}=\C^{\prod_{i\in \Z/k}\#(L_i\cap L_{i+1})}\times  \C^{\prod_{i\in \Z/k}rk(H_1(L_i))}$. This is a slight generalization of the formalism of [KoSo12], where the case of tori only was considered. The fiber bundle is glued from the trivial ones using Stokes isomorphisms, similarly to the loc. cit.  The corresponding Stokes isomorphisms satisfy the so-called $2d-4d$ wall-crossing formulas, differently from [KoSo12], where only $2d$ (a.k.a Cecotti-Vafa) wall-crossing formulas were used. The reason for the appearance of more complicated $4d$ WCF (see [KoSo1]) is the possibility to have for certain directions in $\C_\hbar$ pseudo-holomorphic discs with boundaries on $L_i, i\in \Z/k$. This does not happen for exact Lagrangians in exact symplectic manifolds or for the graphs of closed $1$-forms in the case of cotangent bundles, or for a class of ``good" Lagrangian submanifolds discussed in the next subsection.

\subsection{Normalization problem}\label{normalization}

Normalization problem discussed in this subsection appears when one is trying  to combine the formalism of the Section \ref{wheels of Lagrangians}  with ``real'' examples, e.g. with  WKB expansions as in the example in Section \ref{example  of wkb} above. More generally having two $\hbar$-differential operators $P_1,P_2\in \C\{\hbar\}[x][\hbar\partial_x]$ analyticially depending on $\hbar$ one has two multivalued quantum wave functions $\psi_{L_1}, \psi_{L_2}$ supported on the spectral curves $L_1,L_2\subset T^\ast \C_x$, and obtained via WKB expansions applied to each of the equations $P_1\psi_{L_1}=0, P_2^\ast\psi_{L_2}=0$, where $P_2^\ast$ is the formal adjoint operator. We understand $\psi_{L_1}$ as a left quantum wave function and $\psi_{L_2}$ as a right quantum wave function.These quantum wave functions should be normalized at some points of the  spectral curves (or points ``at infinity" as in the example in Section \ref{example  of wkb}). Then for each point $p$ where $L_1$ and $L_2$ intersect transversally one has the local pairing $\langle \psi_{L_1},\psi_{L_2}\rangle_p:=\langle \psi_{L_1}|\psi_{L_2}\rangle_p\in \C[[\hbar]]$. This pairing can be understood in terms of a certain wheel of five Lagrangian submanifolds $L_i\in T^\ast\C_x, i\in \Z/5$.

 
 Namely we choose  two ``normalizing" Lagrangian submanifolds $L_0, L_3$ as two cotangent fibers for generic $L_1, L_2$ and add the zero section $L_4$. Thus  we obtain a wheel of Lagrangian submanifolds $L_0,L_1,L_2,L_3,L_4$ consisting of five elements.


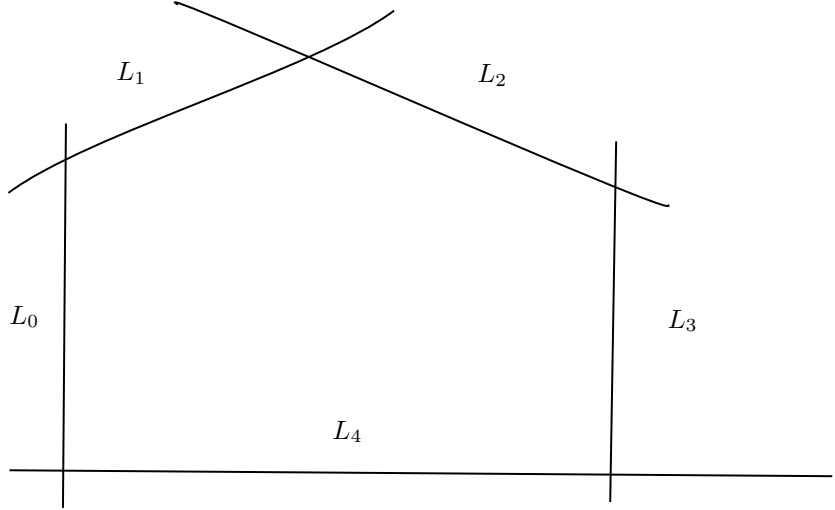
\tikzset{every picture/.style={line width=0.75pt}} 

\begin{tikzpicture}[x=0.75pt,y=0.75pt,yscale=-1,xscale=1]

\draw    (102,129) .. controls (142,99) and (256.5,67) .. (296.5,37) ;
\draw    (187.5,34) .. controls (162.5,19) and (446.5,146) .. (434.5,135) ;
\draw    (131,94) -- (129.5,288) ;
\draw    (408.5,103) -- (405.5,285) ;
\draw    (102.5,269) -- (517.5,272) ;

\draw (101,184.4) node [anchor=north west][inner sep=0.75pt]    {$L_{0}$};
\draw (155,61.4) node [anchor=north west][inner sep=0.75pt]    {$L_{1}$};
\draw (337,62.4) node [anchor=north west][inner sep=0.75pt]    {$L_{2}$};
\draw (433,186.4) node [anchor=north west][inner sep=0.75pt]    {$L_{3}$};
\draw (264,242.4) node [anchor=north west][inner sep=0.75pt]    {$L_{4}$};

\end{tikzpicture}

\vspace{2mm}

One can show that for the corresponding collection $\EE=(E_{L_i})_{i\in \Z/5}$ the element $Tr_{{\bf x},\gamma}(\EE,\mu)$ depends only on the isomorphisms $\EE_{L_0}\simeq \EE_{L_0}^{can}$ and 
$\EE_{L_3}\simeq \EE_{L_3}^{can}$  and does not depend on the similar isomorphisms for other Lagrangians in the wheel. 
We also note that a choice of the fifth Lagrangian submanifold $L_4$ looks artificial, since $L_4$ is chosen with only one purpose, to close the wheel. We will explain below that $L_4$ can be omitted because of a certain ``good" properties of $L_0$ and $L_3$.

The Lagrangian submanifolds $L_0, L_3$ satisfies in the above example the following ``uniqueness'' property: there is a unique (up to isomorphism) rank one  holonomic $DQ$-module supported on each of them. This $DQ$-module corresponds to a $D$-module generated by the delta-function at a point of the zero section. 
This observation suggests that in general we should have a sufficiently large class of ``good" Lagrangian submanifolds $L\subset M$ which satisfy a similar uniqueness property.   ``Good" Lagrangian submanifolds (or, more generally, Lagrangian subvarieties) should be supports of some objects which we call {\it normalizing objects}.  The normalizing objects can be approached via our conjectural Riemann-Hilbert correspondence and should exist both in local and global versions of the relevant Fukaya categories and the categories of holonomic $DQ$-modules. There should be certain compatibility between the four possible versions of the normalizing objects, according to the diagram in Conjecture \ref{diagram of categories}. Although we do not have precise definitions, we are going to list below some desired properties of the normalizing objects:

1) If the category of holonomic $DQ$-modules is defined over the ring of analytic germs $\C\{\hbar\}$ then a normalizing object $E$ is defined over $\C\{\hbar\}$ and the corresponding object $E^{form}$ over $\C[[\hbar]]$ is a holonomic $DQ$-module of rank $1$ supported on a Lagrangian subvariety $L$.

2) Assuming 1) the corresponding object $E^{form}\otimes_{\C[[\hbar]]}\C((\hbar))$ of the local category of holonomic $DQ$-modules supported on $L$ is mapped under the local Riemann-Hilbert functor $RH_{loc}$ into a local system $\rho$ such that the restriction of $\rho$ onto the smooth part $L^{sm}$ is a rank $1$ trivial local system.

3) A normalizing object is endowed with a trivialization of $\rho_{|L^{sm}}$ in the notation of 2).

4) The Stokes functors $g_\theta$ described in Section \ref{analytic families of Fukaya categories} transform $\rho$ from 2) into an isomorphic object.

5) For each Stokes direction $\theta$ we have a choice of the isomorphism  in 4).

6) In the case when we have a global algebra $\OO_{M,\hbar}$   over $\C\{\hbar\}$ together with an isomorphism 
$\OO_{M,\hbar}^{op}\simeq \OO_{M,-\hbar}$, any normalizing object $E$ is self-dual in the sense that the left $\OO_{M,\hbar}^{op}$- module $Ext^n(E, \OO_{M,\hbar})$ is isomorphic to  $E$ after applying the isomorphism $\OO_{M,\hbar}^{op}\simeq \OO_{M,-\hbar}$.

7) In the case when we have a global algebra $\OO_{M,\hbar}$   over $\C\{\hbar\}$ a normalizing object $E$ is endowed with a distinguished generator $vac_E$. 

8) Suppose that $M$ is endowed with a QWFS and the corresponding quantized algebra over $\C[[\hbar]]$ is obtained from the one in 7) by the change of scalars. The the generator $vac_E$ corresponds to a canonical (left) quantum wave function $\psi_L$.

In the example of the WKB expansion  the cotangent fibers  $L_0$ and $L_3$ give rise to normalizing objects as long as we endow them with trivialized local systems of rank $1$. Then properties 2) and 3) are equivalent to a choice of isomorphism $E_{L_i}^{form}\otimes_{\C[[\hbar]]}\C((\hbar))\simeq E_{L_i}^{can}, i=0,3$. The conditions 4) and 5) are redundant.

From the point of view  of our conjectural Riemann-Hilbert correspondence one can try to define a normalizing object supported on $L$ under the condition that the natural functor from the local Fukaya category $\FF_{L, loc,\hbar}(L)$ to the global Fukaya category $\FF_\hbar(M)$ is an embedding for all $\hbar \in \C^\ast$. This is the case for all $\hbar$ which do not belong to the Stokes rays $Arg(\hbar)=\{\int_\gamma\omega^{2,0}), \gamma\in \pi_2(M,L,\Z)\}$, but sometimes this true for all $Arg(\hbar)$.
Namely, we have an embedding of categories (i.e. a faithful functor) in the case of  absence of pseudo-holomorphic discs (for a compatible pseudo-holomorphic structure) with boundary on $L$ for all symplectic forms $Re(\omega^{2,0}/\hbar), \hbar\in \C^\ast$. At the same time there are cases when normalizing objects exist even in the presence of the pseudo-holomorphic discs.

Suppose that instead of a wheel of Lagrangian submanifolds and wheel of objects as before we have a finite sequence of complex Lagrangian submanifolds $L_0, L_1,...,L_k$ and objects $E_0, E_1,..,E_k$ endowed with the data as above (but now $i\in \{0,1,...,k-1\}$ instead of $i\in \Z/k$) and such that $E_0$ and $E_k$ are normalizing objects. Then the collection of elements $\mu_{i,i+1}, 0\le i\le k-1$ and paths $\gamma_i, 1\le i\le k-1$ gives rise to an element $Hol_{{\bf x},\gamma}(\EE,\mu)\in \C[[\hbar]]$ defined in the same way as before.

In particular in the above example with five Lagrangian submanifolds we have only four of them, namely $L_0,L_1,L_2,L_3$, two paths and three elements $\mu_{i,i+1}$. 

In the case $k=1$ assume that we have a symmetric QWFS (i.e. the condition 6) is satisfied). Then we have only two normalizing objects $E_0,E_1$ and one intersection point $x_{0,1}\in L_0\cap L_1$. We do not have to choose any path, but we choose an element $\mu_{0,1}\in Ext^n(E_0,E_1)$. Moreover we expect that there is a natural choice of distinguished element $\mu_{0,1}$. Namely since $Ext^n(E_0,E_1)\simeq E_0^\ast\otimes E_1$ and we have distinguished generators $vac_{E_0^\ast}$ and $vac_{E_1}$ and the corresponding element $\mu_{0,1}^{can}:=vac_{E_0^\ast}\otimes vac_{E_1}$. We expect that the series in $\C[[\hbar]]$ corresponding to $\mu_{0,1}^{can}$ coincides with the pairing $\langle \psi_{L_0}|\psi_{L_1}\rangle_{x_{0,1}}$.

\subsubsection{Some examples of normalizing objects} 

The simplest class of normalizing objects corresponds of such Lagrangian submanifolds $L$ that there are no pseudo-holomorphic discs with boundary on $L$. Example of such $L$ are exact Lagrangians in exact symplectic manifolds or graphs of closed $1$-forms in cotangent bundles.

In this subsection we will describe some other classes of normalizing objects. 

Namely we will see that if $L\subset M$ is a Lagrangian submanifold such the pairing $H_1(L, \Z)\otimes H_{2n-1}(L,\Z)\to \Z$ is trivial then $L$ gives rise to a normalizing object. \footnote{The reader can keep in mind the illustrating example of the curve $z_1+z_2=1$ is the symplectic manifold $((\C^\ast)_{z_1,z_2}^2,{dz_1\over z_1}\wedge{dz_2\over z_2})$. More generally one can consider the case when $L$ is a punctured $\C{\bf P}^1$.} But we start with considerations which do not use this assumption.

Let us fix $\theta=Arg(\hbar)\in \R/2\pi \Z$ and choose a generic  almost complex structure $J_\theta$ compatible with $Re(\omega^{2,0}e^{-i\theta})$.
Let us assume that there exists $\hat{\gamma}\in H_2(M,L,\Z)$ such that $\theta=Arg(\hbar)=Arg(\int_{\hat{\gamma}}\omega^{2,0})$ and such that $\hat{\gamma}$ can be a represented by a $J_\theta$-holomorphic disc with boundary on $L$. Let us assume that $\hat{\gamma}:=\hat{\gamma}_{min}$ is minimal in the sense that $|\int_{\hat{\gamma}}\omega^{2,0}|$ is minimal possible among classes $\hat{\gamma}$ for which such discs exist. It is well-known that for $\hat{\gamma}_{min}$ the corresponding family of $J_\theta$-holomorphic discs in this class forms a smooth compact oriented $(2n-2)$-dimensional manifold. The boundary points of these discs form a $(2n-1)$-dimensional closed manifold which defines a  homology class $\delta\in H_{2n-1}(L,\Z)$.

Let $\gamma\in H_1(L,\Z)$ be the image of  $\hat{\gamma}_{min}$ under the natural map $H_2(M,L,\Z)\to H_1(L,\Z)$.


Let $\rho: H_1(L,\Z)\to \C^\ast$ represents the isomorphism class of a rank one local system on $L$. Then the pair $\gamma,\delta$ as above gives rise to the  partially defined automorphism $g_{\gamma,\delta}$ such that $g_{\gamma,\delta}(\rho):=\rho^{new}: H_1(L,\Z)\to \C^\ast$ is given by the formula
$$\rho^{new}(\mu)=\rho(\mu)(1-e^{-{\int_{\hat{\gamma}}\omega^{2,0}\over \hbar}}\rho(\gamma))^{\langle \mu,\delta\rangle}.$$

The map $g_{\gamma,\delta}$ can be understood as a non-linear automorphism of an open domain in the stack of rank one local systems on $L$. \footnote{Technically speaking we have defined the automorphism on the set of points only. Moreover, one should understand the automorphism as the one for the stack of non-archimedean rank one local systems on $L$.} In the case when all $J_\theta$-holomorphic discs are multiple covers of discs with the class $\hat{\gamma}_{min}$ the Stokes automorphism $g_\theta$ coincides with  $g_{\gamma,\delta}$.

In general e.g. for non-minimal $\hat{\gamma}$ the non-linear automorphism $g_{\gamma,\delta}$ will be a more complicated transformation containing possibly infinite product of cluster transformations of a similar form.

If the pairing  $H_1(L, \Z)\otimes H_{2n-1}(L,\Z)\to \Z$  is trivial then for any choice of $\hat{\gamma}_{min}$ we have $g_{\gamma,\delta}=id$. We expect that  in general (e.g. when the discs with non-minimal $\hat{\gamma}$ appear ) the Stokes isomorphism $g_\theta$ will be trivial. Under this assumption the rank $1$ local systems $\rho$ and $g_\theta(\rho)$ are isomorphic. Nevertheless there is a problem, since there is no preferred isomorphism $\rho\simeq g_\theta(\rho)$.

In order to fix the ambiguity   let us choose a point $x_0\in L$ which does not belong to the boundaries of $J_\theta$-holomorphic discs for all $\theta$. Then instead of the {\it stack} of rank one local systems on $L$ we apply the above considerations to the {\it moduli space} of rank one local systems on $L$ trivialized at $x_0$.

In this way for any $\theta=Arg(\hbar)$ we obtain a trivialization $\tau_{x_0}$ of the trivial torsor with the fiber over a rank one local system $\rho$ equal to $Iso(\rho, g_\theta(\rho))$. 

For a pair of two such trivializations corresponding to the two choices of points $x_0, x_1$ we have a ratio $c_{x_0,x_1}:=\tau_{x_0}\cdot \tau_{x_1}^{-1}$ which  is in general  an invertible function on a non-archimedean open domain in the moduli space of rank one local systems.
In the case when only multiple covers of discs in class $\hat{\gamma}_{min}$ appear there is an explicit formula for $c_{x_0,x_1}$. 

Namely notice that by our assumption on $x_0,x_1$ the class $\delta$ defined above comes naturally from a class $\delta^\prime\in H_{2n-1}(L-\{x_0,x_1\},\Z)$. Therefore we have a well-defined integer $\nu_{x_0,x_1}$ which is the intersection number of any continuous path joining $x_0$ and $x_1$ with $\delta^{\prime}$.


Then 
$$c_{x_0,x_1}(\rho)=(1-e^{-{1\over{\hbar}}\int_{\hat{\gamma}}\omega^{2,0}}\rho(\gamma))^{\nu_{x_0,x_1}}.$$ 


In general we obtain a function $c_{x_0,x_1}$  with values in the series in the variables $exp(-{1\over{\hbar}}\int_{\hat{\gamma}}\omega^{2,0}), \hat{\gamma}\in H_2(M,L,\Z)$ satisfying the $1$-cocycle condition
$$c_{x_0,x_1}c_{x_1,x_2}c_{x_2,x_0}=1.$$
One can also check that the change of a compatible almost complex structure $J_\theta$ does not affect these functions as long as the $x_0,x_1$ do not lie on the boundaries of $J_\theta$-holomorphic discs.

\begin{exa}\label{compactification for quantum dilog}

Let $L\subset (\C^\ast)^2_{z_1,z_2}$ be a curve $z_1+z_2-1=0$, where  $(\C^\ast)^2_{z_1,z_2}$ is endowed with the standard  symplectic form ${dz_1\over{z_1}}\wedge {dz_2\over{z_2}}$. Then $L$ is a Lagrangian submanifold isomorphic to a $\C{\bf P}^1$ without three points $p_1,p_2,p_3$. We think of these points as points ``at infinity'' on $L$. In order to define $\psi_L$ canonically we need to fix three functions $c_{p_1,p_2}, c_{p_2,p_3}, c_{p_3,p_1}$ on the torus satisfying the condition  
$c_{p_1,p_2}c_{p_2,p_3}c_{p_3,p_1}=1$. In this particular case we expect that all $c_{p_i,p_j}=1$. Indeed there is an action of the $\Z/3$ on 
$(\C^\ast)^2_{z_1,z_2}$ which preserves the symplectic structure and the submanifold $L$, and cyclically permutes $p_1,p_2,p_3$. Notice that there are no non-trivial classes in $H_1(L,\Z)$ invariant with respect to $\Z/3$. It follows that the above-defined intersections indices $\nu_{p_i,p_j}=0$. 

The corresponding quantum wave function is related to the quantum dilogarithm.

\end{exa}

These considerations suggest that we can {\it canonically define the quantum wave function $\psi_L$ associated with a normalizing object}. For that we also need to choose a certain partial compactification $\overline{M}$ of $M$ and consider only those $L$ for which the closure $\overline{L}\subset \overline{M}$ is compact. Then we need to define the above function $c_{x,y}$ only for points $x,y\in \overline{M}-M$.

We will explain in the subsequent papers on the project that there is a class of partial compactifications of $M$ called {\it log extensions}. Such a log-extension $M_{log}$ is non-unique. For a Lagrangian submanifold $L$ such that the closure $\overline{L}\subset M_{log}$ is compact, the exceptional divisor of the  blow-up of $M_{log}-M$ at $M_{log}\cap \overline{L}$ provides a good set ``at infinity'' where the functions $c_{x,y}$ should be defined.

\subsection{Resurgence and quantum wave functions}\label{resurgence on wave functions}

Let us discuss  resurgence of perturbative expansions in the case of path integral with boundary conditions. In other words we are interested in the analytic continuations of the Borel transforms of local exponential integrals and their behavior near singular points. 

Recall the considerations of Sections \ref{path integrals}, \ref{wheels of Lagrangians}. Let us  assume that the Hamiltonian $H=0$. Consider the formal power series in $\hbar$ given by the local expansion of the path integral  at the transversal intersection point $x_j\in L_0\cap L_1$:
$I_{x_j}(\hbar)=\langle \psi_{L_{0}}, \psi_{L_1}\rangle_{x_j}$ (i.e.  $I_{x_j}(\hbar)$ corresponds to the {\it modified} exponential integral from Section \ref{WCF}). 

\begin{conj}\label{resurgence for pair of Lagrangians} Assume that $L_0$ and $L_1$ are supports of normalizing objects and $\psi_{L_i}, i=0,1$ are the corresponding canonical quantum wave functions. Then the series $\langle \psi_{L_{0}}, \psi_{L_1}\rangle_{x_j}\in \C[[\hbar]]$ is resurgent.\footnote{For $H\ne 0$ we should modify our considerations using the parallel transport of quantum wave functions. The resulting series are expected to be resurgent.}

\end{conj}

\begin{rmk}\label{resurgent quantum wave functions}
We expect that there is bigger class of {\bf resurgent quantum wave functions} for which the Conjecture \ref{resurgence for pair of Lagrangians} holds.
\end{rmk}

There are two ``explanations" of the Conjecture \ref{resurgence for pair of Lagrangians}. From the point of view of physics it should follow from the interpretation of the local pairing of quantum wave functions in terms of the infinite-dimensional exponential integral. A mathematical explanation should follow from the above-discussed story with generalized local and global Riemann-Hilbert correspondences and analyticity of the arising wall-crossing structures.

\begin{rmk}\label{reminder about normalizing objects}
In relation to the above Conjecture \ref{resurgence for pair of Lagrangians} we remark that we can take as boundary Lagrangian conditions for the path integral normalizing objects  for which   $L_0,L_1$ are either exact Lagrangian submanifolds or graphs of closed holomorphic $1$-forms, or products of varieties of the form $\{x+y=1\}\subset (\C^\ast)^2$, or more generally arbitrary $K_2$-Lagrangians.
\end{rmk}

\

Assume for simplicity that all intersection points of $L_0$ and $L_1$ are transversal, and there are no ``intersection points at infinity". Then we have finitely many formal series $\langle \psi_{L_{0}}, \psi_{L_1}\rangle_{x_j}=\sum_{k\ge 0}c_{k,j}\hbar^k\in \C[[\hbar]]$ labeled by $x_j\in L_0\cap L_1$. Then we can specify the Conjecture \ref{resurgence for pair of Lagrangians}. Namely we expect that singularities of the analytically continued Borel transform $B_j(s):=\sum_{k\ge 0}{c_{k,j}\over {k!}}s^k$ of each of the series belongs to the countable set $\Sigma\subset \C$, which is the image of $H_2(M, L_0\cup L_1, \Z)$ under the map $\gamma\mapsto \int_\gamma\omega^{2,0}\in \C$. Moreover we expect that for any point $s_0\in \Sigma$ the analytic continuation of $B_j(s)$ has singularities at $s_0$ of the form
$${n_{jj^\prime}\over {2\pi i}}\log(s-s_0)B_{j^\prime}(s-s_0)+\text{holomorphic germ at}\, s=s_0.$$

Here $j^\prime$ corresponds to a possibly different intersection point $x_{j^\prime}$, and $n_{jj^\prime}$ is the Stokes index, which  in our case is the virtual number of pseudo-holomorphic $2$-gons $D$ with boundaries on $L_0\cup L_1$ and vertices at  $x_j, x_{j^\prime}$, and such that $s_0=\int_D\omega^{2,0}$.

\begin{rmk}\label{appearance of logarithms}
In general the intersections are not transversal or even isolated. Then we expect that series $\langle \psi_{L_{0}}, \psi_{L_1}\rangle_{x_j}$ will include fractional powers of $\hbar$ as well as polynomials in $\log\hbar$. They will be associated with connected components of $L_0\cap L_1$ and linear functionals on $Hom$-spaces of the corresponding objects of the local Fukaya category of $L_0\cup L_1$.

In case if logarithms are present the Borel transform is defined by:
$$\hbar^\lambda (\log \hbar)^m\mapsto \left({d\over {d\lambda}}\right)^m\left({s^\lambda\over{\Gamma(\lambda+1)}}\right)=
\sum_{i=0}^{m}\binom{m}{i}s^\lambda(\log s)^i\left({d\over {d\lambda}}\right)^{m-i}\left({1\over{\Gamma(\lambda+1)}}\right).$$

\end{rmk}

In the case when $L_0$ and/or $L_1$ are not supports of normalizing objects one can proceed such as follows. For each $L_i, i=0,1$ choose a chain of Lagrangian submanifolds $L_{i,0}, L_{i,1},...,L_{i,m_i}=L_i$ such that $L_{i,0}$ is the support of a normalizing object. After that we choose intersection points, paths etc. as in Section \ref{wheels of Lagrangians}.
}. As we explained in loc.cit. this gives us in the end quantum wave functions $\psi_{L_i}$ supported on $L_i, i=0,1$. In this case we also expect the resurgence of the pairing  $\langle \psi_{L_{0}}, \psi_{L_1}\rangle_{x_j}$,
but the behavior of $B_j(s)$ near singularities will not be as simple as above. The corresponding formula involves the products of local pairings of quantum wave functions with Voros symbols. In general we have countably many resurgent series, and the corresponding analytic wall-crossing structure can be encoded in $2d-4d$ wall-crossing formulas differently from $2d$ wall-crossing formulas above. 

\

Finally, if the Hamiltonian $H\ne 0$ the number of intersection points will be in general infinite. In this case we don't have a precise language of analytic WCS.

\subsection{Conjecture about Fukaya-Seidel category}\label{FS and path space}

Recall that perturbative expansions of finite-dimensional exponential integrals can be interpreted in terms of the pairings of quantum wave functions such as follows.
Let $X$ be a  smooth complex manifold endowed with a holomorphic function $f$ and the holomorphic volume form $vol$. Then $M=T^\ast X$ contains two Lagrangian submanifolds $L_0=X$ (zero section) and $L_1=graph(df)$. We endow $L_i, i=0,1$ with the  quantum wave functions $\psi_i, i=0,1$ corresponding to the half-densities  $(vol)^{1/2}$.
Then the local pairing at the critical point $a\in Crit(f)$, i.e. $\langle \psi_0, \psi_1\rangle_a$ is given by the perturbative expansion of the modified exponential integrals $\int_{\gamma_a}e^{f/\hbar}vol$ over a cycle $\gamma_a$ in the local Betti homology..

On the other hand we can interpret the same pairing in terms of the path integrals.
Namely,  consider the space $\{\phi\in P(L_0, L_1), \phi:[0,1]\to M\}$ of smooth paths in the complex symplectic manifold $M=T^\ast X$ with the endpoints  $\phi(0)\in L_0=X, \phi(1)\in L_1=graph(df)$. The Lagrangian submanifolds $L_i, i=0, 1$ are exact. The  natural closed $1$-form $\eta$ on $P(L_0,L_1)$ is therefore also exact. Then we have the {\it globally} defined potential $W_{L_0,L_1}$ on $P(L_0,L_1)$.

{\it We conclude that the above finite-dimensional exponential integrals are equal to the (ill-defined) infinite-dimensional ones associated with the function $W_{L_0,L_1}$.}

One can wonder if it is possible to develop  infinite-dimensional versions of the mathematical structures discussed previously in this paper, like Betti and de Rham cohomology and their categorifications. Some results concerning the categorification (including a generalization in the framework of $2$-categories) can be found  in [Bou], [DoRez]. In particular one can hope to define a Fukaya-Seidel category $\FF S(P(L_0,L_1), W_{L_0,L_1})$ associated with the potential $W_{L_0,L_1}$.

 The above discussion gives us a hope that the following conjecture holds.

\begin{conj} \label{Equivalence of two definitions of FS} $\FF S(P(L_0,L_1), W_{L_0,L_1})$ is equivalent to $\FF S(X, f)$.

\end{conj}

\subsection{Quantum wave functions and sums associated with polytopes}\label{sums and polytopes}

The meaning of this subsection will be clear in the discussion of Nahm sums.

Let $L\subset \C^{2n}_{q_1,..,q_n, p_1,...,p_n}=T^\ast\C^n_{q_1,...q_n}$ be a (possibly singular) Lagrangian subvariety and $\psi_L$ be a quantum wave function supported on $L$. For each smooth point of $L$ where it is transversal to the vertical foliation $p_i=const_i, 1\le i\le n$ the quantum wave function can be locally written as
$$\psi_L(q)=e^{\FF_0(q)\over \hbar}(\mu_0(q)+\hbar \mu_1(q)+...),$$
where $q=(q_1,...,q_n)$, and $\mu_m=g_m(q)(dq_1\wedge dq_2\wedge...\wedge dq_n)^{1/2}$.

If the quantum wave function belongs to the class mentioned in Remark \ref{resurgent quantum wave functions} then  we can say for every such point the sum $\sum_{n\ge 0}g_m(q)\hbar^m$ is resurgent in $\hbar$ and holomorphically depends on $q$. Notice that this sum can be interpreted as the pairing $\langle\psi_L,\psi_{L^\prime}\rangle_{(q,0)}$, where $L^\prime$ the cotangent fiber $T^\ast_q \C^n_q$.

Let $L^{\prime\prime}$ be the zero section of the cotangent $T^\ast\C^n_{q_1,...q_n}$. If $L$ intersects $L^{\prime\prime}$ at a point $(q^{(0)},0)$ then the pairing $\langle\psi_L,\psi_{L^{\prime\prime}}\rangle_{(q^{(0)},0)}$ is equal to the asymptotic expansion of the local integral $\int_{th_{q^{(0)}}} e^{\FF_0(q)\over \hbar}(g_0(q)+\hbar g_1(q)+...)dq_1\wedge ...\wedge dq_n$, where $th_{q^{(0)}}$ is a local thimble for $\FF_0$ emanating from $th_{q^{(0)}}$.

Recall that we expect a similar answer for non-transversal intersections, and all arising sums in $\hbar$ are expected to be resurgent. Furthermore if there are no critical points at infinite they form a resurgence package (see Section \ref{resurgence package}).

\

Our next goal is to propose a hypothetical picture for an analogous resurgence package in case if we replace integrals over $(q_1,...,q_n)$ by sums over 
the subsets of $\hbar \Z^n$ possibly satisfying a collection of linear inequalities with rational coefficients.

Let $P$ be a rational polyhedral cone  in $\R^n=\R\otimes \Z^n\subset \C^n_{q_1,...,q_n}$ with the vertex at the origin.\footnote{Considerations will be generalized to the case of a rational polyhedral subset, not necessarily conical.}
Assume that $\FF_0, g_i, i\ge 0$ are  well-defined analytic functions on $P$ and and $(g_0(q)+\hbar g_1(q)+...)$ is an asymptotic expansion in $\hbar$ of a $C^\infty$ function $g(q,\hbar)$ defined for $q\in P$ and $\hbar>0$.

We are interested in the resurgence properties of the following sum
$$S_P(\hbar)=\sum_{k=(k_1,...,k_n)\in P\cap \Z^n}e^{\FF_0(k\hbar)\over \hbar}g(k\hbar,\hbar).$$
We claim that under certain convergency conditions (we do not know them precisely at the moment)  the function $S_P(\hbar)$ is well-defined and resurgent. Moreover $S_P(\hbar)$ is a part of a resurgent package (see Section \ref{resurgence package}), and all series which arise in the resurgent package are the local pairings of $\psi_L$ with auxiliary wave function $\psi_P$ associated with the polytope $P$. 

More precisely, we associate with $P$ a complex singular Lagrangian subvariety $L_P$ which is the union of the shifted complexified conormal bundles of all faces of $P$ as well as the zero section $p_1=...=p_n=0$ (the latter corresponds to the open face which is the interior of $P$). The intersection points of $L\cap L_P$ are in one-to-one correspondence with the critical points of the restrictions of $\FF_0$ to all faces of $P$ including the interior.
Without loss of generality we may assume that $P$ is the octant $\R_{\ge 0}^n$ and as a result  the general case can be reduced to the case $n=1$.  For that reason we will discuss only the one-dimensional case.

\

We start with the  case $P=\R$ which is the case when all critical points belong to  the interior of $P$. 
Then we have:
$$S_P(\hbar)=S_\R(\hbar)=\sum_{k\in \Z}e^{\FF_0(k\hbar)\over \hbar}g(k\hbar,\hbar).$$

Assume now that $\FF_0$ restricted to $\R\subset \C_q$ is a complex-valued real analytic function  satisfying the following conditions:

1) $Re(\FF_0)$ has a unique  maximum at $q=q_0\in \R$;

2) $(Re(\FF_0))''(q_0)<0$;

3)  $\FF_0'(q_0)=2\pi i n$ for some $n\in \Z$.

4) $Re(\FF_0)(q)<-C\log |q|$ for some $C>0$ and sufficiently large $|q|$.

For simplicity we also assume that $g(q,\hbar):=g(q)$ is a strictly positive real analytic function in $q$ (in particular it does not depend on $\hbar$), and has at most polynomial growth in $q$ at infinity.

Then the sum $S_P(\hbar)$ is absolutely convergent for $0<\hbar\ll 1$. The assumption 3) means that $q_0$ is in a sense a critical point with respect to index $k$ of the sum defining $S_P(\hbar)$. More precisely, for given $C_1>1$ consider those summands in $S_P(\hbar)$ for which $|k\hbar-q_0|<C_1\hbar$ for all  $0<\hbar\ll 1$.
Then there exists $C_2>0$ such that  the ratio of any two such summands belongs to the disc $|z-1|<C_2\hbar$ for all  $0<\hbar\ll 1$.
Using the Poisson summation formula one can show that the following analog of the stationary phase expansion holds:
$$\hbar S_P(\hbar)\sim\left({2\pi \over -{\FF_0''(q_0)}}\right)^{1/2} e^{\FF_0(q_0)\over \hbar}(\sum_{i\ge 0}c_i\hbar^i),c_i\in \C.$$
and moreover the RHS coincides with the  formal integral
$$\int_{th_{q_0}} e^{\FF_0(q)-\FF_0'(q_0)(q-q_0)\over{\hbar}}g(q)dq,$$
where $th_{q_0}$ is the local Lefschetz thimble emanating from the critical point $q_0$  of the function $\FF_0(q)-\FF_0'(q_0)(q-q_0)=\FF_0(q)-2\pi i \,n(q-q_0)$.

We define the Lagrangian variety $L_P$ as the countable disjoint union of shifts of the zero section $\C_q\subset T^\ast \C_q=\C^2_{q,p}$:
$$L_P=\sqcup_{n\in \Z}(\C_q\times \{p=2\pi i n\}).$$
The quantum wave functions $\psi_P:=\psi_{L_P}$ is given by the collection of ${1\over 2}$-densities $(dq)^{1/2}$ assigned to each line.

The reader can see that the above formal integral can be understood as the pairing of the wave function $\psi_L=e^{\FF_0(q)\over \hbar}g(q)(dq)^{1/2}$ with $\psi_{L_P}$, i.e. as a sum of local pairings over  all intersection points $L\cap L_P$  which are $(q_0,2\pi i n)$.

\

The second case is $P=\R_{\ge 0}$ which corresponds to boundary points of $P$. Then

$$S_P(\hbar)=S_{\R_{\ge 0}}(\hbar)=\sum_{k\ge 0}e^{\FF_0(k\hbar)\over \hbar}g(k\hbar,\hbar).$$

In this case we assume that $\FF_0$ restricted to $\R_{\ge 0}\subset \C_q$ is a complex-valued real analytic function  satisfying the following conditions:

1) $Re(\FF_0)$ has a unique  maximum at $q=0$;

2) $(Re(\FF_0))'(0)<0$;

3) $Re(\FF_0)(q)<-C\log |q|$ for some $C>0$ and sufficiently large $|q|$.

Again we  assume that $g(q,\hbar):=g(q)$ is a  real analytic function in $q$, and has at most polynomial growth in $q$ at infinity.
Then the asymptotic expansion
$\hbar S_P(\hbar)$ looks such as follows

$$\hbar S_P(\hbar)\sim e^{\FF_0(0)\over \hbar}(\sum_{i\ge 0}c_i\hbar^i),c_i\in \C,$$
and the RHS coincides with the pairing of $\psi_L=e^{\FF_0(q)\over \hbar}g(q)(dq)^{1/2}$ with the quantum wave function $\psi_{q=0}$ supported on the vertical line $q=0$ and given by ${(dp)^{1/2}\over {1-e^p}}$.

More precisely, let $p_0=\FF_0^\prime(0)$. Then expanding $\FF_0(k\hbar)$ and $g(k\hbar)$ in $k\hbar$ we obtain:
$$S_P(\hbar)\sim e^{\FF_0(0)\over \hbar}\sum_{k\ge 0}e^{p_0 k}\sum_{0\le i\le 2j}a_{ij}k^i\hbar^j.$$
We see that the RHS is well-defined as a series in $\hbar$ and gives the asymptotic expansion of $S_P(\hbar)$ as $\hbar\to 0$. Changing  the summation order we obtain:
$$S_P(\hbar)\sim e^{\FF_0(0)\over \hbar}\sum_{j\ge 0}\hbar^j\left(\sum_{0\le i\le 2j}a_{ij}\sum_{k\ge 0}k^ie^{p_0 k}\right).$$
Notice that
$$\sum_{k\ge 0}k^ie^{p_0k}=\left({d\over{dp}}\right)^i\left({1\over 1-e^{p}}\right)_{|p=p_0}.$$
It follows that the asymptotic expansion of $\hbar S_P(\hbar)$ coincides with the pairing $\langle \psi_L,\psi_{q=0}\rangle_{(0,p_0)}$. 

In general (i.e. without the conditions 1)-3)) one can expect contributions to the asymptotic expansion of $\hbar S_P(\hbar)$ from the critical points of $\FF_0$ in the interior of the ray $\R_{\ge 0}$. In the latter case the asymptotic expansion coincides  with the pairing of $\psi_L$ with $\psi_{L_P}$, where $P=\R$ as in the previous case. This makes plausible the conjecture that for the summation corresponding to $P=\R_{\ge 0}$ the corresponding Lagrangian subvariety is in fact the union $L_P:=\left(\sqcup_{n\in \Z}(\C_q\times \{p=2\pi i n\}\right)\cup \left(\{q=0\}\times \C_p\right).$ The corresponding quantum wave function $\psi_P$ is given by $1/2$-density which is $(dq)^{1/2}$ on the horizontal lines and ${(dp)^{1/2}\over {1-e^p}}$ on the vertical line. This suggestion is compatible with the Euler-Maclaurin formula. Namely, let $\FF_0=0$ and $g(q)$ be a smooth function on $\R$ with fast decay at infinity. Then the Euler-Maclaurin formula says:
$$\hbar\sum_{k\ge 0}g(k\hbar)\sim \int_{0}^\infty g(x)dx+\left[\hbar{g(0)\over 2}+\sum_{n=1,3,5,...}\hbar^{n+1}{\zeta(-n)\over {n!}}g^{(n)}(0)\right].$$
We can interpret the integral summand as the pairing of $\psi_L:=\psi_{p=0}=g(q)(dq)^{1/2}$ with the quantum wave function $\psi_{\{p=0\}\cup \{q=0\}}$ corresponding to the coordinate cross endowed with  the $1/2$-density which is equal to $(dq)^{1/2}$ on the horizontal line and $-{(dp)^{1/2}\over p}$ on the vertical line. The sum in the square brackets corresponds to the pairing of the same $\psi_L$ with $\psi_{q=0}$, where the latter is equal to 
$({1\over p}+{1\over{1-e^p}})(dp)^{1/2}$. Indeed the Taylor expansion of $({1\over p}+{1\over{1-e^p}})$ at $p=0$ is 
$${1\over 2}+\sum_{n=1,3,5,...}{\zeta(-n)\over {n!}}p^n.$$

\

Let us recall that in the higher-dimensional case we can study sums corresponding to rational polyhedral cones using decompositions into primitive simplicial cones and subsequently reducing to the products of copies of the ray $\R_{\ge 0}$. More generally one can consider similar sums corresponding to {\it rational polyhedra}. 
In this case it is reasonable to assume that $\hbar={1\over N}$ where $N\in \Z_{\ge 1}$, and the corresponding sum is defined by the formula:
$$S_P(\hbar={1\over N}):=\sum_{x=(x_1,...,x_n)\in P\cap (\hbar\Z)^n}e^{\FF_0(x)\over \hbar}g(x,\hbar).$$

One can show that   there exists $M=M(P)\in \Z_{\ge 1}$ such the asymptotic expansion of $S_P({1\over N})$ is well-defined when $N\to \infty$ and the residue $N\,mod\, M\in \Z/M\Z$ is fixed. Hence we get $M$ different asymptotic expansions depending on the residue. Furthermore these asymptotic expansions are equal to the pairings of $\psi_L$ with certain quantum wave functions (each depends on the corresponding residue) similar to those considered above.

\section{Chern-Simons theory}\label{CS theory}

In this section we are going to discuss how our general proposal will look in the case of the  Chern-Simons theory (CS theory for short) of a connected compact oriented $3$-manifold $M^3$.  As in the previous section many parts of the exposition will be rather informal, since we do not aim to give a comprehensive treatment of the  Chern-Simons theory.

We fix a  simple compact group $G_c$ which is the gauge group of the CS theory. Let $G$ denote its complexification. In the complexified CS theory one can consider framed and non-framed connections. To define the former, let us choose a base point  $x_0\in M^3$. We define the infinite-dimensional complex manifold  $\mathcal{A}_\C^{fr}$ as a manifold consisting of pairs $(E,\nabla)$ where $E$ is a principal  $G$-bundle on $M^3$ trivialized at $x_0$ endowed with  connection $\nabla$, where the set of pairs is  taken modulo gauge transformations which are equal to $id$ at $x_0$. We call points of the manifold $\mathcal{A}_\C^{fr}$ {\it framed (at $x_0$) bundles with connections} or simply {\it framed connections}. The group $G$ acts on this manifold by changing the framing at $x_0$. A closely related infinite-dimensional {\it stack} ${\mathcal A}_\C$ is defined similarly but the framing at a point of $M^3$ is not fixed. Explicitly ${\mathcal A}_\C$ is the quotient stack ${\mathcal A}_\C^{fr}/G$. 

Tangent space at a point $(E,\nabla)\in {\mathcal A}_\C^{fr}$ can be identified with the vector space 
$\{\eta\in \Omega^1(M^3,ad(E))/d_\nabla(\Omega^0(M^3,x_0,ad(E)))\}$, where $d_\nabla$ is the differential on the complex 
$\Omega^\bullet(M^3, ad(E))$ arising from the covariant derivative, and $\Omega^0(M^3,x_0,ad(E))$ means the space of sections vanishing at $x_0$. One defines the  {\it complexified Chern-Simons $1$-form} $\alpha_{CS}$ on 
$\mathcal{A}_\C^{fr}$ by the formula
$$\alpha{(E,\nabla,\eta)}=\int_{M^3}Tr(F_\nabla,\eta),$$
where $F_\nabla$ denote the curvature of $\nabla$ and $Tr(\bullet, \bullet)$ is a non-degenerate pairing on the Lie algebra $\g=Lie(G)$ associated with a chosen finite-dimensional representation of $\g$.\footnote{In general the normalization of the $CS$ functional depends on a choice of an element of $H^4(BG, \Z)$.}

For the trivial bundle $E$ it is customary to consider the 
CS functional (a.k.a. $CS$ action) $CS$ such that $d(CS)=\alpha_{CS}$ as the multi-valued functional defined on the space of connections by the formula 
$$CS(A)=\int_{M^3}Tr\left({1\over{2}}A\wedge dA+{1\over{3}}A\wedge A\wedge A\right).$$
The functional $CS$ can be considered as a holomorphic function on the $\Z$-cover $\widehat{\mathcal{A}}_\C^{fr}$ of the manifold $\mathcal{A}_\C^{fr}$. The deck transformation of this universal abelian cover changes the value of $CS$ by $(2\pi i)^2$.

The set of critical points of CS functional coincide with the set of zeros of the holomorphic closed $1$-form $\alpha_{CS}$ which are flat connections trivialized at $x_0$. The latter can be identified with $Rep(\pi_1(M^3,x_0), G)$, which is a scheme of finite type over $\Z$. In the case of ${\mathcal A}_\C$ we have the corresponding algebraic stack $Rep(\pi_1(M^3,x_0), G)/G$.

The above scheme and  the moduli stack carry more structures. 
In particular, by [PT] both the character scheme and the character stack admit  dg versions. In the case of the stack it is  derived (-1)-shifted  symplectic stack. It is probably  part of the folklore that this stack carries a constructible (in fact perverse) sheaf of vanishing cycles $\phi_f(\Z)$, where $f$ is an antiderivative of $\alpha_{CS}$ in a small neighborhood of $\ZZ(\alpha_{CS})$ normalized in such a way that $f_{|\ZZ(\alpha_{CS})}=0$. We may use instead of $f$ the multivalued functional $CS$ thus taking into account critical values of the CS functional. Including the monodromy in the game one can speak about about the perverse sheaf $\phi_{CS}(\Z)$  on $\ZZ(\alpha_{CS})\times S^1$. It can be lifted to the equivariant perverse sheaf on the corresponding character scheme.

When quantizing the CS theory one should give a meaning to the Feynman integral over the totally real subvariety $\Aa_c^{fr}\subset  \Aa_\C^{fr}$ of connections of principal $G_c$-bundles trivialized (framed) at $x_0$
$$I(k)=\int_{\Aa_c^{fr}}e^{kCS(A)/2\pi i}{\mathcal D}A,$$
where $k\in \Z_{>0}$ is the level of the CS theory. 
If understood properly, this integral should give in the end a $3d$ TQFT, coinciding with the one constructed in a more combinatorial fashion (Reshetikhin,Turaev, Viro, Witten and others). We will return to this story in Section \ref{compact gauge group}.

In the case of  complexified CS theory  one considers  local perturbative expansions at the critical points. We will explain below that these local expansions should give an analytic wall-crossing structure, which should be thought of as a mathematical structure underlying analytic continuation of the CS theory (cf. [Wit]).

\subsection{Local asymptotic expansions}\label{local expansions}

Let $\ZZ(\alpha_{CS})=\sqcup_{j\in J}\ZZ_j(\alpha_{CS})$ be the decomposition into the finite union of connected components, i.e.  $J=\pi_0(\ZZ(\alpha_{CS}))$. Since each component $\ZZ_j(\alpha_{CS})$ is a set of $\C$-points of an affine scheme of finite type (representation variety) the cohomology $H_c^0(\ZZ_j(\alpha_{CS}), \phi_{CS}(\Z))$ is an abelian group of finite rank.

Physics predicts that there is a $\Z$-equivariant linear map
$$R_j: H_c^0(\ZZ_j(\alpha_{CS}), \phi_{CS}(\Z))\to \overline{\C((\hbar))}[log\,\hbar],$$
which should be thought of as the formal expansion in $\hbar\to 0$ of
the ill-defined ``local Feynman integral''  $\int_{\gamma_j}e^{CS\over{\hbar}}vol, \hbar=2\pi i/k$ after a certain normalization. Here $\gamma_j$ should be thought of (by Poincar\'e duality) as a ``middle-dimensional integration cycle'', which is a generalization of the local thimble in the case when $\ZZ_j(\alpha_{CS})$ is an isolated Morse critical point. Here $vol$ is a ``complexification of the Feynman measure'' restricted to $\gamma_j$. 

In the RHS of the above formula we have the algebraic closure of the field of Laurent series $\overline{\C((\hbar))}=\cup_{N\ge 1} \C((\hbar^{1/N}))$. The $\Z$-action on the LHS comes from the standard monodromy action on the sheaf of vanishing cycles, while  the $\Z$-action on the RHS is given by $\hbar^{1/N}\mapsto e^{2\pi i/N}\hbar^{1/N}, log\,\hbar\mapsto log\,\hbar +2\pi i, N\ge 1$. The degree of the compactly supported cohomology group (dual to the space of integration cycles) is chosen in such a way that it is compatible with middle perversity, hence the middle-dimensional integration cycle corresponds to the degree zero.

Finally, the normalization mentioned above corresponds to the passing from $I(\hbar)$ to $I^{mod}(\hbar)$ in Section \ref{WCF}: 
$$\int_{\gamma_j}e^{CS\over{\hbar}}vol \mapsto e^{-{CS_j\over \hbar}}(2\pi \hbar)^{-dim_\C \widehat{\mathcal{A}}_\C^{fr}}\int_{\gamma_j}e^{CS\over{\hbar}}vol,$$
where $CS_j$ is the corresponding critical value. The dimension $dim_\C \widehat{\mathcal{A}}_\C^{fr}$ is mathematically equal to infinity, while in physics it is finite and defined via the standard $\zeta$-regularization procedure.

If $\ZZ_j(\alpha_{CS})$ consists of a $G$-orbit of a rigid non-trivial flat irreducible connection $\rho_j$ then $R_j$ is  given by the sum running over $3$-valent graphs of expressions obtained by the standard Feynman rules. All summands are convergent integrals, and the sum does not depend on a choice of propagator. For non-rigid flat $\rho_j$ one can utilize the BV formalism (cf. [Mn]).

In the next few subsections we are going to discuss wall-crossing structure and the associated infinite rank ``constructible" (co)sheaf. In Section \ref{resurgence and functional integral} we are going to relate those questions to the resurgence properties of the series which appear as images of the maps $R_j$.

\subsection{Wall-crossing structures from local perturbative expansions}\label{WCS for CS}

Let us introduce the holomorphic function $f_{CS}:=exp(CS/2\pi i): \Aa_\C^{fr}\to \C^\ast$. 

The WCS for the CS theory is defined by analogy with the one for the finite-dimensional holomorphic closed $1$-forms. It is given by the following data:

1) Lattice $\Gamma=H_1(\C^\ast,Critval(f_{CS}),\Z)$, where $Critval(f)$ denotes the set of critical values of the function $f$.

2) Central charge $Z: \Gamma\to \C, \gamma\mapsto 2\pi i\int_\gamma{dw\over{w}} $.

3) Local system of $\Gamma$-graded Lie algebras $\underline{\g}$ on $\C^\ast_\hbar$ with the fiber given by
$${\g}_\hbar:=\oplus_{\gamma\in \Gamma}\g_{\hbar,\gamma}=$$
$$\bigoplus_{\substack{w_1,w_2\in Critval(f_{CS})\\ \gamma\in \Gamma\,s.t.\,\partial \gamma=[w_1]-[w_2]\\j_1,j_2\in J\,s.t.\,f_{CS|\ZZ_{j_m}(\alpha_{CS})=\hbar w_m, m=1,2}}\\}Hom(H^0_c(\ZZ_{j_1}(\alpha_{CS}),\phi_{CS\over{\hbar}}(\Q)), H^0_c(\ZZ_{j_2}(\alpha_{CS}),\phi_{CS\over{\hbar}}(\Q))\\.$$
(Notice that each graded component is finite-dimensional).

4) For any strict sector $V\subset \C$ with the vertex at the origin ($V$ can be a single ray)  we define a pronilpotent Lie algebra $\widehat{\g}_V=\prod_{{Z(\gamma)\over \hbar}\in V}\g_{V,\gamma}$. Let $G_V=exp(\widehat{\g}_V)$ be the corresponding pronilpotent group. \footnote{Notice that $\widehat{\g}_V$ is well-defined because the corresponding set of $\gamma$ with $\g_{V,\gamma}\ne 0$ belongs to a strict convex cone. Hence the Support Property (see [KoSo1], [KoSo12]) is automatically satisfied.} The Stokes isomorphisms $A_\hbar:=A_{V=l_\hbar}\in G_\hbar:=G_{l_\hbar}$ are not equal to $id$ for at most countable set of rays $l_\hbar=\R_{>0}\cdot\hbar$. Here we use the notation like $G_\hbar$ or $A_\hbar$ for the groups, Stokes isomorphisms etc. associated with the admissible ray $l_\hbar$.

The reader should notice that each graded component $\g_{\hbar,\gamma}$ is a finite-dimensional $\Q$-vector space.

In general the Stokes automorphisms are not well-understood even at the physics level of rigor. 
One can hope  that they should be derived 
from the study of the moduli space of solutions to a generalization of Kapustin-Witten equations 
(see [Wit]) to the case of non-Morse critical points of the $CS$ functional. For isolated Morse critical points  
the Stokes automorphisms are  determined by the Stokes indices  which are the numbers of gradient lines of the function $Re(CS/\hbar)$ 
between the corresponding critical points. 
Then one can propose a description of the Stokes automorphisms analogous to the finite-dimensional case with some modifications. 
In this subsection we will treat the Stokes automorphisms as a black box. \footnote{Sergei Gukov pointed out to us that some information about matrix
elements of the Stokes automorphisms in the case of complexified Chern-Simons theory can be derived from the perturbative expansion of
so-called $\hat{Z}$-invariant which was introduced
and studied in a series of his papers with various collaborators, see e.g. [GuMarPut]. The latter also contains probably the first interpretation of these matrix
elements in terms of pseudo-holomorphic discs.}
Later we will explain that these Stokes automorphisms are the same 
as those coming from quantum wave functions associated with a pair of complex 
Lagrangian submanifolds. In other words we will make a connection with the fully faithful embeddings of local to global Fukaya categories.

First we notice  that the map $f_{CS}: {\mathcal A}_\C^{fr}\to \C^\ast_w$ lifts to the map ${CS\over {2\pi i}}:  \widehat{{\mathcal A}}_\C^{fr}\to \C_z$. Consider now the universal abelian covering $\C_z\to \C^\ast_w$ given by $z\mapsto w=exp(z/2\pi i)$. Then the finite set of critical values of $f_{CS}$ in $\C^\ast_w$ corresponds on the line $\C_z$ to a finite number of  arithmetic series with the step $2\pi i$. The Stokes automorphisms  can be described in terms of straight lines in $\C_z$.
Below we will describe the expected properties of the Stokes automorphism $A_\hbar$. They are different for $\hbar\in \R$ and for $\hbar\in \C-\R$. 

a) Let $\hbar\in \C-\R$.

Then the Stokes automorphism $A_\hbar\in G_\hbar$ is an automorphism of the vector space $\oplus_j H^0_c(\ZZ_j(\alpha_{CS}), \phi_{{CS\over{\hbar}}}(\Z))\otimes \Q$, where the sum is taken over such $j\in J$ that $2\pi i \log(w_j)$ up to a shift by $(2\pi i)^2 m, m\in \Z$ belongs to an oriented line of the slope $\theta=Arg(\hbar)$.

Furthermore, $A_\hbar$ is block upper-triangular with respect to the filtration on the direct sum of cohomologies given by the natural order of points on the oriented line. Moreover, $A_\hbar$ is equal to $id$ on each associated graded vector space. 
Notice that the above description of $A_\hbar$ is invariant with respect to shift of the oriented line by $(2\pi i)^2 m, m\in \Z$.

b) Let $\hbar\in \R$.

Then all corresponding critical values $w_j$ belong to the circles $|w|=const_j$.\footnote{These  $const_j$ are critical values of the imaginary part of the $CS$ functional. Conjecturally they coincide with volumes of $3d$ hyperbolic manifolds.}

The real line consisting of $\hbar\in \R$ contains finitely many arithmetic series with the step $(2\pi i)^2$ given by the real critical values of the $CS$ functional. For each point $x$ of this set we have a $\Q$-vector space $W_x=\oplus_j H^0_c(\ZZ_j(\alpha_{CS}), \phi_{CS}(\Z))\otimes \Q$, where the sum is taken over all connected components $\ZZ_j(\alpha_{CS})$  with the critical value of $f_{CS}$ equal to $x$. There is an invertible action of the deck transformation $T: W_x\to W_{x+(2\pi i)^2}$.  Consider the vector space $W:=\oplus_{x<0}W_x\oplus \prod_{x\ge 0}W_x$. It  carries the action of $T$ and $T^{-1}$ thus  making $W$ into a finite-dimensional  $\Q((T))$-vector space isomorphic to $\oplus_j H^0(\ZZ_j(\alpha_{CS}), \phi_{CS}(\Z))\otimes \Q((T))$. The vector space $W$ carries an  $\R$-filtration $W_{\ge s}=\oplus_{x<0}W_x\oplus \prod_{x\ge s}W_x$. Then $A_\hbar$ is an automorphism of $W$, which preserves the filtration, equals to $id$ on the consecutive quotients.

Finally, in both cases a) and b) the following condition is satisfied:

{\bf Integrality condition}:

{\it $A_\hbar$ preserves the image of $\oplus_j H^0_c(\ZZ_j(\alpha_{CS}), \phi_{{CS\over{\hbar}}}(\Z))$ in $\oplus_j H^0_c(\ZZ_j(\alpha_{CS}), \phi_{{CS\over{\hbar}}}(\Q))$.}

The above collection  of automorphisms $A_\hbar, \hbar \in \C^\ast$ is  equivalent to the following data. Choose two opposite half-planes $P^{\pm}$ with the common boundary line which contains the origin, and moreover which does not contain $z-z^\prime$ for any pair of different critical values $z,z^\prime$ of $CS$. The rays $\pm \R_{>0}\in P^\pm$  are the only accumulation rays for the rays $l_\hbar$ with non-trivial $A_\hbar$.  By our assumptions the common boundary line does not contain Stokes rays.  Then the following clockwise ordered product is well-defined:
$$A_{P^{\pm}}=\prod_{l\subset P^{\pm}}A_l.$$
These are automorphisms of the $\Z((T))$-modules $W_\pm=\oplus H^0_c(\ZZ_j(\alpha_{CS}), \phi_{\pm CS}(\Z))\otimes \Z((T))$. More precisely, $A_{P^+}=id+B_+$, where 
$$B_+\in \prod_{j_1,j_2\in J, k\in \Z}Hom( H^0_c(\ZZ_{j_1}(\alpha_{CS}), \phi_{CS}(\Z)),H^0_c(\ZZ_{j_2}(\alpha_{CS}), \phi_{CS}(\Z))))T^k,$$
and $2\pi i (\log(f_{CS}(\ZZ_{j_2}))-\log(f_{CS}(\ZZ_{j_1}))+(2\pi i)^2k\in P^+$. Here $\log$ means the principal branch of the logarithm function.

Similar condition is satisfied for $A_{P^-}=id+ B_-$, but we should replace $T$ by $T^{-1}$. For any $j\in J$ there are two natural isomorphisms $iso_{j,up}, iso_{j,down}$
$$H^0_c(\ZZ_{j_1}(\alpha_{CS}), \phi_{CS}(\Z)))\to H^0_c(\ZZ_{j_1}(\alpha_{CS}), \phi_{-CS}(\Z))),$$
corresponding to the holonomies along the half-circles in the upper and lower half-planes in $\C_\hbar$.

The pair $A_{P^+}, A_{P^-}$ together with these isomorphisms gives rise to a Riemann-Hilbert problem which completely determines the WCS corresponding to the CS theory.

\begin{conj}\label{analytic WCS for CS}

1) The above-described WCS is analytic.

2) Maps $R_j$ from Section \ref{local  expansions} after multiplication by  sufficiently large powers of $\hbar$ give a formal section of the vector bundle on $\C_\hbar$ with the fiber over $\hbar\ne 0$ given by $\oplus_{j\in J}H^0_c(\ZZ_j(\alpha_{CS}), \phi_{CS/\hbar}(\C))$ and with the Deligne's extension to $\hbar=0$, endowed with the monodromy automorphism.

3) The formal section from 2) is an asymptotic expansion at $\hbar=0$ of an analytic section. 

\end{conj}

By the Resurgence Conjecture from [KoSo12] the formal section lives over $\C\{\hbar\}[\hbar^{-1}]\subset \C((\hbar))$.

\begin{rmk}\label{no Deligne extension}
Our use of Deligne extension in 2) is a bit artificial. There is a more natural extension to $\hbar=0$ analogous to the one we studied in Sections 2,3.

\end{rmk}

Finally, the {\it  analyticity of this  WCS} is equivalent to the following two conditions:

1) $A_{P^+}\in Aut(\oplus_{j\in J}\oplus_{j\in J}H^0_c(\ZZ_j(\alpha_{CS}), \phi_{CS}(\C))\otimes \C\{T\}[T^{-1}].$

2) $A_{P^-}\in Aut(\oplus_{j\in J}\oplus_{j\in J}H^0_c(\ZZ_j(\alpha_{CS}), \phi_{-CS}(\C))\otimes \C\{T^{-1}\}[T].$


\subsection{Perverse (co)sheaf of infinite rank}\label{monodromy}

Let us start with the finite-dimensional motivation. Let $(X, vol, f)$ be a triple consisting of an $n$-dimensional complex manifold endowed with a holomorphic volume form $vol$ and a holomorphic $f: X\to \C^\ast$. We assume that critical values of $f$ is a  finite set $Critval(f)$, and that $f$ defines a locally trivial fiber bundle over the complement of the set of critical values. Furthermore we assume that for any critical value $s\in Critval(f)$ and a sufficiently 
small disc $D$ with center $s$ and a point $p\in \partial D$ the relative homology $H_\bullet(f^{-1}(D), f^{-1}(p), \Z)$ is a finitely generated graded abelian group.

Let $\exp:\C\to \C^\ast$ be the universal abelian covering and
$(X_1,vol_1,f_1)$ be the pullback of the above data to the universal $\Z$-covering $X_1\to X$. Then the set $Critval(f_1))$ of critical values of $f_1$ consists of finitely many arithmetic series, and $f_1$ gives rise to a locally trivial fiber bundle outside of this set.

 For any $s_1\in \C-Critval(f_1)$ we have an isomorphism
 
 $$H_n(X_1,f_1^{-1}(s_1),\Z)\simeq \oplus_{z_i\in Critval(f_1)}H_n(f_1^{-1}(D_\epsilon(z_i)), f_1^{-1}(z_i+\epsilon e^{\sqrt{-1}Arg(\theta(\gamma(s_1,z_i))},\Z),$$
where as before $D_\epsilon(p)$ denote a small disc of radius $\epsilon$ with the center at $p$.

This isomorphism depends on the choice of a collection of paths  $\gamma(s_1,z_i)$ (Gabrielov paths)  from $s_1$ to the critical values $z_i$ which are disjoint outside of $s_1$. Notice that the set of paths is infinite because we have infinitely many critical values of $f_1$. 
In the generic case when $s_1$ does not belong to the countable set of straight lines through different pairs of critical values of $f_1$ there is a canonical choice of the paths consisting of straight intervals. In this case  $\theta(\gamma(s_1,z_i))=Arg(s_1-z_i)$. Otherwise one should take as $\theta(\gamma(s_1,z_i))$ the argument of the tangent vector  to the Gabrielov path at the intersection point of the path with the circle $|s_1-z|=\epsilon$. Let us assume for simplicity that we are in the generic case.

The relative homology groups $H_n(X_1,f_1^{-1}(s_1),\Z)$ form a local system of infinite rank over $\C-Critval(f_1)$. In the case of Morse critical points  the fiber can be identified with $\Z[T^{\pm 1}]\otimes \Z^{Critval(f)}$. In general this local system can be extended as a {\it cosheaf} of abelian groups to the whole complex line $\C$.\footnote{The cosheaf of abelian groups is by definition the sheaf with values in the category opposite to the category of abelian groups. To a cosheaf of vector spaces one can associate the sheaf of vector spaces by passing to the duals.} In our case the stalk of the arising cosheaf $\EE$ at any point $s\in \C$ is defined as $H_n(X_1,f_1^{-1}(s),\Z)$. One can show that $\R\Gamma(\C,\EE)=0$.

By analogy with [KaKoPa1], Section 2.3.2 one can say that for any discrete subset $\Sigma\subset \C$ the category of  cosheaves $\EE$ of abelian groups on $\C$ with $\R\Gamma(\C,\EE)=0$  which are local systems outside of $\Sigma$ is equivalent to the category whose objects are:

$\bullet$  Collection of local systems $\LL_s, s\in \Sigma$ on the standard circle $S^1_\theta$  with the fibers $V_{s,\theta}=H^{-1}(D_\epsilon(s), p,\EE)$, where  for any $s\in \Sigma$ we denote by  $D_\epsilon(s)$ with the center at $s$ a sufficiently small disc with the center at $s$ and a marked point $p=s+\epsilon e^{i\theta}\in \partial D_\epsilon(s)$.

$\bullet$ For a pair $s\ne s^\prime \in \Sigma$ a homomorphism $A_{s,s^\prime}: V_{s,\theta}\to V_{s^\prime,\theta}$, where $\theta=Arg(s^\prime-s)$.

Equivalently, in the notation of Section \ref{Betti cohomology} the cosheaf $\EE$ can be described as a collection of vector spaces $V(B,b)$  subject to the conditions dual to those   from loc.cit. More precisely, let $s\ne s^\prime\in \Sigma$. Let us enumerate all points $z_i, 1\le i\le k$ of $\Sigma$ on the interval $[s,s^\prime]$, so that $z_1=s, z_k=s^\prime$.
Fix  a pair of points $z_-$ and $z_+$ which belong to the opposite half-planes separated by the straight line through the pair of points $s\ne s^\prime\in \Sigma$ and which are close to the straight line. Consider topological discs $B_{\le i,\pm}\subset \C, 1\le i\le k$ such that $z_\pm\in \partial B_{\le i,\pm}$ and each $B_{\le i,\pm}$ contains inside points $z_1,...,z_i$ and no other points of $\Sigma$. Then the natural embeddings $(B_{1,\pm},z_\pm)\subset (B_{2,\pm},z_\pm)\subset...\subset (B_{k,\pm},z_\pm)$ induce increasing filtrations of the vector spaces $V(B_{i,\pm},z_\pm), 1\le i\le k$. The successive quotients can be 
naturally identified with $V(D_i,z_{i,\pm})$ where $D_i$ is a small disc centered at $z_i\in [s,s^\prime]\cap \Sigma,1\le i\le k$ and $z_{i,\pm}\in \partial D_i$. Then one gets an isomorphism $V(B_+,z_+)\to V(B_-,z_-)$ which is upper triangular with respect to the filtrations. Identifying the filtered vector spaces via this isomorphism we interpret the above linear map $A_{s,s^\prime}$ as the matrix element corresponding to the map $V_{D_1,z_{1,+}}=V_{s,\theta}\to V_{D_k,z_{k,+}}=V_{s^\prime,\theta}, \theta=Arg(s^\prime-s)$.

\begin{rmk}\label{two descriptions of WCS}
Suppose we have a set $\Sigma\ne \emptyset$ consisting of finitely many arithmetic series in $\C$ with the step $2\pi i$ and a cosheaf $\EE$ such $\R\Gamma(\C,\EE)=0$, and  which is a local system of $\Q$-vector spaces outside of this set (=set of singularities of $\EE$), covariant with respect to the shifts by $2\pi i \Z$. We also assume that
the relative cohomology groups $V_{s,\theta}=H^{-1}(D_\epsilon(s), p,\EE)$ are   finite-dimensional.

The above data are equivalent to the following wall-crossing structure on $\C^\ast_\hbar$:

a) The  constant local system $\Gamma$ of lattices has fibers $H_1(\C, \Sigma, \Z)_{2\pi i\Z}$, where the notation means the space of coinvariants with respect to the action of the group $2\pi i\Z$. The abelian group $\Gamma$ is naturally isomorphic to $H_1(\C^\ast, \Sigma/2\pi i\Z, \Z)$, which gives a second description of the same lattice. 

b) The central charge $Z:\Gamma\to \C$ is given by the integration over $\gamma$ of the form $dz$ in the first description and of the form $dz/z$ in the second one.

c) The fiber over $\hbar\in \C^\ast$ of the corresponding local system of Lie algebras  is given by 
$\g_\hbar=\oplus_{\gamma\in \Gamma}\g_{\gamma,\hbar}=\left(\oplus_{s,s^\prime\in \Sigma}Hom(V_{s, \theta},V_{s^\prime, \theta})\right)_{2\pi i\Z}, \theta=Arg(\hbar)$.
The summand  $Hom(V_{s, \theta},V_{s^\prime, \theta})_{2\pi i\Z}$ belongs to the graded component $\g_{\gamma,\hbar}$, where $\gamma$ is represented by a path in $\C$ connecting $s$ and $s^\prime$.
The Stokes automorphisms which give the final piece of data in the definition of the WCS are uniquely determined by the maps $A_{s,s^\prime}$ defined as above. Notice that in general the WCS associated with a cosheaf $\EE$ is not analytic.

\end{rmk}

\begin{rmk}\label{rationality}
Under the assumptions of this subsection suppose further that for any $s\in \C^\ast$ and an open disc $D\subset \C^\ast$ the preimage $f_1^{-1}(D)$ has the finite rank cohomology (e.g this holds when $f:X\to \C$ is a regular function on the algebraic variety $X$). Then the matrix-valued formal series $A_{P_\pm}$ from the Remark \ref{two descriptions of WCS} are in fact rational. Hence the WCS described in Remark \ref {two descriptions of WCS} is analytic. Rationality can be shown by considerations similar to those in Section \ref{rationality of WCS}.

\end{rmk}

Motivated by this model example we would like to work out the idea of  Section \ref{infinite-dimensional vanishing cycles} in the case of the complexified Chern-Simons theory.  The $CS$ functional is an analog of the function $f_1$, but now $n=\infty$.

One can hope by analogy with the finite-dimensional case that the {\it CS functional considered as a holomorphic function $\widehat{\mathcal A}^{fr}_\C\to \C$ on the universal abelian cover $\widehat{\mathcal A}^{fr}_\C$ gives rise to a locally trivial infinite-dimensional fiber bundle outside of the set of critical values}.

 Then we replace each summand 
$H_n(f_1^{-1}(D_\epsilon(z_j)), f_1^{-1}(z_j+\epsilon e^{\sqrt{-1}Arg(s_1-z_j)},\Z)$ by the  abelian group $H^0_c((\ZZ_j(\alpha_{CS}), \phi_{CS\cdot e^{-\sqrt{-1}Arg(s_1-z_j)}}(\Z))$. Combining all critical values together we obtain a local system of infinite rank. Its fiber over $s_1\in \C-Critval(CS)$  can be informally thought as ``semi-infinite homology group" $H_{\infty}(\widehat{\mathcal A}_\C^{fr}, CS^{-1}(s_1), \Z)$. Hence the fiber over generic $s_1$ is isomorphic to the direct sum of finitely generated abelian groups
$$\bigoplus_{z_j\in Critval(CS)}H^0_c(\ZZ_j(\alpha_{CS}), \phi_{CS\cdot e^{-\sqrt{-1}Arg(s_1-z_j)} }(\Z)).$$ 
By analogy with the finite-dimensional case we expect that in general the fiber depends on the choice of Gabrielov paths, and moreover the above-defined local system can be extended to the critical values of $CS$ giving rise to a cosheaf $\EE_{CS}$ on $\C$ such that $\R\Gamma(\C,\EE_{CS})=0$.  By the Remark \ref{two descriptions of WCS} these data are equivalent to the data defining the corresponding WCS. Furthermore, they can be uniquely characterized by two matrix-valued series $A_{P_\pm}$.

{\it Let  $\FF_{CS}$ be the dual sheaf. } One can show that $\R\Gamma(\C, \FF_{CS})=0$.\footnote{This fact is not quite trivial, since the Mittag-Leffler property for the restrictions of $\FF_{CS}$ to the discs of increasing radii $R\to \infty$ is not satisfied.}



The above construction is equivariant with respect to the deck transformations. Hence we obtain (co)sheaf on $\C^\ast$ which is a local system outside of finitely many points. 

\begin{rmk} \label{Hodge structure of infinite rank}
a) By analogy with the finite-dimensional case one can speculate that the above-defined abelian group $H_{\infty}(\widehat{\mathcal A}_\C^{fr}, CS^{-1}(s_1), \Z)$ carries a weight filtration of infinite rank, presumably of finite length.
One can hope that  the vector space $H_{\infty}(\widehat{\mathcal A}_\C^{fr}, CS^{-1}(s_1), \Z)\otimes \C$ (maybe after some completion) caries a Hodge filtration of infinite rank.  In this sense we can speak about the Chern-Simons Hodge structure (of infinite rank).

b) The above-described  structure can be thought of as a ``perverse sheaf of infinite rank''. Perverse extension to the singular points is possible because the monodromy matrices are equal to the direct sum of the identity operator in a countably-dimensional space and an operator in a finite-dimensional vector space.

c) This structure should be compared with the previously discussed holomorphic version of the Morse-Novikov theory, where the set $(2\pi i)^2\Z$ is analogous to the set of periods of the holomorphic $1$-form $\alpha$.

\end{rmk}

Let $s$ be a critical value of the $CS$ and $D_\epsilon(s)$ be a small disc of radius $\epsilon>0$ with the center at $s$.

\begin{rmk} \label{cohomology of pair or vanishing cycle}  By analogy with the finite-dimensional case one can hope that there exists a theory of semi-infinite cohomology and a natural isomorphism of graded abelian groups
$$H^{\infty+\bullet}(CS^{-1}(D_\epsilon(s)),CS^{-1}(s+\epsilon),\Z)\to H^{\bullet}(Crit(CS)\cap CS^{-1}(s), \phi_{CS}(\Z)),$$
where $\infty$ should thought of as a ``complex dimension" of the infinite-dimensional manifold $\widehat{\mathcal A}_\C^{fr}$.
Even in this heuristic picture there is a potential caveat, since the $CS$ is not a proper map in any sense, and connected components of $Crit(CS)$ in general are non-compact.
\end{rmk}

\subsection{Interaction of critical points with different stabilizers}\label{interaction of critical points}

Let us now discuss an interesting phenomenon which reflects the ``interaction" of critical points of the CS functional which have different stabilizers.
Recall that for any $j_1,j_2\in J$ and a Stokes direction $\theta=Arg(\hbar)$ (i.e. $\theta=
Arg(\int_{z_1}^{z_2}\alpha_{CS})\ne 0, z_k\in \ZZ_{j_k}(\alpha_{CS}), k=1,2$) we have a natural homomorphism $$A_{\hbar,\gamma}: H_{c}^\bullet(\ZZ_{j_1}(\alpha_{CS}), \phi_{CS/\hbar}(\Z)) \to {H}_{c}^\bullet(\ZZ_{j_2}(\alpha_{CS}), \phi_{CS/\hbar}(\Z)),$$
where $\gamma$ is the homotopy class of a path joining two points in $\ZZ_{j_1}(\alpha_{CS})$ and  $\ZZ_{j_2}(\alpha_{CS})$ respectively such that $\partial \gamma=j_1-j_2\in \pi_0(\ZZ(\alpha_{CS}))$.
This homomorphism is compatible with the natural maps from the equivariant cohomology groups, i.e. we have the following commutative diagram:

\[ \begin{tikzcd}
{H}_{c, G_c}^\bullet(\ZZ_{j_1}(\alpha_{CS}), \phi_{CS/\hbar}(\Z)) \arrow{r} \arrow[swap]{d}{} & {H}_{c, G_c}^\bullet(\ZZ_{j_2}(\alpha_{CS}), \phi_{CS/\hbar}(\Z)) \arrow{d}{} \\%
{H}_{c}^\bullet(\ZZ_{j_1}(\alpha_{CS}), \phi_{CS/\hbar}(\Z)) \arrow{r}{}& {H}_{c}^\bullet(\ZZ_{j_2}(\alpha_{CS}), \phi_{CS/\hbar}(\Z))
\end{tikzcd}
\]

Let us consider this diagram in the case $G=SL(2,\C), G_c=SU(2)$. We take 
$\ZZ_{j_1}(\alpha_{CS})=\{\rho=1\}$ to be the isolated component consisting of the trivial connection. Then the corresponding equivariant cohomology will be $H^\bullet(BSU(2))=\Z[c_2], deg\, c_2=2$ and the non-equivariant cohomology is the group $\Z$ placed in degree $0$. The vertical arrow is the evaluation map at $c_2=0$.

Similarly, take the component $\ZZ_{j_2}(\alpha_{CS})$ which is a free $SL(2,\C)$-oribit of an irreducible rigid connection $\rho$. The non-equivariant cohomology is isomorphic to $H^\bullet(SU(2))=H^\bullet(S^3)$. On the other hand the  equivariant cohomology is isomorphic to $H^\bullet_{c,SU(2)}(SL(2,\C)/SU(2))=H^\bullet_{c,SU(2)}(\R^3)$. We see that there is a homomorphism from the equivariant cohomology of the trivial component $\rho=1$  to the one of the non-trivial rigid $\rho$, but not vice versa.

\begin{rmk}\label{orbits}

a) One can consider a finite-dimensional ``lattice model'' of the CS theory, by decomposing $M^3$ into a finite number of simplices, and assigning the group $G_c$ (or $G$) to edges, etc. Then the above considerations can be made rigorous for such a model.

b) Assume that $H^1(M^3,\Q)=0$. Consider the connected component $\ZZ_{j_0}(\alpha_{CS})$ consisting of the isolated Morse critical point $\rho=1$ corresponding to the trivial local system. Finite-dimensional local model for the above story is the pair $(\C^n, \alpha=d(\sum_{1\le i\le n}z_i^2))$. Finite-dimensional local model for the integration cycle is $\R^n$, which is a thimble for the Morse function $\sum_{1\le i\le n}z_i^2$.

c) Let $G=SL(2,\C)$.Consider a component $\ZZ_{j_1}(\alpha_{CS})\simeq SL(2,\C)$ which is the orbit of a rigid non-trivial flat connection $\rho_{j_1}$. Then the finite-dimensional local model is $(SL(2,\C)\times \C^{n-3}, d(\sum_{1\le i\le n-3}z_i^2))$. Finite-dimensional local model for the integration cycle is now $SU(2)\times \R^{n-3}$. The  ``formal expansion" map $R_j$  corresponds to the integral over the local Lefschetz thimble $\R^{n-3}$.

\end{rmk}

Remarks b) and c) give some sort of finite-dimensional model for the above example which shows why there are gradient trajectories from the trivial local system to non-trivial rigid ones, but not the other way around. One can expect that more generally there is a filtration on the equivariant cohomology group corresponding to the embedding of stabilizers. Moreover one can speculate that this filtration is related to the weight filtration 
in the Remark \ref{Hodge structure of infinite rank}.

\subsection{Resurgence from the functional integral point of view}\label{resurgence and functional integral}

In Section \ref{resurgence on wave functions} we discussed resurgence of the series in $\hbar$ which appear as pairings of quantum wave functions, interpreted as path integrals with boundary conditions. In this subsection we are going to discuss resurgence of the perturbative series arising in complexified Chern-Simons theory from the point of view of the functional integrals. In the next subsection we will see that that these two approaches are in fact equivalent.
\begin{rmk}\label{approaches to resurgence}
The approaches to resurgence based on either path integrals or functional integrals are not well-defined mathematically. In the case of Chern-Simons theory the hypothetical explanation of resurgence is based on the analyticity of the WCS associated with a pair of two complex Lagrangian submanifolds in $(\C^{\ast})^{2d}$. We will explain later what are these complex Lagrangian submanifolds.
\end{rmk}

The approach to resurgence via functional integrals uses  
the local perturbative expansions $R_j$ from Section \ref{local expansions} as well as the perverse sheaf of infinite rank $\FF_{CS}$ from Section \ref{monodromy}.

First, let us consider a special case when critical points of the complexified Chern-Simons functional are either the trivial flat connection $\rho=1$ or rigid  flat connections $\rho\ne 1$ with trivial stabilizers. 

For each rigid flat connection $\rho\ne 1$ with a trivial stabilizer (it is automatically a Morse-Bott critical point of $CS$) consider  the Borel transform ${\mathcal B}(I_\rho(\hbar))$ of the corresponding local asymptotic expansion $I_{\rho}(\hbar)\in \C[[\hbar]]$ of the $CS$ quantum partition function at  $\rho$. Then after the analytic continuation the corresponding multi-valued analytic function ${\mathcal B}(I_\rho(\hbar)):={\mathcal B}(I_\rho(\hbar))(s)$ has poles at the points $s=s_m^\rho$ which belong to the arithmetic series $z_\rho+(2\pi i)^2m, m\in \Z$. Here $z_\rho=CS(\rho)$ is the value of the complexified Chern-Simons functional at the critical point $\rho$. Recall that the image of the set of critical values of $CS$ in $\C/(2\pi i)^2\Z$ coincides with the image of the Beilinson-Borel regulator 
$K_3^{ind}(\overline{\Q})\to \C/(2\pi i)^2\Z$.

\begin{conj}\label{resurgence of the non-trivial inputs}
Local perturbative expansion $I_\rho(\hbar)\in \C[[\hbar]]$ of the CS quantum partition function at the flat connection $\rho\ne 1$ is resurgent. Analytic continuation of  the Borel transform of the germ ${\mathcal B}(I_\rho(\hbar)), \rho\ne 1$ along  paths joining pairs of points $s_m^\rho$ and $s_{\rho^\prime}^l,\rho^\prime\ne 1$ recovers  ${\mathcal B}(I_{\rho^\prime}(\hbar))$.

\end{conj}

Let us comment on the Conjecture.

The analytic continuation of the germ 
$\varphi_{\rho,m}$ of the Borel transform  ${\mathcal B}(I_\rho(\hbar))$ at $s_m^\rho$ to the point $s_{l}^{\rho^\prime}$ has the form $f_{m,l,\rho,\rho^\prime}(s)+n_{m,l}{\log(s-s_l^\rho)\over{2\pi i}}\varphi_{\rho^\prime,l}$, where $f_{m,l,\rho,\rho^\prime}(s)$ is holomorphic in a neighborhood of $s_{l}^{\rho^\prime}$, and $n_{m,l}\in \Z$ is the Stokes index.
Analytic germs $\varphi_{\rho,m}$ enjoy the equivariance property with respect to the action of the group $(2\pi i)^2\Z$ on $\varphi_{\rho, m}$ with fixed $\rho\ne 1$.

In this case the local asymptotic expansion at the trivial connection is $I_{\rho=1}(\hbar)\in \hbar^{3/2}\C[[\hbar]]$. Its Borel transform ``sees" the Borel transforms of those for other rigid $\rho\ne 1$ as logarithmic jumps, as we discussed previously.

More generally, the function $CS$ can have critical points which are non-rigid flat connections or which are flat connections with non-trivial stabilizers. In this case the local perturbative expansions can contain fractional powers of $\hbar$ and positive integer powers of logarithms. Then for the Borel transform we use the formulas from Remark \ref{appearance of logarithms}.

\begin{rmk}\label{infinite derivative}
Recall Remark \ref{relation to Hodge structure}. In the finite-dimensional Morse case the Borel transforms of  local expansion of the modified integral is given by the  $(dim_\C X/2-1)$-th derivative of the volume of the nearby cycle. In the Chern-Simons case the dimension is infinite. \footnote{Physicists use the $\zeta$-regularized dimension which is equal to zero.} If one assumes that the fibers are infinite-dimensional Calabi-Yau varieties  the Borel transforms   ${\mathcal B}(I_\rho(\hbar))(s)$ are the pairings of infinite-dimensional cycles with elements of  the ``middle term" $F_{\infty/2}(H^{\infty-1}(CS^{-1}(s))$ of the above-discussed infinite Hodge filtration. Taking further derivatives with respect to $s$ we will see the part of the putative infinite rank Hodge filtration which corresponds to $F_{\infty/2-i}(H^{\infty-1}(CS^{-1}(s)), 0\le i<\infty$. These considerations suggest that there exists a semi-infinite mixed Hodge structure associated with the above-discussed sheaf $\FF_{CS}\widehat{\otimes}\C$. The weight filtration was already discussed in Remark \ref{Hodge structure of infinite rank}.
\end{rmk}




\subsection{Relation to the integration over the cycle of unitary connections}\label{compact gauge group}

Motivated by the considerations of Section \ref{monodromy} we are going to speculate about the infinite-dimensional generalizations similar to those in the previous subsection. In order to simplify the exposition we will assume that $M^3$ is a homological sphere, $G=SL(2,\C), G_c=SU(2)$, all non-trivial connections are rigid and all critical values are different. We will make a proposal about the integration over the space of unitary connections which can be thought of as the relation between CS wall-crossing structure and Reshetikhin-Turaev invariants.

We have the holomorphic map $f_{CS}=e^{CS/2\pi i}: {\mathcal A}_\C^{fr}\to \C^\ast$ which has finitely many critical values $Critval(f_{CS}):=\{z_j\}_{1\le j\le m}$. On the universal covering they give rise to a finitely many arithmetic series with the step $(2\pi i)^2$. Recall  that outside of the set of critical values we have (conjecturally) a local system of infinite rank $\FF$, whose fiber $\FF_z$ can be informally understood as $H_\infty(\widehat{\mathcal A}_\C^{fr}, CS^{-1}(2\pi i\,log(z)),\Z)$. Our assumption on $M^3$ implies that the locus of critical points $Crit(f_{CS})$ consists of an isolated Morse critical point $x_1=\rho_{triv}$ corresponding to the trivial connection as well as Morse-Bott critical submanifolds corresponding to rigid non-trivial connections. In fact as a set $Crit(f_{CS})= \{\rho_{triv}\}\cup (\sqcup_{2\le j\le m}SL(2,\C))$. The corresponding thimbles are (appropriately defined) unions of gradient lines of the real-valued function $Re(CS/\hbar), \hbar\in i\R_{>0}$. More precisely, in case if there are gradient lines connecting two different critical points we slightly rotate $\hbar$ when defining the corresponding thimble. As a result projections of thimbles will be disjoint rays emanating from $z_j, 1\le j\le m$ with the slope $-{\pi\over 2}+\epsilon$ where $0<\epsilon\ll 1$.

The thimble emanating from $\rho_{triv}$ is an ``infinite-dimensional cell" isomorphic to $\R^\infty$, whereas generalized thimbles emanating from other connected components of $Crit(f_{CS})$ are isomorphic to $S^3\times \R^{\infty-3}$ (cf. Remark \ref{orbits}).
The generic fiber of $\FF$ is non-canonically isomorphic to $\C[T^{\pm 1}]\otimes \C^{\{Critval(f_{CS})\}}$, where $T$ is the monodromy operator. 

Let ${\mathcal A}_{c}^{fr}$ denote the set of unitary (i.e. $SU(2)$) framed connections on $M^3$. By analogy with the finite-dimensional case we can hope that the ``integration cycle''  ${\mathcal A}_{c}^{fr}$ is homologically equivalent to an integer linear combination 
$\sum_\rho n_{\rho_j} \gamma_{\rho_j}$, where $\gamma_{\rho_j}$ are the above-discussed generalized thimbles. 

\begin{conj}\label{unitary cycle and thimbles}

For any $k\in \Z_{\ge 1}$ one has
$$\int_{{\mathcal A}_{unit}^{fr}}e^{k CS\over{2\pi i}}{\mathcal D}A=\sum_{\rho_j,\, s.t.\,|exp(CS(\rho_j)/2\pi i)|\le 1}n_{\rho_j}\int_{\gamma_{\rho_j}}e^{k CS\over{2\pi i}}{\mathcal D}A.$$

Here $n_{\rho_j}\in \Z$. Furthermore $n_{\rho_j}=1$ for unitary flat connections (including the trivial one). For other summands $|exp(CS(\rho_j)/2\pi i)|< 1$, and $n_{\rho_j}$ is the virtual number of the gradient lines of $Re(CS/2\pi i)$ emanating from a unitary connection and terminating at the flat connection $\rho_j$.

\end{conj}

The LHS of the formula in the Conjecture \ref{unitary cycle and thimbles} is understood in terms of the representation theory of the quantum groups at roots of $1$ (Witten-Reshetikhin-Turaev invariants $WRT_k$), and in fact it is an element of a cyclotomic field.\footnote{Later we will discuss a generalization of these invariants in terms of the so-called Nahm sums.} The integrals in the RHS can be interpreted (conjecturally) in terms of the Chern-Simons analytic WCS. According to our general philosophy it gives a holomorphic section of an analytic $m$-dimensional vector bundle over a small disc in $\C_\hbar$ centered at the $\hbar=0$. 
Then the integrals in the RHS of the conjecture should be understood in terms of the evaluation of this section at the points $\hbar=2\pi i/k, k=1,2,...$.

\begin{rmk}\label{non-rigid case} We expect that in the above conjecture the assumptions that $M^3$ is a homological sphere as well as the rigidity of non-trivial flat connections can be dropped. Then $\gamma_{\rho_j}$ will be appropriately defined integration cycles.

\end{rmk}





Consider the generating function
$$G(w)=\sum_{1\le k\le \infty}(\int_{{\mathcal A}_{c}^{fr}}e^{kCS(A)/2\pi i}{\mathcal D}A)w^k.$$

\begin{conj}\label{RT vs WCS}
Assume the Conjecture \ref{analytic WCS for CS}.
Then the generating series $G(w)$ converges in the disc $|w|<1$ and analytically continues to $\C$ with singularities at $\{0\}\cup Critval(CS)$. 
More precisely the above-discussed local system $\FF$ is $(2\pi i)^2 \Z$-equivariant and hence descends to a local system $\FF^\prime$ on $\C^\ast-Critval(f_{CS})$.
The function $G(w)$ gives a morphism of sheaves of abelian groups $\FF^\prime\to \OO_{\C^\ast-Critval(f_{CS})}$.

\end{conj}

The next conjecture makes sense only at the physics level of rigor.

\begin{conj}\label{unitary cycle and Hodge structure}
The  abelian group
$$H_{{\infty}}(\widehat{\mathcal A}_\C^{fr}, Im(CS)\ll 0, \Z)$$
 has finite rank. It is canonically isomorphic to the direct sum of analogous local semi-infinite cohomology groups over the set of critical values of the CS functional. The ``fundamental class" $[{\mathcal A}_{unit}^{fr}]$ defines an element of this abelian group.
\end{conj}

As an illustration of the discussion of this subsection, let us consider a finite-dimensional toy-model example. Namely, the analog of $CS$ is $f=z-log(z)-(1+\epsilon):\widetilde{\C^\ast}\to \C$,
where $\epsilon>0$ is sufficiently small. The analog of $f_{CS}$ is $e^f:\C^\ast\to \C^\ast$. The analog of ${\mathcal A}_{c}^{fr}$ is the cycle $|e^f|=1$. We have:  $Crit(f)=\{z=1\}$ and $Critval(e^f)=\{e^{-\epsilon}\}$.\footnote{We choose $\epsilon\ne 0$ in order to move the critical values of $e^f$ from the unit circle.}
The analog of the integral over the cycle of unitary connections is $I_k=\int_{|e^f|=1}e^{kf(z)}{dz\over{z}}$. By our assumption the cycle $|e^f|=1$ contains a connected component which is a compact cycle homotopy equivalent to the circle $|z|=1$, while the other connected components are contractible.

Then $I_k=e^{k(-1-\epsilon)}\int_{|z|=1}z^{-k}e^{kz}{dz\over{z}}=2\pi i\, e^{k(-1-\epsilon)}{k^k\over{k!}}$. Consider the generating function $G(w)=\sum_{k\ge 1}w^k I_k=2\pi i\sum_{k\ge 0}(we^{-1-\epsilon})^k k^k/k!=2\pi i\sum_{k\ge 0}u^k k^k/k!$, where $u=we^{-1-\epsilon}$. As a function of $w$ it has singularities at the critical values of $e^{-f}$ as well at $w\in \{0,\infty\}$. Hence $G(w)$ has ramification at the set  $w\in \{0, e^{\epsilon}, \infty\}$.  Setting $0^0=1$ we can rewrite $N(w)$ such as follows:
$$G(w=ue^{1+\epsilon})=\sum_{k\ge 0}u^k\int_{|z|=1}(e^z/z)^k{dz\over{z}}=\int_{|z|=1}{u^{-1}\over{u^{-1}-{e^z\over{z}}}}{dz\over{z}}.$$

For a fixed $u\in \C-\{0,e^{-1}\}$ zeros of the denominator $zu^{-1}-e^z$ form a countable discreet subset $S_u\subset \C$. Furthermore the group $H_1(\C-S_u,\Z)$ is spanned by simple loops about points of $S_u$. Let us fix $u_0, 0<u_0\ll 1$. Then the fundamental group $\pi_1(\C-\{0,e^{-1}\},u_0)$ acts by permutations on $S_{u_0}$, and hence on $H_1(\C-S_u,\Z)$ by permutations of the basis. One can show that there is a bijection $S_{u_0}\simeq \Z\cup \{\ast\}$ such that the monodromy about $u=0$ is given by the map $i\mapsto i+1, \ast\mapsto \ast, i\in \Z$, and the monodromy about $u=e^{-1}$ is given by $0\leftrightarrow \ast, i\mapsto i, i\in \Z-\{0\}$.

Let  $z_\ast(u_0)\in S_{u_0}$ denote the unique zero corresponding to  $\ast\in \Z\cup \{\ast\}$ and lying in the disc $|z|\le 1$. By analytic continuation in $u$ we obtain the analytic function $z_\ast(u)=u+O(u)$ for $0<|u|\ll 1$. For such $u$ we have
$$\int_{|z|=1}{u^{-1}\over{u^{-1}-{e^z\over{z}}}}{dz\over{z}}=\int_{|z|=1}{dz\over{z-ue^z}}={2\pi i\over{1-ue^{z_\ast(u)}}}.$$

Hence $G(ue^{1+\epsilon})$ has an analytic continuation which is a multivalued analytic function on $\C_u-\{0,e^{-1}\}$.


\subsubsection{Level shift and resurgence}\label{level shift}

Recall that according to the original computations of Witten (see e.g. [Wit], formula (2.17)) the perturbative expansion of the Chern-Simons partition function at an isolated unitary flat connection $\rho_j$ depends not on the value $k/2$   but on $(k+c_2(G_c))/2$ where $c_2(G_c)$ is the second Chern class of the gauge group $G_c$. E.g. for $G_c=SU(n)$ we have $c_2(G)=2n$, hence for $G_c=SU(2)$ the level $k$ should be replaced by the  shifted level $k\mapsto k+2$. Then in the sum over flat unitary connections we observe the extra factor $e^{2CS(\rho_j)/2\pi i}$ which comes from the $\zeta$-regularized Hessian $(det(\partial^2CS(\rho_j)))^{-1/2}$ (it is also known as $\eta$-invariant). 
When identifying the Chern-Simons functional integral over the cycle of unitary connections with  the Reshetikhin-Turaev  invariant one uses in the latter representations of the quantized enveloping algebra $U_q(sl(2))$ at $q=e^{2\pi i/(k+2)}$, i.e. the quantization parameter $q$ depends on the shifted level. Same shift (a.k.a. dual Coxeter number) appears in the approach to Chern-Simons theory via WZW model or, mathematically, via Kazhdan-Lusztig equivalence of the fusion tensor category of the affine Lie algebra at the integer level with the corresponding modular tensor category of representations of the quantized enveloping algebra.

Comparing with our previous discussion when we used the small parameter $\hbar=2\pi i/k$ we conclude that
we have  another small parameter in the story, namely $\hbar^\prime=2\pi i/(k+2)$. Although the difference in these two choices seems to be minor, it {\it could} affect the resurgent properties of the perturbative expansions.

Nevertheless we claim that this is not the case, i.e. that the change of the small parameter $\hbar$ to $\hbar^\prime$ does not affect substantially the resurgent properties of arising perturbative series. In order to see that, we notice that $\hbar^\prime=\hbar/(1-a\hbar)$, where $a=-2/2\pi i$.
Recall the Borel transform  $\hbar^\lambda\mapsto s^\lambda/\Gamma(\lambda+1)$. It is convenient to use the modified Borel transform
$\hbar^\lambda\mapsto s^\lambda/\Gamma(\lambda)$ (these two transformations are conjugate by the operator $sd/ds$). 
Then applying the modified Borel transform to 
$$(\hbar^\prime)^\lambda=\hbar^\lambda+a\lambda \hbar^{\lambda+1}+a^2{\lambda(\lambda+1)\over 2!}\hbar^{\lambda+2}+...$$
we obtain
$${s^{\lambda}\over{\Gamma(\lambda)}}+a{s^{\lambda+1}\over \Gamma(\lambda)}+{a^2\over 2!}{s^{\lambda+2}\over \Gamma(\lambda)}+...=
e^{as}{s^\lambda\over \Gamma(\lambda)}.$$

Therefore at the level of the modified Borel transform the change of a small parameter $\hbar$ to $\hbar^\prime$ amounts to the multiplication by the entire function $e^{as}$.
Notice that in our case the function $e^{as}$ is invariant under the shift $s\mapsto s+(2\pi i)^2$.

In the next subsection we will consider sums which generalize the WRT-invariants. As a small parameter we will use $\hbar=2\pi i/N$ where $N\in \Z_{\ge 1}$. The reader should think of $N$ as the analog of $k+2$ rather than $k$.

\begin{rmk}\label{MA equation}
We remarked above that the value of the Hessian of the Chern-Simons functional $CS$ at a critical point is equal to the exponent of the appropriately rescaled critical value of $CS$. This coincidence suggests that the functional $CS$ satisfies an infinite-dimensional analog of the Monge-Amp\`ere equation.

\end{rmk}

\subsection{Generalized Nahm sums and quantum wave functions}\label{Nahm sums}

 Let $N,d\in \Z_{\ge 1}$ and $q=e^{2\pi i\over{N}}$ be the corresponding to $N$ primitive root of unity. We will use the following version of the $q$-factorial: $(n)_q!=(1-q)...(1-q^{n})$. 
 
Let us introduce the following data:

a) A rational polyhedron $P\subseteq [0,1]^d$. 

b) A rational symmetric matrix: $b=(b_{i_1i_2})_{1\le i_1,i_2\le d}, b_{i_1i_2}=b_{i_2i_1}, b_{i_1i_2}\in \Q$.

c) A collection of integers $a=(a_i)_{1\le i\le d}, a_i\in \Z$.

d) A collection of rational numbers $c_i\in \Q, 1\le i\le d$.

 e) A character $\chi: \Z^d\to \mu_\infty, \chi(\overline{j})=\prod_{1\le i\le d}\chi_i^{j_i}$, where $\mu_\infty\subset \C^\ast$ is a subgroup of all roots of $1$, i.e. $\mu_\infty=\varinjlim_{n}\mu_n$, and $\overline{j}=(j_1,...,j_d)\in \Z^d$.
 
With these data we  associate the following {\it generalized Nahm sums at roots of unity}\footnote{It was pointed to us by Sergei Gukov that in mathematical physics a closely related object appears under the name ``fermionic sums".}considered as a function in $N$:
$$Z_{N}=\sum_{\substack{0\le j_1,...,j_{d}\le N-1\\
\overline{j}/N:=(j_1/N,...,j_{d}/N)\in P}}\chi(\overline{j})\prod_{1\le i\le d}((j_i)_q!)^{a_{i}}q^{{1\over{2}}\sum_{1\le i_1,i_2\le d}b_{i_1i_2}j_{i_1}j_{i_2}+\sum_i c_ij_i}.$$
 Here by definition we set $q^\lambda=e^{2\pi i\lambda/N},\lambda\in \Q$.
 

Let us define  the following multivalued analytic function 
$$f_{a,b,\chi}({x})=-\sum_{\beta}{a_i}Li_2(x_i)-\sum_{1\le i\le d} log(\chi_i)log(x_i)+{1\over {2}}\sum_{i_1,i_2}b_{i_1i_2}log(x_{i_1})log(x_{i_2}), $$
where $x=(x_1,...,x_d)\in (\C^\ast)^d.$

We would like to compare three objects of different nature: critical points of $Z_N$ as $N\to \infty$, critical points of the function $f_{a,b,\chi}$, and intersection points of certain complex Lagrangian subvarieties in the symplectic manifold $(\C^\ast)^{2d}_{x,y}, x=(x_1,...,x_d),y=(y_1,...,y_d)$ endowed with the symplectic form $\omega^{2,0}=\sum_{1\le i\le d}{dx_i\over{x_i}}\wedge {dy_i\over{y_i}}$.

\begin{rmk}\label{relation to regulator}
Recall that the set of values of the complexified CS functional at critical points belong to the countable set which is the set of values of the Beilinson-Borel regulator map $reg: K_3^{ind}(\overline{\Q})\to \C/(2\pi i)^2\Z$. It is not hard to see that the same set (up to torison) can be described as the set of critical values of $f_{a,b,\chi}$ for some $a,b,\chi$.


\end{rmk}

Here is the idea in the case when $P=[0,1]^d$. Let us write $Z_{N}= \sum_{\overline{j}}Z_{N, \overline{j}}$. If we understand $Z_{N,\overline{j}}$ as a function in $x=(x_1,...,x_d)=(q^{j_i})_{1\le i\le d}\in (S^1)^d\subset (\C^\ast)^d$ then ``critical points" of  the function $\overline{j}\mapsto Z_{N, \overline{j}}$ should correspond to such $x\in (\C^\ast)^d$  that the quotient $Z_{N, (j_1,..., j_i+1,..., j_d)}/Z_{N, (j_1,..., j_i,..., j_d)}$ is equal to $1$ for all $1\le i\le d$, as well as  (possibly) to the points where $x_i=1$. The former conditions give  the following equations:
$$\chi_i(1-qx_i)^{a_i}q^{c_i}\prod_{1\le j\le d}x_j^{b_{ij}}=1, 1\le i\le d.$$

As $N\to \infty$ we have $q\to 1$, and the set of solutions to these equations can be identified  with the intersection of two complex Lagrangian submanifolds in $(\C^\ast)^{2d}_{x,y}$:
$$L_1=\{y_i^{-1}=(1-x_i)^{a_i}, 1\le i\le d\}$$
and
$$L_2=\{\chi_i^{-1}=y_i^{-1}\prod_{1\le j\le d}x_j^{b_{ij}}, 1\le i\le d\}.$$

Furthermore, the same set can be identified with the set of critical points of the multivalued function $f_{a,b,\chi}$.

Notice that this set does not have to belong to $(S^1)^d$, hence the critical points of $f_{a,b,\chi}$ do not correspond literally to the ``critical points" of $Z_N$ discussed above.

\begin{rmk}\label{extended Lagrangians}
In the above considerations we excluded  the possibility $x_i=1$. In order to take it into account, as we will see in examples, we need to modify $L_1$ and consider  the extended Lagrangian variety $L_1^{ext}=\{y_i^{-1}=(1-x_i)^{a_i}1\le i\le k\le d\}\cup\{x_i=1, k+1\le i\le d\}$. Furthermore for general $P$ (which is not necessary of full dimension $d$) we will need also to modify $L_2$ and to consider the Lagrangian subvariety $L_2^{ext}$ which is a union of  Lagrangian tori corresponding to faces of $P$ (multiplicative analogs of conormal bundles). From the point of view of the above discussion  the intersection $L_1^{ext}\cap L_2^{ext}$ corresponds to the set of ``critical points" of  the restriction of the function $\overline{j}\mapsto Z_{N, \overline{j}}$  to  faces of $P$.
\end{rmk}

Next we would like to rewrite $Z_N$  as the pairing of the discrete analogs of quantum wave functions.

Let 
$$\psi_1:=\psi_1(N)=(\psi_1)_{j_1,...,j_{d}}=\prod_{1\le i\le d}((j_i)_q!)^{a_i}, 0\le j_i\le N-1$$
and
$$\psi_2:=\psi_2(N)=(\psi_2)_{j_1,...,j_{d}}={\bf 1}_P({\overline{j}\over{N}})\chi(\overline{j})q^{{1\over{2}}\sum_{1\le i_1,i_2\le d}b_{i_1i_2}j_{i_1}j_{i_2}+\sum_i c_ij_i}, 0\le j_i\le N-1$$
be two vectors in $(\C^{\Z/N\Z})^{d}$. Here ${\bf 1}_P$ is the characteristic function of the polyhedron $P$.
Then $Z_{N}=(\psi_1,\psi_2):=\sum_{(j_1,...,j_d)}(\psi_1)_{j_1,...,j_{d}}(\psi_2)_{j_1,...,j_{d}}$.

In order to relate vectors $\psi_1$ and $\psi_2$ to quantum wave functions associated with Lagrangian subvarieties in $(\C^\ast)^{2d}$ let us look at their asymptotics under the conditions that $N\to \infty, j_i\to \infty, i=1,..., d$, and there is a finite limit  $\overline{j}/N\to \overline{t}:=(t_1,...,t_d)\in (0,1)^d$. One can show that these asymptotic expansions look such as follows:

$$(\psi_1)_{j_1,...,j_d}\sim c_{N}\psi_{L_1}, c_N=(\sqrt{N}e^{{-2\pi {\sf i}\over 24}(N+{1\over N})}e^{2\pi {\sf i}\over 8})^{\sum_ia_{i}}$$
$$(\psi_2)_{j_1,...,j_d}\sim {\bf 1}_P(\overline{t})exp\left({1\over 2\hbar}\sum_{i_1,i_2}b_{i_1,i_2}log(x_{i_1})log(x_{i_2})\right):=\psi_{L_2},$$
where ${\sf i}=\sqrt{-1}, \hbar={2\pi {\sf i}\over N}, x_m=exp({2\pi {\sf i} t_m}), 1\le m\le d$.  

In the RHS of the first formula we use the following notation:
$$\psi_{L_1}(x)=\prod_{1\le i\le d}e^{a_i\left(-{Li_2(x_i)\over{\hbar}}+log(\sqrt{1-x_i})+\sum_{n\ge 0}\hbar^{2n+1}{\zeta(-(2n+1))\over{n!}}R_{2n}(x_i)\right)}=$$
$$e^{a_i\left(-{F_0\over \hbar}+\sum_{1\le j\le d}log(\sqrt{1-x_i})+\sum_{1\le j\le d}\sum_{n\ge 1}\hbar^n{\zeta(-n)\over {n!}}[(x\partial_x)^n(log(1-x))]_{x=x_j}\right)},
$$
where $F_0(x_1,...,x_d)=\sum_{1\le i\le d}Li_2(x_i),$ $R_0(x)={x\over 1-x}, R_m(x)=(x\partial_x)^m({x\over 1-x}), m\ge 1$ and $\zeta(s)$ is the Riemann $\zeta$-function. The notation in the RHS shows that the corresponding expression is the quantum wave function of the Lagrangian subvariety $L_1\subset (\C^\ast)^{2d}$. 


We will see that we need to replace $\psi_{L_1}$ by $\psi_{L_1^{ext}}$ from the Remark \ref{extended Lagrangians}. Then $\psi_{L_1^{ext}}$ coincides with $\psi_{L_1}$ on the subvariety $L_1\subset L_1^{ext}$.

In order to explain this replacement in an example, we will need to recall few facts about quantum tori.
 The algebraic quantum torus $A_q(2d)$ is the unital $\C[q,q^{-1}]$-algebra\footnote{One can work over  $\Z[q,q^{-1}]$ as well.} generated by invertible generators $\hat{x}_1,..., \hat{x}_d, \hat{y}_1,...,\hat{y}_d$ subject to the relations $\hat{x}_i\hat{x}_j=\hat{x}_j\hat{x}_i,\hat{y}_i\hat{y}_j=\hat{y}_j\hat{y}_i, \hat{x}_i\hat{y}_j=q^{\delta_{ij}}\hat{y}_j\hat{x}_i$. 
 
 If $\epsilon=e^{2\pi i/N}$ is a primitive root of unity then the reduction $A_\epsilon(2d)$ of $A_q(2d)$ at $q=\epsilon$ has a big center generated by $\hat{x}_i^{\pm N}, \hat{y}_i^{\pm N}$. The algebra $A_\epsilon(2d)$ becomes the space of global sections of the  bundle of matrix algebras  of the size $N^d$ over the spectrum of the center (Azumaya algebra).


Analogously to the case of $D$-modules one defines the notion of holonomic module over $A_q(2d)$.
The reduction of such a module at a primitive root of $1$ is a module over the corresponding Azumaya algebra, hence a module over its center. The support of the corresponding coherent sheaf on the spectrum of the center is a closed Lagrangian subvariety of the algebraic symplectic manifold 
${\bf G}_m^{2d}$. For a large class of holonomic $A_q(2d)$-modules the support is essentially independent on the order $N$ of the root of $1$ (more precisely it depends on the residue of $N$ modulo some integer, see [Ko1] for the details). This class of holonomic modules seems to include those which appear in relation to the Chern-Simons theory and generalized Nahm sums at roots of unity.

\

Let us now explain in this framework the necessity of enlargement of $L_1$ to $L_1^{ext}$ in the case $d=1$. For $d\ge 1$ one has to consider the $d$-th external tensor power of the one-dimensional case.

We take $q=e^{2\pi i/N}$, and consider as an example the generalized Nahm sum 
$$Z_N=1+\sum_{j=0}^{N-1}((j)_q!)^{-1}q^{bj^2}=(\psi_1,\psi_2),$$
where $\psi_1,\psi_2\in \C^{\Z/N\Z}$ are given by
$$(\psi_1)_j=((j)_q!)^{-1}, (\psi_2)_j=q^{bj^2}, j=0,1,..., N-1.$$

Let $\hat{x},\hat{y}$ be the generators of the quantum torus $A_q(2)$ realized by the $N\times N$ matrices in the standard basis of  $\C^{\Z/N\Z}$
by
$$\hat{x}: e_j\mapsto q^{j}e_j, j=0,1,...,N-1,$$
$$\hat{y}: e_j\mapsto e_{j+1}, j=0,1,..., N-2, e_{N-1}\mapsto e_0.$$
Clearly $\hat{x}\hat{y}=q\hat{y}\hat{x}$. The corresponding matrix algebra is just the fiber of the above-mentioned bundle of matrix algebras at the point $(1,1)\in (\C^\ast)^2$.

Notice that $(1-\hat{x}-\hat{y})\psi_1=-e_0$.
Then 
$$(1-\hat{x})(1-\hat{x}-\hat{y})\psi_1=0.$$

This calculation motivates us to introduce the {\it reducible} left cyclic $A_q(2)$-module $E_1=A_q(2)/A_q(2)(1-\hat{x})(1-\hat{x}-\hat{y})$. 
One checks that the support of the reduction of $E_1$ at any primitive root $\epsilon$ of $1$ is the curve $L_1^{ext}=\{(1-x)(1-x-y)=0\}$ which is a union of two disjoint connected components, each isomorphic to  $\C^\ast$. The latter is also the spectrum of the center of the reduced algebra $A_\epsilon(2)$.
Notice that disconnectedness of $L_1^{ext}$ implies that the reduction of $E_1$ at any root of $1$ (or even in its formal neighborhood) is {\it canonically} isomorphic to a sum of two irreducible modules. We remark that the module $E_1$ together with the cyclic generator is well-defined for any $q$. In particular setting $q=e^{\hbar}$ where $\hbar$ is a formal parameter, we can speak about homomorphisms from $E_1$ to the $\C[[\hbar]]$-module of formal quantum wave functions associated with the smooth part of $L_1^{ext}$. Equivalently we can speak about solutions to the equation
$(1-\hat{x})(1-\hat{x}-\hat{y})\psi=0$ in the space of formal quantum wave functions. They are locally defined up to a factor from $\C[[\hbar]]^\times$. The above choice of $\psi_1$  at $N$-th root of $1$ and its asymptotics as $N\to \infty$ (see above) gives a preferred choice of the quantum wave function.

A similar computation for $\psi_2$ gives us the holonomic right cyclic $A_q(2)$-module $E_2$ such that the support of its reduction at the primitive root  $\epsilon$ is the curve $L_2=\{y=x^{2b}\}\subset (\C^\ast)^2$. Moreover we also get a preferred choice of a quantum wave function.


The quantum wave function $\psi_{L_1^{ext}}$ corresponding to $L_1^{ext}$ coincides with the above-defined $\psi_{L_1}$ when restricted to $L_1\subset L_1^{ext}$.  
Since $L_1^{ext}=Supp(E_1)$ contains also the connected component $x=1$ we should assign to the latter the quantum wave function. It is interesting to observe that this quantum wave function is a non-trivial series $\psi_{x=1}(y,\hbar)=\sum_{l\ge 0}\hbar^lf_l(y)$ rather than the naively expected delta-function supported on the line $x=1$.  In order to find $\psi_{x=1}(y,\hbar)$ one should  solve the following $q$-difference equation

$$(1-\hat{x}-\hat{y})\psi_{x=1}(y,\hbar)=1,$$
where $\hat{y}$ acts by multiplication on the variable $y$ and $\hat{x}(g(y,\hbar))=g(e^{\hbar}y,\hbar)$ for any $g$.

Looking for a WKB expansion $\sum_{l\ge 0}\hbar^lf_l(y)$ of the solution to this equation and taking into account that
$g(e^\hbar y)=e^{\hbar y\partial_y}(g(y))$ we obtain an infinite system of differential equations for $f_l(y), l\ge 0$ which can be solved by induction:

$$ -yf_l=1+\sum_{i=1}^l{(y\partial_y)^i\over{i!}}f_{l-i},$$
where we set $f_{-1}=0$ in order to accommodate the case $l=0$.
Then we find $f_0(y)=-{1\over y}, f_1(y)=-({1\over y}+{1\over y^2}),...$
In particular we find by induction that $f_l(y), l\ge 1$ is a polynomial in $y^{-1}$.


Notice that the finite Fourier transform of the above-defined vector $\psi_1=(\psi_1)_{j_1,...,j_d}\in (\C^{\Z/N\Z})^d$ has an asymptotic expansion as $\hbar=2\pi {\sf i}/N\to 0$  given by $\psi_{x=1}(y,\hbar)=\sum_{l\ge 0}\hbar^lf_l(y)$.

\

One can analyze the reducible module $E_1$ from a different point of view.
Recall the notion of log extension mentioned below Example \ref{compactification for quantum dilog}. Then the above phenomenon is related to the fact that the closures of complex Lagrangian submanifolds $x=1$ and $1-x-y=0$ intersect at a point which belong to a compactifying log-divisor. The quantum wave function corresponding to $Supp(E_1)$ should ``remember" this intersection, since the very definition of the category of holonomic $DQ$-modules  depends on a choice of log extension of the given complex symplectic manifold. In our case the symplectic manifold is  $(\C^\ast)^2$, and the (non-unique) log extension is determined by a choice of ambient Poisson toric surface in which $(\C^\ast)^2$ is an open symplectic leaf. It seems plausible that after making blow-ups at the intersection points of the irreducible components of supports with log-divisors  one can extend the QWFS to the corresponding symplectic manifold.

\

In the higher-dimensional case one has a similar story. In particular one has a  Lagrangian subvariety $L_1^{ext}\supset L_1$ and  the corresponding  reducible left cyclic holonomic module $E_1=E_{L_1^{ext}}$. Furthermore, for a  non-trivial rational polytope $P$ one has a reducible right cyclic holonomic module $E_2=E_{L_2^{ext}}$ over the quantum torus as well as the corresponding Lagrangian subvariety $L_2^{ext}\supset L_2$. The subvariety $L_2^{ext}$ can be described in the following way. First, consider  the union of conormal bundles to all faces of $P$:
$$T^\ast P:=\cup_{\sigma\in \{faces\,of\, P\}}T^\ast_\sigma \C^d\subset \C^{2d}.$$
It gives rise to the union of Lagrangian tori $T^\ast P/\Z^{2d}\subset (\C^\ast)^{2d}$. Then $L_2^{ext}$ is obtained from  $T^\ast P/\Z^{2d}$ by applying  the transformation:
$$x_i\mapsto x_i,$$
$$y_i\mapsto \chi_iy_i\prod_{1\le j\le i}x_j^{b_{ij}}, 1\le i\le d.$$

Recall  now the vector $\psi_2(N)$ which depends on a rational polyhedron $P$, a rational symmetric matrix $(b_{ij})$ and a character $\chi$ of finite order.  Let $\overline{A}_\epsilon(2d)$ denote the finite-dimensional quotient of the reduction  ${A}_{q=\epsilon}(2d)$ of the algebra $A_q(2d)$ at $q=\epsilon$, where $\epsilon=e^{2\pi i/N}$. The quotient is obtained by imposing the additional relations $\hat{x}_i^N=\hat{y}_i^N=1, 1\le i\le d$.

One can show that there exists a positive integer $M$ depending on these data as well as a collection of $M$ cyclic right holonomic $A_q(2d)$-modules $E_{2,l}, l\in \Z/M\Z$ depending rationally on $q$ such that the following holds: 

{\it each $E_{2,l}$ admits a well-defined finite-dimensional reduction at $q=e^{2\pi i/N}$ where $N\in l+M\Z_{\ge 0}$ which is a right cyclic $\overline{A}_\epsilon(2d)$-module with the cyclic vector $\psi_{2,l}(N)$ and moreover this cyclic module coincides with the one generated by $\psi_2(N)$}.


The asymptotics of $\psi_{2}(N)$ as $N\to \infty$ and fixed $l=N\, mod\, M$  gives a preferred choice of a quantum wave function $\psi_{L_2^{ext},l}$ for the right $A_q(2d)$-module associated with $L_2^{ext}$ in agreement with Section \ref{sums and polytopes}. 

\

Recall that for a complex subvariety $Y$ in a complex torus $(\C^\ast)^k$  one can define the corresponding tropical subvariety $trop(Y)\subset \R^k$ by first extending scalars to $\C((z))$ and going to the corresponding non-archimedean subvariety $Y^{an}$ in the non-archimedean torus over the field of Laurent series $\C((z))$, and then taking the image of $Y^{an}$ under the ``tropical map" $log|\bullet|$. We can apply this construction to our situation.

Let us  impose the following   ${\bf Genericity\, Assumption}$:

{\it The tropical subspaces $trop(L_1^{ext})$ and $trop(L_2^{ext})$ considered as subsets of $\R^{2d}$ intersect at $0\in \R^{2d}$ only (this guarantees the absence of intersections of $L_1^{ext}$ and $L_2^{ext}$ ``at infinity"), and the intersection at each of the finitely many points of $L_1^{ext}\cap L_2^{ext}$ is transversal and belongs to the smooth part of $L_2^{ext}$.}

In this case the local pairings $\langle \psi_{L_1^{ext}}, \psi_{L_2^{ext},l}\rangle_p, p\in L_1^{ext}\cap L_2^{ext}$ of quantum wave functions give us  a finite collection of formal series with algebraic coefficients. We expect that these local pairings  form a resurgence package. 
Moreover, we expect that the Borel resummations of the series $\langle \psi_{L_1^{ext}}, \psi_{L_2^{ext},l}\rangle_p$ are analytic functions such that for each $l=N\, mod\, M$ the value of some $\Z$-linear combination of these analytic functions  at $\hbar={2\pi i \over N}$ coincides with the {\it finite} pairing $(\psi_{1}(N),\psi_{2}(N))$.

In case if the Genericity Assumption is not satisfied we expect the following. 

\begin{itemize}
\item There exists an algebraic symplectic manifold $(M,\omega^{2,0})$ (log extension) which contains $(\C^\ast)^{2d}$ as a Zariski open symplectic submanifold such
that the closures $\overline{L_i}^{ext}\subset M, i=1,2$ compact.

\item The standard QWFS on $(\C^\ast)^{2d}$ extends to a QWFS on $M$.

\item Quantum wave functions $\psi_{L_1^{ext}}, \psi_{L_2^{ext},l}$ can be extended to the quantum wave functions $\psi_{\overline{L}_1^{ext}}, \psi_{\overline{L}_2^{ext},l}$ for $\overline{L_i}^{ext}\subset M, i=1,2$.\footnote{Here we use the previously mentioned conjecture that the theory of quantum wave functions exists for singular Lagrangian varieties.}

\item Each rotated symplectic form $Re(e^{i\theta}\omega^{2,0}), \theta\in \R/2\pi \Z$ is ``convex at infinity" (hence the corresponding Fukaya category is well-defined).

\end{itemize}

Under these assumptions we expect that there exist covariantly constant objects $\overline{E}_1$ and $\overline{E}_{2,l}$ in the local system (over $S^1_\theta$) of  the local Fukaya categories of $\overline{L_1}^{ext}\cup\overline{L_2}^{ext}$ with coefficients extended to a cyclotomic field such that the restrictions of them to $L_1^{ext}$ and to the smooth part  $L_{2}^{ext,smooth}\subset L_2^{ext}$ respectively are the trivial trivialized rank one local systems.  Moreover  the conjecturally well-defined pairing $\langle\psi_{\overline{L}_1^{ext}}, \psi_{\overline{L}_2^{ext},l}\rangle\in (Ext^d_{F_{\overline{L_1}^{ext}\cup\overline{L_2}^{ext},loc}}(\overline{E}_1,\overline{E}_{2,l})\otimes {\bf B})^\Z$  (see Remark \ref{pairing of wave functions and Fukaya}) restricted to $L_1^{ext}\cap L_2^{ext,smooth}$ coincides with the corresponding pairing discussed previously.  Similarly to the above we expect that these data gives rise to a resurgence package (depending on $l$) and the finite pairing $(\psi_1(N),\psi_2(N)), N\in l+M\Z$ coincides with a $\Z$-linear combination of the Borel resummations of the corresponding series in fractional powers of $\hbar$ and polynomilas in $log\, \hbar$ evaluated at $\hbar={2\pi i \over N}$.

For each transversal intersection point $p\in L_1^{ext}$ and the smooth part of  $L_2^{ext}$  and for each residue 
$l\in \Z/M\Z$ we have the corresponding local pairing $\langle \psi_{L_1^{ext}}, \psi_{L_2^{ext},l}\rangle_p$ which  is an element of $F[[\hbar]]\subset \C[[\hbar]]$, where $F\subset \overline{\Q}$ is a number field. On the other hand for each $l=N\, mod\, M$  the pairing  $(\psi_{1}(N),\psi_{2}(N))$  belongs to the cyclotomic field $\Q(e^{2\pi i\over{CN}})$ for some universal constant $C\in \Z$ depending on the data for generalized Nahm sum.


\

Finally, we formulate the following analyticity conjecture.

\begin{conj}\label{analytic extension}
Let $Z_N=(\psi_1(N),\psi_2(N))$ defined as above.
Then the generating function $\sum_{N\ge 1}Z_{N}w^N$ 
extends analytically to $\C^\ast_w-\{exp({reg(K_0(\C)\over {2\pi i}})\}$, where $reg$ is the notation for the Beilinson-Borel regulator.

\end{conj}
This conjecture should be compared with the Conjecture \ref{RT vs WCS}.
We remark that although considerations of this subsection are motivated by the CS theory, 
they make sense without it. This remark indicates that there is   a class of holonomic cyclic $q$-$D$-modules 
associated with $K_2$-Lagrangian submanifolds in $(\C^\ast)^{2d}$ which gives rise to wall-crossing 
structure derived from the asymptotic expansions of the local pairings of the cyclic vectors (i.e. quantum wave functions) 
corrected by the Stokes indices given by the virtual numbers of pseudo-holomorphic discs. 
The WCS derived from the CS theory should be a special case of this more general framework.

\subsection{Quantum wave functions  in the Chern-Simons theory}\label{local CS and wave functions}

Let us discuss the relation of the CS theory to our theory of  quantum wave functions. 
In the compact case there are several ways to express this relationship mathematically. 
E.g. one can  construct an invariant of  $M^3$ presenting it as a result of  
a surgery along a knot or link (Witten-Reshetikhin-Turaev invariant) or via a simplicial decomposition of $M^3$ (Turaev-Viro invariant).

In the approach via the surgery we decompose $M^3=M^3_-\cup M^3_0\cup M^3_+$ 
where $M^3_0=\Sigma\times [0,1]$ and $\Sigma$ is an oriented $2$-dimensional closed surface, $M^3_\pm$ 
are oriented manifolds with  boundaries, such that $\partial M^3_-=\Sigma\times \{0\}$ and 
$\partial M^3_+=\Sigma\times \{1\}$.

Recall the space  $\Aa^{fr}_\C:=\Aa_\C^{fr}(M^3)$, the infinite-dimensional 
complex manifolds of framed $G$-bundles with connections on $M^3$ (see Section \ref{CS theory}). 
Consider the submanifold $\Aa_\C^{fr,fl}(M^3)\subset \Aa^{fr}_\C(M^3)$ 
consisting of framed bundles with connections $A$ on $M^3$ such that

1) Both restrictions $A_{|M^3_\pm}$  are flat.

2) The restriction $A_{|\Sigma\times [0,t]}$ is flat for any $0\le t\le 1$.

Let us denote by $L_+$ and  $L_-$  the holomorphic Lagrangian submanifolds in the complex symplectic moduli space $\Aa^{fr,fl}(\Sigma)$ 
of framed bundles with flat $G$-connections on $\Sigma\simeq \partial M^3_\pm$, 
which are defined via the restrictions of flat connections on $M^3_+$  and  $M^3_-$ respectively to the boundary.

Recall the space of paths $P(\Aa^{fr,fl}(\Sigma), L_-,L_+):=P(L_-,L_+)$ from Section \ref{space of paths}. 
By definition we have a canonical projection $\pi: \Aa_\C^{fr,fl}(M^3)\to P(L_-,L_+)$. 
Recall that $P(L_-,L_+)$ carries a holomorphic closed $1$-form $\eta:=\eta_{L_-,L_+}$ (see Section \ref{space of paths}). 
Similarly, the moduli space $\Aa^{fr}_\C(M^3)$ carries the holomorphic closed $1$-form $\alpha_{CS}$ (see Section \ref{CS theory}). 
The pull-backs of $\eta_{L_-,L_+}$ and $\alpha_{CS}$ to $\Aa_\C^{fr,fl}(M^3)$ give holomorphic closed $1$-forms whose sets of zeros coincide. 

This observation makes plausible the idea that {\it the  perturbative expansions for the complexified CS theory as well as the corresponding WCS 
should coincide
with those for the  path integral with holomorphic Lagrangian boundary conditions $L_\pm$ (see Section \ref{infinite-dimensional case}).}
\footnote{Although critical points of both action functionals are naturally identified, 
there might be some delicate issues related to the ``critical points at infinity",
which physicists also have realized although in a different way, see  [GuPut].} 

We also expect that there exist natural compactifications $\overline{L}_-$ and $\overline{L}_+$ obtained as closures of ${L}_-$ and ${L}_+$
in the appropriate partial compactification of $\Aa^{fr,fl}(\Sigma)$
as well as  the quantum wave functions $\psi_{\overline{L}_\pm}$ supported on ${\overline{L}_\pm}$ 
such that the local pairings  $\langle\psi_{\overline{L}_-}|\psi_{\overline{L}_+}\rangle_{\ZZ_\alpha}$ at the connected 
components $\ZZ_{\alpha}\subset \overline{L}_-\cap \overline{L}_+$
coincide with the local perturbative expansions of the above functional (equivalently, path) integral and give rise to a resurgence package.

Conjecturally, we can approach the same problem from the point of view of generalized Nahm sums. As we mentioned previously, Nahm sums
give rise to holomorphic Lagrangian submanifolds in a symplectic torus. The latter can be  partially compactified giving rise to the closures of 
(extended)
Lagrangian subvarieties.
The corresponding quantum wave functions and WCS are interesting by 
themselves without any connection to the complexified CS theory. Nevertheless there are indications that the latter is a special case of the former.



In order to motivate the conjecture below let us recall that after Witten, Turaev and Reshetikhin we know that the 
Chern-Simons partition function for the {\it compact} group $G_c$ can be written in terms of generalized Nahm sums at roots of $1$. 
Recall   $K_2$-Lagrangian subvarieties from Section \ref{Nahm sums}, or more precisely the closures of their extended versions
 $\overline{L_1}^{ext}$ and $\overline{L_2}^{ext}$. 
Let  $\psi_{\overline{L}_1^{ext}}, \psi_{\overline{L}_2^{ext},l}$ be the quantum wave functions defined in the loc.cit. 
We expect that these data give rise to a resurgence package derived from the corresponding WCS.
Then  the conjectural relationship of this story with the complexified Chern-Simons theory can be summarized such as follows.

\begin{conj}\label{CS vs Nahm sums}
1) There is an embedding $\overline{L}_-\cap \overline{L}_+\hookrightarrow \overline{L_1}^{ext}\cap \overline{L_2}^{ext}$ 
with the image which does not depend on a choice of the model producing the 
Lagrangians $L_\pm$ ( e.g.  it does not depend on the surgery along a link, simplicial decomposition of $M^3$, etc.).

2) The resurgence package derived from the pair $\overline{L}_-,\overline{L}_+$ and the quantum wave functions $\psi_{\overline{L}_\pm}$
 is a subset (in the natural sense) of the one derived from the pair $\overline{L_1}^{ext}, \overline{L_2}^{ext}$ and the quantum 
 wave functions $\psi_{\overline{L}_1^{ext}}, \psi_{\overline{L}_2^{ext},l}$ (see  Section \ref{Nahm sums} for the notation).

3) In particular Stokes indices $n_{ij}$ which appear in the  WCS associated with the complexified CS theory are equal to the virtual numbers  of  
holomorphic (in the rotated complex structure) $2$-gons with the 
boundary  on  $K_2$-Lagrangian subvariety  $\overline{L_1}^{ext}\cup \overline{L_2}^{ext}$. 
Degenerating the almost complex structure one can obtain a description of $n_{ij}$ in terms of certain spectral networks.

\end{conj}


\begin{rmk} \label{dg-intersections}

a) A choice of the parameters for the $K_2$-Lagrangian subvarieties as well as 
partial compactification of the symplectic manifold (and hence closures of the Lagrangian subvarieties) in the
 Conjecture \ref{CS vs Nahm sums} was not specified.
 This makes it not very practical from the point of view of computations
 in the complexified CS theory as long as we want to apply it to the flat non-unitary connections.
Perhaps such a   choice can be  tested 
via the interpretation of $\hat{Z}$-invariant of Gukov and others as the pairing of the  quantum wave functions from part 2) of the 
Conjecture \ref{CS vs Nahm sums} . \footnote{We thank to Sergei Gukov 
for discussions on this matter. After our paper was submitted to the arXiv, he and Putrov posted an interesting paper [GuPut] where the previous work on
the relation of Stokes indices and $\hat{Z}$-invariant was summarized and further extended from the physics perspective.}

b) In order to make things mathematically more precise we should understand the 
intersections of Lagrangians as derived intersections of Lagrangian dg-subschemes in the $(-1)$-symplectic stacks 
of all flat connections on $M^3$. Moreover the  intersection locus carries a perverse sheaf of vanishing cycles of an 
analytic function obtained from the complexified CS functional.\footnote{In the recent paper [GunSaf] a similar idea is 
discussed in relation to the approach to skein modules developed in [GunJordSaf].} 
Then the Stokes isomorphisms should be understood as isomorphisms of the cohomology 
groups with coefficients in the sheaves of vanishing cycles.
\end{rmk}

We remark that the above story can be slightly generalized when the $K_2$-Lagrangian subvariety $L_2^{ext}$ 
is replaced by the multiplicative shift $x_i\mapsto x_i, y_i\mapsto c_iy_i, 1\le i\le d$. 
The corresponding holonomic $A_q(2d)$-module corresponds on the Chern-Simons side to the $3$-dimensional manifold with  boundary. 
Typically this boundary is the boundary of the tubular neighborhood of a link in a compact $3$-dimensional manifold. 
This generalization is related to colored link invariants (e.g. to the colored Jones polynomial). 

Regardless of whether $M^3$ has the boundary or not the pairing of the corresponding quantum wave functions 
can be written as a contour integral of the products of quantum dilogarithm functions (see e.g. [AnHan],  [DimGuLeZa]).

\vspace{3mm}

{\bf References}

\vspace{2mm}
[Ab1] M. Abouzaid,  On the wrapped Fukaya category and based loops,  
arXiv:0907.5606.

\vspace{2mm}

[Ab2] M. Abouzaid, A cotangent fiber generates the Fukaya category,  
arXiv:1003.4449.

\vspace{2mm}

[Ab3] M. Abouzaid, Family Floer cohomology and mirror symmetry, arXiv:1404.2659.

\vspace{2mm}

[AnHan] J. Andersen, S. Hansen, Asymptotics of the quantum invariants for surgeries on the figure 8 knot, arXiv:0506456.

\vspace{2mm}

[AgSch] A. D'Agnolo, P. Schapira, Quantization of complex Lagrangian 
submanifolds, arXiv:math/0506064.

\vspace{2mm}

[AgKas] A. D'Agnolo, M. Kashiwara, Riemann-Hilbert correspondence for holonomic
D-modules, arXiv:1311.2374.

\vspace{2mm}

[BarGi] V. Baranovsky, V. Ginzburg, Conjugacy classes in loop groups and 
$G$-bundles on elliptic curves, arXiv:alg-geom/9607008.

\vspace{2mm}

[Bou] P. Bousseau, Holomorphic Floer Theory and Donaldson-Thomas invariants, arXiv:2210.17001.

\vspace{2mm}

[BrDy] C. Brav, T. Dyckerhoff, Relative Calabi-Yau structures, arXiv:1606.00619.

\vspace{2mm}

[BreSo] P. Bressler, Y. Soibelman, Mirror symmetry and deformation quantization, arXiv:hep-th/0202128.

\vspace{2mm}

[CharHurt] B. Charbonneau, J. Hurtubise, Singular Hermitian-Einstein monopoles on the product of a circle and a Riemann surface, arXiv:0812.0221.

\vspace{2mm}

[ChW] S. Cherkis, R. Ward, Moduli of Monopole Walls and Amoebas, arXiv:1202.1294.

\vspace{2mm}

[DKo] J. Dongo, M. Kontsevich, Logarithmic formality theorem, work in progress.
\vspace{2mm}
 
[Di] T. Dimofte, Quantum Riemann surfaces in Chern-Simons theory, arXiv:1102.4847.

\vspace{2mm}

[BeBlDeEs] A. Beilinson, S. Bloch, P. Deligne, H. Esnault, Periods for 
irregular connections on curves, preprint.

\vspace{2mm}
[BezKal] R. Bezrukavnikov, D. Kaledin, Fedosov quantization in algebraic context, arXiv:math/0309290.

\vspace{2mm}

[Del] P. Deligne, letter to Ramis as of 7.01.1986, published in Singularit\'es irr\'eguli\`eres - Correspondance et documents
Pierre Deligne - Bernard Malgrange - Jean-Pierre Ramis,
Documents Math\'ematiques 5 (2007).

\vspace{2mm}
[DelPh] E. Delabaere, F. Pham, Resurgent methods in semi-classical analysis, 
Annales IHP, Section A, 71:1, 1999, 1-94.
\vspace{2mm}

[DimGuLeZa] T. Dimofte, S. Gukov, J. Lenells, D. Zagier, Exact results for perturbative Chern-Simons theory with complex gauge group, arXiv:0903.2472.

\vspace{2mm}

[DimGabGo] T. Dimofte, M. Gabella, A. Goncharov, K-Decompositions and 3d Gauge Theories, arXiv:1301.0192.

\vspace{2mm}

[DiSa] A. Dimca, M. Saito, On the cohomology of a general fiber of a polynomial 
map, Composition Mathematica, 85, 1993, 299-309.
\vspace{2mm}

[DoRez] A. Doan, S. Rezchikov, Holomorphic Floer Theory and the Fueter equation, arXiv:2210.12047.

\vspace{2mm}[DonMar] R. Donagi, E. Markman, Cubics, integrable systems, and Calabi-Yau 
threefolds, arXiv:alg-geom/9408004.

\vspace{2mm}

[El1] Y. Eliashberg, Weinstein manifolds revisited, arXiv:1707.03442.

\vspace{2mm}

[ElGivHof] Y. Eliashberg, A. Givental, H. Hofer, Introduction to Symplectic 
Field Theory,
arXiv:math/0010059.

\vspace{2mm}

[EsKoMoeZo] A. Eskin, M. Kontsevich, M. Moeler, A. Zorich, Lower bounds for Lyapunov exponents of flat bundles on curves,
arXiv: 1609.01170.

\vspace{2mm}

[FreJos] J. Fres\'an, P. Jossen, Exponential motives, book.

\vspace{2mm}

[GV] M.Garay, D. Van Straten, Classical and quantum integrability, Moscow Math. J., 10:3, 2010, p.519-545.

\vspace{2mm}

[GaMoNe1] D. Gaiotto, G. Moore, A. Neitzke,  Four-dimensional wall-crossing via 
three-dimensional field theory, arXiv:0807.4723.

\vspace{2mm}
[GaMoNe2] D. Gaiotto, G. Moore, A. Neitzke,  Wall-crossing, Hitchin Systems, and 
the WKB Approximation, arXiv:0907.3987.

\vspace{2mm}

[GaMoNe3] D. Gaiotto, G. Moore, A. Neitzke,  Wall-crossing in coupled 2d-4d 
systems, arXiv:1103.2598.
\vspace{2mm}

[GaMoNe4] D. Gaiotto, G. Moore, A. Neitzke, Framed BPS states, arXiv:1006.0146.

\vspace{2mm}

[GanParSh-1]  S. Ganatra, J. Pardon, V. Shende, Covariantly functorial Floer 
theory on Liouville sectors, arXiv:1706.03152.

\vspace{2mm}

[GanParSh-2] S. Ganatra, J. Pardon, V. Shende, Microlocal Morse theory of wrapped Fukaya categories, arXiv: 1809.0880.

\vspace{2mm}

[GarZag] S. Garoufalidis, D. Zagier, Asymptotics of Nahm sums at roots of unity, arXiv:1812.07690.

\vspace{2mm}

[GeKazh] I. Gelfand, D. Kazhdan, Some problems of the differential geometry and the calculation of cohomologies of lie algebras of vector fields,
Dokl.Akad.Nauk Ser.Fiz. 200 (1971) 269-272.

\vspace{2mm}

[GuKasSch] S. Guillermou, M, Kashiwara, P. Schapira, Sheaf quantization of 
Hamiltonian isotopies,  arXiv:1005.1517.
\vspace{2mm}

[GuMarPut] S. Gukov, M. Marino, V.Putrov, Resurgence in complex Chern-Simons theory, arXiv:1605.07615.
\vspace{2mm}

[GuPut] S. Gukov, V.Putrov,  On categorification of Stokes matrices in Chern-Simons theory, arXiv:2403.12128.

\vspace{2mm}

[GunSaf] S. Gunningham, P. Safronov, Deformation quantization and perverse sheaves, arXiv:2312.07595.
\vspace{2mm}

[GunJordSaf] S. Gunningham, D. Jordan, P. Safronov, The finiteness conjecture for skein modules, arXiv:1908.05233.

\vspace{2mm}

[Hub] R. Huber, \'Etale cohomology of rigid analytic varieties and adic spaces, Aspects of Mathematics, E30. Friedr. Vieweg and Sohn, Braunschweig, 1996.
\vspace{2mm}
[J] Xin Jin, Holomorphic branes correspond to perverse sheaves, arXiv: 
1311.3756.

\vspace{2mm}

[KaKoPa1] L. Katzarkov, M. Kontsevich, T. Pantev, Hodge theoretic aspects of 
mirror symmetry, arXiv:0806.0107.
\vspace{2mm}

[KaKoPa2] L. Katzarkov, M. Kontsevich, T. Pantev, Bogomologov-Tian-Todorov 
theorems for Landau-Ginzburg models, arXiv:1409.5996.

\vspace{2mm}

[Kap] A. Kapustin, A-branes and Noncommutative Geometry,  	arXiv:hep-th/0502212.
\vspace{2mm}

[KapOr] A. Kapustin, D. Orlov, Remarks on A-branes, Mirror Symmetry, and the Fukaya category,  	arXiv:hep-th/0109098.

\vspace{2mm}

[KasSch1] M. Kashiwara, P. Schapira, Deformation quantization modules, 
arXiv:1003.3304.

\vspace{2mm}

[KasSch2] M. Kashiwara, P. Schapira, Constructibility and duality for simple 
holonomic modules on complex symplectic manifolds,arXiv:math/0512047.

\vspace{2mm} 

[KasSch3] M. Kashiwara, P. Schapira, Regular and irregular holonomic
D-modules, arXiv:1507.00118.

\vspace{2mm}

[KasSch4] M. Kashiwara, P. Schapira, Sheaves on manifolds, book, Springer, 2010.

\vspace{2mm}

[Ko1] M. Kontsevich, Holonomic $D$-modules and positive characteristic, 
arXiv:1010.2908.

\vspace{2mm}

[Ko2] M. Kontsevich, Deformation quantization of algebraic varieties, 
arXiv:math/0106006.
\vspace{2mm}

[Ko3] M. Kontsevich, Geometry in dg-categories, preprint IHP, 2015.
\vspace{2mm}

[Ko4] M. Kontsevich, Symplectic geometry of homological algebra, talk at 
Arbeitstagung, MPI, 2009.
\vspace{2mm}

[Ko5]  M. Kontsevich,Deformation quantization of Poisson manifolds, I, arXiv:q-alg/9709040.

\vspace{2mm}

[Ko6] M. Kontsevich,  XI Solomon Lefschetz Memorial Lecture Series: Hodge structures in non-commutative geometry, arXiv:0801.4760.

\vspace{2mm}
[KoSe] M. Kontsevich, G. Segal,  Wick rotation and the positivity of energy in quantum field theory,  arXiv:2105.10161.

\vspace{2mm}
[KoSo1] M. Kontsevich, Y. Soibelman, Stability structures, motivic 
Donaldson-Thomas invariants and cluster transformations, arXiv:0811.2435.

\vspace{2mm}

[KoSo2] M. Kontsevich, Y. Soibelman, Affine structures and non-archimedean 
analytic spaces, math.AG/0406564.

\vspace{2mm}
[KoSo3] M. Kontsevich, Y. Soibelman, Notes on A-infinity algebras, A-infinity 
categories and non-commutative geometry. I, math.RA/0606241.

\vspace{2mm}

[KoSo4] M. Kontsevich, Y. Soibelman, Deformations of algebras over operads and 
Deligne's conjecture, arXiv:math/0001151.
\vspace{2mm}

[KoSo5] M. Kontsevich, Y. Soibelman, Cohomological Hall algebra, exponential 
Hodge structures and motivic Donaldson-Thomas invariants, arXiv:1006.2706.

\vspace{2mm}

[KoSo6] M. Kontsevich, Y. Soibelman,
Homological mirror symmetry and torus fibrations, arXiv:math/0011041.

\vspace{2mm}

[KoSo7] M. Kontsevich, Y. Soibelman, Wall-crossing structures in 
Donaldson-Thomas invariants, integrable systems and Mirror Symmetry, 
arXiv:1303.3253.

\vspace{2mm}
[KoSo8] M. Kontsevich, Y. Soibelman, Wall-crossing structures and algebraicity, in preparation.

\vspace{2mm}
[KoSo9] M. Kontsevich, Y. Soibelman, Holomorphic Floer theory, in preparation.

\vspace{2mm}

[KoSo10] M. Kontsevich, Y. Soibelman, Deformation theory I. Draft of the book available on
{math.ksu.edu/ $\tilde{}$ {soibel}}.
\vspace{2mm}

[KoSo11] M. Kontsevich, Y. Soibelman, Airy structures and symplectic geometry of topological recursion,
arXiv: 1701.09137.

\vspace{2mm}

[KoSo12] M. Kontsevich, Y. Soibelman,
Analyticity and resurgence in wall-crossing formulas, arXiv: 2005.10651.
\vspace{2mm}

[KoSo13]  M. Kontsevich, Y. Soibelman, Deformations of algebras over operads and
Deligne's conjecture, arXiv:math/0001151.

\vspace{2mm}

[KoZo] M. Kontsevich, A. Zorich, Lyapunov exponents and Hodge theory, arXiv:hep-th/9701164.

\vspace{2mm}

[Mal] B. Malgrange, \'Equations diff\'erentielles \`a coefficients
polynomiaux, Birkh\"auser, 1991.

\vspace{2mm}

[Mn] P. Mnev, Lectures on Batalin-Vilkovisky formalism and its applications in topological quantum field theory, arXiv:1707.08096.

\vspace{2mm}
[Moch1] T. Mochizuki, Holonomic D-modules with Betti structure,
arXiv:1001.2336.
\vspace{2mm}

[Moch2] T. Mochizuki, Doubly periodic monopoles and $q$-difference modules,   arXiv:1902.03551.
\vspace{2mm}

[Moch3] T. Mochizuki, Triply periodic monopoles and difference modules on elliptic curves,
arXiv:1903.03264.
\vspace{2mm}

[NadZas] D. Nadler, E. Zaslow, Constructible Sheaves and the Fukaya Category,  
arXiv:math/0604379.

\vspace{2mm}

[NeTsy] R. Nest, B. Tsygan, Deformations of symplectic Lie algebroids,
deformations of holomorphic symplectic structures,
and index theorems, arXiv:math/9906020.

\vspace{2mm}
 
 [PT] T. Pantev, B. Toen, Poisson geometry of the moduli of local systems on smooth varieties, arXiv:1809.03536.
 
\vspace{2mm}

[RaSauZh]    J.-P. Ramis, J. Sauloy, C. Zhang, Local analytic classification of 
q-difference equations, arXiv:0903.0853.

\vspace{2mm}

[Rai] E. Rains, The noncommutative geometry of elliptic difference equations, 
arXiv:1607.08876.

\vspace{2mm}

[Sab1] C. Sabbah, On the comparison theorem for elementary irregular D-modules, 
Nagoya Math. J. 141, 1996, 107-124.

\vspace{2mm}

[Sab2] C. Sabbah,On a twisted de Rham complex II,, arXiv:1012.3818.

\vspace{2mm}

[Sab3] C. Sabbah, Syst\`emes holonomes d'\'equations aux q-diff\'erences, In:
D-modules and microlocal geometry, 1992, p. 125-147.

\vspace{2mm}

[Sab4] C. Sabbah, Introduction to Stokes structures, arXiv:0912. 2762.

\vspace{2mm}

[SabSai] C. Sabbah, M. Saito, Kontsevich's conjecture on an algebraic formula for vanishing cycles of local systems,
arXiv:1212.0436  [pdf, ps, other]  math.AG

\vspace{2mm}

[ShTrZas1] V.Shende, D. Treumann, E. Zaslow, Legendrian knots and constructible 
sheaves,
arXiv:1402.0490.

\vspace{2mm}   

[ShTrWZas] V.Shende, D. Treumann, H. Williams, E. Zaslow, Cluster varieties from 
Legendrian knots, arXiv:1512.08942.

\vspace{2mm} 

[Se1] P. Seidel, Fukaya categories and Picard-Lefschetz Theory. European 
Mathematical Society, 2008, 326p.

\vspace{2mm}

[SGA7] SGA7, Lecture notes in math. 340 (1973)
\vspace{2mm}

[So1] Y. Soibelman, Remarks on Cohomological Hall algebras and their 
representations,  arXiv:1404.1606.
\vspace{2mm}

[Syl] Z. Sylwan, On partially wrapped Fukaya categories, arXiv:1604.02540.

\vspace{2mm}
[Vo] A. Voros, The return of the quartic oscillator (the complex WKB method), 
Annales Institut H. Poincar\'e, 29:3, 1983, 211-338.

\vspace{2mm}

[Wit] E. Witten, Analytic Continuation Of Chern-Simons Theory, 	arXiv:1001.2933.

\vspace{2mm}

Addresses:

M.K.: IHES, 35 route de Chartres, F-91440, France, {maxim@ihes.fr}

Y.S.: Department of Mathematics, KSU, Manhattan, KS 66506, USA, 
{soibel@math.ksu.edu}

\end{document}